\documentclass[11pt]{article}
\usepackage{latexsym,color,amsmath,amsthm,amssymb,amscd,amsfonts}

\newtheorem*{remark*}{Remark}

\def\R{{\mathbb R}}

\def\upddots{\mathinner{\mkern 1mu\raise 1pt \hbox{.}\mkern 2mu
\mkern 2mu \raise 4pt\hbox{.}\mkern 1mu \raise 7pt\vbox {\kern 7
pt\hbox{.}}} }

\newcommand{\rf}{\overline{\F}}

\newcommand{\mtau}{{\overline{\tau}}}
\newcommand{\spk}{{Sp_{2k}(\F)}}

\newcommand{\soo}{{SO_{2n+1}(\F)}}
\newcommand{\spn}{{Sp_{2n}(\F)}}
\newcommand{\gspn}{{GSp_{2n}(\F)}}
\newcommand{\mspn}{{\overline{Sp_{2n}(\F)}}}
\newcommand{\mspk}{\overline{Sp_{2k}(\F)}}
\newcommand{\tspn}{T_{Sp_{2n}}(\F)}
\newcommand{\bspn}{B_{Sp_{2n}}(\F)}

\newcommand{\zspn}{{Z_{Sp_{2n}}}(\F)}
\newcommand{\pspn}{P_{Sp_{2n}}(\F)}
\newcommand{\wspn}{W_{Sp_{2n}}(\F)}
\newcommand{\wspntag}{W'_{Sp_{2n}}(\F)}
\newcommand{\zspk}{{Z_{Sp_{2k}}}(\F)}
\newcommand{\zspka}{{Z_{Sp_{2k}}}(\A)}
\newcommand {\vt} {{\overrightarrow{t}}}

\newcommand{\tspk}{T_{Sp_{2k}}(\F)}
\newcommand{\tglm}{T_{GL_m}(\F)}
\newcommand{\tgln}{T_{GL_n}(\F)}
\newcommand{\bgln}{B_{GL_n}(\F)}
\newcommand{\bglm}{B_{GL_m}(\F)}
\newcommand{\zgln}{Z_{GL_n}(\F)}
\newcommand{\zglm}{Z_{GL_m}(\F)}
\newcommand{\zglma}{Z_{GL_m}(\A)}
\newcommand{\zglr}{Z_{GL_r}(\F)}
\newcommand{\zgll}{Z_{GL_l}(\F)}

\newcommand{\mspm}{\overline{Sp_{2m}(\F)}}
\newcommand{\mspr}{\overline{Sp_{2r}(\F)}}

\newcommand{\mbk}{\overline{B_{SP_{2k}}(\F)}}
\newcommand{\F}{\mathbb F}
\newcommand{\K}{\mathbb K}
\newcommand{\Of}{\mathbb O_{\F}}
\newcommand{\Pf}{\mathbb P_{\F}}
\newcommand{\Oq}{\mathbb O_{\mathbb Q  _2}}
\newcommand{\Pq}{\mathbb P_{\mathbb Q _2}}
\newcommand{\Kn}{Sp_{2n}(\mathbb O_\F)}
\newcommand{\Kk}{Sp_{2k}(\mathbb O_\F)}
\newcommand{\msigr}{\overline{\sigma_r}}
\newcommand{\msig}{{\overline{\sigma}}}
\newcommand{\spx}{{Sp(X)}}
\newcommand{\vl}{v_{\lambda}}
\newcommand{\tl}{t_{\lambda}}
\newcommand{\mspx}{\overline{Sp(X)}}
\newcommand{\mgspx}{\overline{GSp(X)}}
\newcommand{\gspx}{{GSp(X)}}
\newcommand{\mpt}{\overline{P_{\overrightarrow{t}}(\F)}}
\newcommand{\mmt}{\overline{M_{\overrightarrow{t}}(\F)}}

\newcommand{\mt}{M_{\overrightarrow{t}}(\F)}
\newcommand{\nt}{N_{\overrightarrow{t}}(\F)}
\newcommand{\pt}{{P_{\overrightarrow{t}}}(\F)}
\newcommand{\A}{\mathbb A}
\newcommand{\N}{\mathbb N}
\newcommand{\Z}{\mathbb Z}

\newcommand{\half}{\frac{1}{2}}
\newcommand{\psii}{\widetilde{\psi}}
\newcommand{\ab} {|\!|}
\newcommand{\gl}{{GL_n(\F)}}
\newcommand{\gln}{{GL_n(\F)}}
\newcommand{\gll}{{GL_l(\F)}}
\newcommand{\glr}{{GL_r(\F)}}
\newcommand{\glm}{{GL_m(\F)}}

\newcommand{\gamapsi}{\gamma_\psi^{-1}}
\newcommand{\Q}{\mathbb Q}
\newcommand{\C}{\mathbb C}

\def\>{\rangle}
\def\<{\langle}

\newtheorem{lem}{Lemma}[section]
\newtheorem{thm}{Theorem}[section]
\newtheorem{cor}{Corollary}[section]

\numberwithin{equation}{section}
\newcommand{\chip}{{\chi_\pi}}
\newcommand{\etac}{\chi}

\newcommand{\msl}{\overline{SL_2(\F)}}

\newcommand{\Ofnu}{\mathbb O_{\F_\nu}}

\newcommand{\fnu}{{\F_\nu}}
\newcommand{\tnu}{{ \tau_\nu}}
\newcommand{\pnu}{{ \pi_\nu}}
\newcommand{\glmnu}{{GL_m(\F_\nu)}}
\newcommand{\psinu}{{\psi_\nu}}

\newcommand{\mspknu}{\overline{Sp_{2k}(\F_\nu)}}
\newcommand{\mspnnu}{\overline{Sp_{2n}(\F_\nu)}}
\newcommand{\msignu}{{\overline{\sigma}_\nu}}
\newcommand{\gamapsinu}{\gamma_{\psi{\nu}}^{-1}}\def\dotunion{
\def\dotunionD{\bigcup\kern-9pt\cdot\kern5pt}
\def\dotunionT{\bigcup\kern-7.5pt\cdot\kern3.5pt}
\mathop{\mathchoice{\dotunionD}{\dotunionT}{}{}}} 
\begin {document}

%\author{Dani Szpruch}
%\title {The Langlands-Shahidi method for the metaplectic group, and some applications} \maketitle

Tel Aviv University School of Mathematical Sciences.

{\bf The Langlands-Shahidi Method for the metaplectic group and applications.}

Thesis submitted for the degree "Doctor of Philosophy" by Dani Szpruch.

Prepared under the supervision of Prof. David Soudry.

Submitted to the senate of Tel-Aviv University October 2009.

\newpage
Acknowledgements

I would like to take this opportunity to acknowledge some of the people whose presence helped me significantly through this research.

First, and foremost, I would like to thank my advisor Professor David Soudry. Literally, I could not have found a better teacher. I am forever indebted for his continuous support and for his endless patience. I consider myself extremely privileged for being his student.

Next, I would like to thank Professor Freydoon Shahidi for his kindness. Over several occasions Professor Shahidi explained to me some of the delicate aspects of his method. In fact, Professor Shahidi motivated some of the directions taken in this dissertation.

I would also like to thank Professor James Cogdell for his warm hospitality and for introducing me to the work of Christian A. Zorn which turned out to be quite helpful.

I would like to thank Professor William D. Banks for providing me his PhD thesis which contained important information.

I would like to thank my fellow participants in the representation theory seminar conducted by Professor Soudry, Professor Asher Ben-Artzi, Ron Erez, Eyal Kaplan and Yacov Tanay, for sharing with me many deep mathematical ideas. I would like to thank my office mate Dr. Lior Bary Sorocker for many hours of conversations.

Last, but definitely not least, I would like to thank my wife Norma Musih for her love and support.

\newpage \tableofcontents \newpage
%%%%%%%%%%%%%%%%%%%%%%%%%%%%%%%%%%%%%%%%%%%%%%%%%%%%%%%%%%%%%%%%%%%%%%%%%%%%
\section{Introduction}
Over a period of about thirty years Freydoon Shahidi has developed
the theory of local coefficients and its applications. Nowadays
this method is known as the Langlands-Shahidi method. The
references \cite{Sha78}, \cite{Sha80}, \cite{Sha 1}, \cite{Sha83},
\cite{Sha84}, \cite{Sha85}, \cite{Sha 2}, \cite{Sha 3}, \cite{Sha
4} are among Shahidi`s works from the first half of this period.
These works are used in this dissertation. The applications of
this theory are numerous; see the surveys \cite{GS},
\cite{ShaOxf}, \cite{Sha 96} and \cite{Kim} for a partial list.
Although this theory addresses quasi-split connected reductive
algebraic groups, our aim in this dissertation is to extend this
theory to  $\mspn$, the metaplectic double cover of the symplectic
group over a p-adic field $\F$, which is not an algebraic group,
and present some applications. We shall realize $\mspn$ as the set
$\spn \times \{ \pm 1 \}$ equipped with the multiplication law
$$(g_{1},\epsilon_{1})(g_{2},\epsilon_{2})= \bigl(
g_{1}g_{2},\epsilon_{1}\epsilon_{2}c(g_{1},g_{2})\bigr),$$ where
$c(\cdot,\cdot)$ is Rao's cocycle as presented in \cite{R}.

The properties of $\mspn$ enable the extension of the general
representation theory of quasi-split connected reductive algebraic
groups as presented in \cite{Sil book} and  \cite{Wa03}. A great
part of this extension is already available in the literature; see
\cite{KP}, \cite{F86}, \cite{BanPhD}, \cite{Ban} and \cite{Zo09}
for example.

An analog to Bruhat decomposition holds in $\mspn$. If $\F$ is a
p-adic field, $\mspn$ is an $l$-group in the sense of Bernstein
and  Zelevinsky, \cite{BZ}. Since Rao's cocycle is continuous, it
follows that there exists $U$, an open compact subgroup of $\spn$,
such that $c(U,U)=1$. Thus, a system of neighborhoods of
$(I_{2n},1)$ is given by open compact subgroups of the form
$(V,1)$, where $V \subseteq U$ is an open compact subgroup of
$\spn$. Furthermore, in $\mspn$ the analogs of the Cartan and
Iwasawa decompositions hold as well. If $\F$ is not 2-adic then
$\mspn$ splits over the standard maximal compact subgroup of
$\spn$. Over any local field (of characteristic different than 2)
$\mspn$ splits over the unipotent subgroups of $\spn$.

For a subset $H$ of $\spn$ we denote by $\overline{H}$ its
pre-image in $\mspn$. Let $P=M \ltimes N$ be a parabolic subgroups
of $\spn$. $\overline{P}$ has a "Levi" decomposition.
$\overline{P}=\overline{M} \ltimes \mu(N)$, where $\mu$ is an
embedding of $N$ in $\mspn$ which commutes with the projection
map.

During a course by David Soudry (2008-2009) dedicated to
Waldspurger's {\it La formule de Plancherel pour les groupes
$p$-adiques (d'apr\`{e}s Harish-Chandra)} which is a remake of
Harish-Chandra's theory as presented in \cite{Sil book}, the
author checked that the general theorems regarding Jacquet
modules, $L_2$-representations, matrix coefficients, intertwining
operators, Harish-Chandra's $c$-functions etc. extend to the
metaplectic group. Same holds for Harish-Chandra's completeness
theorem and the Knapp-Stein dimension theorem which follows from
this theorem. In fact many of the geometric proofs that are given
in \cite{BZ} and \cite{BZ77} apply word for word to the
metaplectic group.

We note that many of the properties mentioned in the last
paragraph are common to general $n$-fold covering groups of
classical groups. However, the following property is a special
feature of $\mspn$. Let $g,h \in \spn$. If $g$ and $h$ commute
then the pre-images in $\mspn$ also commute. In particular, the
inverse image of a commutative subgroup of $\spn$ is commutative.
This implies that the irreducible representations of
$\overline{\tspn}$, the inverse image of the maximal torus of
$\spn$, are one dimensional. As noted in \cite{BFH} this is the
reason that a Whittaker model for principal series representation
of $\mspn$ is unique. Furthermore, in Chapter \ref{Uniqueness of
Whittaker model} of this dissertation we prove the uniqueness of
Whittaker model for $\mspn$ in general. We emphasize that this
uniqueness does not hold for general covering groups; see
\cite{GHP} and \cite{BanPhD}. It is the uniqueness of Whittaker
model that enables a straight forward generalization of the
definition of the Langlands-Shahidi local coefficients to the
metaplectic group (see \cite{Bud} for an application of the theory
of local coefficients in the context of non-unique Whittaker
model).

\subsection{Main results}
Let $\F$ be a local field of characteristic 0. Let $\psi$ be a
non-trivial character of $\F$. We regard $\psi$ also as a genuine
non-degenerate character of the inverse image of the unipotent
radical of a symplectic group.

{\large \bf Theorem A.} Let $\pi$ be an irreducible admissible
representation of $\mspn$. Then, the dimension of the space of
$\psi$-Whittaker functionals on $\pi$ is at most 1; see
Chapter \ref{Uniqueness of Whittaker model}.

{\large \bf Theorem B.} Let $\tau$ be an irreducible admissible
generic representation of $\gln$. Let $P(\F)$ be the Siegel
parabolic subgroup of $Sp_{2n}(\F)$ and let $w$ be a particular
representative of the long Weyl element  of $Sp_{2n}(\F)$ modulo
the long Weyl element  of $P(F)$. Let
$C_{\psi}^{\mspn}\bigl(\overline{P(\F)},s,\tau, w \bigr)$ be the
metaplectic analog to the Langlands-Shahidi local coefficient.
Then, there exists an exponential function $c_\F(s)$ such that
$$ C_{\psi}^{\mspn}\bigl(\overline{P(\F)},s,\tau,
w \bigr)=c_\F(s)\frac
{\gamma(\tau,sym^2,2s,\psi)}{\gamma(\tau,s+\half,\psi)}.$$
Furthermore, $c_\F(s)=1$ provided that $\F$ is a p-adic field of
odd residual characteristic, $\tau$ is unramified and $\psi$ is
normalized. See Sections \ref{prin comp with soo} and \ref{true
for cusp}.

Let $\F$ be a p-adic field.

{\large \bf Theorem C.} We keep the notation and assumptions of
Theorem B. Let  $\beta(s,\tau)$ the meromorphic function defined
by the relation
$$A_{w^{-1}}A_{w}=
\beta(s,\tau)Id.$$ Here $A_w$ is the intertwining operator defined
on $\pi$, the representation of $\mspn$ (parabolically) induced
from $\tau$. Then, $\beta(s,\tau)$ has the same analytic
properties as the Plancherel measure attached to $\soo$ and
$\tau$. In particular, if we assume in addition that $\tau$ is
supercuspidal unitary then $\pi$ is irreducible if and only if
$\pi'$ is irreducible. Here $\pi'$ is the representation of $\soo$
(parabolically) induced from $\tau$. See Section \ref{irr max
sec}.

{\large \bf Theorem D.} Let $\pi$ be a principal series
representation of $\mspn$ induced from a unitary character. Then
$\pi$ is irreducible. See Section \ref{irr prin uni sec}.

{\large \bf Theorem E}. For $1 \leq i \leq r$ let $\tau_i$ be an
irreducible admissible supercuspidal unitary  representation of
$GL_{n_i}(\F)$ and let $\msig$ be an irreducible admissible
 supercuspidal generic genuine representation of
$\mspk$. Denote by $\pi$ the corresponding parabolic induction on
$\mspn$. Then, $\pi$ is reducible if and only if there exists $1
\leq i \leq r$ such that $\tau_i$ is self dual and
$$\gamma(\msig \times \tau_i,0, \psi) \gamma(\tau_i,sym^2,0,\psi)
\neq 0.$$ Here $\gamma(\msig \times \tau_i,s, \psi)$ is the
$\gamma$-factor attached to $\msig$ and $\tau_i$ defined by
analogy to the general definition of Shahidi; See section
\ref{General parabolic case}.

{\large \bf Theorem F.} Let $\chi$ be a character of $\F^*$. There
exists a meromorphic function $\widetilde{\gamma}(\chi,\psi,s)$
such that $$
\widetilde{\gamma}(\chi,\psi^{-1},s)\zeta(\phi,\chi,s)=\zeta(\widetilde{\phi},\chi^{-1},1-s)$$
for every  $\phi$, a Schwartz function on $\F$. Here $\zeta$ is
the Mellin Transform and $$\widetilde{\phi}(x)=\int_{\F}
\phi(y)\psi(xy)\gamma_{\psi}^{-1}(xy) \, dy,$$ where $\gamma_\psi$
is the normalized Weil factor attached to a character of second
degree. Furthermore, $\widetilde{\gamma}(\chi,\psi,s)$ has a
relation to the $\msl$ local coefficient which is similar to the
relation that the Tate $\gamma$-factor has with the $SL_2(\F)$
local coefficient, i.e.,
$$\widetilde{\gamma}^{-1}(\chi^{-1},\psi^{-1},1-s)=\epsilon'(\chi,s,\psi)
\frac{\gamma(\chi^2,2s,\psi)}{\gamma(\chi,s+\half,\psi)},$$ where
$\epsilon'(\chi,s,\psi)$ is an exponential factor which equals 1
if $\chi$ is unramified and $\F$ is p-adic field of odd residual
characteristic. See Section \ref{tate remark}.

{\large \bf Theorem G.} Let $\F$ be a field of characteristic
different then 2. The unique extension of Rao`s cocycle from
$\spn$ to $\gspn$ is given by
$$\widetilde{c}(g,h)=v_{\lambda_h}\bigl(p(g) \bigr )c \bigl
(p(g)^{\lambda_h},p(h) \bigr).$$For details see Section \ref{ext
rao}.

\subsection{An outline of
the thesis} Chapters \ref{genral not}-\ref{gen rep meta} are of a
preliminary nature. In Chapter \ref{genral not} we give the
general notation to be used throughout this dissertation. Among
these notations is $(\cdot,\cdot)_\F$, the quadratic Hilbert
symbol.

In Chapter \ref{The metaplectic group} we present the metaplectic
group and prove some of its properties. In Section \ref{The
symplectic group} we introduce some notations and facts related to
symplectic groups. In particular we denote by $\pt$ and $\mt$ a
standard parabolic subgroup and its Levi part.

 In Section \ref{des rao} we
introduce Rao`s cocycle and prove some useful properties. The
metaplectic group is also presented in this section along with its
basic properties. In Section \ref{kubota} we address the rank 1
case of Rao`s cocycle which is identical to Kubota`s cocycle; see
\cite{Kub}. For p-adic fields of odd residual characteristic we
recall the explicit splitting of $SL_2(\Of)$. We use this
splitting to describe some properties of the splitting of
$Sp_{2n}(\Of)$; see Lemma \ref{k split}. For all the p-adic fields
we prove that the cocycle is trivial on small enough open compact
subgroups of $SL_2(\F)$. We describe explicitly the two
isomorphisms between $\overline{SO_2(\R)}$ and $\R / 4\pi \Z$; see
Lemma \ref{so2}. One of the explicit isomorphisms will be used in
Section \ref{real case} where we compute the local coefficients of
$\overline{SL_2(\R)}$. In Section \ref{par mspn} we deal with
parabolic subgroups of $\mspn$ which are defined to be an inverse
image of parabolic subgroups of $\spn$. We prove that these groups
possess an exact analog to Levi decomposition. For p-adic fields
of odd residual characteristic we give in Lemma \ref{p is product
odd f} an explicit isomorphism  between $\mmt$ and $GL_{n_1}(\F)
\times GL_{n_2}(\F) \ldots \times GL_{n_r}(\F) \times \mspk.$ As a
by product of this isomorphism we obtain an explicit splinting of
the Siegel parabolic subgroup. The main ingredient of the proof of
Lemma \ref{p is product odd f} is the existence of functions $\xi:
\F^* \rightarrow \{ \pm 1\}$ with the property
$$\xi(ab)=\xi(a)\xi(b)(a,b)_\F.$$ See Lemma \ref{4split}. In
Section \ref{global} we introduce the global metaplectic group.
This group will appear again only in Chapter \ref{global sec}.

In Section \ref{ext rao} we describe explicitly the unique
extension of Rao`s cocycle to $\gspn$. This extension is
equivalent to a realization of the unique double cover of $\gspn$
which extends the unique non-trivial double cover of $\spn$. For
the rank 1 case, this theorem is proven in \cite{Kub}. The main
ingredient of this extension is a lift of an outer conjugation of
$\spn$ by an element of $\gspn$ to $\mspn$. A particular property
of this lifting (see Corollary \ref{why lift}) will play a crucial
role in Chapter \ref{Uniqueness of Whittaker model} where we prove
the uniqueness of Whittaker model.

 In Chapter \ref{Weil
factor} we introduce $\gamma_\psi:\F^* \rightarrow \{\pm 1, \pm
i\}$, the normalized Weil factor of a character of second degree;
see Theorem 2 of Section 14 of \cite{Weil}. Since the complex case
is trivial and since the real case is clear (see \cite{P}) we
mostly address the p-adic case. As may be expected, the main
difficulty lies in 2-adic fields. We give some formulas (see Lemma
\ref{gamma comp}) that we shall use in Section \ref{padic comp}
where we compute the local coefficients for $\msl$ over p-adic
fields. For the sake of completeness we explicitly compute all the
Weil factors for p-adic fields of odd residual characteristic and
for $\Q_2$. It turns out that if $\F$ is not 2-adic then
$\gamma_\psi$ is not onto $\{\pm 1, \pm i\}$. The Weyl index
defined on $\Q_2$ is onto $\{\pm 1, \pm i\}$. Furthermore, if $\F$
is a p-adic field of odd residual characteristic and if $-1 \in
{\F^*}^2$ then $\gamma_\psi$ equals to one of the $\xi$ functions
presented in Section \ref{par mspn}.

Let $\F$ be a p-adic field. In Chapter \ref{gen rep meta} we
survey some facts from the existing and expected representation
theory of $\mspn$ to be used in this dissertation. This theory is
a straight forward generalization of the theory for algebraic
groups. For $1 \leq i \leq r$ let $\tau_i$ be a smooth
representation of $GL_{n_i}(\F)$ and let $\msig$ be a genuine
smooth representation of $\mspk$. In Section \ref{gen par ind} we
construct a smooth genuine representation of $\mmt$ from these
representations. This is done, roughly speaking, by tensoring
$\tau_1 \otimes \tau_2 \ldots \otimes \tau_r \otimes \msig$ with
$\gamma_\psi$; see Lemma \ref{gen par ind}. Note that this is not
quite the process used for general covering groups; see
\cite{BanPhD} for example. Next we define parabolic induction in
an analogous way to the definition for algebraic groups. In
Section \ref{app bruhat theory} we give a rough condition for the
irreducibility of unitary parabolic induction that follows from
Bruhat theory. Namely, we explain which representations of $\mmt$
are regular; see Theorem \ref{basic gen par}. In Section \ref{sec
io def} we define the intertwining operator $A_{w}$ attached to a
Weyl element $w$. This operator is the meromorphic continuation of
a certain integral; see \eqref{io def}. Its definition and basic
properties are similar to the analogous intertwining operator for
algebraic groups. Chapter \ref{gen rep meta} culminates in Section
\ref{knapp sec} where we give the metaplectic analog to the
Knapp-Stein dimension theorem ;see \cite{Sil} for the p-adic case.
This Theorem describes the (dimension of) the commuting algebra of
a parabolic induction via the properties of the meromorphic
functions $\beta(\overrightarrow{s},\tau_1,\ldots,\tau_r,\msig,w)$
defined by the relation
$$A_{w^{-1}}A_{w}=
\beta^{-1}(\overrightarrow{s},\pi,w)Id.$$ In more details, let
$\pi$ be an irreducible admissible supercuspidal unitary genuine
representation of $\mpt$, let $W(\pi)$ be the subgroup of Weyl
elements which preserve $\pi$, let $\Sigma_{\pt}$ be the set of
reflections corresponding to the roots of $\tspn$ outside $\mt$
and let $W''(\pi)$ be the subgroup of $W(\pi)$ generated by  $w
\in \Sigma_{\pt} \cap W(\pi$) which satisfies
$\beta(\overrightarrow{s},\pi,w)=0$. denote by $I(\pi)$ the
representation of $\mspn$ induced from $\pi$. The Knapp-Stein
dimension theorem states that
$$Dim \bigl(Hom(I(\pi),I(\pi)) \bigr)=[W(\pi) :
W''(\pi)].$$

For connected reductive quasi split algebraic groups
$\beta(\overrightarrow{s},\pi,w)$ is closely related to the
Plancherel measure; see Section 3 of the survey \cite{ShaOxf} for
example. In Chapter \ref{Irreducibility theorems} we shall compute
this function in various cases.

Since $\mspn$ splits over unipotent subgroups of $\spn$ one can
define a genuine $\psi$-Whittaker functionals attached to a smooth
genuine representation $\pi$ of $\mspn$ in a closely related way
to the definition in the linear case. Denote by $W_{\pi,\psi}$ the
space of $\psi$-Whittaker functionals defined on $V_\pi$. In
Theorem \ref{uniquness} of Chapter \ref{Uniqueness of Whittaker
model} we prove that if $\pi$ be an irreducible, admissible
representation of $\mspn$. Then $\dim (W_{\pi,\psi}) \leq 1.$

We emphasize that uniqueness of Whittaker model is not a general
property of covering groups; see \cite{GHP} for the failure of
this uniqueness in the $GL_2(\F)$ case. One of the key reasons
that this uniqueness holds is that $\overline{\tspn}$ is
commutative. We first explain why in the archimedean case the
uniqueness proof is done exactly as in the linear case; see
\cite{H}. Then, in Section \ref{padic whi uni sec} we move to the
non-archimedean case. Our method of proof is a method similar
 to the one in the linear case; see \cite{Sh}, \cite{GK} and
\cite{BZ}. In particular, we adapt the Gelfand-Kazhdan method to
$\mspn$; see Theorem \ref{gk}. We use there $\overline{h} \mapsto
^{\overline{\tau}}\! \! {\overline{h}}$, an involution on $\mspn$
which is a lift of $h \mapsto ^{\tau}\! \! {h}$, the involution
used for the uniqueness proof in the symplectic case; see Lemma
\ref{simple lift}. It is somehow surprising to know that
$\overline{\tau}$ extends $\tau$ in the simplest possible way,
i.e., if $ \overline{h}=(h,\epsilon)$ then $^{\overline{\tau}}\!
\! {\overline{h}}=(^{\tau}\! \! {h},\epsilon)$. As mentioned
before, the explicit computation and crucial properties of
$\overline{\tau}$ follow from the results proven in Section
\ref{ext rao}. The main technical ingredient used for the
uniqueness proof is Theorem \ref{unproven}. As indicated in the
proof itself, provided that the relevant properties of
 $\overline{\tau}$ are proven, Theorem \ref{unproven} is proved exactly as its linear analog which is Lemma \ref{unproven 2}.

Once the uniqueness of Whittaker model is established we define in
Chapter \ref{chapter def lc} the metaplectic Langlands-Shahidi
local coefficient
$$C_{\psi}^{\mspn}(\mpt,\overrightarrow{s},(\otimes_{i=1}^r
\tau_i)\otimes \overline{\sigma},w \bigl)$$ in exactly the same
way as in the linear case; see Theorem 3.1 of \cite{Sha 1}. We
note that the zeros of the local coefficient are among the poles
$A_w$. Furthermore, since by definition
\begin{eqnarray*} && \beta(\overrightarrow{s},\tau_1,\ldots,\tau_n,\msig,w)
\\ &&= C_{\psi}^{\mspn}(\mpt,\overrightarrow{s},(\otimes_{i=1}^r
\tau_i)\otimes \overline{\sigma},w \bigl)
C_{\psi}^{\mspn}(\mpt,\overrightarrow{s^{w}},((\otimes_{i=1}^r
\tau_i)\otimes \overline{\sigma})^{w},w^{-1} \bigl)
\end{eqnarray*} the
importance of the local coefficients for questions of
irreducibility of parabolic induction is clear. Using the local
coefficients we also define $$\gamma(\msig \times \tau,s,
\psi)=\frac {C_{\psi}^{\mspn}
\bigl(\overline{P_{m;k}(\F)},s,\tau\otimes \msig,
j_{m,n}(\omega_m'^{-1}) \bigr)} {C_{\psi}^{\mspm}
\bigl(\overline{P_{m;0}(\F)},s,\tau, \omega_m'^{-1}  \bigr )},$$
the $\gamma$-factor attached to $\overline{\sigma}$, an
irreducible admissible genuine $\psi$-generic representation of
$\mspk$, and $\tau$, an irreducible admissible generic
representation of $\glm$; see \eqref{gama def}. Here $\omega_m'$
and $j_{m,n}(\omega_m'^{-1})$ are appropriate Weyl elements. This
definition of the $\gamma$-factor is an exact analog to the
definition given in Section 6 of \cite{Sha 3} for quasi-split
connected reductive algebraic groups.

Most of Chapter \ref{basic gam chapter} is devoted to the proof of
the multiplicativity properties of $\gamma(\msig \times
\tau,s,\psi).$ See Theorems \ref{mult gama tau} and \ref{mult gama
sigma} of Section \ref{cha gam}. It is the metaplectic analog to
Part 3 of Theorem 3.15 of \cite {Sha 3}. This multiplicativity is
due to the multiplicativity of the local coefficients. The main
ingredient of the proof of the multiplicativity is a certain
decomposition of the intertwining operators; see Lemma \ref{lem
heart}. This decomposition resembles the decomposition of the
intertwining operators in the linear case. The only (small)
difference is that two Weyl elements may carry cocycle relations.
Our choice of Weyl elements is such that this relations are
non-trivial only in the field of real numbers and in 2-adic
fields. In Lemma \ref{princcomp} of Section \ref{prin gama sec} we
compute $\gamma(\msig \times \tau,s,\psi)$ for principal series
representations. Let
$\eta_1,\eta_2,\ldots,\eta_k,\alpha_1,\alpha_2,\ldots,\alpha_m$ be
$n$ characters of $\F^*$. If $\tau$ is induced from
$\alpha_1,\alpha_2,\ldots,\alpha_m$ and $\msig$ is induced from
$\eta_1,\eta_2,\ldots,\eta_k$ (twisted by $\gamma_\psi$) then
there exists $c \in \{ \pm 1\}$ such that
$$\gamma(\overline{\sigma} \times \tau,s,\psi)=
c\prod_{i=1}^k\prod_{j=1}^m \gamma(\alpha_j \times
\eta_i^{-1},s,\psi)\gamma(\eta_i \times \alpha_j,s,\psi).$$ If
$\F$ is a p-adic field of odd residual characteristic, and $\tau$
and $\msig$ are unramified then $c=1$. This computation is an
immediate corollary of Theorems \ref{mult gama tau} and \ref{mult
gama sigma}.

Assume that $\F$ is either $\C$, $\R$ or a p-adic field. In
Sections \ref{padic comp}, \ref{real case} and \ref{easy complex
sec} of Chapter \ref{second paper} we compute the local
coefficients for principal series representations of $\msl$. In
this case $\chi$, the inducing representation is a character of
$\F^*$ and there is only one non-trivial Weyl element. We prove
that

$$C_\psi^{\overline{SL_2(\F)}}\Bigl(\overline{B_{SL_2(\F)}},s,\chi,(\begin{smallmatrix}
& {0}& {-1} \\
& {1} & {0}\end{smallmatrix})\Bigr)=\epsilon'(\chi,s,\psi)
\frac{\gamma(\chi^2,2s,\psi)}{\gamma(\chi,s+\half,\psi)},$$ where
$\epsilon'(\chi,s,\psi)$ is an exponential factor which equals 1
if $\chi$ is unramified, $\psi$ is normalized and $\F$ is p-adic
field of odd residual characteristic. $\epsilon'(\chi,s,\psi)$ is
computed explicitly for p-adic fields and for the field of real
numbers; see Theorems \ref{the formula}, \ref{real thm} and lemma
\ref{easy complex lem}. In this chapter only we write
$C_{\psi_a}(\etac \otimes \gamma_\psi^{-1},s)$ instead of
$C_{\psi_a}^{\overline{SL_2(\F)}}\Bigl(\overline{B_{SL_2(\F)}},s,\chi,(\begin{smallmatrix}
& {0}& {-1} \\
& {1} & {0}\end{smallmatrix})\Bigr),$ where $\psi_a(x)=\psi(ax)$.
This notation emphasizes the dependence of the local coefficient
on two additive characters, rather then on one in the algebraic
case; one is the Whittaker character and the second is the
character defining $\gamma_\psi$.

Section \ref{padic comp} is devoted to the p-adic case. In Section
\ref{red int} we express the local coefficient as the following
"Tate type" integral

$${C_{\psi_a}(\etac \otimes \gamma_\psi^{-1},s)}^{-1}= \int_{\F^*}
\gamma_\psi^{-1}(u) \chi(u) \psi_a(u)\ab u\ab^s d^*u.$$ See Lemma
\ref{reduction}. An immediate corollary of Lemma \ref{reduction}
is Lemma \ref{a to 1} which asserts that
$$
C_{\psi}(\etac \otimes \gamma_\psi^{-1},s)\gamma_\psi(a)\chi(a)\ab
a \ab^s =C_{\psi_a}(\ \chi \cdot (a,\cdot) \otimes
\gamma_\psi^{-1},s).
$$
In light of the computation of $C_{\psi}(\etac \otimes
\gamma_\psi^{-1},s)$ this means that the analytic properties of
$C_{\psi_a}(\etac \otimes \gamma_\psi^{-1},s)$ depend on $a$ if
$\chi^2$ is unramified. This phenomenon has no analog in the
linear case.

After proving some technical lemmas in Section \ref{lemmas} we
compute the "Tate type" integral in Section \ref{int comp}. The
difficult computations appear in the 2-adic case. The parameter
that differentiates between 2-adic fields and p-adic fields of odd
residual characteristic is $e(2,\F)$ which is defined to be the
ramification index of $\F$ over $\Q_2$ if $\F$ is a 2-adic field
and 0 if $\F$ is of odd residual characteristic.

In Section \ref{real case} we compute the local coefficients for
$\overline{SL_2(\R)}$. Since $\overline{SO_2(\R)}$ is commutative
it follows from the Iwaswa decomposition that the
$\overline{SO_2(\R)}$-types are one dimensional. This means that
the intertwining operator maps a function of $\vartheta_n$ type to
a multiple of a function of $\vartheta_n$ type. Thus, we may
compute the local coefficients in a way similar to the computation
of the local coefficients in the $SL_2(\R)$ case. In fact, from a
certain point we use Jacquet computations (see \cite{J}) which
ultimately use the results of Whittaker himself; see \cite{Whi}.

Section \ref{easy complex sec} which follows next is short and
devoted to the complex case. Since $\gamma_\psi(\C^*)=1$ and since
$\overline{SL_2(\C)}=SL_2(\C)\times \{\pm 1 \}$
 the local coefficients of this group are
identical to the local coefficients of $SL_2(\C)$. The $SL_2(\C)$
computation is given in Theorem 3.13 of \cite{Sha85}. The
surprising fact is that the $SL_2(\C)$ computation agrees with the
$\overline{SL_2(\F)}$ computations where $\F$ is either $\R$ or a
p-adic field. This follows from the duplication formula of the
classical $\Gamma$- function.

Section \ref{par remark} is a detailed remark of our choice of
parameterization, i.e., our choice to compute $C_{\psi}(\etac
\otimes \gamma_\psi^{-1},s)$ rather then $C_{\psi_a}(\etac \otimes
\gamma_\psi^{-1},s)$ for some other $a \in \F^*$ (for instance, in
\cite{BFH} $a=-1$ is used). We explain in what sense our
parameterization is the correct one.

Let $\F$ be a p-adic field. In Section \ref{tate remark} we show
that one can define a meromorphic function
$\widetilde{\gamma}(\chi,\psi,s)$ by a similar method to the one
used by Tate (see \cite{T}) replacing to role of the Fourier
transform with $\phi \mapsto \widetilde{\phi}$ defined on the
space of Schwartz functions by
$$\widetilde{\phi}(x)=\int_{\F}
\phi(y)\psi(xy)\gamma_{\psi}^{-1}(xy) \, dy.$$ We Show that
$$C_\psi^{\overline{SL_2(\F)}}\bigl(\overline{B_{SL_2(\F)}},s,\chi,(\begin{smallmatrix}
& {0}& {-1} \\
& {1} &
{0}\end{smallmatrix})\bigr)^{-1}=\widetilde{\gamma}(\chi^{-1},\psi^{-1},1-s).$$
This is an analog to the well known equality
$${C_{\psi}^{SL_2(\F)}\bigl(B_{SL_2}(\F),s,\chi,(\begin{smallmatrix}
& {0}& {-1} \\
& {1} &
{0}\end{smallmatrix})\bigr)}^{-1}=\gamma(\chi^{-1},\psi^{-1},1-s).$$

We conclude this chapter in Section \ref{prin comp with soo} where
we prove that if $\tau$ is a principal series representation of
$\glm$ then there exists an exponential function $c(s)$ such that

\begin{equation} \label{intro prin}
C_{\psi}^{\mspm}\bigl(\overline{P_{m;0}(\F)},s,\tau,
\omega_m'^{-1} \bigr)=c(s)\frac
{\gamma(\tau,sym^2,2s,\psi)}{\gamma(\tau,s+\half,\psi)}.
\end{equation} If $\F$ is a p-adic field of odd residual
characteristic, $\psi$ is normalized and $\tau$ is unramified then
$c(s)=1$ . This computation uses the $\msl$ computations given in
sections \ref{padic comp}, \ref{real case} and \ref{easy complex
sec} and the multiplicativity of the local coefficients proven is
Section \ref{cha gam}; see Theorem \ref{soudry said it}.

In Chapter \ref{global sec} we prove a global functional equation
satisfied by $\gamma(\msig \times \tau,s,\psi)$. See Theorem
\ref{crude} in Section \ref{crude equ}. This theorem is the
metaplectic analog to the crude functional equation proven by
Shahidi in  Theorem 4.1 of \cite{Sha 1} (see also Part 4 of
Theorem 3.15 of \cite{Sha 3}). The proof Theorem \ref{crude}
requires some local computations of spherical Whittaker functions;
see Section \ref{anunramcomp}. These computations use the work of
Bump, Friedberg and Hoffstein; see \cite{BFH}. As a corollary of
the Crude functional equation and the general results presented in
\cite{Sha 2} we prove in Section \ref{true for cusp} that if
$\tau$ is an irreducible  admissible generic representation of
$\glm$ then \eqref{intro prin} holds.

In Chapter \ref{Irreducibility theorems} we use the results of
Chapters \ref{basic gam chapter}-\ref{global sec} combined with
the Knapp-Stein dimension Theorem presented in Chapter \ref{gen
rep meta} to prove certain irreducibility theorems for parabolic
induction on the metaplectic group over p-adic fields. We assume
that the inducing representations are unitary. In Theorem \ref{irr
unitary prinipal} of Section \ref{irr prin uni sec} we prove via
the computation of the local coefficients the irreducibility of
principal series representations of $\mspn$ induced from unitary
characters. The method used in Section \ref{irr prin uni sec} is
generalized in Section \ref{General parabolic case} where we prove
the following.

Let $n_1,n_2,\ldots,n_r,k$ be $r+1$ non-negative integers whose
sum is $n$. For $1 \leq i \leq r$ let $\tau_i$ be an irreducible
admissible supercuspidal unitary representation of $GL_{n_i}(\F)$
and let $\msig$ be  an irreducible admissible supercuspidal
genuine $\psi$-generic representation of $\mspk$. Denote
$\pi=\bigl(\otimes_{i=1}^r (\gamma_{\psi}^{-1}\otimes
{\tau_i})\bigr) \otimes \msig$ and $I(\pi)=Ind_{\mpt}^{\mspn}
\pi$. Then $I(\pi)$ is reducible if and only if there exists $1
\leq i \leq r$ such that $\tau_i$ is self dual and
$$\gamma(\msig \times \tau_i,0,
\psi) \gamma(\tau_i,sym^2,0,\psi) \neq 0.$$ Denote now
$\pi_i=(\gamma_{\psi}^{-1}\otimes {\tau_i}\bigr) \otimes \msig$
and
$I(\pi_i)=Ind_{\overline{P_{n_i;k}(\F)}}^{\overline{Sp_{2(n_i+k)}(\F)}}
\pi_i$. An immediate corollary of the last theorem is that
$I(\pi)$ is irreducible if and only if $I(\pi_i)$ is irreducible
for each $1 \leq i \leq r$.

Let $\soo$ be a split odd orthogonal group and let
$P_{SO_{2n+1}}(\F)$ be a parabolic subgroup of $\soo$ whose Levi
part is isomorphic to $\gln$. In Lemma \ref{lc soo lc meta} of
Section \ref{irr max sec} we prove that if $\tau$ is an
irreducible admissible generic representation of $\gln$ then
$\beta(s,\tau,\omega_n'^{-1})$ has the same analytic properties as
the Plancherel measure attached to $\soo$, $P_{SO_{2n+1}}(\F)$ and
$\tau$. An immediate Corollary is Theorem \ref{metaplectic and so}
which states the following.

Let $\tau$ be an irreducible admissible self dual supercuspidal
representation of $\gln$. Then,
$Ind^{\mspn}_{\overline{P_{n;0}(\F)}} (\gamma_{\psi}^{-1} \otimes
\tau )$ is irreducible if and only if
$Ind^{\soo}_{P_{SO_{2n+1}}(\F)} \tau $ is irreducible.

We then list a few corollaries that follow from Theorem
\ref{metaplectic and so} and from Shahidi`s work, \cite{Sha 4}.

Some of the results presented in this dissertation have already
been published. The main results of Section \ref{ext rao} and
Chapter \ref{Uniqueness of Whittaker model} appeared in \cite{Sz}.
The main results of Sections \ref{padic comp} and \ref{real case}
were published in \cite{Sz09}. Theorems \ref{irr unitary prinipal}
and \ref{metaplectic and so} along with their proofs were outlined
in the introduction of \cite{Sz09}. The main results of Chapters
\ref{basic gam chapter}, \ref{global sec} and \ref{Irreducibility
theorems} were presented in Seminar talks in Tel-Aviv University
(September 2008), Ohio-State University (October 2008) and Purdue
University (October 2008).
\newpage

\section{General notations} \label{genral not}Assume that $\F$ either the field of
real numbers or a finite extensions of $\Q _p$. $\F$ is self dual:
If $\psi$ is a non-trivial (complex) character of $\F$, any
non-trivial character of $\F$ has the form $\psi_a(x)=\psi(ax)$
for a unique $a\in \F^*$. If $\F=\R$ we define $\psi(x)=e^{ix}$.

Let $\F$ be a p-adic field. Let $\Of$ be the ring of integers of
$\F$ and let $\Pf$ be its maximal ideal. Let $q$ be the
cardinality of the residue field $\overline{\F}={\Of}/{\Pf}$. For
$b \in \Of$ denote by $\overline{b}$ its image in $\rf$. Fix
$\pi$, a generator of $\Pf$. Let $\ab \cdot \ab$ be the absolute
value on $\F$ normalized in the usual way: $\ab \pi \ab=q^{-1}$.
Normalize the Haar measure on $\F$ such that $\mu(\Of)=1$ and
normalize the Haar measure on $\F^*$ such that $d^*x=\frac{dx}{\ab
x \ab}$. Define $$e(2,\F)=log_q[\Of:2\Of]=-log_q \ab 2 \ab.$$ Note
that if the residual characteristic of $\F$ is odd then
$e(2,\F)=0$, while if the characteristic of $\rf$ is 2, $e(2,\F)$
is the ramification index of $\F$ over $\Q_2$. We shall write $e$
instead of $e(2,\F)$, suppressing its dependence on $\F$. We also
define
$$\omega=2\pi^{-e} \in \Of^*.$$
We shall often be interested in subgroups of $\Of^*$ of the form
$1+\Pf^n$, where $n$ is a positive integer. By abuse of notation
we write $1+\Pf^0=\Of^*$. For $\chi$, a character of $\F^*$, we
denote by $m(\chi)$ its conductor, i.e, the minimal integer such
that $\chi(1+\Pf^{m(\chi)})=1$. Note that if $\chi$ is unramified
then $m(\chi)=0$. Let $\psi$ be a non-trivial additive character
of $\F$. The conductor of $\psi$ is the minimal integer $n$ such
that $\psi(\Pf^n)=1$. $\psi$ is said to be normalized if its
conductor is 0. Assuming that $\psi$ is normalized, the conductor
of $\psi_a$ is $log_q \ab a ^{-1} \ab$. If the conductor of $\psi$
is $n$, we define $\psi_0(x)=\psi(x\pi^n)$. $\psi_0$ is clearly
normalized.

For any field (of characteristic different then 2) we define
$(\cdot,\cdot)_\F$ to be the quadratic Hilbert symbol of $\F$. The
Hilbert symbol defines a non-degenerate bilinear form on $\F^* /
{\F^*}^2$. For future references we recall some of the properties
of the Hilbert symbol:
\begin{equation} \label{Hilbert properties} 1.\, \, (a,-a)_\F=1 \ \ \ \ \
2.\, \, (aa',b)_\F=(a,b)_\F(a',b)_\F \ \ \ \ \ 3.\,
\,(a,b)_\F=(a,-ab)_\F.
\end{equation}

\newpage
\section{The metaplectic group} \label{The metaplectic group}
Let $\F$ be a local field of characteristic 0. In this chapter we
introduce $\mspn$, the metaplectic group which is the unique
non-trivial double cover of the symplectic group and describe its
basic properties. Through this dissertation, we realize the
metaplectic group via Rao`s cocycle, \cite{R}. In Section \ref{The
symplectic group} we list several notations and standard facts
about the symplectic group and its subgroups. In Section \ref{des
rao} we introduce Rao`s cocycle and list serval of its properties.
In Section \ref{kubota} we describe Kubota`s cocycle which is
Rao`s cocycle for the $SL_2(\F)$ case. Some of the properties of
$\mspn$ may be proven using the properties of $\msl$. Next, in
Section \ref {par mspn}, we describe the properties of parabolic
subgroups of $\mspn$, which are defined to be the inverse images
in $\mspn$ of parabolic subgroups of $\spn$. We show that these
groups have an exact analog to the Levi decomposition. For p-adic
fields of odd residual characteristic we give an explicit
interesting isomorphism between a Levi subgroup and and
$$GL_{n_1}(\F) \times GL_{n_2}(\F) \ldots \times GL_{n_r}(\F) \times \mspk.$$
We do so by proving the existence of functions $\xi: \F^*
\rightarrow \{ \pm 1\}$ with the property
$\xi(ab)=\xi(a)\xi(b)(a,b)_\F.$ In Section \ref{global} we define
the adelic metaplectic group which is the unique non-trivial
double cover of the adelic symplectic group. We conclude this
chapter in Section \ref{ext rao} where we explicitly compute the
unique extension of Rao`s cocycle to $\gspn$. We use this
extension to give a realization of the unique double cover of
$\gspn$ which extends the unique non-trivial double cover of
$\spn$. In Chapter \ref{Uniqueness of Whittaker model} we use this
extension to prove certain properties of an involution on $\mspn$
which plays a crucial role in the proof of the uniqueness of
Whittaker model for irreducible admissible representations of
$\mspn$; see Lemma \ref{simple lift}.
\subsection{The symplectic group} \label{The symplectic group}
Let $\F$ be a field of characteristic different then 2. Let
$X=\F^{2n}$ be a vector space of even dimension over $\F$ equipped
with $<\cdot,\cdot>:X\times X \rightarrow \F$, a non degenerate
symplectic form and let $Sp(X)=\spn$ be the subgroup of $GL(X)$ of
isomorphisms of $X$ onto itself which preserve $<\cdot,\cdot>$.
Following Rao, \cite{R}, we shall write the action of $GL(X)$ on
$X$ from the right. Let
$$E=\{e_{1},e_{2},\ldots,e_{n},e^{*}_{1},e^{*}_{2},\ldots,e^{*}_{n}\}$$
be a symplectic basis of $X$; for $1\leq i,j \leq n$ we have
$<e_{i},e_{j}>=<e^{*}_{i},e^{*}_{j}>=0$ and
$<e_{i},e^{*}_{j}>=\delta_{i,j}.$ In this base $\spx$ is realized
as the group $$\{a \in GL_{2n}(\F) \mid aJ_{2n}a^t=J_{2n} \},$$ where
$J_{2n}=\begin{pmatrix} _{0} & _{I_{n}}\\_{-I_{n} } & _{0},
\end{pmatrix}$. Through this dissertation we shall identify $\spn$
with this realization. For $0 \leq r \leq n$ define $i_{r,n}$ to
be an embedding of $Sp_{2r}(\F)$ in $Sp_{2n}(\F)$ by $$
\begin{pmatrix} _{a} & _{b}\\_{c} & _{d}
\end{pmatrix}\mapsto \begin{pmatrix} _{I_{n-r}} & _{ } & _{ }
\\ _{ } & _{a} & _{ } & _{b}\\ _{ } & _{ }  & _{I_{n-r}} & _{ } \\
_{ } & _{c}  & _{ } & _{d}
\end{pmatrix},$$ where $a,b,c,d \in Mat_{r \times r}(\F)$, and define $j_{r,n}$ to be an embedding of
$Sp_{2r}(\F)$ in $Sp_{2n}(\F)$ by $$ \begin{pmatrix} _{a} &
_{b}\\_{c} & _{d} \end{pmatrix}\mapsto \begin{pmatrix} _{a} & _{ }
& _{b} & _{ } & _{ } \\ _{ } & _{I_{n-r}}  & _{ } & _{ }
\\ _{c} & _{ } & _{d} & _{ } \\
_{ } & _{ }  & _{ } & _{I_{n-r}}
\end{pmatrix}.$$

Let $\tgln$ be the subgroup of diagonal elements of $\gln$, let
$\zgln$ be the group of upper triangular unipotent matrices in
$\gl$ and let $\bgln=\tgln \ltimes \zgln$ be the standard Borel
subgroup of $\gln$. Let $\tspn$ be the subgroup of diagonal
elements of $\spn$ and let $\zspn$ be the following maximal
unipotent subgroup of $\spn$;
$$\Bigl\{\begin{pmatrix} _{z} & _{b}\\_{0} & _{\widetilde{z}}
\end{pmatrix}  \mid z\in Z_{\gl}, b \in Mat_{n \times n}
(\F), \, b^t=z^{-1}b z^t \Bigr \},$$ where for $a \in GL_n$ we
define $\widetilde{a}={^t \!{a} \!^ {-1}}$. The subgroup
$\bspn=\tspn \ltimes \zspn$ of $\spn$ is a Borel subgroup. We call
it the standard Borel subgroup. A standard parabolic subgroup of
$\spn$ is defined to be a parabolic subgroup which contains
$\bspn$. A standard Levi subgroup (unipotent radical) is a Levi
part (unipotent radical) of a standard parabolic subgroup. In
particular a standard Levi subgroup contains $\tspn$ and a
standard unipotent radical is contained in $\zspn$.

Let $n_1,n_2,\ldots,n_r,k$ be $r+1$ nonnegative integers whose sum
is $n$. Put $\overrightarrow{t}=(n_1,n_2,\ldots,n_r;k)$. Let
$M_{\overrightarrow{t}}$ be the standard Levi subgroup of $\spn$
which consists of elements of the form
$$[g_1,g_2,\ldots,g_r,h]=
diag(g_1,g_2,\ldots,g_r,I_k,\widetilde{g_1},\widetilde{g_2},\ldots,\widetilde{g_r},I_k)i_{k,n}(h),$$
where $g_i \in GL_{n_i}(\F), h \in \spk$. When convenient we shall
identify $GL_{n_i}(\F)$ with its natural embedding in $\mt$.
Denote by $\pt$ the standard parabolic subgroup of $\spn$ that
contains $\mt$ as its Levi part. Denote by $\nt$ the unipotent
radical of $\pt$. We denote by $\pspn$ or simply by $P(\F)$, the
Siegel parabolic subgroup of $\spn$:
\begin{equation}\label{parabolic} P_{(n;0)}(\F)= \Bigl\{
\begin{pmatrix} {a} & {b}\\{0} & {\widetilde{a}}
\end{pmatrix} \mid a \in \gl, b \in Mat_{n \times n} (\F),
b^{t}=a^{-1}ba^t \Bigr\}.\end{equation} Note that $M_{(n;0)}(\F)
\simeq \gln$. A natural isomorphism is given by $$g \mapsto
\widehat{g}=
\begin{pmatrix} _{g} & _{ }
\\  & _{ } & _{\widetilde{g}}  \end{pmatrix}.$$
Define $$V=span \{e_{1},e_{2},\ldots,e_{n}\}, \ \ V^{*}=span
\{e^{*}_{1},e^{*}_{2},\ldots,e^{*}_{n}\}.$$ These are two
transversal Lagrangian subspaces of $X$. The Siegel parabolic
subgroup is the subgroup of $\spx$ which consists of elements that
preserve $V^*$. Let $S$ be a subset of $\{1,2,\ldots,n \}$. Define
$\tau_S , a_S$ to be the following elements of $\spx$:
\begin{equation} \label{tau s} e_i \cdot \tau_S=\begin{cases} -e_i^* & i
\in S \\ e_i & $otherwise$ \end{cases},\ \ e_i^*  \cdot \tau_S=
\begin{cases} e_i & i \in S \\ e_i^* & $otherwise$
\end{cases},
\end{equation}
\begin{equation} \label{a s} e_i \cdot a_S=\begin{cases} -e_i & i
\in S \\ e_i & $otherwise$ \end{cases},\ \ e_i^*  \cdot a_S=
\begin{cases} -e_i^* & i \in S \\ e_i^* & $otherwise$
\end{cases}.\end{equation}
The elements $\tau_{S_1}, \,  a_{S_1}, \,  \tau_{S_2}, \, a_{S_2}$
commute. Note that $a_S \in P(\F) , \, a_S^2=I_{2n},$ and that
\begin{equation} \label{t1 mult t2} \tau_{S_1} \tau_{S_2}=\tau_{S_1 \triangle S_2}a_{S_1 \cap
S_2}, \end{equation} where $S_1
 \triangle S_2= S_1 \cup S_2 \backslash S_1 \cap S_2$. In
particular $\tau_S^{2}=a_S$. For $S=\{1,2,\ldots,n\}$ we define
$\tau=\tau_S$, in this case $a_S=-I_{2n}$.

Denote by $\wspntag$ the subgroup of  $Sp_{2n}(\Z)$ generated by
the elements $\tau_S$, and $\widehat{w_\pi}$, where $S \subseteq
\{1,2,\ldots,n \}$, and $w_\pi \in \gl$ is defined by ${w_\pi
}_{i,j}=\delta_{\pi(i),j}$; $\pi$ is a permutation in  $S_n$. If
$\F$ is a p-adic field then $\wspntag$ is a subgroup of $\Kn$.
Note that $\wspntag$ modulo its diagonal elements may be
identified with the Weyl group of $\spn$ denoted by $\wspn$.
Define $W_\pt$ to be the subgroup of $\wspntag$ which consists of
elements $w$ such that
$$\mt ^{w}=w \mt w^{-1}$$ is a standard Levi subgroup and
$$w \bigl(\zspn \cap \mt \bigr)w^{-1} \subset \zspn.$$ This means that up to
conjugation by diagonal elements inside the blocks of $\mt$, we
have
\begin{equation} \label{w action on element}
w[g_1,g_2,\ldots,g_r,s]w^{-1}=
[g^{(\epsilon_1)}_{\pi(1)},g^{(\epsilon_2)}_{\pi(2)},\ldots,g^{(\epsilon_r)}_{\pi(r)},s],\end{equation}
where $\pi$ is a permutation of $\{1,2,\ldots,r \}$, and where for
$g \in \gl$, $\epsilon= \pm 1$ we define
$$g^{(\epsilon)}=\begin{cases} g & \epsilon=1 \\ \omega_n
\widetilde{g} \omega_n & \epsilon=-1
\end{cases},$$
where $$\omega_n=\begin{pmatrix}
_{ } & _{ } & _{ } & _{1} \\
_{ } & _{ } & _{1} & _{ } \\
_{ } & _{\upddots} & _{ } & _{ } \\
_{1} & _{ } & _{ } & _{ }
\end{pmatrix}.$$
We may assume, and in fact do, that $W_\pt$ commutes with $i_{k,n}(\spk).$

For $w \in W_\pt$, let $\pt^{w}$ be the standard parabolic
subgroup whose Levi part is $\mt^{w}$, and let $\nt^{w}$ be its
standard unipotent radical.

\subsection{Rao's cocycle} \label{des rao} In \cite{R},
Rao constructs an explicit non-trivial 2-cocycle $c(\cdot ,\cdot)$
on $Sp(X)$ which takes values in $\{\pm 1 \}$. The set
$\overline{\spx}=\spx \times \{ \pm 1 \}$ is then given a group
structure via the formula
\begin{equation} \label{metaplectic structure}(g_{1},\epsilon_{1})(g_{2},\epsilon_{2})=
\bigl( g_{1}g_{2},\epsilon_{1}\epsilon_{2}c(g_{1},g_{2})\bigr)
.\end{equation} It is called the metaplectic group. For any subset
$A$ of $\spn$ we denote by $\overline{A}$ its inverse image in
$\mspn$. If $\F$ is either $\R$ or a p-adic field, this group is
the unique non-trivial double cover of $\spn$. It is well known
that other classical groups over local fields have more then one
non-trivial double cover. For example, on $\gspn$ one may define a
double cover by the multiplication law
$$(g,\epsilon_1)(h,\epsilon_2)=\bigl(gh,\epsilon_1
\epsilon_2(\lambda_g,\lambda_h)_\F \bigr),$$ where $\lambda_g$ is
the similitude factor of $g$. One may also define a double cover
via the construction given in Section \ref{ext rao}. The
corresponding cocycles are clearly inequivalent. The first is
trivial on $\spn \times \spn$ while the second is an extension of
Rao`s cocycle.

We now describe Rao's cocycle. Detailed proofs can be found in
\cite{R}. Define $$\Omega_j=\{\sigma \in \spx \mid \dim (V^* \cap
V^* \sigma)=n-j\}.$$ Note that $P(\F)=\Omega_0$, $\tau_S \in
\Omega_{\mid \! S \! \mid}$ and more generally, if
$\alpha,\beta,\gamma,\delta \in Mat_{n \times n}(\F)$ and
$\sigma=\begin{pmatrix} _{\alpha} & _{\beta}\\_{\gamma } &
_{\delta}
\end{pmatrix} \in \spx$ then $\sigma \in \Omega_{rank (\gamma)}$.
The Bruhat decomposition states that each $\Omega_j$ is a single
double coset in $_{P(\F)} \diagdown ^\spx \diagup_{P(\F)}$, that
$\Omega_j^{-1}= \Omega_j$ and that $\bigcup _{j=0} ^n
\Omega_j=\spx$. In particular every element of $\spx$ has the form
$p \tau_S p'$, where $p,p' \in P(\F),\, \, S \subseteq
\{1,2,\ldots,n\}$.

Let $p_1, p_2 \in P(\F)$. Rao defines \begin{equation} \label{x}
x(p_1 \tau_S p_2) \equiv \det (p_1p_ {2} \mid_{V^*}) \bigl(mod
(\F^*)^2 \bigr ),
\end{equation} and proves that it is a well defined map from  $\spx$ to
$\F^*/(\F^*)^2$. Note that $x(a_S)\equiv(-1)^{|S|}$. More
generally; if $p=\begin{pmatrix} {a} & {b}\\{0} & {\widetilde{a}}
\end{pmatrix} \in P(\F)$ then $x(p)\equiv\det (a)$. We shall use the notation $\det(p)=\det (a)$.
Also note that $x(\tau_S)\equiv1$ and that for $g\in \Omega_{j}
,\, p_1, p_2 \in P(\F)$,
\begin{equation} \label{x-prop} x(g^{-1})\equiv x(g)(-1)^{j}, \ \
x(p_1gp_2)\equiv x(p_1)x(g)x(p_2).\end{equation} Theorem 5.3 in
Rao's paper states that a non-trivial 2-cocycle on $\spx$ can be
defined by
\begin{equation} \label{cocycle} c(\sigma_1,\sigma_2)=
\bigl(x(\sigma_1),x(\sigma_2)\bigr)_\F
\bigl(-x(\sigma_1)x(\sigma_2),x(\sigma_1\sigma_2)\bigr)_\F
\bigl((-1)^l,d_F(\rho)\bigr)_\F \bigl(-1,-1
\bigr)_\F^{\frac{l(l-1)}{2}}
h_{\F}(\rho), \end{equation}\\
where $\rho$ is the Leray invariant $-q(V^*,V^*\sigma_1,V^*
\sigma_2^{-1})$, $d_\F(\rho)$ and $h_\F(\rho)$ are its
discriminant, and Hasse invariant, and $2l=j_1+j_2-j-\dim (\rho)$,
where $\sigma_1 \in \Omega_{j_1},\, \sigma_2 \in \Omega_{j_2},\,
\sigma_1 \sigma_2 \in \Omega_{j}$. We use Rao's normalization of
the Hasse invariant. (Note that the cocycle formula just given
differs slightly from the one that appears in Rao's paper. There
is a small mistake in Theorem 5.3 of \cite{R}. A correction by
Adams can be found in \cite{Kud}, Theorem 3.1).

An immediate consequence of  Rao's formula is that if $g$ and $h$
commute in $\spn$ then their pre-image in $\mspn$ also commute
(this may be deduced from more general ideas. See page 39 of \cite
{MVW}). In particular, a preimages in $\mspn$ of a commutative
subgroup of $\spn$ is also commutative. This does not hold for
general covering groups; for example, using the cocycle
constructed in Section \ref{ext rao}, the reader may check that
the inverse image of the diagonal subgroup of $\gspn$ in
$\overline{\gspn}$ is not commutative. See \cite{BanPhD} for
the same phenomenon in the $n$-fold cover of $\gln$.

If $\F$ is a local field then $\mspn$ is a locally compact group.
If $\F$ is a p-adic field,  $\mspn$ is an l-group in the sense of
\cite{BZ}; since $c(\cdot,\cdot)$ is continuous, it follows that
there exists $U$, an open compact subgroup of $\spn$, such that
$c(U,U)=1$. Thus, a system of neighborhoods of $(I_{2n},1)$ is
given by open compact subgroups of the form $(V,1)$, where $V
\subseteq U$ is an open compact subgroup of $\spn$. As an example
we describe it explicitly in the $\msl$ case; see Lemma
\ref{example sl2}.

 From \eqref{cocycle} and from previous remarks we obtain the
following properties of $c(\cdot,\cdot)$; for $\sigma,\sigma'\in
\Omega_j$, $p,p' \in P(\F)$ we have
\begin{eqnarray} \label{cocycle prop} \! \! \! \! \! \! \! \! \! \!
&& c(\sigma,\sigma^{-1})=\bigl( x(\sigma),(-1)^jx(\sigma)
\bigr)_\F
(-1,-1)_\F^{\frac{j(j-1)}{2}} \\
\label{coycle prop 3} \! \! \! \! \! \! \! \! \! \! && c(p \sigma,
\sigma' p')=c(\sigma,\sigma') \bigl (x(p),x(\sigma)\bigr)_\F \bigl
(x(p'),x(\sigma')\bigr)_\F \bigl (x(p),x(p')\bigr)_\F \bigl
(x(pp'),x(\sigma\sigma')\bigr)_\F.\end{eqnarray} As a consequence
of \eqref{coycle prop 3} we obtain

\begin{equation} \label{cocycle prop 2}  c(p,\sigma)=c(\sigma,p)=\bigl(
x(p),x(\sigma) \bigr)_\F. \end{equation} Another property of the
cocycle noted in \cite{R} is that
\begin{equation} \label{cocycle tau}
c(\tau_{S_1},\tau_{S_2})=(-1,-1)_\F^{\frac {j(j+1)}{2}},
\end{equation} where $j$ is the cardinality of $S_1 \cap S_2$.
From \eqref{t1 mult t2}, \eqref{coycle prop 3} and \eqref{cocycle
tau} we conclude that if $S$ and $S'$ are disjoint then for $p,p'
\in P(\F)$ we have

\begin{equation} \label{cocycle prop 4}
c(p \tau_S,\tau_{S'}p')=\bigl(x(p),x(p') \bigr)_\F.\end{equation}
It follows from \eqref{cocycle prop 2} and \eqref{Hilbert
properties} that
\begin{equation} \label{pgp-1}
(p,\epsilon_1)(\sigma,\epsilon)(p,\epsilon_1)^{-1}=(p \sigma
p^{-1},\epsilon),
\end{equation}
for all $\sigma \in \spx , \, p \in P(\F), \, \epsilon_1 ,
\epsilon \in \{ \pm 1 \}$. Furthermore, assume that $p \in P(\F)$,
$\sigma \in \spx$ satisfy $\sigma p \sigma ^{-1} \in P(\F)$. Then
\begin{equation} \label{gpg-1}
(\sigma,\epsilon_1)(p,\epsilon)(\sigma,\epsilon_1)^{-1}=(\sigma p
\sigma^{-1} ,\epsilon). \end{equation} Indeed, due to
\eqref{pgp-1} and the Bruhat decomposition we only need to show
that if $\tau_S p {\tau_S} ^{-1} \in P(\F)$ then
$$c(p,\tau_S)c(\tau_S p ,\tau_S^{-1})c(\tau_S,{\tau_S}^{-1})=1.$$
Define $j$ to be the cardinality of $S$. From \eqref{cocycle prop}
it follows that $c(\tau_S,{\tau_S}^{-1})=(-1,-1)_\F^{\frac
{j(j-1)}{2}}$. From \eqref{cocycle prop 2} it follows that
$c(p,\tau_S)=1$. It is left to show that if $V^*\tau_S p
\tau_S^{-1}=V^*$ then \begin{equation}\label{an coc com}c(\tau_S p
,{\tau_S}^{-1})=(-1,-1)_\F^{\frac {j(j-1)}{2}}.\end{equation}
Recall that the Leray invariant is stable under the action of
$\spx$ on Lagrangian triplets; see Theorem 2.11 of \cite{R}.
Therefore,
$$q(V^*,V^* \tau_S p, V^*\tau_S)=q(V^*{\tau_S}^{-1},V^* \tau_S p
{\tau_S}^{-1}, V^*)=q(V^*{\tau_S}^{-1},V^*, V^*)$$ is an inner
product defined on the trivial space. \eqref{cocycle} implies now
\eqref{an coc com}.

We recall Corollary 5.6 in Rao's paper. For $S \subset
\{1,2,\ldots,n\}$ define $$X_S=span \{e_i,e_i^* \mid i \in S\}.$$
We may now consider $x_S$ and $c_{X_S}(\cdot,\cdot)$ defined by
analogy with $x$ and $c(\cdot,\cdot)$. Let $S_1$ and $S_2$ be a
partition of $\{ 1,2,\ldots,n \}$. Suppose that
$\sigma_1,\sigma_1' \in Sp(X_{S_1})$ and that $\sigma_2,\sigma_2'
\in Sp(X_{S_2})$. Put $\sigma=diag(\sigma_1,\sigma_2)$,
$\sigma'=diag(\sigma_1',\sigma_2')$. Rao proves that
$c(\sigma,\sigma')$ equals \begin{eqnarray} \label{multiplicative}
c_{S_1}(\sigma_1,\sigma_1')c_{S_2}(\sigma_2,\sigma_2') \bigl
(x_{S_1}(\sigma_1),x_{S_2}(\sigma_2) \bigr )_\F \bigl
(x_{S_1}(\sigma_1'),x_{S_2}(\sigma_2') \bigr )_\F \bigl
(x_{S_1}(\sigma_1 \sigma_1'),x_{S_2}(\sigma_2 \sigma_2')
\bigr)_\F.
\end{eqnarray}
From \eqref{multiplicative} it follows that $(s,\epsilon) \mapsto
\bigl(i_{r,n}(s), \epsilon \bigr)$ and $(s,\epsilon) \mapsto
\bigl(j_{r,n}(s), \epsilon \bigr)$ are two embeddings of
$\overline{Sp_{2r}(\F)}$ in $\overline{Sp_{2n}(\F)}$. We shall
continue to denote these embeddings by  $i_{r,n}$ and $j_{r,n}$
respectively. Note that the map $g \mapsto (\widehat{g},1)$ is not
an embedding of $\gln$ in $\mspn$, although, by \eqref{cocycle
prop 2}, its restriction to $\zgln$ is an embedding.
\subsection{Kubota's cocycle} \label{kubota} For n=2 Rao's cocycle reduces to Kubota's
cocycle, \cite{Kub}:
\begin{equation}\label{rao}c(g_1,g_2)=\bigl(x(g_1),x(g_2)\bigr)_\F
\bigl(-x(g_1)x(g_2),x(g_1 g_2)\bigr)_\F,\end{equation} where
$$x \begin{pmatrix} _{a} & _{b}\\_{c} &
_{d}\end{pmatrix}=\begin{cases} c & c \neq 0 \\ d & c=0
\end{cases}.$$
For $\F$, a p-adic field of odd residual characteristic it is
known (see page 58 of \cite{M}) that $\msl$ splits over
$SL_2(\Of)$, the standard maximal compact subgroup of $SL_2(\F)$
and that
$$\iota_{2}:SL_2(\Of) \rightarrow \{\pm 1\}$$ defined by
$$\iota_2 \begin{pmatrix} _{a} & _{b} \\
_{c} & _{d}
\end{pmatrix}=\begin{cases} (c,d)_\F & 0<\ab c \ab <1 \\ 1 &
$otherwise$ \end{cases}$$ is the unique map such that the map $$k
\mapsto \kappa_{2}(k)= \bigl(k,\iota_2(k) \bigr)$$ is an embedding
of $SL_2(\Of)$ in $\msl$. More generally, it is known, see
\cite{M}, that if $\F$ is a p-adic field of odd residual
characteristic then $\mspn$ splits over $\Kn$; there exists a map
$$\iota_{2n}:\Kn \rightarrow \{\pm 1\}$$ such that the map $$k
\mapsto \kappa_{2n}(k)= \bigl(k,\iota_{2n}(k)\bigr)$$ is an
embedding of $\Kn$ in $\mspn$. Since $\kappa_{2}$ is the unique
splitting of $SL_2(\Of)$ in $\overline{SL_2(\F)}$ and since $\Kn$
is generated by various embedding of $SL_2(\Of)$ it follows that
$\iota_{2n}$ is also unique.

\begin{lem} \label{k split}The restrictions of $\iota_{2n}$ to $P(\F) \cap \Kn$
and to $\wspntag$ are trivial.
\end{lem}
\begin{proof}
Since for odd residue characteristic $(\Of^*,\Of^*)_\F=1$ (see
Lemma \ref{usefull} for this well known fact), we conclude, using
\eqref{cocycle prop 2}, that $\iota_{2n}$ restricted to $P(\F)
\cap \Kn$ is a quadratic character and hence has the form $p
\mapsto \chi(\det \, p)$, where $\chi$ is a quadratic character of
$\Of^*$. By the inductivity property of Rao`s cocycle and by the
formula of $\iota_2$ we conclude that
$$\iota_{2n}\big(i_{1,n}(SL_2(\Of)\bigr)\cap P(\F)=1.$$ Thus $\chi=1$. We now
move to the second assertion. We note that the group generated by
the elements of the form $\tau_S$ is the group of elements of the
form $\tau_{S_1}a_{S_2}$, where $S,S_1,S_2 \subseteq
\{1,2,\ldots,n)$. The group $\{\widehat{w_\pi}\mid \pi \in S_n \}$
is disjoint from that group and normalizes it. Hence we need only
to show that for all $S_1,S_2 \subseteq \{1,2,\ldots,n \}$, and
for all $\pi \in S_n$ $$\iota_{2n}(\widehat{w_\pi} a_{S_1}
\tau_{S_2})=1.$$ The fact that $\iota_{2n}(\widehat{w_\pi}
a_{S_2})$=1 was proved already. We now show that
$\iota_{2n}(\tau_{S_1})=1$: The fact that for $|S_1|=1$:
$\iota_{2n}( \tau_{S_1})=1$ follows from the properties of
$\iota_2$ and the inductivity properties of Rao's cocycle. We
proceed by induction on the cardinality of $S_1$; suppose that if
$\mid S_1 \mid \leq l$ then $\iota_{2n}( \tau_{S_1})=1$. Assume
now that $\mid S_1 \mid=l+1$. Write $S_1=S'\cup S''$, where $S'$
and $S''$ are two non-empty disjoint sets. By \eqref{cocycle tau}
we have
$$\iota_{2n}(\tau_{S_1})=\iota_{2n}(\tau_{S'})\iota_{2n}(\tau_{S''})
c(\tau_{S'},\tau_{S''})=1.$$ Finally, $$\iota_{2n}(\widehat{w_\pi}
a_{S_2} \tau_{S_1})=
\iota_{2n}(\widehat{w_\pi})\iota_{2n}(a_{S'})\iota_{2n}(\tau_{S})c(\widehat{w_\pi},
a_{S'})c(\widehat{w_\pi} a_{S'},\tau_{S})=1.$$
\end{proof} For all the p-adic fields we have the following lemma; define
$$K_n=\{\begin{pmatrix} _{1+a} & _{b}\\_{c} & _{1+d}\end{pmatrix}
\in SL_2(\F) \mid a,b,c,d \in \Pf^n \}$$
\begin{lem} \label{example sl2}$c(K_n,K_n)=1$ for all $n \geq 2e+1$.
\end{lem}
\begin{proof}
In Lemma \ref{usefull} we prove that $1+\Pf^n \in {\F^*}^2$ for
all $n \geq 2e+1$. We shall use this fact here. Suppose that
$g=\begin{pmatrix} _{a} & _{b}\\_{0} & _{a^{-1}}\end{pmatrix}$,
where $a \in 1+\Pf^n, \, n \geq 2e+1$. Then, for all $h \in
SL_2(\F)$, since $a \in {\F^*}^2$ we have
$$c(g,h)=c(h,g)=\bigl(a,x(h)\bigr)_{\F}=1.$$ It is left to show
that for $g=\begin{pmatrix} _{a} & _{\frac{ad-1}c}\\_{c} &
_{d}\end{pmatrix} \in K_n, \, \, \, h=\begin{pmatrix} _{x} &
_{\frac{xw-1}z}\\_{z} & _{w}\end{pmatrix}\in K_n$, where $c,z \neq
0$, we have $c(g,h)=1$. Note that
$$x(gh)=\begin{cases} cx+dz & cx+dz \neq 0 \\ax+(ad-1)\frac z c &
cx+dz=0 \end{cases}.$$ If $cx+dz=0$ then
$$c(g,h)=(c,z)_\F(-cz,ax+(ad-1)\frac z c)_\F=(c,z)_\F(-cz,\frac z
c) _\F=(-c,-z)_\F.$$ The fact that $c=-\frac d x z$ and that
$\frac x d \in {\F^*}^2$ implies now that $c(g,h)=1$. Suppose now
that $cx+dz \neq 0$. In this case
$$c(g,h)=(c,z)_\F(-cz,cx+dz)_\F=(-cz,x+\frac {dz} c)_\F=(-cz,1+\frac {dz} {cx})_\F.$$
Since for all $p,m \in \F^* , p \neq 1$, we have
$(p,m)_\F=\bigl(p,p(1-m)\bigr)_\F$, we finally get
$$c(g,h)=(1+\frac {dz} {cx}, \frac d xz^2)_\F=1.$$

\end{proof}
Assume now that $\F=\R$. The group $SO_2(\R)$ is a maximal compact
subgroup of $SL_2(\R)$ and it is commutative. Every $k \in
SO_2(\R)$ can be written uniquely as
$$k(t)=\begin{pmatrix} _{cos(t)} & _{sin(t)}\\_{-sin(t)} &
_{cos(t)}\end{pmatrix}$$ for some $0 \leq t < 2 \pi.$ We have:

$$c \bigl(k(t),k(t)^{-1}\bigr)=\begin{cases} 1 & t \neq \pi
 \\-1 & t=\pi  \end{cases}.$$ If $k(t) \neq \pm I_2$ then

$$c \bigl(k(t),-k(t^{-1})\bigr)=\begin{cases} 1 & sin(t)>0 \\-1 & sin(t)<0
\end{cases},$$and
$$c \bigl(k(t),-I_2\bigr)^{-1}=\begin{cases} 1 & sin(t)>0 \\-1 & sin(t)<0 \end{cases}.$$
If none of $k(t_1), k(t_2), k(t_1)k(t_2)$ equals $\pm I_2$ then
$$c\bigl((k(t_1),k(t_2)\bigr)=\begin{cases} -1 & sign\bigl(sin(t_1)\bigr)=sign\bigl(sin(t_2)\bigr)
\neq sign\bigl(sin(t_1+t_2)\bigr) \\1 &  $otherwise$\end{cases}.$$
$\overline{SL_2(\R)}$ does not split over $SO_2(\R)$. Indeed, Let
$H$ be any subgroup of  $SL_2(\R)$ which contains $-I_2$. Then,
there is no $\beta:H \rightarrow \{\pm 1 \}$ such that
$$\beta(ab)=c(a,b)\beta(a)\beta(b)$$ since this function must
satisfy $$1=\beta\bigl((-I_2)^2\bigr)=c(-I_2,-I_2) \beta^2(-I_2)$$
or, equivalently, $ \beta^2(-I_2)=-1$. It is well known, see page
74 of \cite{G} for example, that $\overline{SO_2(\R)} \simeq \R /
4\pi \Z$. We determine all the isomorphisms between these two
groups:
\begin{lem} \label{so2} Realize $\R / 4\pi \Z$ as the segment $[0,4\pi)$.
There are exactly two isomorphisms between  $\R / 4\pi \Z$ and
$\overline{SO_2(\R)}$. One is $t \mapsto
\phi(t)=\bigl(k(t),\theta(t)\bigr)$, where $\theta:\R / 4\pi \Z
\rightarrow \{\pm 1\}$ is the unique function that satisfies
\begin{equation} \label{theta prop}\theta(t_1+t_2)=\theta(t_1)\theta(t_2)c\bigl(k(t_1),k(t_2)\bigr).\end{equation}
The second is $t \mapsto \phi(-t)$. $\theta$ is given by
\begin{equation} \label{theta}\theta(t)=\begin{cases} 1 & 0\leq t \leq \pi$ or $
3\pi<t<4\pi
\\-1 &  \pi<t\leq 3\pi\end{cases}.\end{equation}
\end{lem}
\begin{proof} Using the cocycle formulas given above, one may check directly that $\theta$ as defined in
\eqref{theta} satisfies \eqref{theta prop}. The uniqueness of
$\theta$ follows from the fact that any function from  $\R / 4\pi
\Z$ to $\{\pm 1\}$ that satisfies \eqref{theta prop} must satisfy
$\theta(t)=c \bigl(k(\frac t 2),k (\frac t 2) \bigr)$. It is now
clear that $\phi$ is indeed an isomorphism. Let $\phi':\R / 4\pi
\Z \rightarrow \overline{SO_2(\R)}$ be an isomorphism. Since
$\phi^{-1}\phi'$ is an automorphism of $\R / 4\pi \Z$ it follows
that either $\phi^{-1}\phi'(t)=t$ or $\phi^{-1}\phi'(t)=-t$. This
implies $\phi'(t)=\phi(t)$ or $\phi'(t)=\phi(-t)$.
\end{proof}

\subsection{Some facts about parabolic subgroups of $\mspn$}
\label{par mspn} Let $\F$ be a local field. By a parabolic subgroup of $\mspn$ we mean an
inverse image of a parabolic subgroup of $\spn$.

\begin{lem}
Let $Q$ a parabolic subgroup of $\spn$. Write $Q=M \ltimes N$, a
Levi decomposition. Then, there exists a unique function $\mu':N
\rightarrow \{\pm 1 \}$, such that $n \mapsto
\mu(n)=\bigl(n,\mu'(n)\bigr)$ is an embedding of $N$ in $\mspn$.
Furthermore: $\overline{Q}=\overline{M} \ltimes \mu(N)$. By abuse
of language we shall refer to the last equality as the Levi
decomposition of $\overline{Q}$.
\end{lem}
\begin{proof} Suppose first that $Q$ is standard. From the fact $\mspn$ splits
over $N$ via the trivial section it follows that $\mu'$ is a
quadratic character of $N$. Since $N=N^2$ we conclude that $\mu'$
is trivial. Using \eqref{gpg-1} we get $\overline{Q}=\overline{M}
\ltimes \mu(N)$. Assume now that $Q$ is a general parabolic
subgroup. Then, $\overline{Q}=(w,1)\overline{Q'}(w,1)^{-1}$ for
some $w \in \spn$, and a standard parabolic subgroup $Q'$. For $n
\in Q$ define $n'=w^{-1}nw$. From the proof in the standard case
it follows that
$$\mu'(n)=c(w,n')c(wn',w^{-1})c(w^{-1},w)=c(nw,w{-1})c(w^{-1},w)$$
 is the unique function mentioned
in the lemma. The fact that $\overline{Q}=\overline{M} \ltimes
\mu(N)$ follows also from the standard case. \end{proof}

\begin{lem} Let $\F$ be a p-adic field. $\overline{\Kk}$ is a maximal
open compact subgroup of $ \mspn$. For any parabolic subgroup $Q$
of $\spn$ we have
$$\mspn=(Q,1)\overline{\Kk}=\overline{\Kk}(Q,1).$$ If $\F$ is a
p-adic field of odd residual characteristic then
$\mspn=\overline{Q}\kappa_{2n}\bigl(Sp_{2n}(\Of)\bigr)$.

\end{lem}
By an abuse of natation we call the last decomposition an Iwasawa
decomposition of $\mspn$.
\begin{proof} This follows immediately from the analogous lemma in the algebraic
case. \end{proof}

Let $\F$ be a p-adic field of odd residual
characteristic.  Pick $\alpha \in \{\pm 1\}$ and $\chi$, a
quadratic character of $\Of^*$. Define $\xi_{\alpha,\chi} \! \!
:\F^* \rightarrow \{\pm 1\}$ by
\begin{eqnarray*} \xi_{\alpha,\chi}(\epsilon
\pi^{2n})&=& \chi(\epsilon)(\pi,\pi)_\F^n  \\
\xi_{\alpha,\chi}(\epsilon \pi^{2n+1}) &=& \alpha
\chi(\epsilon)(\pi,\pi)_\F^n(\epsilon,\pi)_\F, \end{eqnarray*}
where $\epsilon \in \Of^*$ and $n \in \Z$.

\begin{lem} \label{4split} Let $\F$ be a p-adic field of odd residual
characteristic. There are exactly four maps $\xi: \F^* \rightarrow
\{ \pm 1\}$ such that
\begin{equation} \label{h split}\xi(ab)=\xi(a)\xi(b)(a,b)_\F.
\end{equation} Each of them has the form $\xi_{\alpha,\chi}$.
\end{lem}
\begin{proof} The reader can check that
$\xi_{\alpha,\chi}$ satisfies \eqref{h split}. We now show that
each $\xi$ satisfying \eqref{h split} has this form. Indeed,
define $\alpha=\xi(\pi)$ and define $\chi$ to be the restriction
of $\xi$ to $\Of^*$. Since $(\Of^*, \Of^*)_\F=1$ it follows that
$\chi$ is  a quadratic character. We shall show that
$\xi=\xi_{\alpha,\chi}$:
\begin{eqnarray*} \xi(\epsilon
\pi^{2n}) &=&
\xi(\epsilon)\xi(\pi^{2n})(\epsilon,\pi^{2n})_\F=\chi(\epsilon)\xi(\pi^{2n}),
\\
\xi(\epsilon \pi^{2n+1}) &=&
\xi(\epsilon)\xi(\pi^{2n+1})(\epsilon,\pi)_\F= \alpha
\chi(\epsilon)\xi(\pi^{2n})(\epsilon,\pi)_\F. \end{eqnarray*} It
remains to show that
\begin{equation} \label{second cor}\xi(\pi^{2n})=(\pi,\pi)_\F^n\ \end{equation}
for all $n \in \N$.  Since
$\xi(1)=1$ and since $$\xi(\pi^2)=\xi^2(\pi)(\pi,\pi)_\F=\alpha^2(\pi,\pi)_\F=(\pi,\pi)_\F,$$ it
follows that for all $n\in \Z$
\begin{eqnarray*}
\xi(\pi^{2n+2}) &=& \xi(\pi^{2n})\xi(\pi^2)=\xi(\pi^{2n})(\pi,\pi)_\F.\end{eqnarray*}
Similarly, $\xi(\pi^{2n-2})=\xi(\pi^{2n})(\pi,\pi)_\F$ for all $n\in \Z$.
This completes the proof of \eqref{second cor}.

\end{proof}
{\bf Remark}: Since $\F^* \simeq \Of^* \times \Z$ it follows that
$\F^*/{\F^*}^2 \simeq \Of^* / {\Of^*}^2 \times \{\pm 1 \}$. Thus,
there is an isomorphism between the group of quadratic characters
of $\F^*$ and the group of quadratic characters of $\Of^*$ $\times
\{\pm 1\}$. The isomorphism is $\eta \mapsto \bigl
(\eta|_{\Of^*},\eta(\pi) \bigr)$. Since every quadratic character
of $\F^*$ has the form $\eta(b)=\eta_a(b)=(a,b)_\F$, we have
proven that every $\xi$ satisfying \eqref{h split} has the form
$\xi=\xi_a$, where $a \in \F^*$ and $\xi_a$ is defined as follows $$\xi_a(\epsilon
\pi^{2n})=(a,\epsilon)(\pi,\pi)^n,$$
$$\xi_a(\epsilon \pi^{2n+1})=\xi_a(\epsilon
\pi^{2n})(\pi,a\epsilon).$$

\begin{lem} \label{p is product odd f}
If $\F$ is a p-adic field of odd residual characteristic then
$$GL_{n_1}(\F) \times GL_{n_2}(\F) \ldots \times GL_{n_r}(\F)
\times \mspk \simeq \mmt.$$ An isomorphism is given by
$$(g_1,g_2,\ldots,g_r,\overline{h})  \mapsto
\bigl(j_{n-k,n}(\widehat{g}),1)(i_{k,n}(h),\epsilon\xi_{\alpha,\chi}(\det
g)\bigr),$$ where for $1 \leq i \leq r$, $g_i \in GL_{n_i}(\F)$,
$g=diag(g_1,g_2,\ldots,g_r) \in GL_{n-k}(\F)$,
$\overline{h}=(h,\epsilon) \in \mspk$. In particular, $\mspn$
splits over the Siegel parabolic subgroup via the map $$p \mapsto
\bigl(p,\xi_{\alpha,\chi}(\det p)\bigr).$$

\end{lem}
\begin{proof} Once the existence of a map $\xi: \F^*
\rightarrow \{ \pm 1\}$ satisfying \eqref{h split} is established,
this lemma reduces to a straight forward computation. One uses the
fact, following from \eqref{cocycle prop 2}, that for $p,p' \in
P(\F)$ we have
$$\xi_{\alpha,\chi}\bigl(\det(p)\bigr)
\xi_{\alpha,\chi}\bigl(\det(p')\bigr)c(p,p')=
\xi_{\alpha,\chi}(\det(pp')\bigr).$$ A closely related argument is
presented in Lemma \ref{def rep levi}.
\end{proof}
We shall not use this lemma.
\subsection{The global metaplectic group} \label{global}
Let $\F$ be a number field, and let $\A$ be its
Adele ring. For every place $\nu$ of $\F$ we denote by $\fnu$ its
completion at $\nu$. We denote by
$\widehat{Sp_{2n}(\A)}$ the restricted product $\prod'_\nu
Sp_{2n}(\fnu)$ with respect to $$\Bigl \{
\kappa_{2n}\bigl(Sp_{2n}(\Ofnu) \bigr) \mid \nu \, \, is\, \,
finite\, \, and\, \, odd \Bigr \}.$$ $\widehat{Sp_{2n}(\A)}$ is
clearly not a double cover of $Sp_{2n}(\A)$. Put $$C'=\Bigr \{
\prod_\nu(I,\epsilon_\nu) \mid \prod_\nu\epsilon_\nu=1 \Bigl \}.$$
We define
$$\overline{Sp_{2n}(\A)}=C'\setminus \widehat{Sp_{2n}(\A)}$$ to be the
metaplectic double cover of $Sp_{2n}(\A)$. It is shown in page 728
of \cite{JS 07} that $$k \mapsto C'\prod_\nu(k,1)$$ is an
embedding of $\spn$ in $\overline{Sp_{2n}(\A)}$.

\subsection{Extension of Rao's cocycle to GSp(X).} \label{ext rao}
Let $\F$ be a local field. Let $\gspx$ be the similitude group of $\spx$. This is the subgroup of $GL_{2n}(\F)$
which consists of elements which satisfy
$<v_1g,v_2g>=\lambda_g<v_1,v_2>$ for all $v_1,v_2 \in X$, where
$\lambda_g \in \F^*$. $\F^*$ is embedded in $\gspx$ via $$\lambda
\mapsto i(\lambda)=\begin{pmatrix} _{I_n} & _{0}\\_{0} & _{\lambda
I_n}\end{pmatrix}.$$ Using this embedding we define an action of
$\F^*$ on $\spx$: $$(g,\lambda) \mapsto
g^{\lambda}=i(\lambda^{-1})g \, i(\lambda).$$ Let $\F^* \ltimes
\spx$ be  the semi-direct product corresponding to this action.
For $g \in \gspx$ define $$p(g)=\begin{pmatrix} _{\alpha} &
_{\beta}\\_{\lambda_g^{-1}\gamma}  & _{\lambda_g^{-1}\delta}
\end{pmatrix} \in \spx,$$
where $g=\begin{pmatrix} _{\alpha} & _{\beta}\\_{\gamma } &
_{\delta}
\end{pmatrix} \in \gspx$, $\alpha,\beta,\gamma, \delta \in Mat_{n \times n}(\F)$.
Note that $g=i(\lambda_g)p(g)$. The map $$g \mapsto
\iota(g)=\bigl(\lambda_g,p(g)\bigr)$$ is an isomorphism between
$\gspx$ and $\F^* \ltimes \spx.$

We know that we can lift the outer conjugation $g \mapsto g^\lambda$
of $\spx$ to $\mspx$ (see page 36 of \cite{MVW}), namely we can define
a map $v_\lambda:\spx \rightarrow \{ \pm 1 \}$ such that
$$(g,\epsilon) \mapsto
(g,\epsilon)^{\lambda}=\bigl(g^{\lambda},\epsilon
v_{\lambda}(g)\bigr)$$ is an automorphism of $\mspx$. In \ref {my
calc} we compute $v_\lambda$. We shall also show there that
$$\bigl(\lambda,(g,\epsilon)\bigr) \mapsto
(g,\epsilon)^{\lambda}$$ defines an action of $\F^*$ on $\mspx$.
Let us show how this computation enables us to
extend $c(\cdot,\cdot)$ to a 2-cocycle
$\widetilde{c}(\cdot,\cdot)$ on $\gspx$ which takes values in
$\{\pm 1\}$ and hence write an explicit multiplication formula of
$\mgspx$, the unique metaplectic double cover of $\gspx$ which
extends a non-trivial double cover of $\spx$. We define the group
$\F^* \ltimes \mspx$ using the multiplication formula
$$\bigl(a,(g,\epsilon_1)\bigr) \bigl
(b,(h,\epsilon_2)\bigr)=\bigl(ab,(g,\epsilon_1)^b(h,\epsilon_2)\bigr).$$
We now define a bijection from the set $\gspx \times \{ \pm 1 \}$
to the set $\F^* \times \mspx$ by the formula
$$\overline{\iota}(g,\epsilon)=\bigl(\lambda_g,(p(g),\epsilon)\bigr),$$
whose inverse is given by
$$\overline{\iota}^{-1}\bigl(\lambda,(h,\epsilon)\bigr)=\bigl(i(\lambda)h,\epsilon
\bigr ).$$ We use $\overline{\iota}$ to define a group structure
on $\gspx \times \{ \pm 1 \}$. A straightforward computation shows
that the multiplication in $\mgspx$ is given by
$$(g,\epsilon_1)(h,\epsilon_2)=\bigl(gh,v_{\lambda_h}(p(g)\bigr)c \bigl(p(g)^{\lambda_h},p(h)\epsilon_1\epsilon_2 \bigr).$$
Thus,
\begin{equation} \label{extension}
\widetilde{c}(g,h)=v_{\lambda_h}\bigl(p(g) \bigr )c \bigl
(p(g)^{\lambda_h},p(h) \bigr)
\end{equation}
serves as a non-trivial 2-cocycle on $\gspx$ with values in $\{\pm
1\}$. We remark here that Kubota, \cite{Kub} (see also \cite{G}),
used a similar construction to extend a non-trivial double cover
of $SL_2(\F)$ to a non-trivial double cover of $GL_2(\F)$. For
$n=1$ our construction agrees with Kubota's.

\subsubsection{Computation of $v_\lambda(g)$} \label{my calc} In
\cite{Bar} Barthel extended Rao's unnormalized cocycle to $\gspx$.
One may compute $v_\lambda(g)$ using Barthel's work and Rao's
normalizing factors. Instead, we compute $v_\lambda(g)$ using
Rao's (normalized) cocycle. Fix $ \lambda \in \F^*$. Since
$(g,\epsilon) \mapsto (g,\epsilon)^\lambda$ is an automorphism,
$v_{\lambda}$ satisfies:
\begin{equation} \label{property of v}
v_\lambda(g)v_\lambda(h)v_\lambda(gh)=\frac
{c(g^{\lambda},h^{\lambda})} {c(g,h)}. \end{equation} We shall
show that this property determines $v_{\lambda}$.

We first note the following: \begin{equation} \label{on p}
\begin{pmatrix} _{a} & _{b}\\_{0} & _{\widetilde{a}}
\end{pmatrix} ^{\lambda}=\begin{pmatrix} _{a} & _{\lambda b}\\_{0} & _{\widetilde{a}}
\end{pmatrix},\ \ x\Biggl(\begin{pmatrix} _{a} & _{b}\\_{0} &
_{\widetilde{a}}
\end{pmatrix} ^{\lambda}\Biggr) \equiv x \begin{pmatrix} _{a} & _{b}\\_{0} & _{\widetilde{a}}
\end{pmatrix}.\end{equation}
For $S \subseteq \{1,2,\ldots,n\}$ define $a_S(\lambda) \in P(\F)$
by
\begin{equation} \label{a_S lambda} e_i \cdot a_S(\lambda)=\begin{cases} \lambda^{-1}e_i & i
\in S \\ e_i & $otherwise$ \end{cases}, \,     e_i^* a_S(\lambda)=
\begin{cases} \lambda e_i^* & i \in S \\ e_i^* & $otherwise$
\end{cases}.\end{equation}
Note that $a_S=a_S(-1)$. One can verify that \begin{equation}
\label{tau lambda} \tau_S^\lambda=a_S(\lambda) \tau_S=\tau_S
a_S(\lambda^{-1}). \end{equation}Since $x \bigl (a_S(\lambda)
\bigr) \equiv \lambda^{|S|}$, we conclude, using \eqref{on p},
\eqref{tau lambda}, the Bruhat decomposition and the properties of
$x$ presented in Section \ref{des rao} that
$\Omega_j^{\lambda}=\Omega_j$, and that for $g \in \Omega_j$
\begin{equation} \label{x lambda} x(g^{\lambda}) \equiv \lambda^jx(g).
\end{equation}
\begin{lem} \label{v identity on p,g} For $p \in P(\F), \  g \in
\Omega_j$ we have
\begin{equation} \label{v id on p,g}
\vl(p)\vl(g)\vl(pg)=\bigl(x(p),\lambda^j \bigr)_\F
\end{equation} and \begin{equation} \label{v id on p,g 2}
\vl(g)\vl(p)\vl(gp)=\bigl(x(p),\lambda^j \bigr)_\F
\end{equation}\end{lem}

\begin{proof} We prove \eqref{v id on p,g} only.
\eqref{v id on p,g 2} follows in the same way. We use
\eqref{property of v}, \eqref{cocycle prop 2}, \eqref{x lambda}
and \eqref{Hilbert properties}: $$\vl(p)\vl(g)\vl(pg)=\frac
{c(p^{\lambda},g^{\lambda})} {c(p,g)}=\frac{\bigl
(x(p^{\lambda}),x(g^{\lambda})\bigr)_\F}{\bigl(x(p),x(g)\bigr)_\F}=
\frac{\bigl(x(p),x(g)\lambda^j \bigr
)_\F}{\bigl(x(p),x(g)\bigr)_\F}=\bigl(x(p),\lambda^j \bigr)_\F.$$
\end{proof}
\begin{lem} \label{v is ch on p}
There exists a unique $\tl \in \F^*/ {\F^*}^2$ such that for all
$p \in P(\F)$: $$\vl(p)=\bigl(x(p),\tl \bigr)_\F.$$
\end{lem}
\begin{proof} Substituting $p'\in P(\F)$ instead of $g$ in \eqref{v id on
p,g} we see that ${\vl}|_{_{P(\F)}}$ is a quadratic character.
Since $N_{n;0}(\F)$, the unipotent radical of $P(\F)$ is
isomorphic to a vector space over $\F$, it follows that
$N_{n;0}(\F)^2=N_{n;0}(\F)$. Thus, $\vl \mid _{_{N_{n;0}(\F)}}$ is
trivial. We conclude that ${\vl}|_{_{P(\F)}}$ is a quadratic
character of $\gl$ extended to $P(\F)$. Every quadratic character
of $\gl$ is of the form $g \mapsto \chi \bigl(\det(g) \bigr)$,
where $\chi$ is a quadratic character of $\F^*$. Due to the
non-degeneracy of the Hilbert symbol, every quadratic character of
$\F^*$ has the form $\chi(a)=(a,t_{\chi})_\F$, where $t_{\chi} \in
\F^*$ uniquely determined by $\chi$ up to multiplication by
squares.
\end{proof}
\begin{lem}\label{skeleton of formula} For $\sigma \in \Omega_j$ we have
$$\vl(\sigma)=\bigl(x(\sigma),\tl \lambda^j \bigr)_\F
\vl(\tau_S),$$ where $S \subseteq \{1,2,\ldots,n\}$ is such that
$|S|=j$. In particular, if $|S|=|S'|$ then
$\vl(\tau_S)=\vl(\tau_{S'})$.
\end{lem}
\begin{proof}
An element $\sigma \in \Omega_j$ has the form $\sigma=p\tau_S p'$
where $p,p' \in P(\F),\, \mid \! S \! \mid=j$. Substituting
$g=\tau_S p'$ in \eqref{v id on p,g} yields $$\vl(p\tau_S
p')=\vl(p)\vl(\tau_S p')\bigl(x(p),\lambda^j \bigr)_\F.$$
Substituting $g=\tau_S$ and $p=p'$ in \eqref{v id on p,g 2} yields
$$\vl(\tau_S p')=\vl(\tau_S)\vl(p')\bigl(x(p'),\lambda^j
\bigr)_\F.$$ Using the last two equalities, Lemma \ref{v is ch on
p} and \eqref{Hilbert properties} we obtain $$\vl(p\tau_S p')=
\bigl(x(pp'),\tl \lambda^j \bigr)_F\vl(\tau_S).$$ Since $|S|=|S'|$
implies $p\tau_S p^{-1}=\tau_{S'}$, for some $p \in P(\F)$, the
last argument shows that $\vl(\tau_S)=\vl(\tau_{S'})$.
\end{proof}
It is clear now that once we compute $\tl$ and $\vl(\tau_S)$ for
all $S\ \subseteq \{1,2,\ldots,n \}$ we will find the explicit
formula for $\vl$.
\begin{lem} \label{factors}$\tl=\lambda$ and $\vl(\tau_S)=(\lambda,\lambda)_\F^{\frac{|S|(|S|-1)}{2}}.$
\end{lem}
\begin{proof} Let $k$ be a symmetric matrix in $\gl$. Put $$p_k=\begin{pmatrix} _{k} & _{-I_n}\\_{0} &
_{k^{-1}} \end{pmatrix} \in P(\F), \, \, n_k=\begin{pmatrix}
_{I_n} & _{k}\\_{0} & _{I_n}\end{pmatrix} \in N_{n;0}(\F),$$ and
note that $x(n_k) \equiv 1$, $x(p_k) \equiv \det(k)$, and that
\begin{equation} \label {nice mult} \tau n_k \tau=n_{-k^{-1}}\tau
p_k.\end{equation} We are going to compute $\vl(\tau)\vl(n_k
\tau)\vl(\tau n_k \tau)$ in two ways: First, by Lemma \ref
{skeleton of formula} and by $\eqref{nice mult}$ we have
$$\vl(\tau)\vl(n_k \tau)\vl(\tau n_k
\tau)=\vl(\tau)\vl(\tau)\vl(n_{-k^{-1}}\tau
p_k)=\vl(n_{-k^{-1}}\tau p_k).$$ Since \begin{equation}
\label{nice mult x} x(\tau n_k \tau) \equiv x(n_{-k^{-1}}\tau p_k)
\equiv \det (k),
\end{equation} we obtain, using Lemma \ref {skeleton of formula} again,
\begin{equation} \label{side 1}  \vl(\tau)\vl(n_k \tau)\vl(\tau
n_k \tau)=\bigl(\det(k),\tl \lambda^n \bigr)_\F
\vl(\tau).\end{equation} Second, by \eqref{property of v} we have
\begin{equation} \label{side 2} \vl(\tau)\vl(n_k \tau)\vl(\tau n_k
\tau)=\frac {c \bigl(\tau,n_k\tau \bigr) }{c \bigl (\tau
^{\lambda},(n_k \tau)^{\lambda} \bigr)}.\end{equation} Recall that
$\tau=\tau_{\{1,2\ldots,n \}}$. We shall compute the two terms
on the right side of \eqref{side 2}, starting
with $c(\tau,n_k\tau)$: Let $\rho$ and $l$ be the factors
in \eqref{cocycle}, where $\sigma_1=\tau$, $\sigma_2=n_k \tau$.
Since,
$$q \bigl (V^*,V^* \tau, V^*
(n_k \tau)^{-1} \bigr)=q \bigl( V^*,V,V^*(-I_{2n} \tau
n_{-k})\bigr )=q(V^*,V,V n_{-k}).$$ We conclude that $\rho=k$,
$l=0$. Using \eqref{cocycle}, and \eqref{nice mult x} we observe
that
\begin{equation} \label{c(tau,nk tau)} c(\tau,n_k
\tau)=\bigl (-1,\det (k) \bigr)_\F h_\F(k). \end{equation}

We now turn to $c \bigr(\tau ^{\lambda},(n_k
\tau)^{\lambda}\bigl)$: Let $\rho$ and $l$ be the factors in
\eqref{cocycle}, where $\sigma_1=\tau^{\lambda}$, $\sigma_2=(n_k
\tau)^{\lambda}$. Note that by \eqref{on p} and \eqref{tau
lambda}, $(n_k \tau)^{\lambda}=n_{\lambda k}\lambda I_{2n} \tau$,
hence $$q\Bigl(V^*,V^* \tau ^{\lambda},V^*{\bigl((n_k \tau
)^{\lambda}\bigr)}^{-1} \Bigr)=q(V^*,V,Vn_{-\lambda k}).$$ Thus,
$\rho=mk$, $l=0$, and we get $$c\bigl(\tau ^{\lambda},(n_k
\tau)^{\lambda}\bigr)=\bigl(x(\lambda I_{2n}),x(\lambda
I_{2n})\bigr)_\F \bigl(-1,x(\tau^{\lambda}(n_k
\tau)^{\lambda})\bigr)_\F h_\F(\lambda k).$$ We recall \eqref{nice
mult} and note now that
$$\tau^{\lambda}(n_k \tau)^{\lambda}=(\tau n_k
\tau)^{\lambda}=(n_{-k^{-1}} \tau p_k)^\lambda=n_{-\lambda k^{-1}}
{\lambda}I_{2n} \tau \begin{pmatrix} _{k} & _{\lambda k}\\_{0} &
_{k^{-1}} \end{pmatrix}.$$ Hence, $c\bigl(\tau ^{\lambda},(n_k
\tau)^{\lambda}\bigr)=(\lambda ^n,\lambda ^n)_\F
\bigl(-1,\lambda^n \det (k)\bigr)_\F h_\F(\lambda k),$ or, using
\eqref{Hilbert properties}:
\begin{equation} \label{second c} c \bigl(\tau
^{\lambda},(n_k \tau)^{\lambda}\bigr)=\bigl(-1, \det (k)\bigr)_\F
h_\F(\lambda k).
\end{equation}
Using  \eqref{side 1}, \eqref{side 2},\eqref {c(tau,nk tau)} and
\eqref{second c} we finally get \begin{equation} \label{useful}
\vl(\tau)\bigl(\det (k), \tl \lambda^n \bigr
)_\F=\frac{h_\F(\lambda k)}{h_\F(k)}. \end{equation} By
substituting $k=I_n$ in \eqref{useful}, we get
$\vl(\tau)=(\lambda,\lambda)_\F^{\frac {n(n-1)}{2}}.$ We can now
rewrite \eqref{useful} as
\begin{equation} \label{useful 2}\bigl(\det (k),\tl \lambda^n \bigr)_\F=
\frac{h_\F(\lambda k)}{h_\F(k)}(\lambda,\lambda)_\F^{\frac
{n(n-1)}{2}}.\end{equation} In order to find $\tl$ we note that
for any $y \in \F^*$ we can substitute $k_y=diag(1,1,\ldots,1,y)$
in \eqref{useful 2} and obtain $$(y,\tl \lambda^n)_\F=
(\lambda,\lambda)_\F^{\frac {(n-1)(n-2)}{2}}(\lambda,\lambda
y)_\F^{n-1}(\lambda,\lambda)_\F^{\frac {n(n-1)}{2}}.$$ For both
even and odd $n$ this is equivalent to
$(y,\lambda)_\F=(y,\tl)_\F$. The validity of the last equality for
all $y \in \F^*$ implies that $\tl \equiv \lambda \bigl(mod
(\F^*)^2 \bigr )$.

We are left with the computation of $\vl(\tau_S)$ for $S
\subsetneq \{1,2,\ldots,n\}$. For such $S$, define $_S \tau \in
Sp(X_S)$ by analogy with $\tau \in \spx$. We can embed $Sp(X_S)$
in $\spx$ in a way that maps $_S \tau$ to $\tau_S$. We may now use
\eqref{multiplicative} and repeat the computation of $\vl (\tau)$.
\end{proof}
Joining Lemma \ref{skeleton of formula} and Lemma \ref{factors} we
write the explicit formula for $\vl$: For $g \in \Omega_{j}$ we
have \begin{equation} \label{vl}
\vl(g)=\bigl(x(g),\lambda^{j+1}\bigr)_\F(\lambda,\lambda)^{\frac{j(j-1)}{2}}.
\end{equation} One can easily check now that $\vl(g)
v_{\eta}(g^{\lambda})=v_{\lambda \eta}(g)$, and conclude that the
map $\bigl( \lambda,(g,\epsilon)\bigr) \mapsto
(g,\epsilon)^{\lambda}$ defines an action of $\F^*$ on $\mspx$,
namely that
${\bigl((g,\epsilon)^{\lambda}\bigr)}^{\eta}=(g,\epsilon)^{\lambda
\eta}.$
\begin{cor} \label{why lift}Comparing \eqref{vl} and \eqref {cocycle prop},
keeping \eqref{Hilbert properties} in mind, we note that
\begin{equation} \label{for w}
v_{-1}(g)=v_{-1}(g^{-1})=c(g,g^{-1}). \end{equation} \end{cor}This fact will
play an important role in the proof of the uniqueness of Whittaker
models for $\mspx$, see Chapter \ref{Uniqueness of Whittaker
model}.
\newpage
\section{Weil factor attached to a character of a second degree}
\label{Weil factor}
Let $\F$ be a local field. More accurately, assume that $\F$ is
either $\R$, $\C$ or a finite extension of $\Q _p$. In this chapter
we introduce $\gamma_{\psi}$, the normalized Weil
factor associated with a character of second degree of $\F$ which
 takes values in $\{ \pm i,\pm 1 \} \subset \C^*$ . This factor is an
ingredient in the definition of genuine parabolic induction on
$\mspn$. We prove certain properties of the Weil factor that will
be crucial to the computation the local coefficients in Chapter
\ref{second paper}. We also compute it for p-adic fields of odd
residual characteristic and for $\Q_2$. The Weil factor is trivial
in the field of complex numbers. In the real case it is computed
in \cite{P}. It is interesting to know that unless $\F$ is a
2-adic field $\gamma_\psi$ is not onto $\{\pm 1, \pm i \}$.

Let $\psi$ be a non-trivial character of $\F$. For $a\in \F^*$ let
$\gamma_\psi(a)$ be the normalized Weil factor associated with
the character of second degree of $\F$ given by $x \mapsto
\psi_a(x^2)$ (see Theorem 2 of Section 14 of \cite{Weil}). It is
known that
\begin{equation} \label{gammaprop}
\gamma_{\psi}(ab)=\gamma_{\psi}(a)\gamma_{\psi}(b)(a,b)_\F
\end{equation} and that
\begin{equation} \label{gammaprop square} \gamma_\psi(b^2)=1
,\gamma_{\psi}(ab^2)=\gamma_{\psi}(a),\, \gamma_\psi^4(a)=1.
\end{equation}
From Appendix A-1-1 of \cite{P} it follows that for $\F=\R$ we
have
\begin{equation} \label{perrin}\gamma_{\psi_a}(y)=\begin{cases} 1 & y>0
 \\ -sign(a)i  & y<0  \end{cases}.\end{equation}
(recall that for  $\F=\R$ we use the notation
$\psi_b(x)=e^{ibx}$). Unless otherwise stated, until the end of
this chapter we assume that $\F$ is a p-adic field. It is known
that
\begin{equation} \label{gammadef 1} \gamma_\psi(a)= \ab a \ab
^\half \frac{\int_\F \psi_a(x^2) \, dx}{\int_\F \psi(x^2) \,
dx},\end{equation} see page 383 of \cite{BFH} for example. These
integrals, as many of the integrals that will follow, should be
understood as principal value integrals. Namely,
\begin{equation}
\label{gammadef 2} \gamma_\psi(a)= \ab a \ab ^\half
 \frac{lim_{n \rightarrow
\infty}\int_{\Pf^{-n}} \psi(ax^2) \, dx}{lim_{n \rightarrow
\infty}\int_{\Pf^{-n}} \psi(x^2) \, dx}.\end{equation} For $a \in
\F^*$ define $c_\psi(a)=\int_\F \psi(ax^2) \, dx$. With this
notation:
$$\gamma_\psi(a)=\ab a \ab ^\half c_\psi(a)c^{-1}_\psi(1).$$
Recall the following definitions and notations from Chapter
\ref{genral not}: $e=e(\F)=[\Of:2\Of]$, $\omega=2\pi^{-e}\in \Of^*$.
If $x \in \Of$ then $\overline{x}$ is its image in the residue field
$\overline{\F}$.

\begin{lem} \label{2e} $1+\Pf^{2e+1} \subseteq {\F^*}^2$,
$1+\Pf^{2e} \nsubseteq {\F^*}^2õ$
\end{lem}

\begin{proof}  The lemma is known for $\F$ of odd residual
characteristic, see Theorem 3-1-4 of \cite{W} for example. We now
assume that $\F$ is of even residual characteristic. The proof of
the first assertion in this case resembles the proof for the case
of odd residual characteristic. We must show the for every $b \in
1+\Pf^{2e+1}$ the polynomial $f_b(x)=x^2-b \in \Of[x]$ has a root
in $\F$. This follows from Newton`s method: By Theorem 3-1-2 of
\cite{W}, it is sufficient show that $\ab \frac
{f_b(1)}{f_b'^2(1)} \ab <1$, and it is clear. We now prove the
second assertion. Pick $b\in \Of^*$ such that the polynomial
$$p_b(x)=x^2+\overline{\omega} x-\overline{b} \in \rf [x]$$ does not
have a root in $\rf$. Such a polynomial exists: The map $x \mapsto
\phi(x)=x^2+\overline{\omega}x$ from $\rf$ to itself is not
injective since $\phi(0)=\phi(\overline{\omega})$ and since
$\overline{\omega} \neq 0$. Therefore, it is not surjective. Thus,
there exists $b \in \Of^*$ such that $x^2+\overline{\omega}x \neq
\overline{b}$ for all $x \in \rf$. Define now $y=1+b \pi^{2e} \in
1+\Pf^{2e}$. We shall see that $y \notin {\F^*}^2$. It is
sufficient to show that there is no $x_0 \in \Of^*$ such that
$x_0^2=y$. Suppose that such an $x_0$ exists. Since $x_0^2-1
\equiv 0 \, (mod \, \Pf)$ it follows that $\overline{x_0}$ is a
root of $$f(x)=x^2-1=(x-1)^2 \in \rf[x].$$ This polynomial has a
unique root. Thus, we may assume $x_0=1+a\pi^k+c\pi^{k+1}$ for
some $k \geq 1, \, a\in \Of^*, \, c\in \Of$. If $k<e$ then
$$x_0^2=1+a^2\pi^{2k}+c'\pi^{2k+1}$$ for some $c'\in \Of$. Hence
$x_0^2 \neq y$. If $k>e$ then $$x_0^2=1+ \omega a\pi^{k+e}+c''
\pi^{k+e+1}$$ for some $c''\in \Of$. Again $x_0^2 \neq y$. So,
$k=e$ and $$x_0^2=1+(a^2+\omega a)\pi^{2e}+c'''\pi^{2e+1}$$ for
some $c'''\in \Of$. Since we assumed that $x_0^2=y$ it follows
that $a^2+\omega a=b$. This contradicts the fact that $p_b(x)$
does not have root in $\rf$.
\end{proof}
\begin{lem} \label{gamma comp} Assume that $\psi$ is normalized. For $a \in \Of^*$ we have
\begin{eqnarray}
 \label{gamma on of*} \gamma_\psi(a)&=&c^{-1}_\psi(1) \Bigl(1+ \sum_{n=1}^e
q^n\int_{\Of^*} \psi(\pi^{-2n}x^2a) \, dx \Bigr) \\
 \label{gamma on pf-1}
\gamma_\psi(\pi a)&=&q^{-\half}c^{-1}_\psi(1) \Bigl(1+
\sum_{n=1}^{e+1} q^n\int_{\Of^*} \psi(\pi^{1-2n}x^2a) \,  dx \Bigr)\\
 \label{gamma-1 on of*} \gamma_\psi^{-1}(a)&=&c^{-1}_\psi(-1) \Bigl(1+ \sum_{n=1}^e
q^n\int_{\Of^*} \psi^{-1}(\pi^{-2n}x^2a) \, dx \Bigr) \\
\label{gamma-1 on pf-1} \gamma_\psi^{-1}(\pi
a)&=&q^{-\half}c^{-1}_\psi(-1) \Bigl(1+ \sum_{n=1}^{e+1}
q^n\int_{\Of^*} \psi^{-1}(\pi^{1-2n}x^2a) \, dx \Bigr)
\end{eqnarray}
\end{lem}
\begin{proof}
In order to prove \eqref{gamma on of*} and \eqref{gamma on pf-1},
it is sufficient to show that
\begin{equation} \label{cpsi on of*} c_\psi(a)=1+ \sum_{n=1}^e
q^n\int_{\Of^*} \psi(\pi^{-2n}x^2a)\, dx
\end{equation} and that
\begin{equation} \label{cpsi on pf-1} c_\psi(a\pi)=1+
\sum_{n=1}^{e+1} q^n\int_{\Of^*} \psi(\pi^{1-2n}x^2a) \, dx.
\end{equation} We prove only the first assertion. The
proof of the second is done in the same way. For $a \in \Of^*$ we
have \begin{eqnarray*} c_\psi(a)&=&\int_{\Of} \psi(ax^2)\,
dx+\int_{\ab x \ab
>1} \psi(ax^2) \,dx\\ &=&1+\sum_{n=1}^\infty \int_{\ab x \ab =q^n}
\psi(ax^2) \,dx=1+\sum_{n=1}^\infty q^n \int_{\Of^*} \psi(a \pi^
{-2n}x^2) \, dx.\end{eqnarray*} It remains to show that for $n>e$
we have $\int_{\Of^*} \psi(a \pi^ {-2n}x^2) \, dx=0$. Indeed, by
Lemma \ref{2e}, for every $u \in \Pf^{2e+1}$ we may change the
integration variable $x \mapsto x\sqrt{1+u}$ and obtain
$$\int_{\Of^*} \psi(a \pi^ {-2n}x^2)dx=\int_{\Of^*} \psi(a \pi^
{-2n}x^2)\psi(a \pi^ {-2n}x^2u)dx.$$ Hence
\begin{eqnarray*} \int_{\Of^*} \psi(a \pi^
{-2n}x^2)dx &=& \mu(\Pf^{2e+1})^{-1}\int_{\Pf^{2e+1}}\int_{\Of^*}
\psi(a \pi^ {-2n}x^2)\psi(a \pi^ {-2n}x^2u)dx \, du \\
&=& q^{2e+1}\int_{\Of^*} \psi(a \pi^ {-2n}x^2)
\biggl(\int_{\Pf^{2e+1}} \psi(a \pi^ {-2n}x^2u) \, du
\biggr)dx.\end{eqnarray*} For $n>e$ the last inner integral
vanishes. To prove \eqref{gamma-1 on of*} and \eqref{gamma-1 on
pf-1} note that $|\gamma_\psi(a)|=1$ implies
$\gamma_\psi(a)^{-1}=\overline{\gamma_\psi(a)}$ and since $c_\psi$
is defined via integrations over compact sets we have:
$$\gamma_\psi^{-1}(a)=\ab a \ab ^\half
c_\psi(-a)c^{-1}_\psi(-1).$$
\end{proof}
Since $\gamma_\psi$ is determined by its values on $\Of^*$ and
$\pi \Of^*$ the last lemma provides a description of $\gamma_\psi$
for normalized $\psi$. Furthermore, there is no loss of generality
caused by the normalization assumption; recalling the definition
of $\gamma_\psi$ and \eqref{gammaprop} we observe that
\begin{equation} \label{norgama}
\gamma_\psi(a)=\frac
{\gamma_{\psi_0}(a\pi^{-n})}{\gamma_{\psi_0}(\pi^{-n})}=\gamma_{\psi_0}(a)(a,\pi^n)_\F.
\end{equation}
(recall that in Chapter \ref{genral not} we have defined
$\psi_0(x)=\psi(x \pi^n)$, where $n$ is the conductor of $\psi$. $\psi_0$ is normalized)

\begin{lem} \label{usefull} Assume that $\psi$ is normalized.\\
1. $\gamma_\psi(uk)=\gamma_\psi(k)$ for all $k \in
\Of^*, u\in 1+\Pf^{2e}$.\\
2. $\gamma_\psi(u)=1$ for all $u\in 1+\Pf^{2e}$.\\
3. $(1+\Pf^{2e},\Of^*)_\F=1$.\\
4. The map $x \mapsto (x,\pi)_\F$ defined on $1+\Pf^{2e}$ is a non-trivial character.\\
\end{lem}
{\bf Remark}: A particular case of the first assertion of this
lemma is the well known fact that that $(\Of^*, \Of^*)_\F=1$ for
all the p-adic fields of odd residual characteristic.

\begin{proof} Due to \eqref{gamma on of*}, to prove the first
assertion  it is sufficient to show that
$$\psi(\pi^{-2n}x^2k)=\psi(\pi^{-2n}x^2uk)$$ for all $k,x \in
\Of^*, \, u\in 1+\Pf^{2e}, \, n \leq e$ and this is clear. The
second assertion is a particular case of the first. The third
assertion follows from the first and from \eqref{gammaprop}. We
prove the forth assertion. By Lemma \ref{2e} we can  pick $u \in
1+\Pf^{2e}$ which is not a square. By the non-degeneracy of the
Hilbert symbol, there exists $x \in \F^*$ such that $(x,u)_\F=-1$.
Write $x=\pi^n k$ for some $n \in \Z, \, k\in \Of^*$. Then, by the
third assertion of this lemma we have
$-1=(x,u)_\F={(\pi,u)_\F}^n$. Thus, $n$ is odd and the assertion
follows.
\end{proof}
The last two lemmas will be sufficient for our purposes, i.e.,
the computation of the local coefficients in Chapter \ref{second
paper}. However, for the sake of completeness we prove the
following:
\begin{lem} \label{gamma complete f odd}Suppose $\F$ is of odd residual characteristic.
If the the conductor of $\psi$ is even then for $a \in \Of^*$ we
have:
\begin{eqnarray} \label{f odd cond even} \gamma_\psi(a) &=& 1 \\ \label{f odd cond even pi} \gamma_\psi(\pi a) &=&
\begin{cases} \gamma_\psi(\pi) & a \in {\F^*}^2  \\ -\gamma_\psi(\pi) & a \notin {\F^*}^2
\end{cases}. \end{eqnarray}
If the the conductor of $\psi$ is odd then for $a \in \Of^*$ we
have:
\begin{eqnarray} \label{f odd cond odd} \gamma_\psi(a) &=& \begin{cases} 1 & a \in {\F^*}^2  \\ -1 & a \notin {\F^*}^2
\end{cases}  \\ \gamma_\psi(\pi a) &=&
\ \gamma_\psi(\pi).\end{eqnarray} Also
\begin{equation} \label{f odd gamapi }\gamma_\psi(\pi )=q^{-\half}\sum_{\overline{x} \in
\overline{\F}} \overline{\psi}(\overline{x}^2) \times
\begin{cases} -1 & -1 \notin {\F^*}^2 $ and the conductor of
$\psi$ is odd$ \\ 1  & $otherwise$
\end{cases} \in \begin{cases} \{\pm 1\} & -1 \in {\F^*}^2  \\  \{\pm i\} & -1 \notin {\F^*}^2
\end{cases},\end{equation} where $\overline{\psi}$
is the character of $\overline{\F}$ defined by
$\overline{\psi}(\overline{x})=\psi_0(\pi^{-1}x)$.
\end{lem}
\begin{proof}
The fact that $(\Of^*, \Of^*)_\F=1$ for all the  p-adic fields of
odd residual characteristic combined with the non-degeneracy of
the Hilbert symbol implies that for $a\in \Of^*$ we have
\begin{equation} \label {(of,pi)}  (a,\pi)_\F= \begin{cases} 1 & a \in {\F^*}^2  \\ -1 &
a \notin {\F^*}^2 \end{cases}. \end{equation} It follows that
$$\gamma_\psi(\pi)^2=\gamma_\psi(\pi^2)(\pi,\pi)_\F=(\pi,-1)_\F=\begin{cases} 1 & \, -1 \in {\F^*}^2   \\ -1 &
-1 \notin {\F^*}^2 \end{cases}.$$ This implies that
$$\gamma_\psi(\pi) \in \begin{cases} \{\pm 1\} & -1 \in {\F^*}^2  \\  \{\pm i\} & -1 \notin {\F^*}^2
\end{cases}.$$
By \eqref{(of,pi)} and \eqref{norgama} it is sufficient to prove
the rest of the assertions mentioned in this lemma for normalized
characters. In this case \eqref{f odd cond even} is simply the
second part of Lemma \ref{usefull}. \eqref {f odd cond even pi}
follows now from \eqref{gammaprop}, \eqref{f odd cond even} and
\eqref{(of,pi)}. It remains to prove that if $\psi$ is normalized
then
$$\gamma_\psi(\pi )=q^{-\half} \sum_{\overline{x} \in \overline{\F}}
\overline{\psi}(\overline{x}^2).$$ Indeed, by \eqref{gamma on
pf-1} and by \eqref{cpsi on of*} we have
$$\gamma_\psi(\pi)=q^{-\half} \Bigl(1+ q \int_{\Of^*}
\psi(\pi^{-1}x^2) dx \Bigr)=q^{-\half} \Bigl(1+ q \sum_{x \in
\Of^* / 1+\Pf} \int_{1+\Pf} \psi(\pi^{-1}x^2y^2) dy \Bigr)$$
Note that
$$\psi(\pi^{-1}x)=\psi(\pi^{-1}xy)$$ for all $x \in \Of, \, \, y\in 1+\Pf$.
This shows that $\overline{\psi}$ is well defined and  that
$$\gamma_\psi(\pi)=q^{-\half} \Bigl(1+ \sum_{x \in \Of^* /
1+\Pf} \psi(\pi^{-1}x^2) \Bigr)=q^{-\half} \sum_{\overline{x} \in
\overline{\F}} \overline{\psi}(\overline{x}^2).$$

\end{proof}

\begin{lem} \label{gamma complete q2}Let $\psi$ be a non-trivial character of $\Q_2$.
Denote its conductor by $n$. Assume that $a \in \Oq^*$. We have
$\psi(a2^{n-2}) \in \{ \pm i \}$. Also, if $n$ is even then
\begin{eqnarray} \label{q2 cond even of} \gamma_\psi(a) &=&
\begin{cases} 1 & a=1 $ mod $ 4   \\ \psi(-2^{n-2}) & a= -1 $ mod $ 4
\end{cases} \\ \label{q2 cond even 2of} \gamma_\psi(2a)
&=&\gamma_\psi(2)
\begin{cases} 1 & a=1 $ mod $ 8   \\ \psi(-2^{n-2}) & a= -1 $ mod $
8
 \\ -1 & a= 5 $ mod $8
 \\ -\psi(-2^{n-2}) & a= -5$ mod $ 8
\end{cases}, \end{eqnarray}
while if $n$ is odd we have
\begin{eqnarray} \label{q2 cond odd of} \gamma_\psi(a) &=&
\begin{cases} 1 & a=1 $ mod $ 8   \\ \psi(-2^{n-2}) & a= -1 $ mod $
8
 \\ -1 & a= 5 $ mod $8
 \\ -\psi(-2^{n-2}) & a= -5$ mod $ 8
\end{cases}  \\ \label{q2 cond odd 2of} \gamma_\psi(2a) &=&
\gamma_\psi(2) \begin{cases} 1 & a=1 $ mod $ 4   \\ \psi(-2^{n-2})
& a= -1 $ mod $ 4
\end{cases}. \end{eqnarray}
Finally, \begin{eqnarray} \label{gamma2}\gamma_\psi(2)=\frac
{\sqrt{2}\psi(2^{n-3})} {1+\psi(2^{n-2})} \in \{ \pm 1
\}\end{eqnarray}

\end{lem}
\begin{proof}
The computation of the Hilbert symbol for $\Q_2$ is well known
(see Theorem 1 of Chapter 3 of \cite{Ser} for example). For $a \in
\Oq^*$ we have
\begin{eqnarray} \label{hilbert on q2}(2,a)_{\Q_2} &=&  \begin{cases} 1 & a= \pm 1 $
mod $ 8
\\-1 & a= \pm5 $ mod $ 8
\end{cases}  \end{eqnarray}
Note that from \eqref{Hilbert properties} it follows that for any
field of characteristic different then 2 we have
\begin{eqnarray} \label{genral fact}(2,2)_\F=(2,-1)_\F=\bigl(2,-(1-2) \bigr)_\F=1 \end{eqnarray}
These facts, combined with \eqref{norgama} imply that it is
sufficient to prove this lemma with the additional assumption that
$\psi$ is normalized. From \eqref{cpsi on of*} it follows that if
$a\in \Oq^*$ then \begin{equation}
\label{cq2of}c_\psi(a)=1+2\int_{\Of^*} \psi(\frac 1 4 x^2a)\,
dx.\end{equation} Next note that for any $x\in \Oq^*$
\begin{equation} \label{psicont}\psi(\frac 1 4x^2a)=\psi(\frac a 4).\end{equation}
Indeed, $[\Oq^*:1+\Pq]=\mid {\overline{\Q_2}}^* \mid =1.$ Thus,
$\Oq^*=1+\Pq$. Therefore, for  any $x\in \Oq^*$ there exists $t
\in \Pq$ such that $x=1+t$. This implies that $$\psi(\frac 1
4x^2a)=\psi(\frac a 4)\psi(\frac {at} 2 + \frac {t^2}
4)=\psi(\frac a 4).$$ \eqref{cq2of} and \eqref{psicont} imply
now that

$$\gamma_\psi(a)=\frac {1+\psi(\frac a 4)}  {1+\psi(\frac 1 4)}$$

By Lemma \ref{usefull}, $\gamma_\psi(\Oq^*)$ is determined on
$\Oq^* / 1+\Pq^2$. Since $\{ \pm 1 \}$ is a complete set of
representatives of $\Oq^* / 1+\Pq^2$, \eqref{q2 cond even of} will
follow once we show that $\psi(\frac 1 4)$ is a primitive fourth
root of 1. Since $\psi^2(\frac 1 4)=\psi(\half)$ it is sufficient
to show that $\psi(\half)=-1$. Since $\psi^2(\half)=1$ and since
$\psi$ is normalized we only have to show that for $x \in \Oq^*$
we have $\psi(\frac x 2)=\psi(\half)$. Write $x=1+t$, where $t\in
\Pq$. We have: $\psi(\frac x 2)=\psi(\frac 1 2)(\frac t
2)=\psi(\frac 1 2)$. \eqref{q2 cond even 2of} follows now from
\eqref{q2 cond even of} and from \eqref{hilbert on q2}.

We now prove \eqref{gamma2}. From \eqref{genral fact} it follows
that in any field of characteristic different then 2,
$\gamma_\psi^2(2)=1.$ From \eqref{gamma on pf-1} and \eqref{cq2of}
it follows that $$\gamma_\psi(2)=  \frac{1+ \sum_{n=1}^{2}
2^n\int_{\Of^*} \psi(2^{1-2n}x^2) dx} {\sqrt{2} \bigl(1+\psi(\frac
1 4) \bigr )}.$$ We have already seen that $\psi(\half \Oq^*)=-1$.
This implies
$$\gamma_\psi(2) = \frac{
4\int_{\Of^*} \psi(\frac 1 8 x^2) dx} {\sqrt{2} \bigl(1+\psi(\frac
1 4) \bigr )}.$$ We now show that $\psi(\frac 1 8 x^2)=\psi(\frac
1 8)$ for all $x\in \Oq^* $ and conclude that $$\gamma_\psi(2) =
\frac{ \sqrt{2} \psi(\frac 1 8)} { \bigl(1+\psi(\frac 1 4) \bigr
)}.$$ Indeed, for $x \in \Oq^*$ we write again $x=1+m$, where $m
\in \Pq$. With this notation: $\psi(\frac 1 8 x^2)=\psi(\frac 1
8)\psi(\frac {2m} 8)\psi(\frac {m^2} 8)$. We finally note that if
$\ab m \ab=\half$ then $\psi(\frac {2m} 8)=\psi(\frac {m^2} 8)=-1$
and if $\ab m \ab < \half$ then $\psi(\frac {2m} 8)=\psi(\frac
{m^2} 8)=1$.
\end{proof}

The last two lemmas give a complete description of $\gamma_\psi$
for p-adic field of odd residual characteristic and for $\F=\Q_2$.
From Lemma \ref {gamma complete f odd} it follows that for $\F$, a
p-adic field of odd residual characteristic we have
\begin{eqnarray}
\label{less4}\mid \gamma_\psi(\F^*) \mid=\begin{cases} 2 & $if $ \,-1 \in {\F^*}^2  \\
3 & $if $ \, -1 \notin {\F^*}^2
\end{cases}. \end{eqnarray}
Note that in the case where $-1 \in {\F^*}^2$ we must conclude
that $\gamma_\psi$ is one of the four functions $\xi_{a,\chi}$
described in Lemma \ref{4split}. Combining \eqref{less4} with
\eqref{perrin} we observe that if $\F$ is either $\R$ or a p-adic
field of odd residual characteristic then $\gamma_\psi$ is not
onto $\{\pm 1, \pm i \}$. Lemma \ref{gamma complete q2} shows that
there exists a 2-adic field such that the cardinality of
$\gamma_\psi(\F^*)$ is 4.
\newpage
\section{Some facts from the representation theory of $\mspn$}
\label{gen rep meta} Let $\F$ be a p-adic field. Throughout this
section we list several general results concerning the
representation theory of $\mspn$. Although $\mspn$ is not a linear
algebraic group, its properties enable the extension of the theory
presented in \cite{Sil book}, \cite{Wa03}, \cite{BZ}, and
\cite{BZ77}. Among these properties are the obvious analogs of
Bruhat, Iwasawa, Cartan and Levi decompositions and the fact that
$\mspn$ is an l-group in the sense of \cite{BZ}. These properties
are common for $n$-fold covering groups.

We first note that since $\mspn$ is an l-group we have the same
definitions of smooth and admissible representations as in the
linear case; see page 18 of \cite{BanPhD} or page 9 of \cite{Zo09}
for example. We shall mainly be interested in smooth, admissible,
genuine representations of $\mspn$. A representation
$(V,\overline{\sigma})$ of $\mspn$ is called genuine if
$$\overline{\sigma}(I_{2k},-1)=-Id_V.$$
It means that $\msig$ does not factor through the projection map
$Pr:\mspn \rightarrow \spn$. Same definition applies to
representations of $\mmt$.
\subsection{Genuine parabolic induction} \label{gen par ind}
For a representation $(\tau,V)$ of $\gln$ and a complex number $s$
we denote by  $\tau_{(s)}$ the representation of $\gln$ in $V$
defined by
$$g \mapsto \ab \det(g) \ab^s \tau(g).$$
Put $\overrightarrow{t}=(n_1,n_2,\ldots,n_r;k)$, where
$k+\sum_{i=1}^r n_i=n$. Let $(\tau_1,V_{\tau_1}),
(\tau_2,V_{\tau_2}),\ldots,(\tau_r,V_{\tau_r})$ be $r$
representations of $GL_{n_1}(\F),GL_{n_2}(\F),\ldots,GL_{n_r}(\F)$
respectively. Let $(\overline{\sigma},V_{\overline{\sigma}})$ be a
genuine representation of $\mspk$. We shall now describe a
representation of $\mpt$ constructed from these representations.
We cannot repeat the algebraic construction since generally
$$\overline{\mt} \not \simeq GL_{n_1}(\F) \times GL_{n_2}(\F)
\ldots \times GL_{n_r}(\F) \times \mspk$$ (these groups are
isomorphic in the case of p-adic fields of odd residual
characteristic, see Lemma \ref{p is product odd f}). Instead we
define $$\bigl(\otimes_{i=1}^r (\gamma_{\psi}^{-1}\otimes
{\tau_i}_{(s_i)}) \bigr)\otimes \overline{\sigma}:\mmt \rightarrow
GL \bigl((\otimes_{i=1}^{i=r} V_{\tau_i}) \otimes
V_{\overline{\sigma}} \bigr)$$ by \begin{equation} \label{mta ind
def}
 \bigl(\otimes_{i=1}^r (\gamma_{\psi}^{-1}\otimes
{\tau_i}_{(s_i)}) \bigr)\otimes \overline{\sigma}
\bigl(j_{n-k,n}(\widehat{g}),1)(i_{k,n}(h),\epsilon)\bigr)=
\gamma_{\psi}^{-1}\bigl(\det(g)\bigr)\Big(\otimes_{i=1}^r
{\tau_i}_{(s_i)}(g_i)\Bigr)\otimes \overline{\sigma}(h,\epsilon),
\end{equation}
where for $1 \leq i \leq r$, $g_i \in GL_{n_i}(\F)$,
$g=diag(g_1,g_2,\ldots,g_r) \in GL_{n-k}(\F)$, $h \in \spk$ and
$\epsilon \in \{\pm 1 \}$. When convenient we shall use the
notation
\begin{eqnarray} \label{ind short} \pi(\overrightarrow{s}) &=&\bigl(\otimes_{i=1}^r
(\gamma_{\psi}^{-1}\otimes {\tau_i}_{(s_i)}) \bigr)\otimes
\overline{\sigma} \\ \nonumber  \pi&=&\pi(0).\end{eqnarray}
\begin{lem} \label{def rep levi} $\pi(\overrightarrow{s})$ is
a representation of $\mmt$.
\end{lem}

\begin{proof} Let $\alpha=\bigl(j_{n-k,n}(\widehat{g}),1 \bigr)\bigl(i_{k,n}(h),\epsilon\bigr)$
and $\alpha'=\bigl(j_{n-k,n}(\widehat{g'}),1
\bigr)\bigl(i_{k,n}(h'),\epsilon'\bigr)$, where for $1 \leq i \leq
r$, $g_i,g_i' \in GL_{n_i}(\F)$, $g=diag(g_1,g_2,\ldots,g_r), \,
g'=diag(g_1',g_2',\ldots,g_r') \in GL_{n-k}(\F)$, $h,h' \in \spk$
and $\epsilon, \epsilon' \in \{\pm 1 \}$ be two elements in
$\mmt$. It is sufficient to show that if $v=(\otimes_{i=1}^r
v_i)\otimes w$ is a pure tensor in $(\otimes_{i=1}^{i=r}
V_{\tau_i}) \otimes V_{\overline{\sigma}}$ then $\pi(\alpha'
\alpha)v= \pi(\alpha')\bigl(\pi (\alpha)v\bigr).$ Indeed,
$$\pi(\alpha')\bigl(\pi (\alpha)v\bigr)=
\gamma_{\psi}^{-1}\bigl(\det(g)\bigr)\gamma_{\psi}^{-1}\bigl(\det(g')\bigr)
\Big(\otimes_{i=1}^r {\tau_i}(g_i'){\tau_i}(g_i)v_i\Bigr)\otimes
\overline{\sigma}(h',\epsilon')\overline{\sigma}(h,\epsilon)w.$$
Due to \eqref{cocycle prop 2}, for $p,p' \in P(\F)$ we have
\begin{equation} \label{coc gam} \gamma_{\psi}^{-1}\bigl(\det(p)\bigr)
\gamma_{\psi}^{-1}\bigl(\det(p')\bigr)c(p,p')=
\gamma_{\psi}^{-1}(\det(pp')\bigr).
 \end{equation}
Recalling that $\msig$ is genuine we now see that
\begin{equation} \label
{right side} \pi(\alpha')\bigl(\pi (\alpha)v\bigr)=\epsilon
\epsilon' \gamma_{\psi}^{-1}\bigl(\det(gg')\bigr)
c\bigl(j_{n-k,n}(\widehat{g'}),j_{n-k,n}(\widehat{g})\bigr)c(h',h)
\Big(\otimes_{i=1}^r {\tau_i}(g_i'g_i)v_i\Bigr)\otimes
\overline{\sigma}(h'h,1)w. \end{equation} Next, we note that since
$\bigl(j_{n-k,n}(\widehat{g'}),1 \bigr)$ and
$\bigl(i_{k,n}(h),\epsilon \bigr)$ commute we have
$$\alpha'\alpha=\Bigl(j_{n-k,n}(\widehat{gg'}),1 \Bigr)
\Bigl(i_{k,n}(h'h),\epsilon
\epsilon'c(h',h)c\bigl(j_{n-k,n}(\widehat{g'}),j_{n-k,n}(\widehat{g})\bigr)\Bigr).$$
\eqref{right side} implies now that $\pi(\alpha' \alpha)v=
\pi(\alpha')\bigl(\pi (\alpha)v\bigr).$
\end{proof}
As in the linear case we note that if
$\tau_1,\tau_2,\ldots,\tau_r$ and $\msig$ are smooth (admissible)
representations then $\pi(\overrightarrow{s})$ is also smooth
(admissible). Due to \eqref{gpg-1} it is possible to extend
$\pi(\overrightarrow{s})$ to a representation of $\pt$ by letting
$(\nt,1)$ act trivially.

Assuming that $\tau_1,\tau_2,\ldots,\tau_r$ and $\msig$ are smooth
we define smooth induction
\begin{equation} \label{induced notation}I\bigl(\pi(\overrightarrow{s})\bigr)=I({\tau_1}_{(s_1)},{\tau_2}_{(s_2)},\ldots,{\tau_r}_{(s_r)},\msig)=
Ind^{\mspn}_{\mpt} \pi(\overrightarrow{s}) \end{equation} and
\begin{equation} \label{induced notation no s} I(\pi)=I(\tau_1,\tau_2,\ldots,\tau_r,\msig)
=Ind^{\mspn}_{\mpt} \pi.\end{equation} All the induced
representations in this dissertation are assumed to be normalized,
i.e., if $(\pi,V)$ is a smooth representation of $H$, a closed
subgroup of a locally compact group $G$, then $Ind^G_H \pi$ acts
in the space of all right-smooth functions on $G$ that take values
in $V$ and satisfy $f(hg)=\sqrt{\frac
{\delta_H(h)}{\delta_G(h)}}\pi(h)f(g),$ for all $h \in H, g \in
G$. Whenever we induce from a parabolic subgroup (a pre-image of a
parabolic subgroup in a metaplectic group) we always mean that the
inducing representation is trivial on its unipotent radical (on
its embedding in the metaplectic group).

We claim that if the inducing representations are admissible, then
$\pi(\overrightarrow{s})$ is also admissible. Indeed, the proof of
Proposition 2.3 of \cite{BZ77} which is the p-adic analog of this
claim applies to $\mspn$ as well since it mainly uses the
properties on an $l$-group; see Proposition 4.7 of \cite{Zo09} or
Proposition 1 of \cite{BanPhD}.

Next we note that similar to the algebraic case we can define the
Jacquet functor, replacing the role of $\zspn$ with $(\zspn,1)$.
The notion of a supercuspidal representation is defined via the
vanishing of Jacquet modules along unipotent radicals of parabolic subgroups. The (metaplectic) Jacquet functor
has similar properties to those of the Jacquet functor in the
linear case; see Section 4.1 of \cite{BanPhD} or Proposition 4.7
of \cite{Zo09}. This similarity follows from the fact that
$(\zspn,1)$ is a limit of compact groups.

\subsection{An application of Bruhat theory} \label{app bruhat
theory} Let $\overrightarrow{t}=(n_1,n_2,\ldots,n_r;k)$ where
$n_1,n_2,\ldots,n_r,k$ are $r+1$ non-negative integers whose sum
is $n$. Recall that $W_{\pt}$ is the subgroup of Weyl group of $\spn$ which consists of elements which
maps $\mt$ to a standard Levi subgroup and commutes with $i_{k,n}(\spk)$; see Section \ref{The symplectic group}.
\begin{lem} \label{soudry asked}
For $1 \leq i \leq r$ let $g_i$ be an element of $GL_{n_i}(\F)$ and
let $(h,\epsilon)$ be an element of $\mspk$. Denote again $g=diag(g_1,g_2,\ldots,g_r)
\in GL_{n-k}(\F)$. For $w \in W_\pt$ we have
$$(w,1)\bigl(j_{n-k,n}(\widehat{g}),1 \bigr)\bigl(i_{k,n}(h),\epsilon\bigr)(w,1)^{-1}=
\bigl(wj_{n-k,n}(\widehat{g})w^{-1},1 \bigr)\bigl(i_{k,n}(h),\epsilon\bigr)$$
\end{lem}
\begin{proof} This follows from \eqref{gpg-1} and from the fact that $w$ commutes with $h$.
\end{proof}

For $1 \leq i \leq r$ let $\tau_i$ be an irreducible admissible supercuspidal
  representation of $GL_{n_i}(\F)$. Let
$\msig$ be an irreducible admissible supercuspidal  genuine
representation of $\mspk$. Then,
$$\pi=\bigl(\otimes_{i=1}^r (\gamma_{\psi}^{-1}\otimes
{\tau_i}) \bigr)\otimes \overline{\sigma}$$ is an irreducible
admissible supercuspidal genuine representation of
$\overline{\mt}$. For $w \in W_{\pt}$ denote by $\pi^{w}$ the
representation of $\overline{\mt^w}$ defined by $$(m,\epsilon)
\mapsto \pi\bigl((w,1)^{-1}(m,\epsilon)(w,1) \bigr).$$ Exactly as
in the algebraic case, see Section 2 of \cite{BZ77}, $I(\pi)$ and
$I(\pi^{w})$ have the same Jordan Holder series.

Recalling \eqref{w action on element}, we note that due to
Lemma \ref{soudry asked} and the fact that
$\gamma_{\psi}(a)=\gamma_{\psi}(a^{-1})$ we have
\begin{eqnarray} \label{w action} \Bigl((\otimes_{i=1}^r \gamma_{\psi}^{-1}\otimes {\tau_i}_{(s_i)})\otimes
\overline{\sigma}\Bigr)\Bigl({(w,1
)}^{-1}\bigl(j_{n-k,n}(\widehat{g}),1)(i_{k,n}(h),\epsilon)\bigr){(w,1
)}\Bigr)\\ \nonumber =\gamma_{\psi}^{-1}(\det(g))(\otimes_{i=1}^r
 \ab \det(g_i) \ab ^{\epsilon_i s_i}
{\tau_i}^{(\epsilon_i)}(g_i)
\overline{\sigma}(h,\epsilon),\end{eqnarray} where for $1 \leq i
\leq r$, $g_i \in GL_{n_i}(\F)$, $g=diag(g_1,g_2,\ldots,g_r) \in
GL_{n-k}(\F)$, $h \in \spk$, $\epsilon \in \{\pm 1 \}$ and where
$\tau_i^{(\epsilon_i)}(g_i)=\begin{cases} \tau_i (g_i) &
 \epsilon_i=1 \\  \tau_i \bigl((\omega_n \widetilde{g_i}\omega_n\bigr)
& \epsilon_i=-1 \end{cases}$, where $t_{w}^{-1}\bigl(\omega_n
\widetilde{g_i}\omega_n\bigr) t_{w}$ is to be understood via the
identification of $GL_{n_i}(\F)$ with its image in the relevant
block of $M_{\overrightarrow{t}}$. Hence, it makes sense to denote
by $\pi^w$ the representation $$(\otimes_{i=1}^r \gamma_{\psi}^{-1}\otimes
{(\tau_i^{(\epsilon_i)})}_{ (\epsilon_i s_i)})\otimes
\overline{\sigma}.$$ Note that $\tau_i^{(-1)} \simeq
\widehat{\tau_i}$, where $\widehat{\tau_i}$ is the dual
representation of $\tau$; see Theorem 4.2.2 of \cite{B} for
example.

Define $W(\pi)$ to be the following subgroup of $W_{\pt}$:
\begin{equation} \label{wpi def} W(\pi)=\{w \in W_{\pt} \mid
\pi^{w} \simeq \pi \}.\end{equation} $\pi$ is called {\it regular}
if $W(\pi)$ is trivial and { \it singular} otherwise. Bruhat
theory, \cite{Bru61}, implies the following well known result :
\begin{thm} \label{bruhat} If $\pi$ is regular then
\begin{equation} \label{bruhat dimension}
Hom_{\mspn} \bigl(I(\pi),I(\pi)\bigr) \simeq \C. \end{equation}
\end{thm}
See the Corollary in page 177 of \cite{Har73} for this theorem in
the context of connected algebraic reductive p-adic groups, and
see Proposition 6 in page 61 of \cite{BanPhD} for this theorem in
the context of an $r$-fold cover of $\gln$. This theorem follows
immediately from the description of the Jordan-Holder series of a
Jacquet module of a parabolic induction; see Proposition Theorem 5
in page 49 of \cite{BanPhD}. This Jordan-Holder series is an exact
analog to the the Jordan-Holder series of a Jacquet module of a
parabolic induction in the linear case. The proof of Proposition
Theorem 5 in page 49 of \cite{BanPhD} is done exactly as in the
linear case. It uses Bruhat decomposition and a certain filtration
(see Theorem 5.2 of \cite{BZ77}) which applies to both linear and
metaplectic cases. An immediate corollary of Theorem \ref{bruhat}
is the following.
\begin{thm} \label{basic gen par}
Let $\overrightarrow{t}=(n_1,n_2,\ldots,n_r;k)$ where
$n_1,n_2,\ldots,n_r,k$ are $r+1$ non-negative integers whose sum
is $n$. For $1 \leq i \leq r$ let $\tau_i$ be an irreducible
admissible supercuspidal representation of $GL_{n_i}(\F)$. Let
$\msig$ be an irreducible admissible supercuspidal genuine
representation of $\mspk$. Denote $\pi=\bigl(\otimes_{i=1}^r
(\gamma_{\psi}^{-1}\otimes {\tau_i}) \bigr)\otimes
\overline{\sigma}$. If $\tau_i \neq \tau_j$ for all $1 \leq i < j
\leq n$ and $\tau_i \neq \widehat{\tau_j}$ for all $1 \leq i \leq
j \leq n$ then \eqref{bruhat dimension} holds. If we assume in
addition that $\tau_i$ is unitary for each $1 \leq i \leq r$ then
$I(\pi)$ is irreducible.
\end{thm}
\begin{proof} From \eqref{w action} it follows that $\pi$ is regular. Thus,
\eqref{bruhat dimension} follows immediately from Theorem
\ref{bruhat}. Note that since the center of $\mspk$ is finite
$\msig$ is unitary. Therefore, the assumption that $\tau_i$ is
unitary for each $1 \leq i \leq r$ implies that $\pi$ is unitary.
Hence $I(\pi)$ is unitary. The irreducibility of $I(\pi)$ follows
now from \eqref{bruhat dimension}.
\end{proof}

\subsection{The intertwining operator} \label{sec io def}Let $n_1,n_2,\ldots,n_r,k$
be $r+1$ nonnegative integers whose sum is $n$. Put
$\overrightarrow{t}=(n_1,n_2,\ldots,n_r;k)$. For $1 \leq i \leq r$
let  $(\tau_i,V_{\tau_i})$ be an irreducible admissible
representation of $GL_{n_i}(\F)$ and let
$(\overline{\sigma},V_{\overline{\sigma}})$ be an irreducible
 admissible genuine representation of $\overline{Sp_{2k}(\F)}.$
Fix now $w\in W_\pt$. Define
$$N_{\vt,w}(\F)=\zspn \cap(w \nt^- w^{-1}),$$ where $\nt^-$ is the
unipotent radical opposite to $\nt$, explicitly in the $\spn$
case:
$$\nt^-=J_{2n} \nt J_{2n}^{-1}.$$
For $g \in \mspn$ and $f\in I\bigl( \pi(\overrightarrow{s})\bigr)$
define
\begin{equation} \label{io def} A_{w}f(g)=
\int_{N_{\vt,w}} f \bigl( ((w t_w),1)^{-1}(n,1)g \bigr) \, dn,
\end{equation}
 where $t_{w}$ is a particular diagonal element in $\mt \cap
\wspntag$ whose entries are either 1 or -1. The exact definition of $t_{w}$ will be given in Chapter
\ref{chapter def lc}; see \eqref{tw def}. The appearance of $t_{w}$ here is technical.

As in the algebraic case (see Section 2 of \cite{Sha 1}) the last
integral converges absolutely in some open set of $\C ^r$ and has
a meromorphic continuation to $\C ^r$. In fact, it is a rational
function in $q^{s_1},q^{s_2},\ldots,q^{s_r}$. See Chapter 5 of
\cite{BanPhD} for the context of a p-adic covering group. We
define its continuation to be the intertwining operator
$$A_{w}:I({\tau_1}_{(s_1)},{\tau_2}_{(s_2)},\ldots,{\tau_r}_{(s_r)},\msig)
\rightarrow I\bigl({(\tau_{\pi(1)}^{(\epsilon_1)})}_{(\epsilon_1
s_1)},{(\tau_{\pi(2)}^{(\epsilon_2)})}_{(\epsilon_2
s_2)},\ldots,{(\tau_{\pi(r)}^{(\epsilon_r)})}_{(\epsilon_r
s_r)},\sigma \bigr)$$  Denote $\overrightarrow{s}^{w}=(\epsilon_1
s_1,\epsilon_2 s_2,\ldots,\epsilon_r s_r)$ and ${(\otimes_{i=1}^r
\tau_i)}^{w}=(\otimes_{i=1}^r \tau_{\pi(i)}^{(\epsilon_i)})$.

\subsection{The Knapp-Stein dimension theorem} \label{knapp sec} We keep the notation
and assumptions of the first paragraph of Section \ref{app bruhat
theory}. From Theorem \ref{bruhat} it follows that outside a
Zariski open set of values of $(q^{s_1},\ldots,q^{s_r}) \in
{(\C^*)}^r$,
$$Hom_\mspn \Bigl(I\bigl(\pi(\overrightarrow{s})\bigr),
I\bigl(\pi(\overrightarrow{s})\bigr) \Bigr)\simeq \C.$$ This
implies that for every $w_0 \in W_\pt$ there exists a meromorphic
function
$\beta(\overrightarrow{s},\tau_1,\ldots,\tau_r,\msig,w_0)$ such
that
\begin{equation} \label{def beta}A_{w_0^{-1}}A_{w_0}=
\beta^{-1}(\overrightarrow{s},\tau_1,\ldots,\tau_n,\msig,w_0)Id.
\end{equation}
{\bf Remark}: In the case of connected reductive quasi-split
algebraic group this function differs from the Plancherel measure
by a well understood positive function; see Section 3 of
\cite{ShaOxf} for example. For connected reductive quasi-split
algebraic groups it is known that if (with addition to all the
other assumptions made here) we assume that
$\tau_1,\ldots,\tau_r,\msig$ are unitary then
$\beta(\overrightarrow{s},\tau_1,\ldots,\tau_r,\msig,w_0)$, as a
function of $\overrightarrow{s}$ is analytic on the unitary axis.
This is (part of) the content of Theorem 5.3.5.2 of \cite{Sil
book} or equivalently Lemma V.2.1 of \cite{Wa03}. The proof of
this property has a straight forward generalization to the
metaplectic group.

Let $\Sigma_{\pt}$ be the set of reflections corresponding to the
roots of $\tspn$ outside $\mt$. $W_{\pt}$ is generated by
$\Sigma_{\pt}$. Following \cite{Sil} we denote by $W''(\pi)$ the
subgroup of $W(\pi)$ generated by the elements $w \in \Sigma_{\pt}
\cap W(\pi$) which satisfy
$$\beta(\overrightarrow{0},\tau_1,\ldots,\tau_n,\msig,w)=0.$$
The Knapp-Stein dimension theorem states the following:
\begin{thm} \label{silberger} Let $\overrightarrow{t}=(n_1,n_2,\ldots,n_r;k)$ where
$n_1,n_2,\ldots,n_r,k$ are $r+1$ non-negative integers whose sum
is $n$. For $1 \leq i \leq r$ let $\tau_i$ be an irreducible admissible supercuspidal
unitary representation of $GL_{n_i}(\F)$.
Let $\msig$ be an irreducible admissible supercuspidal genuine
representation of $\mspk$. We have:
 $$Dim \bigl(Hom(I(\pi),I(\pi) \bigr)=[W(\pi) :
W''(\pi)].$$
\end{thm}
The Knapp-Stein dimension theorem was originally proved for
real groups, see \cite{K}. Harish-Chandra and Silberger proved it
for algebraic p-adic groups; see \cite{Sil}, \cite{Sil79} and
\cite{Sil book}. It is a consequence of Harish-Chandra`s
completeness theorem; see Theorem 5.5.3.2 of \cite{Sil book}.
Although Harish-Chandra`s completeness theorem was proved for
algebraic p-adic groups its proof holds for the metaplectic case
as well; see the remarks on page 99 of \cite{F86}. Thus, Theorem
\ref{silberger} is a straight forward generalization of the
Knapp-Stein dimension theorem to the metaplectic group.  The
precise details of the proof will appear in a future publication
of the author. We note here that the theorem presented in
\cite{Sil} is more general; it deals with square integrable
representations. The version given here will be sufficient for our
purposes.
\newpage

\section{Uniqueness of Whittaker model} \label{Uniqueness of Whittaker
model} Let $\F$ be a local field and let $\psi$ be a non-trivial
character of $\F$. We shall continue to denote by $\psi$ the
non-degenerate character of $\zspn$ given by
$$\psi(z)=\psi \bigr(z_{(n,2n)}+\Sigma_{i=1}^{n-1} z_{(i,i+1)}\bigl).$$ From
\ref{cocycle prop 2} it follows that $\overline{\zspn} \simeq
\zspn \times \{\pm 1\}$. We define a character of
$\overline{\zspn}$ by $(z,\epsilon) \mapsto \epsilon \psi(z),$ and
continue to denote it by $\psi$. We shall also denote by $\psi$
the characters of $Z_{\glm}$ identified with
$\bigl(i_{k,n}\bigl(j_{m,k}(\widehat{Z_{\glm}})\bigr),1 \bigr)$
and of $\overline{\zspk}$ identified with
$i_{k,n}\bigl(\overline{\zspk} \bigr)$, obtained by restricting
$\psi$ as a character of $\overline{\zspn}$.

Let $(\pi,V_{\pi})$ be a smooth representation of $\mspn$ (of
$\gl$). By a $\psi$-Whittaker functional on $\pi$ we mean a linear
functional $w$ on $V_\pi$ satisfying $$w \bigl (\pi(z)v
\bigr)=\psi(z)w(v)$$ for all $v \in V_{\pi}$, $ z \in
\overline{\zspn}(\zgln)$. Define $W_{\pi,\psi}$ to be the space of
Whittaker functionals on $\pi$ with respect to $\psi$.  If $\F$ is
archimedean  we add smoothness requirements to the definition of a
Whittaker functional, see \cite{H} or \cite{Sha80}. $\pi$ is
called $\psi$-generic or simply generic if it has a non-trivial
Whittaker functional with respect to $\psi$. If $w$ is a
non-trivial $\psi$-Whittaker functional on $(\pi,V)$ then one may
consider $W_w(\pi,\psi)$, the space of complex functions on
$\mspn$ (on $\gln$) of the form $$g \mapsto w\bigl( \pi(g)v
\bigr),$$ where $v \in V$. $W_w(\pi,\psi)$ is a representation
space of $\mspn$( of $\gln$). The group acts on this space by
right translations. It is clearly a subspace of
$Ind^{\mspn}_{\overline{\zspn}} \psi$ (of $Ind^{\gln}_{\zgln}
\psi$).  From Frobenius reciprocity it follows that if $\pi$ is
irreducible and $\dim (W_{\pi,\psi})=1$ then $W_w(\pi,\psi)$ is
the unique subspace of $Ind^{\mspn}_{\overline{\zspn}} \psi$ (of
$Ind^{\gln}_{\zgln} \psi$) which is isomorphic to $\pi$. In this
case we drop the index $w$ and we write $W(\pi,\psi)$. One can
identify $\pi$ with $W(\pi,\psi)$ which is called the Whittaker
model of $\pi$.

For quasi-split algebraic groups, uniqueness of Whittaker
functional for irreducible admissible representations is well
known, see \cite{Sh}, \cite{GK} and \cite{BZ} for the p-adic case,
and see \cite{H} and \cite{Sh} for the archimedean case.
Uniqueness of Whittaker models for irreducible admissible
representations of $\mspn$ is a well known fact to the experts.
Although this fact had been used many times, there is no written
proof of this uniqueness in general (a uniqueness theorem for
principal series representations in the case of p-adic fields can
be found in \cite{BFH}). In this chapter we correct the situation
and prove
\begin{thm} \label{uniquness} Let $\pi$ be an irreducible
admissible representation of $\mspn$. Then, $$\dim (W_{\pi,\psi}) \leq
1.$$
\end{thm}
We note that uniqueness of Whittaker model does not hold in
general for covering groups. See the introduction of \cite{Ban}
for a $k$ fold cover of $\gln$, see \cite{Bud} for an application
of a theory of local coefficients in a case where Whittaker model
is not unique and see \cite {GHP} for a unique model for genuine
representations of a double cover of $GL_2(\F)$.

 Uniqueness for $\overline{SL_2(\R)}$ was proven in \cite{Wa}. To
prove uniqueness for $\overline{Sp_{2n}(\R)}$ for general $n$, it
is sufficient to consider principal series representations. The
proof in this case follows from Bruhat Theory. In fact, the proof
goes almost word for word as the proof of Theorem 2.2 of \cite{H}
for minimal parabolic subgroups. The proof in this case is a
heredity proof in the sense of \cite{Rod}. The reason that this,
basically algebraic, proof works for $\overline{Sp_{2n}(\R)}$ as
well is the fact that if the characteristic of $\F$ is not 2 then
the inverse image in $\mspn$ of a maximal torus of $\spn$ is
commutative. This implies that its irreducible representations are
one dimensional. Uniqueness for $\overline{Sp_{2n}(\C)}$ follows
from the uniqueness for $Sp_{2n}(\C)$ since
$\overline{Sp_{2n}(\C)}=Sp_{2n}(\C) \times \{ \pm 1\}$. However,
the p-adic case is not as easy. Although it turns out that one may
use similar methods to those used in~\cite{GK}, \cite{Sh}, and
\cite{BZ} for quasi-split groups defined over $\F$, one has to use
a certain involution on $\mspn$ whose crucial properties follow
from the results of Section \ref{ext rao}.

\subsection{Non-archimedean case} \label{padic whi uni sec}
Until the rest of this Chapter, $\F$ will denote a p-adic field.
We may assume that the representations of $\mspn$ in discussion
are genuine, otherwise the proof is reduced at once to the well
known $\spn$ case.
\begin{thm} \label{u of w} If $(\pi,V_\pi)$ is  an irreducible admissible
representation of $\mspn$ then, $$\dim (W_{\pi,\psi}) \cdot \dim
(W_{\widehat{\pi},\psi^{-1}}) \leq 1.$$
\end{thm}

The proof of Theorem \ref{u of w} will show:

\begin{thm} \label{u of w 2}
Suppose $(\pi,V_\pi)$ is  an irreducible admissible representation
of $\mspn$. If from the existence of a non-trivial Whittaker
functional on $\pi$ with respect to $\psi$ one can deduce the
existence of a non-trivial Whittaker functional on $\widehat{\pi}$
with respect to $\psi^{-1}$, then $ \dim (W_{\pi,\psi}) \leq 1$.
\end{thm}

\begin{cor} \label{cor for uni}
 If $(\pi,V_\pi)$ is  an irreducible
admissible unitary representation of $\mspn$ then $ \dim
(W_{\pi,\psi}) \leq 1$.
\end{cor}
\begin{proof}
We show that the conditions of Theorem \ref{u of w 2} hold in this
case. Indeed, if $(\pi,V_\pi)$ is an irreducible admissible
unitary representation of $\mspn$ one can realize the dual
representation in the space $\widetilde{V_\pi}$ which is identical
to $V_\pi$ as a commutative group. The scalars act on
$\widetilde{V_\pi}$ by $\lambda \cdot v=\overline{\lambda}v$. The
action of $\widehat{\pi}$ in this realization is given by
$\widehat{\pi}(g)=\pi(g)$. It is clear now that if $L$ is a
non-trivial Whittaker functional on $\pi$ with respect to $\psi$
then $L$, acting on $\widetilde{V_\pi}$, is a non-trivial
Whittaker functional on $\widehat{\pi}$ with respect to
$\psi^{-1}$.
\end{proof}

Since every supercuspidal representation $\pi$ is unitary, it
follows from Corollary \ref{cor for uni} that $ \dim
(W_{\pi,\psi}) \leq 1$. Furthermore, assume now that $\pi$ is an
irreducible admissible representation of $\mspn$. Then $\pi$ is a
subquotient of a representation induced from a supercuspidal
representation of a parabolic subgroup: Let $H$ be a parabolic
subgroup of $\spn$. We may assume that $H=\pt$ where
$\overrightarrow{t}=(n_1,n_2,\ldots,n_r;k)$. Suppose that for
$1\leq i\leq r$, $\sigma_i$ is a supercuspidal representation of
$GL_{n_i}(\F)$, and that $\pi'$ is a supercuspidal representation
of $\overline{Sp_{2k}(\F)}$. Denote by $\psi_i$ and $\psi'$ the
restriction of $\psi$ to $Z_{GL_{n_i}}(\F)$ and $\zspk$
respectively, embedded in $\overline{M_H}$. Note that $ \dim
(W_{\sigma_i,\psi_i}) \leq 1$ and that $\dim (W_{\pi',\psi'}) \leq
1$. Let $\tau$ be the representation of $\overline{M_H}$ defined
by
$$\bigl(diag(g_1,g_2,\ldots,g_r,h,{^t g_r^{-1}},{^t
g_{r-1}^{-1}},\ldots,{^t g_1^{-1}}),\epsilon \bigr) \mapsto
\otimes_{i=1}^{i=r}\sigma_i(g_i) \gamma_\psi(\det g_i)\otimes
\pi(h,\epsilon),$$ We extend $\tau$
 from $M_H$ to $\overline{H}$, letting the unipotent radical act trivially. Define
$Ind_{\overline{H}} ^{\mspn} \tau $ to be the corresponding
induced representation. One may use the methods of Rodier, \cite
{Rod}, extended by Banks, \cite{Ban}, to a non-algebraic setting
and conclude that $$\dim (W_{{Ind_{\overline{H}}^{\mspn}
\tau},\psi})=\dim(W_{\pi',\psi'}) \Pi_{i=1}^{i=r} \dim
(W_{\sigma_i,\psi_i}).$$ Now, if $V_2 \subseteq V_1 \subseteq
Ind_{\overline{H}}^{\mspn} \tau$ are two $\mspn$ modules then
clearly the dimension of the space of Whittaker functionals on
$V_1$ and $V_2$ with respect to $\psi$ is not greater then $\dim
(W_{{Ind_{\overline{H}}^{\mspn} \tau},\psi})$. It follows now that
$\dim(W_{\pi,\psi}) \leq 1$. This is due to the exactness of
twisted Jacquet functor (this is a basic property of twisted
Jacquet functor. It is proven for the algebraic case in 5.12 of
\cite{BZ}. The proof given there continues to hold in the
metaplectic case; see page 72 of \cite{Wang} for example). Thus,
we proved Theorem \ref{uniquness}.
\subsubsection {Proof of Theorem \ref{u of w}} \label{proof of u
of w}Define the following map on $\spn$: \begin{equation}
\label{tau def} g \mapsto {^\tau g}=\sigma_0 ({^t g})
\sigma_0^{-1},
\end{equation} where $\sigma_0=\begin{pmatrix} _{0} & _{\epsilon_n}\\_{\epsilon_n} &
_{0} \end{pmatrix}, \ \ \epsilon_n=diag
\bigl(1,-1,1\ldots,(-1)^{n+1}\bigr) \in \gl.$ \\ We note that
$\sigma_0^{-1}={^t\sigma_0}=\sigma_0$, and that $\sigma_0 \in
GSp_{2n}(\F)$, with similitude factor -1. Hence, $g \mapsto {^\tau
\! g}$ is an anti-automorphism of order 2 of $\spn$. We now extend
$\tau$ to $\mspn$. A similar lifting was used in \cite{GHP} for
$\overline{GL_2(F)}$.

\begin{lem} \label{simple lift} The map $(g,\epsilon) \mapsto
 {^{\overline{\tau}}(g, \epsilon)}=(^\tau \! g, \epsilon)$ is an anti-automorphism of order 2 of
$\mspn$. It preserves both $\overline{\zspn}$ and $\psi$ and
satisfies $^{\overline{\tau}}\bigl(\spn,\epsilon
\bigr)=\bigl(\spn,\epsilon \bigr).$

\end{lem}
\begin{proof} We note that if $g \in \spn$ then $^tg=Jg^{-1}(-J).$ Hence, $${^\tau
g}=\begin{pmatrix}
  \epsilon_n & 0 \\
  0 & \epsilon_n
\end{pmatrix} \begin{pmatrix}
  -I_n & 0 \\
  0 & I_n
\end{pmatrix} g^{-1} \begin{pmatrix}
  -I_n & 0 \\
  0 & I_n
\end{pmatrix} \begin{pmatrix}
  \epsilon_n & 0 \\
  0 & \epsilon_n
\end{pmatrix}^{-1}.$$
Thus, the map $$(g,\epsilon) \mapsto \bigl(^\tau g, \epsilon
c(g,g^{-1})v_{-1}(g^{-1})c(p_\epsilon,\widetilde{g})c (p_\epsilon
\widetilde{g},p_\epsilon^{-1})c(p_\epsilon,p_\epsilon^{-1})\bigr),$$
where $ p_\epsilon=\begin{pmatrix}
  \epsilon_n & 0 \\
  0 & \epsilon_n \\
\end{pmatrix} \in P, \, \, \widetilde{g}=\begin{pmatrix}
  -I_n & 0 \\
  0 & I_n
\end{pmatrix} g^{-1} \begin{pmatrix}
  -I_n & 0 \\
  0 & I_n
\end{pmatrix}$
is an anti-automorphism of $\mspn$. We now show that
$$c(g,g^{-1})v_{-1}(g^{-1})c(p_\epsilon,\widetilde{g})c (p_\epsilon
\widetilde{g},p_\epsilon^{-1})c(p_\epsilon,p_\epsilon^{-1})=1.$$
Indeed, the fact that $c(p_\epsilon,\widetilde{g})c (p_\epsilon
\widetilde{g},p_\epsilon^{-1})c(p_\epsilon,p_\epsilon^{-1})=1$ is
a property of Rao's cocycle and is noted in \eqref{pgp-1}. The
fact that $c(g,g^{-1})v_{-1}(g^{-1})=1$ is a consequence of the
calculation of $\vl(g)$ and is noted in \ref{for w}. The rest of
the assertions mentioned in this Lemma are clear.
\end{proof}

 Let $S\bigl (\mspn \bigr)$ be the space of Schwartz functions on $\mspn$.
For $h \in \mspn$, $\phi \in S \bigl(\mspn \bigr)$ we define
$\lambda (h)\phi$, $\rho(h) \phi$, and $^\mtau \! \phi$ to be the
following elements of $S\bigl (\mspn \bigr)$:
\begin{equation} \label{lambda rho} \bigl(\rho(h)\phi \bigr)(g)=\phi(gh), \ \
\ \ \bigl(\lambda(h)\phi \bigr)(g)=\phi({h^{-1}}g), \ \ \ \ ^\mtau
\phi(g)=\phi(^\mtau g).
\end{equation} $S\bigl (\mspn \bigr)$ is given an algebra structure, called the Hecke
Algebra, by the convolution. Given $(\pi,V_{\pi})$, a smooth
representation of $\mspn$, we define, as usual a representation of
this algebra in the space $V_{\pi}$ by
\begin{equation} \label {hecke op} \pi(\phi)v=\int_{\mspn}
\phi(g)\pi(g)v \, dg. \end{equation}

The following theorem, known as Gelfad-Kazhdan Theorem, in the
context of $\gl$ (see \cite{GK} or \cite{BZ}), will be used in the
proof of Theorem \ref{u of w}.
\begin{thm} \label{gk} Suppose that $D$ is a functional on $S\bigl (\mspn \bigr)$ satisfying
$$D \bigl(\lambda(z_1) \rho(z_2) \phi
\bigr)=\psi(z_2z_1^{-1})D(\phi)$$ for all $\phi \!  \in \! S\bigl (\mspn \bigr), \, z_1, z_1 \! \in \! \overline{\zspn}$.
Then $D$ is $\mtau$ invariant, i.e.,
$D(^{\mtau} \! \phi)=D(\phi)$ for all $\phi \in S\bigl (\mspn \bigr)$.\\
\end{thm}
We will give the proof of this Theorem in \ref{proof of gk}. Here
we shall use this theorem to prove Theorem \ref{u of w}.
\begin{proof}
Since any irreducible admissible representation of $\mspn$ may be
realized as a dual representation, the proof of Theorem \ref{u of w}
amounts to showing that if $W_{\pi,\psi} \neq {0}$ then $\dim
(W_{\widehat{\pi},\psi^{-1}})\leq 1$. We shall use a similar
argument to the one used in \cite {So1}, Theorem 2.1. Let $w$ be a
non-trivial Whittaker functional on $(\pi,V_\pi)$ with respect to
$\psi$. Suppose $\widehat{w_1}$ and $\widehat{w_2}$ are two
non-trivial Whittaker functionals on $\widehat{\pi}$ with respect
to $\psi^{-1}$. The proof will be achieved once we show that
$\widehat{w_1}$ and $\widehat{w_2}$ are proportional.

For $\phi \in S\bigl (\mspn \bigr)$, let $\pi^*(\phi)w$ be a
functional on $V_\pi$ defined by: \begin{equation}
\label{pi*(phi)w} \bigl(\pi^*(\phi)w \bigr)v=\int _{\mspn}
\phi(g)w \bigl (\pi(g^{-1})v \bigr) \ dg.
\end{equation} Note that since $\phi \in S\bigl (\mspn \bigr)$, $\pi^*(\phi)w$ is smooth even if $w$ is
not. Thus, $\pi^*(\phi)w \in V_{\widehat{\pi}}$. By a change of
variables we note that \begin{equation} \label {hat *}
\widehat{\pi}(h) \bigl(\pi^*(\phi)w \bigr)=\pi^* \bigl
(\lambda(h)\phi \bigr)w.
\end{equation} Define now  $R_1$ and $R_2$ to be the following two functionals on
$S\bigl (\mspn \bigr)$:
\begin{equation} \label{Ri} R_i(\phi)=\widehat{w_i} \bigl(\pi^*(\phi)w \bigr),\ \ \ i=1,2. \end{equation}

Using \eqref{hat *}, the facts that $w,\, \widehat{w_1}, \,
\widehat{w_2}$ are Whittaker functionals, and by changing
variables we observe that for all $z \in \overline{\zspn}$
\begin{equation}
\label{lambda Ri} R_i \bigl(\lambda(z)\phi \bigr
)=\psi^{-1}(z)R_i(\phi), \, \, \, R_i \bigl (\rho(z)\phi \bigr
)=\psi(z)R_i(\phi). \end{equation} From Theorem \ref{gk} it
follows now that $R_i(\phi)=R_i(^{\mtau} \! \phi).$ Hence,
\begin{equation} \label {first of page 3} \widehat{w_i} \bigl(\widehat{\pi}(h)\pi^*(\phi)w \bigr)=
\widehat{w_i} \biggl(\pi^*\Bigl(^\mtau \! \bigl(\lambda(h)\phi
\bigr  ) \Bigr) w \biggr).
\end{equation} Using a change of variables again we also note
that \begin{equation} \label {second of page 3} \pi^*\Bigl (^\mtau
\! \bigl (\lambda(h)\phi \bigr) \Bigr)w=\pi^* (^{\mtau} \! \phi)
\bigl ( \pi^* (^\mtau \! h)w \bigr ).
\end{equation} Joining \eqref{first of page 3} and \eqref{second
of page 3} we obtain \begin{equation} \label{third of page 3}
\widehat{w_i} \bigl (\widehat{\pi}(h)\pi^*(\phi)w \bigr
)=\widehat{w_i} \Bigl(\pi^* \bigl(^{\mtau} \phi \bigl ) \bigr (
\pi^*(^\mtau h)w \bigr ) \Bigr ).
\end{equation} In particular, if for some $\phi \in S\bigl (\mspn \bigr)$, $\pi^*(\phi)w$ is the zero functional, then for all $h
\in \mspn$: $$\widehat{w_i} \Bigl(\pi^* \bigl(^{\mtau} \phi \bigl
) \bigr ( \pi^*(^\mtau h)w \bigr ) \Bigr )=0.$$ In this case, for
all $f \in S\bigl (\mspn \bigr)$:
\begin{equation} \label {second of page 4} 0=\widehat{w_i}
\biggl( \int_{\mspn} f(h) \widehat{w_i} \Bigl(\pi^* \bigl(^{\mtau}
\phi \bigl ) \bigr ( \pi^*(^\mtau h)w \bigr ) \Bigr ) \biggr).
\end{equation}
By the definition of $\pi^*(\phi)w$ and by changing the order of
integration, we have, for all $v \in V_\pi$:
\begin{equation} \label {third of page 4} \int_{\mspn} f(h) \widehat{w_i} \Bigl(\pi^* \bigl(^{\mtau} \phi \bigl ) \bigr (
\pi^*(^\mtau h)w \bigr ) \Bigr )(v) dh=\widehat{\pi}(^\mtau \!
\phi) \bigl( \pi^*(f)w \bigr)(v).
\end{equation} \eqref{second of page 4} and \eqref{third of page
4} yield that if $\pi^*(\phi)w=0$ then for all $v \in V_\pi,$ $f
\in S\bigl (\mspn \bigr)$:
\begin{equation} \label{first of page 5} \widehat{w_i} \Bigl(\widehat{\pi} \bigl(^\mtau
\phi \bigr ) \bigl(\pi ^*(f)w \bigr) \Bigr )(v)=0. \end{equation}
From the fact $\pi$ is irreducible one can conclude that
$$\pi^* \Bigl (S\bigl (\mspn \bigr)\Bigr)w =V_{\widehat{\pi}}.$$
Indeed, since $\pi^* \Bigl (S\bigl (\mspn \bigr)\Bigr)w $ is an
$\mspn$ invariant subspace we only have to show that $\pi^* \Bigl
(S\bigl (\mspn \bigr)\Bigr)w  \neq \{0 \}$, and this is clear.

Hence, using a change of variables once more we see that for all
$\xi \in V_{\widehat{\pi}}$ we have
$$0=\widehat{w_i} \bigl(\widehat{\pi}(^\mtau \phi) \xi \bigr)=
\int_{\mspn} {^\mtau \phi({g^{-1}})} \widehat{w_i} \bigl
(\widehat{\pi}(g^{-1} \xi) \bigr) dg .$$ For $g\in \mspn$ define
$^\omega g=^\mtau \! \! g^{-1}$, and for $\phi \in S\bigl (\mspn \bigr)$ define $^\omega \phi(g)=\phi(^\omega g).$ We have just
seen that if $\pi^*(\phi)w=0$ then
\begin{equation} \label{enable Si} \bigl( (\widehat{\pi})^* (^\omega
\phi) \bigr ) \widehat{w_i}=0.
\end{equation} This fact and the fact that
$\pi^*\Bigl ( S\bigl (\mspn \bigr)\Bigr )w =V_{\widehat{\pi}}$
show that the following linear maps are well defined: For $i=1,2$
define $S_i:V_{\widehat{\pi}} \rightarrow
V_{\widehat{\widehat{\pi}}}$ by $\label{Si} S_i \bigl(\pi^*(\phi)w
\bigr)= \bigl((\widehat{\pi})^*(^\omega \phi) \bigr)
\widehat{w_i}$. One can easily check that $S_1$ and $S_2$ are two
intertwining maps from $\widehat{\pi}$ to $h \mapsto
\widehat{\widehat{\pi}}( {^{\overline{\tau}} h ^{-1}})$. The last
representation is clearly isomorphic to $h \mapsto \pi(
{^{\overline{\tau}} h ^{-1}})$, which due to the irreducibility of
$\pi$ is irreducible. Schur's Lemma guarantees now the existence
of a complex number $c$ such that $S_2=cS_1$. So, for all $\phi
\in S\bigl (\mspn \bigr)$ and for all $\xi \in V_{\widehat{\pi}}$:
$$\int_{\mspn}
\phi(g)(\widehat{w_2}-c\widehat{w_1})\bigl(\widehat{\pi}
(g^{-1})\xi \bigr)dg=0.$$ We can now conclude that $\widehat{w_1}$
and $\widehat{w_2}$ are proportional.
\end{proof}
{\bf Remark}: Let $\pi$ be a genuine generic admissible
irreducible representation of $\mspn$. Assume that $\widehat{\pi}$
is also generic (by Corollary \ref{cor for uni} this assumption
holds if $\pi$ is unitary, in particular if $\pi$ is
supercuspidal). Then, from the last step in the proof just given,
it follows that
$$\widehat{\pi }\simeq h\mapsto \pi( ^{\overline{\tau}} h
^{-1}).$$ In page 92 of \cite{MVW} a similar result is proven for
$\pi$, an irreducible admissible representation of $\mspn$
provided that its character is locally integrable. Recently, Sun
has proved this result for any genuine admissible irreducible
representation of $\mspn$, see \cite{SB}.

\subsubsection {Proof of Theorem \ref{gk}} \label{proof of gk} Put
$G=\mspn$ and define $H=\overline{\zspn} \times \overline{\zspn}$.
Denote by $\psii$ be the character of $H$ defined by
$\psii(n_1,n_2)=\psi({n_1^{-1}}n_2)$. $H$ is acting on $G$ by
$(n_1,n_2)\cdot g=n_1 g n_2^{-1}.$ For $g \in G$ we denote by
$H_g$ the stabilizer of $g$ in $H$. It is clearly a unimodular
group. Let $Y$ be an $H$ orbit, that is, a subset of $G$ of the
form $H \cdot g=\zspn g\zspn$, where $g$ is a fixed element in
$G$. Let $S(Y)$ be the space of Schwartz functions on $Y$. $H$
acts on $S(Y)$ by
$$(h \cdot \phi)(k)=\phi(h^{-1} \cdot k) \psii^{-1}(h).$$

With this notation, the proof of Theorem \ref{gk} goes almost
word for word as Soudry's proof of Theorem 2.3 in \cite{So1}. The
main ingredient of that proof was Theorem 6.10 of \cite{BZ} which
asserts that the following four conditions imply Theorem \ref{gk}:\\
1. The set $\{(g,h \cdot g) \mid g \in G, h \in H \}$ is a union
of finitely many locally closed subsets of $G \times G$.\\
2. For  each $h \in H$ there exists $h_\mtau \in H$ such that for
all $g \in G$: $h \cdot {^\mtau g}={^\mtau(h_\tau \cdot g)}$.\\
3. $\mtau$ is of order 2. \\
4. Let $Y$ be an $H$ orbit. Suppose that there exists a non zero
functional on S(Y), satisfying $D(h \cdot \phi)=D(\phi)$ for all
$\phi \in S(Y), h \in H$,  then $^\mtau Y=Y$ and for all $\phi \in
S(Y)$, $D(^\mtau \phi)=D(\phi)$.

Of these four conditions, only the forth requires some work. In
order to make Soudry's proof work in our context we only have to
change Theorem 2.2 of \cite{So1} to

\begin{thm} \label{unproven} Fix $g \in G$. If for all $h \in
H_g$ we have $\psii^{-1}(h)=1$ then, there exists $h^g \in H$ such
that $h^g \cdot g={^\mtau g}$ and $\psii^{-1}(h)=1$.
\end{thm}
Before we prove this theorem  we state and prove its analog for
$\spn$. After the proof we give a short argument which completes
the proof of Theorem \ref{unproven}

\begin{lem} \label{unproven 2}
 For a fixed $g \in
\spn$ one of
the following holds:\\
A. There exist $n_1,n_2 \in \zspn$ such that $n_1gn_2=g$ and
$\psi(n_1n_2) \neq 1$.\\
B. There exist $n_1,n_2 \in \zspn$ such that $n_1gn_2={^\tau \!
g}$ and $\psi(n_1n_2)=1$.
\end{lem}
(in fact this lemma gives the uniqueness of Whittaker model in the linear case)

\begin{proof} Due to the fact that $\tau$ preserves both $\zspn$ and
$\psi$, it is enough to prove this lemma only for a complete set
of representatives of $_\zspn \diagdown ^\spn \diagup _\zspn$. We
recall the Bruhat decomposition: $$\spn=\bigcup_{w \in \wspn}
\zspn \tspn w \zspn.$$ We realize the set of Weyl elements as $$
\bigl\{ \overline{w}_\sigma \tau_S \mid \sigma \in S_n , \, S
\subseteq \{1,2,\ldots,n \} \bigr\},$$ where for $\sigma \in S_n$
we define $w_\sigma \in \gl$ by
$(w_\sigma)_{i,j}=\delta_{i,\sigma(j)}$, and
$\overline{w}_\sigma=\begin{pmatrix}
  w_\sigma &0 \\
  0 & w_\sigma
\end{pmatrix} \in \spn$. Thus, as a complete set of
representatives of $_\zspn \diagdown ^\spn \diagup _\zspn$ we may
take $$\bigl\{ d^{-1} \overline{w}_\sigma^{-1} \varphi_S \mid d
\in T, \, \sigma \in S_n, \, S \subseteq \{1,2,\ldots,n
\}\bigr\},$$ where $\varphi_S=\tau_{S^c}a_{S^c}=\begin{pmatrix}
 M_{S} &M_{S^c} \\
  -M_{S^c}& M_{S}
\end{pmatrix} $, where for $S \subseteq \{1,2,\ldots,n \}
$, $M_{S} \in Mat_{n \times n}(F)$ is defined by
${(M_S)}_{i,j}=\delta_{i,j} \delta_{i \in S}$.

Denote by $w_k$ the $k \times k$ invertible matrix defined by
$(w_k)_{i,j}=\delta_{i+j,n+1}$. Suppose that
$k_1,k_2,\ldots,k_p,k$ are non negative integers such that
$k+\sum_{i=1}^p k_i=n.$ Suppose also that $a_1,a_2,\ldots,a_p \in
\F^*$, and $\eta \in \{\pm 1 \}$. For
\begin{equation} \label{wtypeb}
\overline{w}_\sigma=diag(w_{k_1},w_{k_1},\ldots,w_{k_p},I_k,w_{k_1},w_{k_1},\ldots,w_{k_p},I_k),\end{equation}

\begin{equation} \label{dtypeb} d=diag(a_1\epsilon_{k_1},a_2\epsilon_{k_2},\ldots,a_p\epsilon_{k_p},\eta
I_k,a_1^{-1}\epsilon_{k_1},a_2^{-1}\epsilon_{k_2},\ldots,a_p^{-1}\epsilon_{k_p},\eta
I_k), \end{equation} and $S=\{n-k+1,n-k,\ldots,n\}$ one checks
that
\begin{equation}  ^\tau(d^{-1} \overline{w}_\sigma^{-1} \varphi_S)=d^{-1}
\overline{w}_\sigma^{-1} \varphi_S. \end{equation} Thus $d^{-1}
\overline{w}_\sigma^{-1} \varphi_S $ is of type $B$.

We shall show that all other representatives are of type $A$: We
will prove that in all the other cases one can find $n_1,n_2 \in
\zspn$ such that  \begin{equation}
\label{typea1}\overline{w}_\sigma
dn_1d^{-1}\overline{w}_\sigma^{-1}=\varphi_S
n_2^{-1}\varphi_S^{-1}, \end{equation} and \begin{equation}
\label{typea2} \psi(n_1n_2) \neq 1.\end{equation} We shall use the
following notations and facts: Denote by $E_{p,q}$ the $n \times
n$ matrix defined by $(E_{p,q})_{i,j}=\delta_{p,i} \delta_{q,j}$.
For $i,j \in \{1,2,\ldots,n \},\, i\neq j$ we define the root
subgroups of $\spn$:
\begin{eqnarray} && U_{i,j}=\{ u_{i,j}(x)=\begin{pmatrix}
  I_n+xE_{i,j} &0 \\
  0 & I_n-xE_{j,i}
\end{pmatrix} \mid x\in F \} \simeq F,\\ &&V_{i,j}=\{ v_{i,j}(x)=\begin{pmatrix}
 I_n &xE_{i,j}+xE_{j,i} \\
  0 & I_n \end{pmatrix} \mid x\in F \} \simeq F, \\ && V_{i,i}=\{ v_{i,i}(x)=\begin{pmatrix}
 I_n &xE_{i,i} \\
  0 & I_n \end{pmatrix} \mid x\in F \} \simeq F.\end{eqnarray}
if $i<j$ we call $U_{i,j}$ a positive root subgroup. If $j=i+1$ we
call $U_{i,j}$ a simple root subgroup, and if $j>i+1$ we call
$U_{i,j}$ a non-simple root subgroup. We call
$U_{i,j}^{t}=U_{j,i}$ the negative of $U_{i,j}$. $S_n$ acts on the
set $\bigr\{U_{i,j} \mid \, i,j \in \{1,2,\ldots,n \},\, i\neq j
\bigl\}$ by

\begin{equation}  \label{sigma on u} \overline{w}_\sigma u_{i,j}(x)
\overline{w}_\sigma^{-1}=u_{\sigma(i),\sigma(j)}(x),\end{equation}
and on the set $\bigl\{V_{i,j} \mid \, i,j \in \{1,2,\ldots,n \}
\bigr\}$ by
\begin{equation} \label{sigma on v} \overline{w}_\sigma v_{i,j}(x)
\overline{w}_\sigma^{-1}=v_{\sigma(i),\sigma(j)}(x).\end{equation}
$T$ acts on each root subgroup via rational characters:
\begin{equation} \label{t on u}
du_{i,j}(x)d^{-1}=u_{i,j}(xd_id_j^{-1}),
\end{equation} and \begin{equation} \label{t on v}
dv_{i,j}(x)d^{-1}=v_{i,j}(xd_id_j), \end{equation} where
$d=diag(d_1,d_2,\ldots,d_n,d_1^{-1},d_2^{-1},\ldots,d_n^{-1})$. We
also note the following : If $i \in S$ then
\begin{equation} \label{var on vii} \varphi_S
v_{i,i}(x)\varphi_S^{-1}=v_{i,i}(x),
\end{equation} if $i\in S, j \notin S$ then \begin{equation} \label{var on
vij} \varphi_S v_{i,j}(x)\varphi_S^{-1}=u_{i,j}(x), \end{equation}
if $i \in S, j \in S, i \neq j$ then \begin{equation} \label{var
on uijs} \varphi_S u_{i,j}(x)\varphi_S^{-1}=u_{i,j}(x)
\end{equation} and if $i \notin S, j \notin S, i \neq j$ then
\begin{equation} \label{var on uij ns} \varphi_S
u_{i,j}(x)\varphi_S^{-1}=u_{j,i}(-x)={^t u_{i,j}(x)}^{-1}.
\end{equation}

Consider the representative
$d^{-1}\overline{w}_\sigma^{-1}\varphi_S$. Assume first that $S$
is empty. Suppose that there exists a simple root subgroup
$U_{k,k+1}$ that $\sigma$ takes to the negative of a non-simple
root subgroup. Then we choose $n_1=u_{k,k+1}(x), \, n_2={^t
u_{\sigma(k),\sigma(k+1)}(d_kd_{k+1}^{-1}x)}.$ For such a choice,
by \eqref{sigma on u}, \eqref{t on u} and \eqref{var on uij ns},
\eqref{typea1} holds. Also, since
$\psi(n_1n_2)=\psi(n_1)=\psi(x),$ it is possible by choosing $x$
properly, to satisfy \eqref{typea2}. Suppose now that there exists
a non-simple positive root subgroup $U_{i,j}$ that $\sigma$ takes
to the negative of a simple root subgroup. Then, we choose
$n_1=u_{i,j}(x), \, n_2={^t
u_{\sigma(i),\sigma(j)}(d_id_{j}^{-1}x)},$ and repeat the previous
argument. Thus, $d^{-1}\overline{w}_\sigma^{-1}\varphi_\emptyset$
is of type A unless $\sigma$ has the following two properties: 1)
If $\sigma$ takes a simple root to a negative root, then it is
taken to the negative of a simple root. 2) If $\sigma$ takes a
non-simple positive root to a negative root, then it is taken to
the negative of a non-simple root. An easy argument shows that if
$\sigma$ has these two properties, $\overline{w}_\sigma$ must be
as in \eqref{wtypeb} with $k=0$. We assume now that
$\overline{w}_\sigma$ has this form. To finish the case
$S=\emptyset$ we show that unless $d$ has the form \eqref{dtypeb},
with $k=0$ and $k_1,k_2,\ldots, k_p$ corresponding to
$\overline{w}_\sigma$, then
$d^{-1}\overline{w}_\sigma^{-1}\varphi_\emptyset$ is of type A.
Indeed, suppose that there exist $d_k$ and $d_{k+1}$ that belong
to the same block in $\overline{w}_\sigma$, such that $d_k \neq
-d_{k+1}.$ Then, we choose $n_1=u_{k,k+1}(x), \, n_2=
{u_{\sigma(k),\sigma(k+1)}(d_k d_{k+1}^{-1}x)}^t.$ For such a
choice, as before, \eqref{typea1} holds, and $\psi(n_1n_2)=\psi
\bigl (x(1+d_k d_{k+1}^{-1}) \bigr).$ Therefore, it is possible,
by choosing $x$ properly, to satisfy \eqref{typea2}.

We may now assume $\mid \! S \! \mid \geq 1$. We show that if $n
\notin S$ then $d^{-1}\overline{w}_\sigma^{-1}\varphi_S$ is of
type A. Indeed, if $\sigma(n) \in S$, in particular $\sigma(n)
\neq n$, then for all $x\in \F$, if we choose $n_1=v_{n,n}(x),\,
n_2=v_{\sigma(n),\sigma(n)}(-xd_n^2)$, by \eqref{sigma on v},
\eqref{t on u}, and \eqref{var on vii}, \eqref{typea1} holds.
Clearly, there exists $x \in \F$ such that
$\psi(n_1n_2)=\psi(n_1)=\psi(x) \neq 1.$ Suppose now that $n
\notin S \neq \emptyset$, and that $\sigma(n) \notin S$. In this
case we can find $1 \leq k \leq n-1$ , such that
\begin{equation} \label{k for nin s} \sigma(k) \in S, \ and \
\sigma(k+1) \notin S.
\end{equation} We choose
\begin{equation} \label{n1n2 for nin
s}n_1=u_{k,k+1}(x), \, n_2=v_{\sigma(k),\sigma(k+1)}(-xd_k d
_{k+1}^{-1}).\end{equation} By \eqref{sigma on u}, \eqref{t on u}
and \eqref{var on vij}, \eqref{typea1} holds, and since
$\psi(n_1n_2)=\psi(n_1)=\psi(x),$ we can satisfy \eqref{typea2} by
properly choosing $x$. We assume now $n \in S$. We also assume
$\sigma(n) \in S$, otherwise we use the last argument. Fix
$n_1=v_{n,n}(x), \, n_2=v{\sigma(n),\sigma(n)}(-xd_n^2).$ One can
check that by \eqref{sigma on v}, \eqref{t on v}, and \eqref{var
on vii}, \eqref{typea1}
holds. Note that $$\psi(n_1n_2)=\begin{cases} x& \sigma(n)\neq n \\
x(d_n^2-1) & \sigma(n)=n \end{cases},$$ hence unless $\sigma(n)=n,
d_n=\pm 1$, $d^{-1}\overline{w}_\sigma^{-1}\varphi_S$ is of type
A. We now assume $\sigma(n)=n, d_n=\pm 1$. If $S=\{n \}$, we use a
similar argument to the one we used for $S=\emptyset$, analyzing
the action of $\sigma$ on $\{1,2,3,\ldots,n-1 \}$ this time.

We are left with the case $\sigma(n)=n, d_n=\pm 1, \, S\supsetneq
\{n \}$. If $\sigma(n-1) \notin S$ we repeat an argument we used
already: We choose $1 \leq k \leq n-2,\, n_1,\, n_2$ as in
\eqref{k for nin s} and \eqref{n1n2 for nin s}. We now  assume
$\sigma(n-1) \in S$. We choose $n_1=u_{n-1,n}(x),
\,n_2=u_{\sigma(n-1),n}(-xd_{n-1} d_n^{-1}).$ Using \eqref{sigma
on u}, \eqref{t on u} and \eqref{var on uijs} we observe that
\eqref{typea1} holds. Also,
$$\psi(n_1n_2)=\begin{cases} x& \sigma(n-1)\neq n-1 \\
x(d_{n-1}d_n^{-1}-1) & \sigma(n-1)=n-1 \end{cases}.$$ Thus, unless
$\sigma(n-1)=n-1$ and $d_{n-1}=d_n=\pm 1$,
$d^{-1}\overline{w}_\sigma^{-1}\varphi_S$ is of type A. Therefore,
we should only consider the case $\sigma(n)=n, \,
\sigma(n-1)=n-1,\,  d_{n-1}=d_n=\pm 1, \, \{n-1,n\} \subseteq S$.
We continue in the same course: If $S=\{n-1,n \}$ we use similar
argument we used for $S=\emptyset$, analyzing the action of
$\sigma$ on $\{1,2,3,\ldots,n-2 \}$. If $S\supsetneq \{n-1,n \}$
we show that unless $\sigma(n-2)=n-2 \in S$ and
$d_{n-2}=d_{n-1}=d_n=\pm 1$ we are in type A etc`.
\end{proof}
We now complete the proof of Theorem \ref{unproven}. We Define
types $\overline{A}$ and $\overline{B}$ for $\mspn$ by analogy
with the definitions given in Lemma \ref{unproven 2} and show that
each element of $$_{\overline{\zspn}} \diagdown ^{\mspn} \diagup
_{\overline{\zspn}}$$ is either of type $\overline{A}$ or of type
$\overline{B}$. Given $\overline{g}=(g,\epsilon) \in \mspn$, if
$g$ is of type A then there are $n_1,n_2 \in \zspn$ such that
$n_1gn_2=g$ and $\psi(n_1n_2) \neq 1$. Let
$\overline{n_i}=(n_i,1)$. Clearly
$\overline{n_1}\overline{g}\overline{n_2}=\overline{g}$ and
$\psi(\overline{n_1}\overline{n_2})=\psi(n_1n_2) \neq 1$. If $g$
is not of type A then by Lemma \ref {unproven 2} it is of type B:
There are $n_1,n_2 \in \zspn$ such that $n_1gn_2={^\tau g}$ and
$\psi(n_1n_2)=1$. Define $\overline{n_i}$ as before. Note that
$\psi(\overline{n_1}\overline{n_2})=\psi(n_1n_2)=1$. From Lemma
\ref{simple lift} it follows that $\overline{n_1} \overline{g}
\overline{n_2}= ^{\overline{\tau}} \overline{g}$. This proves
Lemma \ref{unproven 2} for $\mspn$, which is Theorem
\ref{unproven}.
\newpage
\section{Definition of the local coefficient and of $\gamma(\msig \times \tau,s,
\psi$)} \label{chapter def lc} Unless otherwise is mentioned,
through this chapter $\F$ will denote a p-adic field. We shall
define here the metaplectic analog of the local coefficients
defined by Shahidi in Theorem 3.1 of \cite{Sha 1} for connected
reductive quasi split algebraic groups. Let $n_1,n_2,\ldots,n_r,k$
be $r+1$ nonnegative integers whose sum is $n$. Put
$\overrightarrow{t}=(n_1,n_2,\ldots,n_r;k)$. Let
$(\tau_1,V_{\tau_1}),
(\tau_2,V_{\tau_2}),\ldots,(\tau_r,V_{\tau_r})$ be $r$ irreducible
admissible generic representations of
$GL_{n_1}(\F),GL_{n_2}(\F),\ldots,GL_{n_r}(\F)$ respectively. It
is clear that for $s_i \in \C$, ${\tau_i}_{(s_i)}$ is also
generic. In fact, if $\lambda$ is a $\psi$-Whittaker functional on
$(\tau_i,V_{\tau_i})$ it is also a $\psi$-Whittaker functional on
$({\tau_i}_{(s_i)},V_{\tau_i})$. Let
$(\overline{\sigma},V_{\overline{\sigma}})$ be an irreducible
 admissible $\psi$-generic genuine representation of
$\overline{Sp_{2k}(\F)}.$ Let
$$I\bigl(\pi(\overrightarrow{s})\bigr)=I({\tau_1}_{(s_1)},{\tau_2}_{(s_2)},\ldots,{\tau_r}_{(s_r)},\msig)$$
be the parabolic induction defined in Section \ref{gen par ind}.
Since the inducing representations are generic, then, by a theorem
of Rodier, \cite{Rod}, extended to a non algebraic setting in
\cite{Ban}, $I\bigl(\pi(\overrightarrow{s})\bigr)$ has a unique
$\psi$-Whittaker functional. Define $\lambda_{\tau_1, \psi},
\lambda_{\tau_2, \psi}, \ldots, \lambda_{\tau_r, \psi}$ to be
non-trivial $\psi$-Whittaker functionals on
$V_{\tau_1},V_{\tau_2},\ldots,V_{\tau_r}$ respectively and fix
$\lambda_{\overline{\sigma},\psi}$, a non-trivial $\psi$-Whittaker
functional on $V_{\overline{\sigma}}$. Define
$$\epsilon(\vt)=j_{n-k}(diag(\epsilon_{n_1},\epsilon_{n_2},\ldots,\epsilon_{n_r},
\epsilon_{n_1},\epsilon_{n_2},\ldots,\epsilon_{n_r}),$$ where as
in Section \ref{proof of u of w},  $\epsilon_n=diag
\bigl(1,-1,1\ldots,(-1)^{n+1}\bigr) \in \gl.$ We fix $J_{2n}$ as a
representative of the long Weyl element of $\spn$ and
$$\omega_n=\begin{pmatrix}
_{ } & _{ } & _{ } & _{1} \\
_{ } & _{ } & _{1} & _{ } \\
_{ } & _{\upddots} & _{ } & _{ } \\
_{1} & _{ } & _{ } & _{ }
\end{pmatrix}$$
as a representative of the long Weyl element of $\gln$. We now fix
$$w_l(\overrightarrow{t})=
j_{n-k,n}\bigl(diag(\epsilon_{n_1}\omega_{n_1},\epsilon_{n_2}\omega_{n_2},\ldots,\epsilon_{n_r}\omega_{n_r},
\epsilon_{n_1}\omega_{n_1},\epsilon_{n_2}\omega_{n_2},\ldots,\epsilon_{n_r}\omega_{n_r})\bigr)i_{k,n}(-J_{2k})$$
as a representative of the long Weyl element of $\pt$. We also
define
\begin{equation} \label{good long rep}
w_l'(\overrightarrow{t})=w_l(\overrightarrow{t})J_{2n}.
\end{equation}
Note that $w_l'(\vt)$ is a representative of minimal length of the
longest Weyl element of $\mspn$ modulo the Weyl group of $\mt$. It
maps the positive roots outside $\mt$ to negative roots and maps
the positive roots of $\mt$ to positive roots. The presence of
$\epsilon(\vt)$ in the definition of $w_l(\vt)$ and $w_l'(\vt)$ is
to ensure that \begin{equation} \label{why epsilon}
\psi\Bigl(\bigl(w_l'(\vt),1\bigr)^{-1}n
\bigl(w_l'(\vt),1\bigr)\Bigr)=\psi(n)\end{equation} for all $n \in
\overline{\mt^{w_l'(\vt)}} \cap \overline{ \zspn}$ (the reader may
verify this fact by \eqref{gpg-1} and by a matrix multiplication).
Consider the integral
\begin{equation} \label{whittake1}\lim_{r \rightarrow \infty}
\int_{N_{\overrightarrow{t}}(\Pf^{-r})^{w}}\Bigl((\otimes_{i=1} ^r
\lambda_{\tau_i, \psi^{-1}})\otimes
\lambda_{\overline{\sigma},\psi}\Bigr) \Bigl(f\bigl(
w_l'(\overrightarrow{t}),1\bigr)^{-1}(n,1)\Bigr) \psi^{-1}(n) \,
dn.
\end{equation} By abuse of notations we shall write
\begin{equation} \label{whi on ind}
\int_{\nt^w}\Bigl((\otimes_{i=1} ^r \lambda_{\tau_i,
\psi^{-1}})\otimes \lambda_{\overline{\sigma},\psi}\Bigr)
\Bigl(f\bigl( w_l'(\overrightarrow{t}),1\bigr)^{-1}(n,1)\Bigr)
\psi^{-1}(n) \, dn.
\end{equation} Since $\zspn$ splits in $\mspn$ via the trivial section
the integral converges exactly as in the algebraic case, see
Proposition 3.1 of \cite{Sha 1} (and see Chapter 4 of
\cite{BanPhD} for the context of a p-adic covering group). In
fact, it converges absolutely in an open subset of $\C^r$. Due to
\eqref{pgp-1}, it defines, again as in the algebraic case, a
$\psi$-Whittaker functional on
$I({\tau_1}_{(s_1)},{\tau_2}_{(s_2)},\ldots,{\tau_r}_{(s_r)},\msig)$.
We denote this functional by
$$\lambda \bigl( \overrightarrow{s}, (\otimes_{i=1}^r
\tau_i)\otimes \overline{\sigma},\psi \bigr),$$ where
$\overrightarrow{s}=(s_1,s_2,\ldots,s_r)$. Also, since the
integral defined in \eqref{whittake1} is stable for a large enough
$r$ it defines a rational function in
$q^{s_1},q^{s_2},\ldots,q^{s_r}$.

{\bf Remark}: One can check that $w_l(\vt)\mt w_l(\vt)=\mt$. Thus,
in \eqref{whittake1} and \eqref{whi on ind} we could have written
$\nt$ instead of $\nt^w$. However, in other groups this does not
always happen and we find it appropriate to describe this
construction so that it will fit the general case (see for example a
$\glm$ maximal parabolic case in Section \ref{cha gam}).

Fix now $w\in W_\pt$. Let $t_{w}$ be the unique diagonal element
in $\mt \cap \wspntag$ whose first entry in each block is $1$ such
that  \begin{equation} \label{tw def} \psi \bigl({(w
t_{w},1)}^{-1} n (w t_{w},1) \bigr)=\psi(n) \end{equation} for
each $n \in \overline{\mt^{w}} \cap \overline{\zspn}$.
 \eqref{tw def} assures that  $$\lambda  \bigl(
\overrightarrow{s}^w, (\otimes_{i=1}^r \tau_i)^{w}\otimes
\overline{\sigma},\psi \bigr) \circ A_{w}$$ is another
$\psi$-Whittaker functional on $I\bigl(
\pi(\overrightarrow{s})\bigr)$, (here $A_{w}$ is the intertwining
operator defined in Section \ref{sec io def}). It now follows from
the uniqueness of Whittaker functional that there exists a complex
number
$$C_{\psi}^{\mspn}(\mpt,\overrightarrow{s},(\otimes_{i=1}^r
\tau_i)\otimes \overline{\sigma},w \bigl)$$ defined by the
property

\begin{equation} \label{local c def} \lambda  \bigl( \overrightarrow{s}, (\otimes_{i=1}^r
\tau_i)\otimes \overline{\sigma},\psi
\bigr)=C_{\psi}^{\mspn}(\mpt,\overrightarrow{s},(\otimes_{i=1}^r
\tau_i)\otimes \overline{\sigma},w \bigl)\lambda  \bigl(
\overrightarrow{s}^w, (\otimes_{i=1}^r \tau_i)^{w}\otimes
\overline{\sigma},\psi \bigr) \circ A_{w}. \end{equation} It is
called the {\bf local coefficient} and it clearly depends only on
$\overrightarrow{s}, w$ and the isomorphism classes of
$\tau_1,\tau_2,\ldots,\tau_r$ (not on a realization of the
inducing representations nor on $\lambda_{\tau_1, \psi},
\lambda_{\tau_2, \psi}, \ldots, \lambda_{\tau_r, \psi}$
$\psi^{-1},\lambda_{\overline{\sigma},\psi}$). Also it is clear by
the above remarks that $$\overrightarrow{s} \mapsto
C_{\psi}^{\mspn}(\mpt,\overrightarrow{s},(\otimes_{i=1}^r
\tau_i)\otimes \overline{\sigma},w \bigl)$$ defines a rational
function in function in $q^{s_1},q^{s_2},\ldots,q^{s_r}$. Note
that \eqref{local c def} implies that the zeros of the local
coefficient are among the poles of the intertwining operator.

{\bf Remark}: Assume that the residue characteristic of $\F$ is
odd. Then, by Lemma \ref{k split} $\kappa_{2n}(w)=(wt_w,1)$ for
all $w \in W_\pt$ and
$\kappa_{2n}\bigl(w_l'(\overrightarrow{t})\bigr)=\bigl(w_l'(\overrightarrow{t}),1
\bigr ).$ Keeping the Adelic context in mind, whenever one
introduces local integrals that contain a pre-image of $w \in
\wspntag \subset \Kn$ inside $\mspn$ one should use elements of
the form $\kappa_{2n}(w)$.

Let $\overrightarrow{t}=(n_1,n_2,\ldots,n_r;k)$ where
$n_1,n_2,\ldots,n_r,k$ are $r+1$ non-negative integers whose sum
is $n$. For each $1 \leq i \leq r$ let $\tau_i$ be an irreducible
admissible generic representation of $GL_{n_i}(\F)$ and let
$\msig$ be an irreducible admissible generic genuine
representation of $\mspk$. From the definition of the local
coefficients it follows that
\begin{eqnarray} \label{beta is lc} && \lefteqn{\beta(\overrightarrow{s},\tau_1,\ldots,\tau_n,\msig,w_0)}
\\ \nonumber && = C_{\psi}^{\mspn}(\mpt,\overrightarrow{s},(\otimes_{i=1}^r
\tau_i)\otimes \overline{\sigma},w_0 \bigl)
C_{\psi}^{\mspn}(\mpt,\overrightarrow{s^{w_0}},((\otimes_{i=1}^r
\tau_i)\otimes \overline{\sigma})^{w_0},w_0^{-1} \bigl),
\end{eqnarray}
where
$\beta(\overrightarrow{s},\tau_1,\ldots,\tau_n,\msig,w_0^{-1})$ is
the function defined in \eqref{def beta}. Recalling Theorem
\ref{silberger}, the significance of the local coefficients for
questions of irreducibility of a parabolic induction is clear. In
Chapter \ref{Irreducibility theorems} we shall compute
$\beta(\overrightarrow{s})$ in various cases via the computations
of the local coefficients.

Let $\overline{\sigma}$ be an irreducible admissible generic
irreducible admissible generic genuine representation of $\mspk$.
Let $\tau$ be an irreducible admissible $\psi$-generic
representation of $\glm$. Put $n=m+k$. We define:
\begin{equation} \label{gama def} \gamma(\msig \times \tau,s,
\psi)=\frac {C_{\psi}^{\mspn}
\bigl(\overline{P_{m;k}(\F)},s,\tau\otimes \msig,
j_{m,n}(\omega_m'^{-1}) \bigr)} {C_{\psi}^{\mspm}
\bigl(\overline{P_{m;0}(\F)},s,\tau, \omega_m'^{-1}  \bigr )},
\end{equation} where $\omega_m'=\begin{pmatrix}
& _{ }& _{-\omega_m } \\
& _{\omega_m} & _{ }\end{pmatrix}$. It is clearly a rational
function in $q^s$. Note that if $k=0$ then $\msig$ is the
non-trivial character of the group of 2 elements and $\gamma(\msig
\times \tau,s, \psi)=1$.

This definition of the $\gamma$-factor is an exact analog to the
definition given in Section 6 of \cite{Sha 3} for quasi-split
connected algebraic groups over a non-archimedean field. We note
that similar definitions hold for the case $\F=\R$. In this case
the local coefficients are meromorphic functions.
\newpage

\section{Basic properties of $\gamma(\overline{\sigma} \times \tau,s,\psi)$} \label{basic gam chapter}
 In this chapter $\F$ can be any of the local
fields in discussion. Most of this chapter is devoted to the proof
of a certain multiplicativity property of the $\gamma$-factor; see
Theorems \ref{mult gama tau} and \ref{mult gama sigma} of Section
\ref {cha gam}. This properly in an analog of Part 3 of Theorem
3.15 of \cite {Sha 3}. Aside from technicalities, the proof of the
multiplicativity property follows from a decomposition of the
intertwining operators; see Lemma \ref{lem heart}. We remark that
the proof of this property for $\gamma$-factors defined via the
Rankin-Selberg integrals is harder; see Chapter 11 of \cite{So93}
for example. In Section \ref{prin gama sec}, as an immediate
corollary of the multiplicativity of the $\gamma$-factor, we
compute $\gamma(\overline{\sigma} \times \tau,s,\psi)$ in the case
where $\msig$ and $\tau$ are principal series representations; see
Corollary \eqref{princcomp}.

\subsection{Multiplicativity of $\gamma(\overline{\sigma} \times \tau,s,\psi)$}
\label{cha gam} Let $\overline{\sigma}$ be a genuine irreducible admissible
$\psi$-generic representation of $\mspk$. Let $\tau$ be an
irreducible admissible $\psi$-generic representation of $\glm$.
For two nonnegative integers $l,r$ such that $l+r=m$ denote by
$P^0_{l,r}(\F)$ the standard parabolic subgroup of $\glm$ whose
Levy part is $$M^0_{l,r}(\F)=\begin{pmatrix}
& _{\gll}& _{ } \\
& _{ } & _{\glr}\end{pmatrix}.$$ Denote its unipotent radical by
$N^0_{l,r}(\F)$. Let $\tau_l, \tau_r$ be two irreducible
admissible $\psi$-generic representations of $\gll$ and $\glr$
respectively. In the p-adic case; see \cite{Sha84}, Shahidi
defines
\begin{equation} \label{gama def gl} \gamma(\tau_l \times \tau_r
\,s,\psi)= \chi_{\tau_r}(-I_r)^l
C_{\psi}^{\glm}\bigl(P_{l,r}^0(\F),(\frac{s}{2},\frac{-s}{2}),\tau_l
\otimes \widehat{\tau_r}, \varpi_{r,l} \bigr),  \end{equation}
where $\varpi_{r,l}=\begin{pmatrix}
& _{ }& _{I_r} \\
& _{I_l} & _{ } \end{pmatrix}$, $\chi_{\tau_r}$ is the central
character of $\tau_r$ and $C_{\psi}^{\glm}(
\cdot,\cdot,\cdot,\cdot)$, the $\glm$ local coefficient in the
right-hand side defined via a similar construction to the one
presented in Chapter \ref{chapter def lc}. In the same paper the
author proves that the $\gamma$-factor defined that way is the
same arithmetical factor defined by Jacquet, Piatetskii-Shapiro
and Shalika via the Rankin-Selberg method, see \cite{JPS}. Due to
the remark given in the introduction of \cite{Sha85}, see page 974
after Theorem 1, we take \eqref{gama def gl} as a definition in
archimedean fields as well. The archimedean $\gamma$-factor
defined in this way agrees also with the definition given via the
Rankin-Selberg method, see \cite{J09}.

For future use we note the following:
\begin{equation} \label{varpi action}
\varpi_{r,l}^{-1}=\varpi_{r,l}^{t}=\varpi_{l,r}, \, \, \, \,
\varpi_{r,l}M^0_{l,r}(\F)\varpi_{l,r}=M^0_{r,l}(\F),
\end{equation} and
\begin{equation} \label{varpi action on n}\psi(\varpi_{l,r}n\varpi_{r,l})=\psi(n), \end{equation}
for all $n \in \zglm \cap  M^0_{l,r}(\F)$. $\varpi_{r,l}$ is a
representative of the long Weyl element of $\gln$ modulo the long
Weyl element of $M^0_{l,r}(\F)$.

\begin{thm} \label{mult gama tau}
Assume that $\tau=Ind^{\glm}_{P^0_{l,r}(\F)} \tau_l \otimes
\tau_r$,where $\tau_l, \tau_r$ are two irreducible admissible
generic representations of $\gll$ and $\glr$ respectively,
where $l+r=m$, then
\begin{equation} \label{gama mult tau}
 \gamma(\msig \times \tau,s, \psi)=
 \gamma(\msig \times \tau_l,s, \psi) \gamma(\msig \times \tau_r,s,
 \psi). \end{equation}\end{thm}
 \begin{thm} \label{mult gama sigma}
Assume that $\msig=Ind^{\mspk}_{P_{l;r}(\F)} (\gamapsi \otimes
\tau_l) \otimes \overline{\sigma_r}$, where $\tau_l$ is an
irreducible admissible generic representation of $\gll$ and
$\overline{\sigma_r}$ is an irreducible admissible genuine
$\psi$-generic representation of $\mspr$, where $l+r=k$. Let
$\tau$ be an irreducible admissible generic representation of
$\glm$. Then,
\begin{equation} \label{gama mult sigma}
\gamma(\msig \times \tau,s,
\psi)=\chi^l_{\tau}(-I_m)\chi^m_{\tau_l}(-I_l)(-1,-1)_\F^{ml}\gamma(\msigr
\times \tau,s, \psi) \gamma( \widehat{\tau_l} \times \tau,s, \psi)
\gamma(\tau_l\times \tau,s, \psi). \end{equation}
\end{thm}

We start with proving Theorem \ref{mult gama tau}. We proceed
through the following lemmas:
\begin{lem} \label{ind by parts} With notations in Theorem \ref{mult gama tau} we have:
\begin{equation} C_{\psi}^{\mspn}\bigl(\overline{P_{m;k}(\F)},s,\tau \otimes \msig ,
j_{m,n}(\omega_m'^{-1}) \bigr )=C_{\psi}^{\mspn} \bigl
(\overline{P_{l,r;k}(\F)},(s,s),\tau_l \otimes \tau_r \otimes
\msig , j_{m,n}(\omega_m'^{-1})\bigr).
\end{equation} In particular, for $m=n$, that is when $k=0$ and
$\msig$ is the non trivial representation of the group of two
elements, we have \begin{equation} \label{ind by parts p}
C_{\psi}^{\mspm}(\overline{P_{m;0}(\F)},s,\tau,
\omega_m'^{-1})=C_{\psi}^{\mspm}(\overline{P_{l,r;0}(\F)},(s,s),\tau_l
\otimes \tau_r, \omega_m'^{-1}). \end{equation}
\end{lem}
\begin{proof} We find it convenient to assume that the inducing representations
$\tau_l$, $\tau_r$ and $\msig$ are given in their $\psi$-Whittaker
model. In this realization $f \mapsto f(I_l)$, $f \mapsto f(I_r)$
and  $f \mapsto f(I_{2n},1)$ are $\psi$-Whittaker functionals on
$\tau_l$, $\tau_r$ and $\msig$ respectively. We realize the space
on which $\tau$ acts as a space of functions
$$f:\glm \times \gll \times \glr \rightarrow \C$$ which are smooth
from the right in each variable and which satisfies
$$f(\begin{pmatrix}
& _{a}& _{*} \\
& _{ } & _{b}\end{pmatrix}g,nl_o,n'r_0)= \ab \det(a) \ab ^{\frac
{r} {2}} \ab \det(b) \ab ^{\frac {-l} {2}} \psi(n) \psi(n')
f(g,l_0a,r_0b),$$ for all $g \in \glm$, $l_0,a \in \gll$, $r_0,b
\in \glr, \, n \in \zgll, \, n' \in \zglr$. In this realization
$\tau$ acts by right translations of the first argument (see pages
11 and 65 of \cite{So93} for similar realizations). According to
the $\glm$ analog to the construction given in Chapter
\ref{chapter def lc}, i.e., Proposition 3.1 of \cite{Sha 1} for
$\glm$, and due to \eqref{varpi action} and \eqref{varpi action on
n}, a $\psi$-Whittaker functional on $\tau$ is given by
\begin{equation} \label {whit on tau} \lambda_{\tau,\psi}(f)=
\int_{N^0_{r,l}(\F)} f(\varpi_{l,r}n,I_l,I_r) \psi^{-1}(n) \,
dn.\end{equation} We realize the representation space of
$I(\tau_{(s)} ,\msig )$ as a space of functions $$f:\mspn \times
\mspk \times \glm \times \gll \times \glr \rightarrow \C$$ which
are smooth from the right in each variable and which satisfy
\begin{eqnarray} \label{real i1} \lefteqn{f \bigl((j_{m,n}(\widehat{m}),1)
i_{k,n}(h)u s,ny,\begin{pmatrix}
& _{a}& _{*} \\
& _{ } & _{b}\end{pmatrix}g,n'l_o,n''r_0
\bigr)=f(s,yh,gm,l_0a,r_0b)}\\ \nonumber & & \gamapsi \bigl(\det
(m)\bigr) \ab \det(m) \ab ^{\frac {2k+m+1} {2}+s} \ab \det(a) \ab
^{\frac {r} {2}} \ab \det(b) \ab ^{\frac {-l} {2}} \psi(n)
\psi(n') \psi(n''),
\end{eqnarray} for all $s \in \mspn$, $h,y \in \mspk $, $\begin{pmatrix}
& _{a}& _{*} \\
& _{ } & _{b}\end{pmatrix} \in P^0_{l,r}(\F)$, $m,g \in \glm$,
$l_0 \in \gll$, $r_0 \in \glr$, $u \in (N_{m;k}(\F),1), \, n \in
\overline{\zspk}, \, n' \in \zgll, \, n'' \in \zglr$. In this
realization $\mspn$ acts by right translations of the first
argument. Due to \eqref{whit on tau}, we have\begin{equation}
\lambda (s,\tau \otimes \msig,\psi \bigr)(f)=
 \int_{N_{m;k}(\F)}
\lambda_{\tau,\psi} \Bigl( f \bigl(w_l'(m;k)
n,1),(I_{2k},1),I_m,I_l,I_r \bigr) \biggr) \psi^{-1}(n) dn=$$ $$
 \int_{n \in N_{m;k}(\F)} \int_{n' \in N^0_{r,l}(\F)} f
\bigl((j_{m,n}(-\widehat{\epsilon_m}
\omega_m')n,1),(I_{2k},1),\varpi_{l,r}n',I_l,I_r \bigr)
\psi^{-1}(n') \psi^{-1}(n) \, dn' \, dn.
\end{equation}
Note that $\omega_m'=\begin{pmatrix}
& _{ }& _{-I_m} \\
& _{I_m} & _{ }\end{pmatrix} \begin{pmatrix}
& _{\omega_m}& _{ } \\
& _{ } & _{\omega_m}\end{pmatrix}$. Thus, $x(-\widehat{\epsilon_m}
\omega_m')=(-1)^m$. Due to \eqref{real i1} and \eqref{cocycle prop
2}, we observe that for $n \in N_{m;k}(\F), n' \in M^0_{l,r}(\F)$
\begin{eqnarray*} \lefteqn{f \Bigl( \bigl(j_{m,n}(-\widehat{\epsilon_m} \omega_m')
n,1),(I_{2k},1  \bigr),\varpi_{l,r} n',I_l,I_r \Bigr) =}\\ &&
(-1,-1)_\F^{mrl} \gamapsi(-1^{rl}) f
\biggl(\bigl(j_{m,n}(-\widehat{\varpi_{l,r}} \widehat{n'}
\widehat{\epsilon_m} \omega_m')n ,1 \bigr ),(I_{2k},1),I_m,I_l,I_r
\biggr) .\end{eqnarray*} We shall write
$$-\widehat{\varpi_{l,r}} \widehat{n'}
\widehat{\epsilon_m
}\omega_m'=-\widehat{\varpi_{l,r}}\widehat{\epsilon_m} \omega_m'
\widehat{\breve{n'}},$$ where for $g \in \glm$ we define
$$\breve{g}=( \epsilon_m \omega_m )^{-1} \widetilde{g}(\epsilon_m
\omega_m).$$ Since $n \mapsto \breve{n}$ maps $ N^0_{r,l}(\F)$ to
$N^0_{l,r}(\F)$, we get by \eqref{why epsilon} and \eqref{real
i1}:
\begin{eqnarray}\label{fact lambda on i1}  \lefteqn{\lambda
(s,\tau \otimes \msig,\psi)(f)=} \\ \nonumber &&(-1,-1)_\F^{mrl}
\gamapsi(-1^{lr})
  \int_{N_{l,r;k}(\F)} f \bigl(
j_{m,n}( -\widehat{\varpi_{l,r}}\widehat{\epsilon_m} \omega_m')n
,1),(I_{2k},1),I_m,I_l,I_r \bigr) \psi^{-1}(n) \,dn.
\end{eqnarray} (Clearly the change of integration variable does not
require a correction of the measure). We turn now to
$I({\tau_l}_{(s)},{\tau_r}_{(s)},\msig)$. We realize the space of
this representation as a space of functions
$$f:\mspn \times \mspk\times \gll \times \glr \rightarrow \C$$
which are smooth from the right in each variable and which satisfy

\begin{eqnarray} \label{real i2} \lefteqn{f((j_{m,n}\widehat{ \begin{pmatrix}& _{a}& _{ } \\
& _{ } & _{b}\end{pmatrix}},1) i_{k,n}(h)u
s,ny,n'l_o,n''r_0)=f(s,yh,l_0a,r_0b)} \\ \nonumber &&
\gamapsi(\det \begin{pmatrix}& _{a}& _{ } \\
& _{ } & _{b}\end{pmatrix}) \ab \det(\begin{pmatrix}& _{a}& _{ } \\
& _{ } & _{b}\end{pmatrix}) \ab ^{\frac {2k+m+1} {2}+s} \ab
\det(a) \ab ^{\frac {r} {2}} \ab \det(b) \ab ^{\frac {-l} {2}}
\psi(n) \psi(n') \psi(n''),
\end{eqnarray} for all $s \in \mspn$, $h,y \in \mspk $, $\begin{pmatrix}
& _{a}& _{ } \\
& _{ } & _{b}\end{pmatrix} \in M^0_{l,r}(\F)$,  $l_0 \in \gll$,
$r_0 \in \glr$, $u \in \bigl(N_{l,r;k}(\F),1 \bigr), \, n \in
\overline{\zspk}, \, n' \in \zgll, \, n'' \in \zglr$. In this
realization $\mspn$ acts by right translations of the first
argument. Recall that in \eqref{good long rep} we have defined
$w_l'(\vt)$ to be a particular representative of minimal length of
the longest Weyl element of $\mspn$ modulo the Weyl group of
$\mt$. Since
$$w_l'(l,r;k)=j_{m,n}(
-\widehat{g_0}\widehat{\varpi_{l,r}}\widehat{\epsilon_m}
\omega_m'),$$ where $$g_0=\begin{pmatrix}
& _{(-I_l)^r}& _{ } \\
& _{ } & _{I_r}\end{pmatrix},$$ We have: \begin{equation}
\label{f214} \lambda \bigl((s,s),\tau_l \otimes \tau_r \otimes
\msig,\psi \bigr)(f)= \int_{N_{l,r;k}(\F)} f \biggl( \bigl(
j_{m,n}( -\widehat{g_0}\widehat{\varpi_{l,r}}\widehat{\epsilon_m}
\omega_m')n ,1 \bigr),(I_{2k},1),I_l,I_r \biggr) \psi^{-1}(n) \,
dn.\end{equation}   Note that \begin{eqnarray*} \lefteqn{f \biggl(
\bigl( j_{m,n}(
-\widehat{g_0}\widehat{\varpi_{l,r}}\widehat{\epsilon_m}
\omega_m')n ,1 \bigr),(I_{2k},1),I_l,I_r \biggr)} \\ &=& \nonumber
(-1,-1)_\F^{mlr}\gamma_\psi(-1^{lr})f \biggl( \bigl( j_{m,n}(
-\widehat{g_0}\widehat{\varpi_{l,r}}\widehat{\epsilon_m}
\omega_m')n ,1 \bigr),(I_{2k},1),(-I_l)^r,I_r \biggr)\\ &=&
\nonumber
(-1,-1)_\F^{mlr}\gamma_\psi(-1^{lr})\chi_{\tau_l}^r(-I_l)(f
\biggl( \bigl( j_{m,n}(
-\widehat{\varpi_{l,r}}\widehat{\epsilon_m} \omega_m')n ,1
\bigr),(I_{2k},1),I_l,I_r \biggr).\end{eqnarray*} Thus,
\begin{eqnarray} \label{lamda on l,r,0}
\lefteqn{\lambda \bigl((s,s),\tau_l \otimes \tau_r \otimes
\msig,\psi
\bigr)(f)=(-1,-1)_\F^{mlr}\gamma_\psi(-1^{lr})\chi_{\tau_l}^r(-I_l)}
\\ \nonumber && \int_{N_{l,r;k}(\F)}f \biggl( \bigl( j_{m,n}(
-\widehat{\varpi_{l,r}}\widehat{\epsilon_m} \omega_m')n ,1
\bigr),(I_{2k},1),I_l,I_r \biggr) \psi^{-1}(n)
\end{eqnarray}
For $f \in I(\tau_{(s)} , \msig )$ define $$\widetilde{f}:\mspm
\times \mspk \times \gll \times \glr \rightarrow \C$$ by
$\widetilde{f}(s,y,r_0,l_0)=f(s,y,I_m,r_0,l_0)$. The map $f
\mapsto \widetilde{f}$ is an $\mspm$ isomorphism from
$I(\tau_{(s)} , \msig )$ to $I({\tau_l}_{(s)},{\tau_r}_{(s)},
\msig )$. Comparing  \eqref{fact lambda on i1} and \eqref{lamda on
l,r,0}we see that
\begin{equation} \label{lambdas the same up} \lambda
(s,\tau \otimes \msig,\psi)(f)=  \chi_{\tau_l}^r(-I_l)\lambda
\bigl((s,s),\tau_l \otimes \tau_r \otimes \msig,\psi
\bigr)(\widetilde{f}). \end{equation} We now introduce an
intertwining operator
$$A_{j_{m,n}(\omega_m'^{-1})}:I(\tau_{(s)},\msig) \rightarrow
I(\widehat{\tau}_{(-s)},\msig),$$ defined by
$$A_{j_{m,n}(\omega_m'^{-1})}(s,y,m,l_0,r_0)=\int_{j_{m,n}(N_{(m;0)}(\F))} f
\bigl((j_{m,n}({\widehat{\epsilon_m}}\omega_m')n,1)s,y,m,l_0,r_0
\bigr) dn.$$ Note that for $f \in I(\tau_{(s)},\msig)$ we have
$$A_{j_{m,n}(\omega_m'^{-1})}f:\mspm \times \mspk \times \glm \times
\gll \times \glr \rightarrow \C$$ is smooth from the right in each
variable and satisfies
\begin{equation} \label{real Ai1} A_{j_{m,n}(\omega_m'^{-1})}(f)((j_{m,n}(\widehat{m}),1)
i_{k,n}(h)u s,ny,\begin{pmatrix}
& _{a}& _{*} \\
& _{ } & _{b}\end{pmatrix}g,n'l_o,n''r_0)=$$ $$ \gamapsi(\det (m))
\ab \det(m) \ab ^{\frac {2k+m+1} {2}-s} \ab \det(a) \ab ^{\frac
{r} {2}} \ab \det(b) \ab ^{\frac {-l} {2}} \psi(n) \psi(n')
\psi(n'')\psi f(s,yh,g \breve{m},l_0a,r_0b),
\end{equation} for all $s \in \mspn$, $h,y \in \mspk $, $\begin{pmatrix}
& _{a}& _{*} \\
& _{ } & _{b}\end{pmatrix} \in P^0_{l,r}(\F)$, $m,g \in \glm$,
$l_0 \in \gll$, $r_0 \in \glr$, $u \in \bigl(N_{m;k}(\F),1 \bigr),
\, n \in \overline{\zspk}, \, n' \in \zgll, \, n'' \in \zglr$.
Since $\breve{\breve{g}}=g$ and since
$\breve{\varpi_{l,r}}=h_0\varpi_{r,l}$, where
$$h_0=\begin{pmatrix}
& _{(-I_r)^l}& _{ } \\
& _{ } & _{(-I_l)^r}\end{pmatrix},$$ we have, \begin{eqnarray*} &&
\lambda (-s,\breve{\tau} \otimes \msig,\psi \bigr)(f)\\ &&
=\int_{n \in N_{m;k}(\F)} \int_{n' \in N^0_{r,l}(\F)} f
\bigl((j_{m,n}(-\widehat{\epsilon_m} \omega_m')n,1),(I_{2k},1),h_0
h_0\varpi_{l,r}n',I_l,I_r \bigr) \psi^{-1}(n') \psi^{-1}(n) \, dn'
\, dn \\  && = \chi_{\tau_l}^r(-I_l) \chi_{\tau_r}^l(-I_r)
(-1,-1)_\F^{mrl} \gamapsi(-1^{lr})\\ && \int_{n \in N_{m;k}(\F)}
\int_{n' \in N^0_{r,l}(\F)} f
\bigl((j_{m,n}(-\varpi_{r,l}\breve{n'}\widehat{\epsilon_m}
\omega_m')n,1),(I_{2k},1),I_m,I_l,I_r \bigr) \psi^{-1}(n')
\psi^{-1}(n) \, dn' \, dn\\ &&=\chi_{\tau_l}^r(-I_l)
\chi_{\tau_r}^l(-I_r) (-1,-1)_\F^{mrl} \gamapsi(-1^{lr}) \\
&& \int_{n \in N_{m;k}(\F)} \int_{n' \in N^0_{r,l}(\F)} f
\bigl((j_{m,n}(-\varpi_{r,l}\widehat{\epsilon_m}
\omega_m')n'n,1),(I_{2k},1),I_m,I_l,I_r \bigr) \psi^{-1}(n')
\psi^{-1}(n) \, dn' \, dn .\end{eqnarray*} We have shown:
\begin{eqnarray} \label{f1} \lefteqn{\lambda (-s,\breve{\tau}
\otimes \msig,\psi \bigr)(f)= \chi_{\tau_l}^r(-I_l)
\chi_{\tau_r}^l(-I_r) (-1,-1)_\F^{mrl} \gamapsi(-1^{lr})} \\
\nonumber &&
 \int_{n \in N_{r,l,k}(\F)}
 f \bigl((j_{m,n}(-\varpi_{r,l}\widehat{\epsilon_m}
\omega_m')n'n,1),(I_{2k},1),I_m,I_l,I_r \bigr) \psi^{-1}(n) \, dn.
\end{eqnarray} Consider now
$$\widetilde{A}_{j_{m,n}(\omega_m'^{-1})}:I({\tau_l}_{(s)},{\tau_r}_{(s)},\msig) \rightarrow
I(\widehat{\tau_r}_{(-s)},\widehat{\tau_l}_{(-s)},\msig),$$
defined by
$$\widetilde{A}_{j_{m,n}(\omega_m'^{-1})}(s,y,l_0,r_0)=\int_{j_{m,n}(N_{m;0})} f
\bigl((j_{m,n}(\widehat{\epsilon_m}\omega_m'^{-1})n,1)s,y,l,r
\bigr) dn.$$ Note that for $f \in
I({\tau_l}_{(s)},{\tau_r}_{(s)},\msig)$, we have
$$\widetilde{A}_{j_{m,n}(\omega_m'^{-1})}f:\mspn \times \mspk
\times \gll \times \glr \rightarrow \C$$ is smooth from the right
in each variable and satisfies
\begin{eqnarray} \label{real Ai1} \lefteqn{\widetilde{A}_{j_{m,n}(\omega_m'^{-1})}(f)((j_{m,n}\widehat{\begin{pmatrix}
& _{b}& _{ } \\
& _{ } & _{a}\end{pmatrix}},1) i_{k,n}(h)u
s,ny,n'l_o,n''r_0)=f(s,yh,l_0\breve{a},r_0\breve{b})} \\
\nonumber && \gamapsi(\det \begin{pmatrix}
& _{b}& _{ } \\
& _{ } & _{a}\end{pmatrix}) \ab \det(\begin{pmatrix}
& _{b}& _{ } \\
& _{ } & _{a}\end{pmatrix}) \ab ^{\frac {2k+m+1} {2}-s} \ab
\det(a) \ab ^{\frac {r} {2}} \ab \det(b) \ab ^{\frac {-l} {2}}
\psi(n) \psi(n') \psi(n''),
\end{eqnarray} for all $s \in \mspn$, $h,y \in \mspk $, $\begin{pmatrix}
& _{b}& _{ } \\
& _{ } & _{a}\end{pmatrix} \in M^0_{r,l}$,  $l_0 \in \gll$, $r_0
\in \glr$, ${u \in \bigl(N_{r,l;k}(\F),1 \bigr)}, \, n \in
\overline{\zspk}, \, n' \in \zgll, \, n'' \in \zglr$. Similar to
\eqref{lamda on l,r,0} we have:
\begin{equation} \label{lamda on A l,r,0} \lambda
\bigl((-s,-s),\breve{\tau_r} \otimes\breve{ \tau_l},\psi
\bigr)(f)=(-1,-1)_\F^{mlr}\gamma_\psi(-1^{lr})\chi_{\tau_r}^l(-r_l)
$$ $$\int_{N_{r,l;k}(\F)}f \biggl( \bigl( j_{m,n}(
-\widehat{g_0}\widehat{\varpi_{r,l}}\widehat{\epsilon_m}
\omega_m')n ,1 \bigr),(I_{2k},1),I_l,I_r \biggr) \psi^{-1}(n)
\end{equation} For $f \in I(\breve{\tau}_{(-s)} \otimes \msig )$
define $$\widehat{f}:\mspm \times \mspk \times \gll \times \glr
\rightarrow \C$$ by
$$\widehat{f}(s,y,r_0,l_0)=f(s,y,I_m,r_0,l_0).$$ The map $f \mapsto
\widehat{f}$ is an $\mspn$ isomorphism from $I(\breve{\tau}_{(-s)}
, \msig )$ to $I(\breve{\tau_r}_{(-s)},\breve{\tau_l}_{(-s)},
\msig )$. By \eqref{f1} and \eqref{lamda on A l,r,0} we have,
\begin{equation} \label{lambdas the same down} \lambda
(-s,\breve{\tau} \otimes \msig,\psi)(f)=  \chi_{\tau_l}^r(-I_l)
\lambda \bigl((-s,-s),\breve{\tau}_r \otimes \breve{\tau}_l
\otimes \msig,\psi \bigr)(\widehat{f}). \end{equation} We use
\eqref{lambdas the same up}, \eqref{lambdas the same down} and the
fact that for all $f \in I(\tau^s \otimes \msig )$ , we have
$$\widetilde{A}_{j_{m,n}(\omega_m'^{-1})}(\widetilde{f})=\widehat{A_{j_{m,n}(\omega_m'^{-1})}(f)},$$
to complete the lemma:
$$C_{\psi}^{\mspn}\bigl(\overline{P_{m,k;0}(\F)},s,\tau \otimes \msig ,
j_{m,n}(\omega_m'^{-1})\bigr)= \frac {\lambda (s,\tau \otimes
\msig,\psi)(f)} {\lambda (-s,\breve{\tau} \otimes
\msig,\psi)\bigl(A_{j_{m,n}(\omega_m`^{-1})}(f)\bigr)}=$$
$$ \frac {\lambda
\bigl((s,s),\tau_l \otimes \tau_r \otimes \msig,\psi
\bigr)(\widetilde{f})}{\lambda \bigl((-s,-s),\breve{\tau}_r
\otimes \breve{\tau}_l \otimes \msig,\psi
\bigr)\bigl(\widetilde{A}_{j_{m,n}(\omega_m`^{-1})}(\widetilde{f})\bigr)}=
C_{\psi}^{\mspn}\bigl(\overline{P_{l,r;k}(\F)},(s,s),\tau_l
\otimes \tau_r \otimes \msig , j_{m,n}(\omega_m')\bigr).$$
\end{proof}

The heart of the proof of Theorem \ref{mult gama tau} is the
following lemma which is  a slight modification of the following
argument, originally proved for algebraic groups (see \cite {Sha
1} for example): If $w_1, w_2$ are two Weyl elements such that
$l(w_1w_2)=l(w_1)+l(w_2)$, where $l(\cdot)$ is the length function
in the Weyl group, then $A_{w_1}\circ A_{w_2}=A_{w_1w_2}$. See
Lemma 1 of Chapter VII of \cite{BanPhD} for a proof of this
factorization in the case of an $r$-fold cover of $\gln$.
\begin{lem} \label{lem heart}
$$ C_{\psi}^{\mspn}\bigl(\overline{P_{l,r;k}(\F)},(s,s),\tau_l
\otimes \tau_r \otimes \msig ,
j_{m,n}(\omega_m'^{-1})\bigr)=\phi^{-1}_{\psi}(r,l,\tau_r)c_1(s)c_2(s)c_3(s),$$
where

\begin{eqnarray*} &&
c_1(s)=C_{\psi}^{\mspn}\bigl(\overline{P_{l,r;k}(\F)},(s,s),\tau_l
\otimes \tau_r \otimes \msig , j_{m,n}(w_1^{-1})\bigr), \\ &&
c_2(s)=C_{\psi}^{\mspn}\bigl(\overline{P_{l,r;k}(\F)},(s,-s),\tau_l
\otimes \widehat{\tau_r} \otimes \msig , j_{m,n}(w_2^{-1})\bigr),
\\ &&c_3(s)=C_{\psi}^{\mspn}\bigl(\overline{P_{l,r;k}(\F)},(-s,s),
\widehat{\tau_r} \otimes \tau_l \otimes \msig ,
j_{m,n}(w_3^{-1})\bigr)\end{eqnarray*} and where
$$w_1=\begin{pmatrix}
& _{I_l}& _{ } &_{ } &_ { } \\
& _{ } & _{ }&_{ } &_ {-\omega_r } \\ &_{ } &_ { } & _{I_l}& _{ }\\
&_{ } &_ {\omega_r} &_{ } &_{ }
\end{pmatrix},\quad w_2=\widehat{\varpi_{l,r}}, \, \, w_3=\begin{pmatrix}
& _{I_r}& _{ } &_{ } &_ { } \\
& _{ } & _{ }&_{ } &_ {-\omega_l } \\ &_{ } &_ { } & _{I_r}& _{ }\\
&_{ } &_ {\omega_l} &_{ } &_{ }
\end{pmatrix},$$ and $$\phi_{\psi}(r,l,\tau_r)=(-1,-1)_\F^{\frac{l^2(l-1)}{2}} \chi_{\tau_r}(-I_r)^l
\gamapsi \ \bigl((-1)^{rl} \bigr).$$ In particular : $$
C_{\psi}^{\mspm}\bigl(\overline{P_{l,r;0}(\F)},(s,s),\tau_l
\otimes \tau_r  , \omega_m'^{-1}
\bigr)=\phi^{-1}_{\psi}(r,l,\tau_r)c'_1(s)c'_2(s)c'_3(s),$$ where
\begin{eqnarray*}
&&c'_1(s)=
C_{\psi}^{\mspm}\bigl(\overline{P_{l,r;0}(\F)},(s,s),\tau_l
\otimes \tau_r , w_1^{-1}\bigr)
\\
&&
c'_2(s)=C_{\psi}^{\mspm}\bigl(\overline{P_{l,r;0}(\F)},(s,-s),\tau_l
\otimes \widehat{\tau_r}  , w_2^{-1}\bigr) \\
&& c'_3(s)=C_{\psi}^{\mspm}\bigl(\overline{P_{l,r;0}(\F)},(-s,s),
\widehat{\tau_r} \otimes \tau_l  , w_3^{-1}\bigr).\end{eqnarray*}
\end{lem}
\begin{proof}
We keep the realizations used in Lemma \ref{ind by parts} and most
of its notations. Consider the following three intertwining
operators:
\begin{eqnarray} \label{3 io} && A_{j_{m,n}(w_1^{-1})} : I({\tau_l}_{(s)},{\tau_r}_{(-s)},\msig) \rightarrow
I({\tau_l}_{(s)},{\breve{\tau_r}}_{(-s)},\msig)\\ \nonumber &&
A_{j_{m,n}(w_2^{-1})} :
I({\tau_l}_{(s)},{\breve{\tau_r}}_{(-s)},\msig)
 \rightarrow I({\breve{\tau_r}}_{(-s)},{\tau_l}_{(s)},\msig)\\ \nonumber
&& A_{j_{m,n}(w_3^{-1})} :
I({\breve{\tau_r}}_{(-s)},{\tau_l}_{(s)},\msig) \rightarrow
I({\breve{\tau_r}}_{(-s)},{\breve{\tau_l}}_{(-s)},\msig)
\end{eqnarray}
Suppose that we show that \begin{equation} \label{heart}
A_{j_{m,n}(w_3^{-1})} \circ A_{j_{m,n}(w_2^{-1})} \circ
A_{j_{m,n}(w_1^{-1})}=\phi_{\psi}(r,l,\tau_r)A_{j_{m,n}(\omega_m'^{-1})},
\end{equation}
This will finish the lemma at once since

\begin{eqnarray*} && \! \! \! \! \! \! \! C_{\psi}^{\mspn} \bigl(\overline{P_{l,r;k}(\F)},(s,s),\tau_l
\otimes \tau_r \otimes \msig , j_{m,n}(\omega_m')\bigr)= \frac
{\lambda \bigl((s,s),\tau_l \otimes \tau_r \otimes \msig,\psi
\bigr)(f)}{\lambda \bigl((-s,-s),\breve{\tau}_r \otimes
\breve{\tau}_l \otimes \msig,\psi
\bigr)\bigl(A_{j_{m,n}(\omega_m`^{-1})}(f)\bigr)}=\\
&& \! \! \! \! \! \! \! \frac {\lambda \bigl((s,s),\tau_l \otimes
\tau_r \otimes \msig,\psi \bigr)(f)}{\lambda \bigl((s,-s),\tau_l
\otimes \breve{\tau}_r \otimes \msig,\psi
\bigr)\bigl(A_{j_{m,n}(w_1^{-1})}(f)\bigr)} \frac {\lambda
\bigl((s,-s),\tau_l \otimes \breve{\tau}_r \otimes \msig,\psi
\bigr)\bigl(A_{j_{m,n}(w_1^{-1})}(f)\bigr)} {\lambda \bigl((-s,s),
\breve{\tau}_r \otimes  \tau_l\otimes \msig,\psi
\bigr)\bigl(A_{j_{m,n}(w_2^{-1}\bigr)} \circ
A_{j_{m,n}(w_1^{-1})}(f)\bigr)} \\
&& \! \! \! \! \! \! \! \frac {\lambda \bigl((-s,s),
\breve{\tau}_r \otimes  \tau_l\otimes \msig,\psi
\bigr)\bigl(A_{j_{m,n}(w_2^{-1})} \circ A_{j_{m,n}(w_1^{-1})}(f)
\bigr )} {\lambda \bigl((-s,-s),  \breve{\tau}_r \otimes \breve{
\tau}_l\otimes \msig,\psi \bigr)\bigl(A_{j_{m,n}(w_3^{-1})} \circ
A_{j_{m,n}(w_2^{-1})} \circ A_{j_{m,n}(w_1^{-1})}(f)\bigr)}.
\end{eqnarray*}
By definition the right-hand side equals
$$\phi^{-1}_{\psi}(r,l,\tau_r)c_1(s)c_2(s)c_3(s).$$ Thus, we prove
\eqref{heart}. It is sufficient to prove it for $Re(s)>>0$ where
$A_{j_{m,n}(\omega_m'^{-1})}$ is given by an absolutely convergent
integral. Our argument will use Fubini's theorem, whose use will
be justified by \eqref{comp io 2}.
\begin{eqnarray} \label{comp io} && \! \! \! \! \! \! \! \lefteqn{A_{j_{m,n}(w_3^{-1})} \circ A_{j_{m,n}(w_2^{-1})} \circ
A_{j_{m,n}(w_1^{-1})}(f)(s,y,l,r)=}\\
\nonumber && \! \! \! \! \! \! \! \int_{N_{w_3}(\F)}
A_{j_{m,n}(w_2^{-1})} \circ
 A_{j_{m,n}(w_1^{-1})}(f)\bigl(\bigl(j_{m,n}(t_{w_3}w_3)n_3,1 \bigr)s,y,l,r \bigr)
 dn_3=\\
 \nonumber && \! \! \! \! \! \! \! \int_{N_{w_3}(\F)} \int_{N_{w_2}(\F)} A_{j_{m,n}(w_1^{-1})}(f)
 \bigl(\bigl(j_{m,n}(t_{w_2}w_2)n_2,1 \bigr) \bigl(j_{m,n}(t_{w_3}w_3)n_3,1 \bigr)s,y,l,r \bigr)
 dn_2 dn_3=\\
  \nonumber && \! \! \! \! \! \! \! \int_{N_{w_3}(\F)} \! \! \int_{N_{w_2}(\F)}   \! \! \int_{N_{w_1}(\F)}f
 \bigl(\bigl(j_{m,n}(t_{w_1}w_1)n_1,1 \bigr) \bigl(j_{m,n}(t_{w_2}w_2)n_2,1 \bigr)
 \!  \! \bigl(j_{m,n}(t_{w_3}w_3)n_3,1 \bigr)s,y,l,r \bigr)
 dn_1 dn_2 dn_3, \end{eqnarray}
where, as a straight forward computation will show, $N_{w_1}(\F)$
is the group of elements of the form
$$\begin{pmatrix}
&_{I_l} & ~  & ~ &_{0_l} ~  & ~ & ~\\
& ~ &_{I_r}   & ~ & ~  &_{z}  & ~ \\
& ~  & ~ &_{I_k}   & ~ & ~ &_{0_k}\\
& ~  & ~ & ~ &_{I_l} & ~  & ~\\
& ~  & ~  & ~ & ~ &_{I_r} & ~ \\
& ~ & ~  & ~  & ~ & ~ &_{I_k} \end{pmatrix},\quad z \in
Mat^{sym}_{r \times r}(\F), $$ $N_{w_2}(\F)$ is the group of
elements of the form $$\begin{pmatrix}
&_{I_l} & _{z}  & ~ & ~  & ~ & ~\\
& ~ &_{I_r}   & ~ & ~  & ~ & ~ \\
& ~  & ~ &_{I_k}   & ~ & ~ &~\\
& ~  & ~ & ~ &_{I_l} & ~  & ~\\
& ~  & ~  & ~ & _{-z^t} &_{I_r} & ~ \\
& ~ & ~  & ~  & ~ & ~ &_{I_k} \end{pmatrix},\quad z \in Mat_{l
\times r}(\F), $$ $N_{w_3}(\F)$ is the group of elements of the
form
$$\begin{pmatrix}
&_{I_r} & ~  & ~ &_{0_r} ~  & ~ & ~\\
& ~ &_{I_l}   & ~ & ~  &_{z}  & ~ \\
& ~  & ~ &_{I_k}   & ~ & ~ &_{0_k}\\
& ~  & ~ & ~ &_{I_r} & ~  & ~\\
& ~  & ~  & ~ & ~ &_{I_l} & ~ \\
& ~ & ~  & ~  & ~ & ~ &_{I_k} \end{pmatrix},\quad z \in
Mat^{sym}_{l \times l}(\F),$$ and where
$t_{w_1}=\widehat{\begin{pmatrix} &_{I_l} & ~\\ &~ &_{\epsilon_r}
\end{pmatrix}}, \qquad t_{w_2}=I_m,\qquad t_{w_3}=\widehat{\begin{pmatrix}
&_{I_r} & ~\\ &~ &_{\epsilon_l} \end{pmatrix}}. $ We consider the
first argument of $f$ in \eqref{comp io}: By \eqref{cocycle prop
2} and \eqref{gpg-1} we have: \begin{eqnarray} \label{arg for
decomp} && \bigl(j_{m,n}(t_{w_1}w_1)n_1,1 \bigr)
\bigl(j_{m,n}(t_{w_2}w_2)n_2,1 \bigr)
\bigl(j_{m,n}(t_{w_3}w_3)n_3,1 \bigr)\\ \nonumber &&=
\bigl(j_{m,n}(t_{w_1}w_1),1 \bigr)(n_1,1)
\bigl(j_{m,n}(t_{w_2}w_2),1 \bigr) (n_2,1)
\bigl(j_{m,n}(t_{w_3}w_3),1 \bigr)(n_3,1)\\
\nonumber && =\bigl(j_{m,n}(t_{w_1}w_1),1 \bigr)
\bigl(j_{m,n}(t_{w_2}w_2),1 \bigr) \bigl(j_{m,n}(t_{w_3}w_3),1
\bigr)(n_1'n_2'n_3,1)\\ \nonumber &&=(I_{2n},\epsilon)
\bigl(j_{m,n}(t_{w_1}w_1 t_{w_2}w_2 t_{w_3}w_3),1 \bigr)
(n_1'n_2'n_3,1) \\ \nonumber &&
=(I_{2n},\epsilon)\bigl(j_{m,n}\widehat{(\begin{pmatrix} &_{I_l} &
~\\ &~ &_{I_r(-1)^l}
\end{pmatrix}} \widehat{\epsilon_m} \omega_m'),1
\bigr) (n_1'n_2'n_3,1) ,\end{eqnarray}where
$$n_1'=j_{m,n}(w_3t_{w_3}w_2t_{w_2})^{-1}(n_1)j_{m,n}(w_2t_{w_2}w_3t_{w_3})^{-1},
\quad n_2'=j_{m,n}(w_3t_{w_3})^{-1}(n_2)j_{m,n}(w_3t_{w_3}),$$ and
where
$$\epsilon=c(t_{w_1}w_1,t_{w_2}w_2)c(t_{w_1}w_1t_{w_2}w_2,t_{w_3}w_3).$$
We compute:
$$j_{m,n}(w_3t_{w_3}w_2t_{w_2})^{-1}\begin{pmatrix}
&_{I_l} & ~  & ~ &_{0_l} ~  & ~ & ~\\
& ~ &_{I_r}   & ~ & ~  &_{z}  & ~ \\
& ~  & ~ &_{I_k}   & ~ & ~ &_{0_k}\\
& ~  & ~ & ~ &_{I_l} & ~  & ~\\
& ~  & ~  & ~ & ~ &_{I_r} & ~ \\
& ~ & ~  & ~  & ~ & ~ &_{I_k}
\end{pmatrix}j_{m,n}(w_2t_{w_2}w_3t_{w_3})^{-1}=\begin{pmatrix}
&_{I_r} & ~  & ~ &_{z} ~  & ~ & ~\\
& ~ &_{I_l}   & ~ & ~  &~  & ~ \\
& ~  & ~ &_{I_k}   & ~ & ~ &~\\
& ~  & ~ & ~ &_{I_r} & ~  & ~\\
& ~  & ~  & ~ & ~ &_{I_l} & ~ \\
& ~ & ~  & ~  & ~ & ~ &_{I_k} \end{pmatrix},$$

$$j_{m,n}(w_3t_{w_3})^{-1}\begin{pmatrix}
&_{I_l} & _{z}  & ~ & ~  & ~ & ~\\
& ~ &_{I_r}   & ~ & ~  & ~ & ~ \\
& ~  & ~ &_{I_k}   & ~ & ~ &~\\
& ~  & ~ & ~ &_{I_l} & ~  & ~\\
& ~  & ~  & ~ & _{-z^t} &_{I_r} & ~ \\
& ~ & ~  & ~  & ~ & ~ &_{I_k}
\end{pmatrix}j_{m,n}(w_3t_{w_3})=\begin{pmatrix}
&_{I_l} & ~  & ~ & ~ & _{z`}   & ~\\
& ~ &_{I_r}   & ~ & _{z`^t}  & ~ & ~ \\
& ~  & ~ &_{I_k}   & ~ & ~ &~\\
& ~  & ~ & ~ &_{I_l} & ~  & ~\\
& ~  & ~  & ~ & ~ &_{I_r} & ~ \\
& ~ & ~  & ~  & ~ & ~ &_{I_k} \end{pmatrix},$$where
$z`=z\omega_rt_{w_3}.$ Hence we can change the three integrals in
\eqref{comp io} to a single integration on $j_{m,n}
\bigl(N_{m;0}(\F)\bigr)$ without changing the measure and obtain
\begin{eqnarray} \label{comp io 2} && A_{j_{m,n}(w_3^{-1})} \circ
A_{j_{m,n}(w_2^{-1})} \circ A_{j_{m,n}(w_1^{-1})}(f)(s,y,l,r)
\\ \nonumber &&= \epsilon \int_{j_{m,n}(N_{m;0})} f\biggl(\bigl(j_{m,n}\widehat{(\begin{pmatrix}
&_{I_l} & ~\\ &~ &_{I_r(-1)^l}
\end{pmatrix}} \widehat{\epsilon_m} \omega_m')n,1
\bigr)s,y,l,r \biggr) dn \\ \nonumber &&=\epsilon
\chi_{\tau_r}(-I_r)^l \gamapsi \ \bigl((-1)^{rl} \bigr)
A_{j_{m,n}(\omega_m'^{-1})}(f)(s,y,l,r).
\end{eqnarray}
It is left to show that
$\epsilon=(-1,-1)_\F^{\frac{l^2(l-1)}{2}}$. Indeed, we have
$$t_{w_1}w_1=\begin{pmatrix}
& _{I_l}& _{ } &_{ } &_ { } \\
& _{ } & _{\epsilon_r\omega_r}&_{ } &_ { } \\ &_{ } &_ { } & _{I_l}& _{ }\\
&_{ } &_ { } &_{ } &_{\epsilon_r\omega_r }
\end{pmatrix} \begin{pmatrix}
& _{I_l}& _{ } &_{ } &_ { } \\
& _{ } & _{ }&_{ } &_ {-I_r } \\ &_{ } &_ { } & _{I_l}& _{ }\\
&_{ } &_ {I_r} &_{ } &_{ }
\end{pmatrix}, \quad t_{w_2}w_2=\begin{pmatrix}
& _{ }& _{I_l } &_{ } &_ { } \\
& _{I_r} & _{ }&_{ } &_ { } \\ &_{ } &_ { } & _{}& _{I_l }\\
&_{ } &_ { } &_{I_r} &_{ }\end{pmatrix}.$$ Thus, by
$\eqref{cocycle prop 2}$ we have
$c(t_{w_1}w_1,t_{w_2}w_2)=(-1,-1)_\F^{\frac{r^3l(r-1)}{2}}$. Since
\begin{eqnarray*} &&t_{w_1}w_1t_{w_2}w_2=\begin{pmatrix}
& _{I_l}& _{ } &_{ } &_ { } \\
& _{ } & _{\epsilon_r\omega_r}&_{ } &_ { } \\ &_{ } &_ { } & _{I_l}& _{ }\\
&_{ } &_ { } &_{ } &_{\epsilon_r\omega_r }\end{pmatrix}
 \begin{pmatrix}
& _{ }& _{I_l } &_{ } &_ { } \\
& _{I_r} & _{ }&_{ } &_ { } \\ &_{ } &_ { } & _{}& _{I_l }\\
&_{ } &_ { } &_{I_r} &_{ }\end{pmatrix}
\begin{pmatrix}
& _{ }& _{ } &_{-I_r } &_ { } \\
& _{ } & _{I_l}&_{ } &_ { } \\ &_{-I_r} &_ { } & _{}& _{ }\\
&_{ } &_ { } &_{ } &_{I_l}
\end{pmatrix},\\ \nonumber && t_{w_3}w_3=\begin{pmatrix}
& _{I_r}& _{ } &_{ } &_ { } \\
& _{ } & _{ }&_{ } &_ {-I_l } \\ &_{ } &_ { } & _{I_r}& _{ }\\
&_{ } &_ {I_l} &_{ } &_{ }
\end{pmatrix}\begin{pmatrix}
& _{I_r}& _{ } &_{ } &_ { } \\
& _{ } & _{\epsilon_l\omega_l}&_{ } &_ { } \\ &_{ } &_ { } & _{I_r}& _{ }\\
&_{ } &_ { } &_{ } &_{\epsilon_r\omega_r }
\end{pmatrix}, \end{eqnarray*}we conclude, using \eqref{cocycle prop 4}, that
$c(t_{w_1}w_1t_{w_2}w_2,t_{w_3}w_3)=(-1,-1)_\F^{\frac{r^3l(r-1)}{2}+\frac{l^2(l-1)}{2}}.$
\end{proof}

\begin{lem} \label{id with c} Keeping the notations of the previous lemmas we have:
\begin{equation} \label{same 1}C_{\psi}^{\mspn}\bigl(\overline{P_{l,r;k}(\F)},(s,s),\tau_l
\otimes \tau_r \otimes \msig ,
j_{m,n}(w_1^{-1})\bigr)=C_{\psi}^{\overline{Sp_{2(r+k)}(\F)}}\bigl(\overline{P_{r;k}(\F)},s,
\tau_r \otimes \msig , j_{r,r+k}(\omega_r'^{-1})\bigr).
\end{equation}
\begin{equation} \label{same 2}C_{\psi}^{\mspn}\bigl(\overline{P_{l,r;k}(\F)},(s,-s),\tau_l \otimes
\widehat{\tau_r} \otimes \msig ,
j_{m,n}(w_2^{-1})\bigr)=C_{\psi}^{\glm}\bigl(P_{l,r}^0(\F),(s,-s),\tau_l
\otimes \widehat{\tau_r}, \varpi_{l,r}^{-1}\bigr).\end{equation}
\begin{equation} \label{same 3}C_{\psi}^{\mspn}\bigl(\overline{P_{r,l;k}(\F)},(-ss,s),\widehat{\tau}_r \otimes
\tau_l \otimes \msig ,
j_{m,n}(w_3^{-1})\bigr)=C_{\psi}^{\overline{Sp_{2(l+k)}(\F)}}\bigl(\overline{P_{l;k}(\F)},s,
\tau_l \otimes \msig , j_{l,l+k}(\omega_r'^{-1})\bigr).
\end{equation}

In particular:
\begin{eqnarray*} && C_{\psi}^{\mspm}\bigl(\overline{P_{l,r;0}(\F)},(s,s),\tau_l
\otimes \tau_r,
w_1^{-1}\bigr)=C_{\psi}^{\overline{Sp_{2r}(\F)}}\bigl(\overline{P_{r;0}(\F)},s,
\tau_r  , \omega_r'^{-1}\bigr).\\ &&
C_{\psi}^{\mspm}\bigl(\overline{P_{l,r;0}(\F)},(s,-s),\tau_l
\otimes \widehat{\tau_r} ,
w_2^{-1}\bigr)=C_{\psi}^{\glm}\bigl(P_{l,r}^0(\F),(s,-s),\tau_l
\otimes \widehat{\tau_r}, \varpi_{l,r}^{-1}\bigr). \\ &&
C_{\psi}^{\mspm}\bigl(\overline{P_{r,l;0}(\F)},(-s,s),\widehat{\tau}_r
\otimes \tau_l \otimes \msig ,
w_3^{-1}\bigr)=C_{\psi}^{\overline{Sp_{2l}(\F)}}\bigl(\overline{P_{l;0}(\F)},s,
\tau_l , \omega_r'^{-1}\bigr). \end{eqnarray*}
\end{lem}
\begin{proof}
We prove \eqref{same 1} and \eqref{same 2} only. \eqref{same 3} is
proven exactly as \eqref{same 1}. We start with \eqref{same 1}: As
in Lemma \ref{ind by parts} we realize $I({\tau_r}_{(s)},\msig)$
as a space of functions $$f:\overline{Sp_{2(r+k)}(\F)} \times
\mspk \times \glr \rightarrow \C$$ which are smooth from the right
in each variable and which satisfies

\begin{equation} \label{real i10} f(j_{r,r+k}(\widehat{b}) i_{k,n}(h)u s,ny,n'r_0)=
\gamapsi(\det (b)) \ab \det(b) \ab ^{k+r+\half} \psi(n) \psi(n')
\psi f(s,yh,r_0b),
\end{equation} for all $s \in \overline{Sp_{2(r+k)}(\F)}$, $h,y \in \mspk $, $b,r_0 \in \glr$,
$u \in (N_{l;k}(\F),1)$, $n \in \overline{\zspk}$ $n'\in \zglm$.
For $f \in I({\tau_l}_{(s)},{\tau_r}_{(s)},\msig)$, $g \in \mspn$
we define
$$f_g:\overline{Sp_{2(r+k)}(\F)} \times \mspk \times \glr \rightarrow
\C$$ by $$f_g(s,y,r)=f(i_{r+k,n}(g)s,y,I_l,r).$$ Recalling
\eqref{real i2} we note that $f_g \in I({\tau_r}_{(s)},\msig)$. We
want to write the exact relation between $\lambda
\bigl((s,s),\tau_l \otimes \tau_r \otimes \msig,\psi \bigr)$ and
$\lambda \bigl(s,\tau_r \otimes \msig,\psi \bigr)$. To do so we
consider the left argument of $f$ in  \eqref{f214}: We decompose
$n \in N_{l,r;k}(\F)$ as $n=n'n'',$ where $$n' \in
i_{r+k,n}(N_{r,k}(\F)), \, n'' \in U_0(\F)= \Bigg \{
u=\begin{pmatrix}
&_{I_l} & *  & * & *  & * & *\\
& ~ &_{I_r}   & _{0{r \times k}} & *  & _{0{r \times r}} &  _{0{r \times k}} \\
& ~  & ~ &_{I_k}   & * &_{0{k \times r}} &_{0{k \times k}}\\
& ~  & ~ & ~ &_{I_l} & ~  & ~\\
& ~  & ~  & ~ & * &_{I_r} & ~ \\
& ~ & ~  & ~  & * & ~ &_{I_k}
\end{pmatrix} \mid u \in \spn \Bigg \}. $$  We have \begin{eqnarray*}
&& \lefteqn{( w_l'(l,r,k)n,1)} \\ && =\bigr(\begin{pmatrix}
&_{I_l} & ~  & ~ & ~  & ~ & ~\\
& ~ & ~   & ~ & ~  & _{\epsilon_r\omega_r} & ~ \\
& ~  & ~ &_{I_k}   & ~ & ~ &~\\
& ~  & ~ & ~ &_{I_l} & ~  & ~\\
& ~  & _{-\epsilon_r\omega_r}  & ~ & ~ & ~ & ~ \\
& ~ & ~  & ~  & ~ & ~ &_{I_k}
\end{pmatrix} \begin{pmatrix}
&~ & ~  & ~ & _{\epsilon_l\omega_l} & ~ & ~\\
& ~ & _{I_r}   & ~ & ~  & ~ & ~ \\
& ~  & ~ &_{I_k}   & ~ & ~ &~\\
& _{-\epsilon_l\omega_l}  & ~ & ~ & ~ & ~  & ~\\
& ~  & ~  & ~ & ~ & _{I_r} & ~ \\
& ~ & ~  & ~  & ~ & ~ &_{I_k}
\end{pmatrix}n'n'',1 \bigl) \\ && =\biggl(i_{r+k,n}\bigr(j_{r,r+k}(-\widehat{\epsilon_r}\omega_r')\bigr)n'\begin{pmatrix}
&~ & ~  & ~ & _{\epsilon_l\omega_l} & ~ & ~\\
& ~ & _{I_r}   & ~ & ~  & ~ & ~ \\
& ~  & ~ &_{I_k}   & ~ & ~ &~\\
& _{-\epsilon_l\omega_l}  & ~ & ~ & ~ & ~  & ~\\
& ~  & ~  & ~ & ~ & _{I_r} & ~ \\
& ~ & ~  & ~  & ~ & ~ &_{I_k}
\end{pmatrix}n'',1 \biggr)\\ \end{eqnarray*}
Thus, for an appropriate $\epsilon$ independent of $n'$ and $n''$,
we have:

\begin{eqnarray*}
&& \lefteqn{( w_l'(l,r,k)n,1)=} \\ &&
\biggl(i_{r+k,n}\bigr(j_{r,r+k}(-\widehat{\epsilon_r}\omega_r')\bigr)n',1
\biggr) \biggl(\begin{pmatrix}
&~ & ~  & ~ & _{\epsilon_l\omega_l} & ~ & ~\\
& ~ & _{I_r}   & ~ & ~  & ~ & ~ \\
& ~  & ~ &_{I_k}   & ~ & ~ &~\\
& _{-\epsilon_l\omega_l}  & ~ & ~ & ~ & ~  & ~\\
& ~  & ~  & ~ & ~ & _{I_r} & ~ \\
& ~ & ~  & ~  & ~ & ~ &_{I_k}
\end{pmatrix}n'',\varepsilon \biggr).\end{eqnarray*}
We denote the right element of the last line by $g(n'')$, and we
see that
\begin{eqnarray} && \lefteqn{\label{reduction up 1}\lambda \bigl((s,s),\tau_l
\otimes \tau_r \otimes \msig,\psi \bigr)(f)} \\ \nonumber &&=
\int_{U_0(\F)} \int_{N_{r;k}(\F)}
f_{g(n'')}\bigl(j_{r,r+k}(-\widehat{\epsilon_r}\omega_r')n',(I_{2k},1),I_r
\bigr) \psi^{-1}(n') \psi^{-1}(n'') \, dn' \, dn''\\ \nonumber
&&=\int_{U_0(\F)} \lambda \bigl(s, \tau_r \otimes \msig,\psi
\bigr)(f_{g(n'')}) \psi^{-1}(n'') \, dn''. \end{eqnarray} For $f
\in I({\tau_l}_{(s)},{\breve{\tau_r}}_{(-s)},\msig)$, $s \in
\mspn$ we define $f_g$ as we did for
$I({\tau_l}_{(s)},{\tau_r}_{(s)},\msig)$. In this case $f_g \in
I({\breve{\tau_r}}_{(-s)},\msig)$. Exactly as \eqref{reduction up
1} we have \begin{equation} \label{reduction up 2} \lambda
\bigl((s,-s),\tau_l \otimes \breve{\tau}_r \otimes \msig,\psi
\bigr)(f)=\int_{U_0(\F)} \lambda \bigl(-s, \breve{\tau}_r \otimes
\msig,\psi \bigr)(f_{g(n'')}) \psi^{-1}(n'')  dn''.
\end{equation}

 Let
$$A_{j_{m,n}(w_1^{-1})}:I({\tau_l}_{(s)},{\tau_r}_{(s)},\msig)
\rightarrow I({\tau_l}_{(s)},{\breve{\tau_r}}_{(-s)},\msig)$$ be
as in lemma \ref{lem heart} and let
$$A_{j_{r,r+k}({\omega'_r}^{-1})}:I({\tau_r}_{(s)},\msig)
\rightarrow I({\breve{\tau_r}}_{(-s)},\msig)$$ be the intertwining
operator defined by
$$A_{j_{r,r+k}({\omega'_r}^{-1})}(f)=\int_{j_{r,r+k}(N_{r;0}(\F))}
f(j_{r,r+k}(\epsilon_r\omega'_r n),(I_{2k},1),I_r) \psi^{-1}(n)
dn$$ Using the fact that
$A_{j_{r,r+k}(\epsilon_r\omega'_r)^{-1}}(f_g)=\bigl(A_{j_{m,n}(w_1^{-1})}(f)\bigr)_g$
we prove \eqref{same 1}:
\begin{eqnarray*} && \! \! \! \! \! \! \! C_{\psi}^{\mspn}\bigl(\overline{P_{l,r;k}(\F)},(s,s),\tau_l
\otimes \tau_r \otimes \msig , j_{m,n}(w_1^{-1})\bigr)\\
&& \! \! \! \! \! \! \! = \frac {\lambda \bigl((s,s),\tau_l
\otimes \tau_r \otimes \msig,\psi \bigr) \bigl(f \bigr )}{\lambda
\bigl((s,-s),\tau_l \otimes \breve{\tau}_r \otimes \msig,\psi
\bigr)\bigr(A_{j_{m,n}(w_1^{-1})} (f) \bigl )}\\
&& \! \! \! \! \! \! \! =\frac { \int_{U_0(\F)} \lambda \bigl(s,
\tau_r \otimes \msig,\psi \bigr)(f_{g(n'')}) \psi^{-1}(n'')  dn''}
{\int_{U_0(\F)} \lambda \bigl(-s, \breve{\tau}_r \otimes
\msig,\psi \bigr) \bigl(
A_{j_{r,r+k}(\epsilon_r\omega'_r)^{-1}}(f_{g(n'')})\bigr)
\psi^{-1}(n'') dn''}\\
&& \! \! \! \! \! \! \!= \frac { \int_{U_0(\F)}
C_{\psi}^{\overline{Sp_{2(r+k)}(\F)}}\bigl(\overline{P_{r;k}(\F)},s,
\tau_r \otimes \msig , j_{r,r+k}(\omega_r'^{-1})\bigr) \lambda
\bigl(-s, \breve{\tau}_r \otimes \msig,\psi \bigr) \bigl(
A_{j_{r,r+k}(\epsilon_r\omega'_r)^{-1}}(f_{g(n'')})\bigr)
\psi^{-1}(n'')  dn''} {\int_{U_0(\F)} \lambda \bigl(-s,
\breve{\tau}_r \otimes \msig,\psi \bigr) \bigl(
A_{j_{r,r+k}(\epsilon_r\omega'_r)^{-1}}(f_{g(n'')})\bigr)
\psi^{-1}(n'') dn''}\\
&& \! \! \! \! \! \! \! =
C_{\psi}^{\overline{Sp_{2(r+k)}(\F)}}\bigl(\overline{P_{r;k}(\F)},s,
\tau_r \otimes \msig ,
j_{r,r+k}(\omega_r'^{-1})\bigr).\end{eqnarray*}

 To prove \eqref{same 2} one uses similar arguments. The key point is that
for $f \in I({\tau_l}_{(s)},{\breve{\tau_r}}_{(-s)},\msig), \, g
\in \mspn$, the function $$f_g:\glm \times \gll \times \glr
\rightarrow \C$$ defined by $$f_g(a,l,r)=\ab \det a
\ab^{-\frac{2k+m+1}{2}} \gamma_\psi(a)f\biggl(
\bigl(j_{m,n}(a),1\bigr)g,(I_{2k},1),l,r\biggr),$$ lies in
$I({\tau_l}_{(s)},{\breve{\tau_r}}_{(-s)})$.
\end{proof}
These three lemmas provide the proof of Theorem \ref{mult gama
tau}
\begin{eqnarray*} && \! \! \! \! \! \! \! \lefteqn{ \gamma(\msig \times \tau,s, \psi)} \\
&& \! \! \! \! \! \! \! =\frac
{C_{\psi}^{\mspn}(\overline{P_{m;k}(\F)},s,\tau\otimes \msig,
j_{m,n}(\omega_m'^{-1}))}{C_{\psi}^{\mspm}(\overline{P_{m;0}(\F)},s,\tau,
\omega_m'^{-1} )} \\
&& \! \! \! \! \! \! \! =\frac
{C_{\psi}^{\mspn}(\overline{P_{l,r;k}(\F)},(s,s),\tau_l \otimes
\tau_r \otimes \msig , j_{m,n}(\omega_m'^{-1}))}{
C_{\psi}^{\mspm}(\overline{P_{l,r;0}(\F)},(s,s),\tau_l \otimes
\tau_r, \omega_m'^{-1})} \\
&& \! \! \! \! \! \! \! =\frac
{C_{\psi}^{\mspn}\bigl(\overline{P_{l,r;k}(\F)},(s,s),\tau_l
\otimes \tau_r \otimes \msig ,
j_{m,n}(w_1^{-1})\bigr)}{C_{\psi}^{\mspm}\bigl(\overline{P_{l,r;0}(\F)},(s,s),\tau_l
\otimes \tau_r , w_1^{-1}\bigr)}
\frac{C_{\psi}^{\mspn}\bigl(\overline{P_{l,r;k}(\F)},(s,-s),\tau_l
\otimes \widehat{\tau_r} \otimes \msig ,
j_{m,n}(w_2^{-1})\bigr)}{C_{\psi}^{\mspm}\bigl(\overline{P_{l,r;0}(\F)},(s,-s),\tau_l
\otimes \widehat{\tau_r}  , w_2^{-1}\bigr) }\\
 && \! \! \! \! \! \! \!
\frac{C_{\psi}^{\mspn}\bigl(\overline{P_{l,r;k}(\F)},(-s,s),
\widehat{\tau_r} \otimes \tau_l \otimes \msig ,
j_{m,n}(w_3^{-1})\bigr)}{C_{\psi}^{\mspm}\bigl(\overline{P_{l,r;0}(\F)},(-s,s),
\widehat{\tau_r} \otimes \tau_l  , w_3^{-1}\bigr)} \\
&& \! \! \! \! \! \! \!
=\frac{C_{\psi}^{\overline{Sp_{2(r+k)}(\F)}}\bigl(\overline{P_{r;k}(\F)},s,
\tau_r \otimes \msig ,
j_{r,r+k}(\omega_r'^{-1})\bigr)}{C_{\psi}^{\overline{Sp_{2r}(\F)}}\bigl(\overline{P_{r;0}(\F)},s,
\tau_r  , \omega_r'^{-1}\bigr)}
\frac{C_{\psi}^{\glm}\bigl(P_{l,r}^0(\F),(s,-s),\tau_l \otimes
\widehat{\tau_r},
\varpi_{l,r}^{-1}\bigr)}{C_{\psi}^{\glm}\bigl(P_{l,r}^0(\F),(s,-s),\tau_l
\otimes \widehat{\tau_r}, \varpi_{l,r}^{-1}\bigr)} \\
&& \! \! \! \! \! \! \! \frac
{C_{\psi}^{\overline{Sp_{2(l+k)}(\F)}}\bigl(\overline{P_{l;k}(\F)},s,
\tau_l \otimes \msig ,
j_{l,l+k}(\omega_r'^{-1})\bigr)}{C_{\psi}^{\overline{Sp_{2l}(\F)}}\bigl(\overline{P_{l;0}(\F)},s,
\tau_l , \omega_r'^{-1}\bigr)} =\gamma(\msig \times \tau_l,s,
\psi) \gamma(\msig \times \tau_r,s,
 \psi). \end{eqnarray*}

The proof of Theorem \ref{gama mult sigma} is achieved through
similar steps to those used in the proof of Theorem \ref{mult gama
tau}. We outline them: First one proves an analog to Lemma
\ref{ind by parts} and shows that
\begin{equation} \label{ind by parts 2}
C_{\psi}^{\mspn} \bigl(\overline{P_{m;k}(\F)},s,\tau \otimes \msig
,
j_{m,n}(\omega_m'^{-1})\bigr)=C_{\psi}^{\mspn}\bigl(\overline{P_{m,l;r}(\F)},(s,0),\tau
\otimes \tau_l \otimes \msigr , j_{m,n}(\omega_m'^{-1}) \bigr).
\end{equation}
Then one gives an analog to Lemma \ref{lem heart}: Using a
decomposition of $A_{j_{m,n}(\omega_m'^{-1})}$, one shows that
\begin{equation} \label{heart 2}
C_{\psi}^{\mspn}\bigl(\overline{P_{m,l;r}(\F)},(s,0),\tau \otimes
\tau_l \otimes \msigr ,
j_{m+l,n}(\omega_m'^{-1})\bigr)=(-1,-1)_\F^{ml}k_1(s)k_2(s)k_3(s),
\end{equation} where
\begin{eqnarray}
\nonumber &&
k_1(s)=C_{\psi}^{\mspn}\bigl(\overline{P_{m,l;r}(\F)},(s,0),\tau
\otimes \tau_l \otimes \msigr , j_{m+l,n}(w_4^{-1})\bigr), \\
\nonumber &&
k_2(s)=C_{\psi}^{\mspn}\bigl(\overline{P_{l,m;r}(\F)},(0,s),\tau_l
\otimes \tau \otimes \msigr , j_{m+l,n}(w_5^{-1})\bigr),
\\ \nonumber && k_3(s)=C_{\psi}^{\mspn}\bigl(\overline{P_{m,l;r}(\F)},(0,-s),
\tau_l  \otimes \widehat{\tau}  \otimes \msigr ,
j_{m+l,n}(w_6^{-1})\bigr)\end{eqnarray}  and where
$$w_4=\widehat{\varpi_{m,l}}, \quad w_5=\begin{pmatrix}
& _{I_l}& _{ } &_{ } &_ { } \\
& _{ } & _{ }&_{ } &_ {-\omega_m } \\ &_{ } &_ { } & _{I_l}& _{ }\\
&_{ } &_ {\omega_m} &_{ } &_{ }
\end{pmatrix}, \quad w_6=\widehat{\varpi_{l,m}}.$$ The third step
is, an analog to Lemma \ref{id with c}: \begin{equation} \label{id
with c 2}
C_{\psi}^{\mspn}\bigl(\overline{P_{m,l;r}(\F)},(s,0),\tau \otimes
\tau_l \otimes \msigr , j_{m+l,n}(w_4^{-1})\bigr)=
C_{\psi^{-1}}^{GL_{m+l}(\F)}\bigl(P^0_{m,l}(\F)(\frac{s}{2},\frac{-s}{2}),\tau
\otimes \tau_l,\varpi_{m,l}^{-1} \bigr), $$ $$
C_{\psi}^{\mspn}\bigl(\overline{P_{l,m;r}(\F)},(0,s),\tau_l
\otimes \tau \otimes \msigr ,
j_{m+l,n}(w_5^{-1})\bigr)=C_\psi^{\overline{Sp_{2(m+r)}(\F)}}\bigl(\overline{P_{l;r}(\F)},s,\tau
\otimes \msigr,j_{m,m+k}(\omega_m'^{-1}) \bigr), $$ $$
C_{\psi}^{\mspn}\bigl(\overline{P_{l,m;r}(\F)},(0,-s),\tau_l
\otimes \widehat{\tau }\otimes \msigr ,
j_{m+l,n}(w_6^{-1})\bigr)=C_{\psi}^{GL_{m+l}(\F)}\bigl(P^0_{m,l}(\F)(\frac{s}{2},\frac{-s}{2}),\tau_l
\otimes \widehat{\tau},\varpi_{l,m}^{-1} \bigr).
\end{equation}
Combining \eqref{ind by parts 2}, \eqref{heart 2}, \eqref{id with
c 2} we have: \begin{eqnarray*} && \gamma(\msig \times \tau,s,
\psi)=(-1,-1)_\F^{ml}
 \gamma(\msigr \times \tau,s, \psi) \\ &&
C_{\psi^{-1}}^{GL_{m+l}(\F)}\bigl(P^0_{m,l}(\F)(\frac{s}{2},\frac{-s}{2}),\tau
\otimes \tau_l,\varpi_{m,l}^{-1} \bigr)
C_{\psi^{-1}}^{GL_{m+l}(\F)}\bigl(P^0_{m,l}(\F)(\frac{s}{2},\frac{-s}{2}),\tau_l
\otimes \widehat{\tau},\varpi_{l,m}^{-1} \bigr).\end{eqnarray*}
With \eqref{gama def gl} we finish.

\subsection{Computation of $\gamma(\overline{\sigma} \times \tau,s,\psi)$
for principal series representations} \label{prin gama sec} Assume
that $\F$ is either $\R$, $\C$ or a p-adic field. Let
$\eta_1,\eta_2,\ldots,\eta_k$ be $k$ characters of $\F^*$ and let
$\gamapsi \otimes \chi=(\gamapsi \circ \det) \otimes \chi$ be the
character of $\overline{\tspk}$ defined by
$$\bigl(diag(t_1,t_2,\ldots,t_k,t_1^{-1},t_2^{-1},\ldots,t_k^{-1}),\epsilon
\bigr)\mapsto \epsilon
\gamma_\psi^{-1}(t_1t_2\ldots,t_k)\prod_{i=1}^{k} \eta_i(t_i).$$
Let $\alpha_1,\alpha_2,\ldots,\alpha_m$ be m characters of $\F^*$
and let $\mu$ be the character of $\tglm$ defined by
$$Diag(t_1,t_2,\ldots,t_m) \mapsto \prod_{i=1}^{m}\alpha_i(t_i).$$
Define $\msig$ and $\tau$ to be the corresponding principal series
representations:
$$\msig=I(\chi)=Ind^{\mspk}_{\mbk} \gamapsi \otimes \chi, \quad
\tau=I(\mu)=Ind^{\glm}_{B_{\glm}} \mu.$$
\begin{lem} \label{princcomp} There exists $c \in \{ \pm 1 \}$
such that
\begin{equation} \label{mult cor prin}\gamma(\overline{\sigma} \times \tau,s,\psi)= c\prod_{i=1}^k\prod_{j=1}^m
\gamma(\alpha_j \times \eta_i^{-1},s,\psi)\gamma(\eta_i \times
\alpha_j,s,\psi). \end{equation}
\end{lem}
\begin{proof} Note that $\tau \simeq Ind^{\glm}_{P^0_{1,m-1}(\F)}
\alpha_1 \otimes \tau'$, where
$\tau'=Ind^{GL_{m-1}(\F)}_{B_{GL_{m-1}}(\F)}\otimes_{i=1}^{m-1}\alpha_i$.
Theorem \ref{mult gama tau} implies that
$$\gamma(\overline{\sigma} \times \tau,s,\psi)=\gamma(\msig \times
\alpha_1,s, \psi) \gamma(\msig \times \tau',s,
 \psi).$$
Repeating this argument $m-1$ more times we observe that
\begin{equation} \label{use multtau}\gamma(\overline{\sigma} \times
\tau,s,\psi)=\prod_{j=1}^m\gamma(\msig \times \alpha_j,s,
\psi).\end{equation} Next we note that
$\msig=Ind^{\mspk}_{\overline{P_{1;k-1}(\F)}} (\gamapsi \otimes
\eta_1) \otimes \overline{\sigma'}$, where
$\overline{\sigma'}=Ind^{\overline{Sp_{2(k-1)}(\F)}}_{\overline{B_{Sp_{2(k-1)}}(\F)}}\gamma_\psi^{-1}
\otimes \bigl(\otimes_{j=1}^{k-1}\eta_j \bigr)$. By using Theorem
\ref{mult gama sigma} we observe that for all $1\leq j \leq m$.
There exists $c' \in \{ \pm 1 \}$ such that \begin{equation}
\label{use multsgima}\gamma(\msig \times \alpha_j,s,
\psi)=c'\gamma(\overline{\sigma'} \times \alpha_j,s, \psi)
\gamma(\alpha_j \times \eta_1^{-1},s, \psi) \gamma(\eta_1 \times
{\alpha_j},s, \psi).\end{equation} By Repeating this argument
$k-1$ more times for each $1 \leq j \leq m$ and by using
\eqref{use multtau} we finish.
\end{proof}

\newpage

\section{Computation of the local coefficients for principal series representations of $\msl$}
\label{second paper} Assume that $\F$ is either
$\R$ , $\C$ or a finite extension of $\Q _p$. Let $\psi$ be a non-trivial
character of $\F$ and let $\chi$ be a character of $\F^*$. For $s
\in \C$ let
$$I(\chi \otimes \gamma_\psi^{-1},s)=Ind^{\msl}_{\overline{B_{SL_2}(\F)}}
\gamma_{\psi}^{-1}\otimes \chi_{(s)} $$ be the corresponding
principal series representation of $\msl$. Its space consists of
smooth complex functions on $\msl$ which satisfy:$$f\Bigl(
\bigl(\begin{pmatrix} _{a} & _{b}\\_{0} & _{a^{-1}}\end{pmatrix},
\epsilon \bigr) g \Bigr)=\epsilon \ab a \ab^{s+1} \chi(a)
\gamma_{\psi}^{-1}(a)f(g),$$ for all $a\in \F^*$, $b\in \F$, $g
\in \overline{SL_2(\F)}$. $\msl$ acts on this space by right
translations. $\lambda \bigl(s,\chi,\psi_a)$, the $\psi_a$
Whittaker functional on $I(\chi \otimes \gamma_\psi^{-1},s)$,
defined in Chapter \ref{chapter def lc}, is the analytic
continuation of
\begin{equation}\label{whi func} \int_\F f
\Biggl(\begin{pmatrix} _{0} & _{1}\\_{-1} &
_{0}\end{pmatrix}\begin{pmatrix} _{1} & _{x}\\_{0} &
_{1}\end{pmatrix},1 \Biggr) \psi_a^{-1}(x) dx.\end{equation} The
intertwining operator corresponding to the unique non-trivial Weyl
element of $SL_2(\F)$,
$$A(s):I(\chi \otimes \gamma_\psi^{-1},s) \rightarrow
I(\chi^{-1}\otimes \gamma_\psi^{-1},-s)$$ is defined by the
meromorphic continuation of
\begin{equation} \label{io int}\bigl(A(s)(f)\bigr)(g)=\int_\F f
\Biggl( \Bigl(\begin{pmatrix} _{0} & _{1}\\_{-1} &
_{0}\end{pmatrix}\begin{pmatrix} _{1} & _{x}\\_{0} &
_{1}\end{pmatrix},1\Bigr)g\Biggr) dx.\end{equation}

We shall prove in this chapter, see Theorems \ref{the formula},
\ref{real thm} and Lemma \ref{easy complex lem}, that there exists
an exponential function, $\widetilde{\epsilon}(\chi,s,\psi)$ such
that
\begin{equation} \label{res l}C_\psi^{\overline{SL_2(\F)}}\Bigl(\overline{B_{SL_2(\F)}},s,\chi,(\begin{smallmatrix}
& {0}& {-1} \\
& {1} &
{0}\end{smallmatrix})\Bigr)=\widetilde{\epsilon}(\chi,s,\psi)\frac{L(\chi,s+\half)}
{L(\chi^{-1},-s+\half)} \frac{L(\chi^{-2},-2s+1)}
{L(\chi^2,2s)}.\end{equation} Furthermore, if $\F$ is a p-adic
field of odd residual characteristic, $\psi$ is normalized and
$\chi$ is unramified then $\widetilde{\epsilon}(\chi,s,\psi)=1$.
Recall that
$$\gamma(\chi,s,\psi)=\epsilon(\chi,s,\psi)
\frac{L(\chi^{-1},1-s)}{L(\chi,s)},$$ where
$\epsilon(\chi,s,\psi)$ is an exponential factor. Thus, \eqref
{res l} can be written as
\begin{equation} \label{res g}C_\psi^{\overline{SL_2(\F)}}\Bigl(\overline{B_{SL_2(\F)}},s,\chi,(\begin{smallmatrix}
& {0}& {-1} \\
& {1} & {0}\end{smallmatrix})\Bigr)=\epsilon'(\chi,s,\psi)
\frac{\gamma(\chi^2,2s,\psi)}{\gamma(\chi,s+\half,\psi)},\end{equation}
where $\epsilon'(\chi,s,\psi)$ is an exponential factor which
equals 1 if $\chi$ is unramified and $\F$ is p-adic field of odd
residual characteristic.

In this chapter only, we shall write
$$C_{\psi_a}(\etac \otimes \gamma_\psi^{-1},s)$$
Instead of
$$C_{\psi_a}^{\overline{SL_2(\F)}}\Bigl(\overline{B_{SL_2(\F)}},s,\chi,(\begin{smallmatrix}
& {0}& {-1} \\
& {1} & {0}\end{smallmatrix})\Bigr).$$ This notation emphasizes
the dependence of the local coefficient on two additive
characters, rather than on one in the algebraic case.

This chapter is organized as follows. In Section \ref{padic comp}
we present the p-adic computation. Our computations include the
often overlooked case of 2-adic fields. In Section \ref{real case}
the computation for the real case is given. In both two sections
the exponential factor is computed explicitly. The complex case is
addressed in Section \ref{easy complex sec}. Since
$\overline{SL_2(\C)}=SL_2(\C) \times \{\pm 1 \}$ and
$\gamma_\psi(\C^*)=1$ the local coefficients for this group are
identical to the local coefficients of $SL_2(\C)$. The $SL_2(\C)$
computation is given in Theorem 3.13 of \cite{Sha85}. It turns out
that the $SL_2(\C)$ computation of the local coefficients agrees
with \eqref{res g}. Next we give two detailed remarks; one on the
chosen parameterizations and the other on the existence of a
non-archimedean metaplectic $\widetilde{\gamma}$-factor defined by
a similar way to the definition of the Tate $\gamma$-factor (see
\cite{T}) whose relation to the computed local coefficients is
similar to the relation that the Tate $\gamma$-factor has with the
local coefficients of $SL_2(\F)$. These are Sections \ref{par
remark} and \ref{tate remark}. We conclude this chapter in Section
\ref{prin comp with soo} where we show that if $\tau$ is a
principal series representation then
$$C_{\psi}^{\mspm}\bigl(\overline{P_{m;0}(\F)},s,\tau,
\omega_m'^{-1} \bigr)=c(s)\frac
{\gamma(\tau,sym^2,2s,\psi)}{\gamma(\tau,s+\half,\psi)},$$ where
$c(s)$ is an exponential factor which equals 1 if $\F$ is a p-adic
field of odd residual characteristic, $\psi$ is normalized and
$\tau$ is unramified. We do so using the $\msl$ computations
presented in Sections \ref{padic comp}, \ref{real case} and
\ref{easy complex sec},
 the multiplicativity of the local coefficients proved in Chapter
\ref{cha gam} and the theory of local coefficients and $\gamma$-
factors attached $\gln$; see \cite{Sha84}. In Section \ref{true
for cusp} we shall show that the same formula holds if $\tau$ is
supercuspidal.

\subsection{Non-archimedean case} \label{padic comp}
For a non-trivial $\psi$ and ramified $\chi$ we define
$$G(\chi,\psi)=\int_{\Of^*} \psi_0(\pi^{-m(\chi)}u)\,
\chi(u) \, du.$$ An easy modification of the evaluation of
classical Gauss sums, see Proposition 8.22 of \cite{I} for
example, shows that $\mid G(\chi,\psi)\mid
=q^{-\frac{m(\chi)}{2}}$. If $\chi$ is unramified we define
$G(\chi,\psi)=1$. We now state the main theorem of this section.
\begin{thm} \label{the formula} There exist $k_\chi \in \C^*$ and $d_\chi \in \Z$
such that \begin{equation} \label{the formula l} C_{\psi}(\etac
\otimes \gamma_\psi^{-1},s) =k_\chi q^{d_\chi s}
\frac{L(\chi,s+\half)} {L(\chi^{-1},-s+\half)}
\frac{L(\chi^{-2},-2s+1)} {L(\chi^2,2s)}.
\end{equation} Furthermore, if $\chi^2$ is unramified we
have:$$k_\chi=\gamma_{\psi_0}^{-1}(\pi^n) \frac
{c_{\psi_0}(-1)\chi(-\pi^{2e-m(\chi)-n})}{G(\chi,\psi)}q^{-\frac
{m(\chi)}{2}}, \, \, \, \, d_\chi=m(\chi)-2e+n,$$ where $n$ is the
conductor of $\psi$.

\end{thm}

Recall that for $\F$, a p-adic field,
$$L_{\F}(\chi,s\bigr)=\begin{cases} \frac 1 {1-\chi(\pi)q^{-s}} &
\chi$ is unramidied$
\\ 1 & $otherwise$  \end{cases}$$
For unramified characters and $\F$ of odd residual characteristic
\eqref {the formula l} can be proved by using the existence of a
spherical function, see \cite{BFH}.

We proceed as follows: In Subsection \ref {red int} we reduce the
computation of the local coefficient to a computation of a certain
"Tate type" integral; see Lemma \ref{reduction}. In Subsection
\ref{lemmas} we collect some facts that will be used in Subsection
\ref{int comp} where we compute this integral.
\subsubsection{The local coefficient expressed as an integral}
\label{red int}
\begin{lem}\label{reduction} For $Re(s)>>0$:
$${C_{\psi_a}(\etac \otimes \gamma_\psi^{-1},s)}^{-1}= \int_{\F^*}
\gamma_\psi^{-1}(u) \chi_{(s)}(u) \psi_a(u) d^*u.$$

\end{lem}
This integral should be understood as a principal value integral,
i.e., $$\int_{\Pf^{-N}} \gamma_\psi^{-1}(u) \chi_{(s)}(u)
\psi_a(u) \, d^*u,$$ for $N$ sufficiently large.
\begin{proof}
We recall that the integral in the right side of \eqref{whi func}
converges in principal value for all $s$, and furthermore, for all
$f \in I(\chi \otimes \gamma_\psi^{-1},s)$ there exists $N_f$ such
that for all $N>N_f$:
$$
\lambda (s,\chi,\psi_a)(f)=\int_{\Pf^{-N}} f
\Biggl(\begin{pmatrix} _{0} & _{1}\\_{-1} &
_{0}\end{pmatrix}\begin{pmatrix} _{1} & _{x}\\_{0} &
_{1}\end{pmatrix},1 \Biggr) \psi_a^{-1}(x) dx.
$$
Assume that $Re(s)$ is so large such that the integral in the
right hand side of \eqref{io int} converges absolutely for all $f
\in I(\chi \otimes \gamma_\psi^{-1},s)$. Define $f \in I(\chi
\otimes \gamma_\psi^{-1},s)$ to be the following function:
$$
f(g)=\begin{cases} 0 & g \in \overline{B} \\ \epsilon
\chi_{s+1}(b) \gamma_\psi^{-1}\bigr(b) \phi(x) & g= \Bigr(
\begin{pmatrix} _{b} & _{c}\\_{0} & _{b^{-1}}\end{pmatrix}
\begin{pmatrix} _{0} & _{1}\\_{-1} & _{0}\end{pmatrix}\begin{pmatrix} _{1} & _{x}\\_{0}
& _{1}\end{pmatrix},\epsilon\bigl)  \end{cases},
$$
where $\phi \in S(\F)$ is the characteristic function of $\Of$. It
is sufficient to show that
\begin{equation}
\label{whatwewant}\lambda (-s,\chi^{-1},\psi_a)(A(s)f)=\lambda
(s,\chi,\psi_a)(f)\int_{\F^*} \gamma_\psi^{-1}(u) \chi_{(s)}(u)
\psi_a(u) d^*u.
\end{equation}
Indeed,
\begin{eqnarray*}
\lefteqn{ \lambda (-s,\chi^{-1},\psi_a)(A(s)f)}&&\\
&=&\int_\F \int_{\F^*} f \Biggl(\Bigl(\begin{pmatrix} _{0} &
_{1}\\_{-1} & _{0}\end{pmatrix}\begin{pmatrix} _{1} & _{u}\\_{0} &
_{1}\end{pmatrix},1)(\begin{pmatrix} _{0} & _{1}\\_{-1} &
_{0}\end{pmatrix}\begin{pmatrix} _{1} & _{x}\\_{0} &
_{1}\end{pmatrix},1\Bigr) \Biggr) \psi_a^{-1}(x) \, du \, dx.
\end{eqnarray*}
By matrix multiplication and by \eqref{rao} we have:
\begin{eqnarray*}
\lefteqn{ \Biggl(\begin{pmatrix} _{0} & _{1}\\_{-1} &
_{0}\end{pmatrix}\begin{pmatrix} _{1} & _{u}\\_{0} &
_{1}\end{pmatrix},1 \Biggr)\Biggl(\begin{pmatrix} _{0} &
_{1}\\_{-1} & _{0}\end{pmatrix}\begin{pmatrix} _{1} & _{x}\\_{0} &
_{1}\end{pmatrix},1 \Biggr)}&&\\ &=&\Biggl(\begin{pmatrix}
_{-u^{-1}} & _{1}\\_{0} & _{-u}\end{pmatrix},1 \Biggr)
\Biggl(\begin{pmatrix} _{0} & _{1}\\_{-1} &
_{0}\end{pmatrix}\begin{pmatrix} _{1} & _{x-u^{-1}}\\_{0} &
_{1}\end{pmatrix},1\Biggl).\end{eqnarray*} Hence, for $N$
sufficiently large:
\begin{eqnarray*} \lefteqn{\lambda (-s,\chi^{-1},\psi_a)(A(s)f \bigr)} && \\ &=& \int_{\Pf^{-N}} \int_{\F^*}
\chi_{(s)}(-u^{-1}) \ab u \ab^{-1}
 \gamma_\psi(-u)^{-1} f \bigl(\begin{pmatrix} _{0} & _{1}\\_{-1}
& _{0}\end{pmatrix}\begin{pmatrix} _{1} & _{x-u^{-1}}\\_{0} &
_{1}\end{pmatrix},1  \bigr) \psi_a^{-1}(x) \, du \, dx \\
&=&\int_{\Pf^{-N}} \int_{\F^*} \chi_{(s)}(-u) \ab u
\ab^{-1}\gamma_\psi(-u)^{-1} f \bigl(\begin{pmatrix} _{0} &
_{1}\\_{-1} & _{0}\end{pmatrix}\begin{pmatrix} _{1} & _{x-u}\\_{0}
& _{1}\end{pmatrix},1  \bigr) \psi_a^{-1}(x) \, du \,
dx.\end{eqnarray*} Since the map $y \mapsto f
\bigl(\begin{pmatrix} _{0} & _{1}\\_{-1} &
_{0}\end{pmatrix}\begin{pmatrix} _{1} & _{y}\\_{0} &
_{1}\end{pmatrix},1  \bigr)$ is supported on $\Of$ and since we
may assume that $N>0$ we have \begin{eqnarray*} \lefteqn{\lambda
(-s,\chi^{-1},\psi_a)(A(s)f \bigr)} && \\&=& \int_{\ab x \ab <
q^N} \int_{0<\ab u \ab< q^N} \chi_{(s)}(-u) \gamma_\psi(-u)^{-1} f
\bigl(\begin{pmatrix} _{0} & _{1}\\_{-1} &
_{0}\end{pmatrix}\begin{pmatrix} _{1} & _{x-u}\\_{0} &
_{1}\end{pmatrix},1  \bigr) \psi_a^{-1}(x) \, d^*u \,
dx.\end{eqnarray*} By changing the order of integration and by
changing $x \mapsto x+u$ we obtain
\begin{eqnarray*} \lefteqn {\lambda (-s,\chi^{-1},\psi_a)(A(s)f \bigr)} && \\ &=&\int_{\ab x \ab <q^N} f
\Biggl(\begin{pmatrix} _{0} & _{1}\\_{-1} &
_{0}\end{pmatrix}\begin{pmatrix} _{1} & _{x}\\_{0} &
_{1}\end{pmatrix},1 \Biggr) \psi_a^{-1}(x) dx \int_{0<\ab u \ab
<q^N} \gamma_\psi^{-1}(-u) \chi_{(s)}(-u) \psi_a^{-1}(u)d^*u
.\end{eqnarray*} \eqref{whatwewant} is proven once we change $u
\mapsto -u$.
\end{proof}
The next lemma reduces the proof of  Theorem \ref{the formula} to
case where $\psi$ is normalized.
\begin{lem} \label{reduction psi} Define $n$ to be the conductor of $\psi$. We have:

$$C_{\psi}(\etac \otimes \gamma_\psi^{-1},s) =
\gamma_{\psi_0}^{-1}(\pi^n) \chi(\pi^{-n})q^{ns}C_{\psi_0}(\etac
\otimes \gamma_{\psi_0}^{-1},s).$$
\end{lem}
\begin{proof}
Recalling the definition of $\gamma_\psi$ we observe that
$$
\gamma_\psi(a)=\frac
{\gamma_{\psi_0}(a\pi^{-n})}{\gamma_{\psi_0}(\pi^{-n})}.
$$
Thus, by Lemma \ref{reduction} we have
\begin{eqnarray*}
C_\psi(\etac \otimes \gamma_\psi^{-1},s)^{-1} &=&
\int_{\F^*} \gamma_\psi^{-1}(u) \chi_{(s)}(u) \psi(u) \, d^*u\\
&=& \gamma_{\psi_0}(\pi^n)\int_{\F^*}
\gamma_{\psi_0}^{-1}(u\pi^{-n}) \chi_{(s)}(u) \psi^{-1}(u) \,
d^*u\\
&=&\gamma_{\psi_0}(\pi^n)\int_{\F^*} \gamma_{\psi_0}^{-1}(u)
\chi_{(s)}(u \pi^n) \psi^{-1}(u\pi^n)\, d^*u\\
&=& \gamma_{\psi_0}(\pi^n)\chi_{(s)}(\pi^n)\int_{\F^*}
\gamma_{\psi_0}^{-1}(u) \chi_{(s)}(u ) \psi_0^{-1}(u) \,
d^*u\\
&=&\gamma_{\psi_0}(\pi^n)q^{-ns}\chi(\pi^n) C_{\psi_0}(\etac
\otimes \gamma_{\psi_0}^{-1},s)^{-1}.
\end{eqnarray*}

\end{proof}
We note that in Theorem \ref{the formula}, we compute
$$C_{\psi}(\etac \otimes \gamma_\psi^{-1},s)$$
rather then $$C_{\psi_a}(\etac \otimes \gamma_\psi^{-1},s),$$ that
is, we use the same additive character in the definition of
$I(\chi \otimes \gamma_\psi^{-1},s)$ and in the Whittaker
functional. The advantage of this choice is explained in Section
\ref{par remark}. We conclude this section with a lemma that
describes the relation between these two local coefficients.
\begin{lem} \label{a to 1}
$$
C_{\psi}(\etac \otimes
\gamma_\psi^{-1},s)\gamma_\psi(a)\chi_{(s)}(a)=C_{\psi_a}(\ \chi
\cdot  (a,\cdot) \otimes \gamma_\psi^{-1},s).
$$
\end{lem}
\begin{proof} By Lemma \ref{reduction}
\begin{eqnarray*}C_{\psi_a}(\etac \otimes \gamma_\psi^{-1},s)^{-1}&=&
\int_{\F^*} \gamma_\psi^{-1}(u) \chi_{(s)}(u) \psi(au)
d^*u\\
&=&\int_{\F^*} \gamma_\psi^{-1}(a^{-1}u) \chi_{(s)}(a^{-1}u)
\psi(u) d^*u\\
&=& \gamma_\psi^{-1}(a)\chi_{(s)}(a^{-1})\int_{\F^*}
\gamma_\psi^{-1}(u) \chi_{(s)}(u)(a,u)_\F \psi(u)
d^*u\\
&=&\gamma_\psi^{-1}(a)\chi_{(s)}(a^{-1})C_{\psi}(\etac \cdot
(a,\cdot)\otimes \gamma_\psi^{-1},s)^{-1}.
\end{eqnarray*}
\end{proof}

\subsubsection{Some lemmas.} \label{lemmas}
\begin{lem} \label{m chi quad}Assume that $\chi^2$ is unramified and that
$m(\chi) \leq 2e$. Then $m(\chi)$ is even.
\end{lem}
\begin{proof} The lemma is trivial for $\F$ of odd residual
characteristic. Assume that $\F$ is of even residual
characteristic. Since $\Of^* /(1+\Pf) \simeq \overline{\F}^*$ is a
cyclic group of odd order it follows that if $\chi(1+\Pf)=1$ then
$\chi$ is unramified. Thus, $m(\chi) \neq 1$ and it is left to
show that for any $2 \leq k \leq e-1$ if $\chi(1+\Pf^{2k+1})=1$
then $\chi(1+\Pf^{2k})=1$. For $l \in \N$, define $G_l=\Of^* /
1+\Pf^l$.  The number of quadratic characters of $G_l$ is
$b_l=[G_l : {G_l^2}]$. We want to prove that for any $2 \leq k
\leq e-1$, $b_{2k}=b_{2k+1}$. Note that for any $l>1$, $G_l \simeq
\rf^* \times H_l$, where $H_l=1+\Pf / 1+\Pf^{l}$ is a commutative
group of order $2^{fl}$ where $f=[\rf: \Z_2]$ is the residue class
degree of $\F$ over $\Q_2$. Thus, there exists $c_l \in \N$ such
that $H_l=\prod_{j=1}^{c_l} \Z /
  2^{d_l(j)}\Z$ and $G_l / {G_l^2}=H_l / {H_l^2}={\bigl(\Z / 2 \Z
  \bigr)}^{c_l}.$ Hence, the proof is done once we show that for
any $2 \leq k \leq e-1$, $c_{2k}=c_{2k+1}$. We shall prove that
for any $2 \leq k \leq e-1$, $c_{2k}=c_{2k+1}=fk$. Note that $$
\mid \{x \in H_l \mid x^2=1 \} \mid =2^{c_l}.$$ Also note that
$H_l$ may be realized as $$\{1+\sum_{j=1}^{l-1}a_j \pi^j \mid a_j
\in A \},$$ where $A$ is a set of representatives of $\Of / \Pf$
which contains 0. Hence, it is sufficient to show that for
$x=1+\sum_{j=1}^{2k-1}a_j \pi^j$, where $a_j \in A$ we have $x^2
\in 1+\Pf^{2k}$ if and only if $a_j=0$ for all $1 \leq j \leq
k-1$, and that for $x=1+\sum_{j=1}^{2k}a_j \pi^j$, where $a_j \in
A$ we have $x^2 \in 1+\Pf^{2k+1}$ if and only if $a_j=0$ for all
$1 \leq j \leq k$. We prove only the first claim, the second is
proven in the same way. Suppose $x \neq 1$ and that
$x=1+\sum_{j=1}^{2k-1}a_j \pi^j$, where $a_j \in A$. Note that
$$\bigl(1+\sum_{j=1}^{2k-1}a_j \pi^j \bigr)^2=1+\sum_{i=1}^{2k-1}
a_i^2\pi^{2i} +\sum_{i=1}^{2k-1}  \omega a_i\pi^{e+i}+ \! \! \!
\sum_{1\leq i<j
 \leq 2k-1} \omega a_ia_j\pi^{i+j+e}.$$ (recall that $\omega=2\pi^{-e} \in \Of^*$).
Since $k<e$ we have $\ab x^2-1 \ab=q^{-2r}$, where $r$ is the
minimal index such that $a_r \neq0$. The assertion is now proven.
\end{proof}
The following sets will appear in the computations of $\int_{\F^*}
\gamma_\psi^{-1}(u) \chi_{(s)}(u) \psi(u) d^*u$. For $n\in \N$
define
\begin{eqnarray}H(n,\F) &=&  \bigl\{x \in \Of^* \mid \ab1-x^2 \ab \leq
q^{-n}\bigr\},\\
D(n,\F) &=& \bigl\{x \in \Of^* \mid \ab1-x^2 \ab=q^{1-n}\bigr\}.
\end{eqnarray}
\begin{lem} \label{H D comp}
Suppose $\F$ is of odd residual characteristic. Then
\begin{eqnarray} \label{D 1 Odd}D(1,\F) &=& \bigl\{x \in \Of \mid \overline{x} \notin
\{0,\pm 1\}\bigr\},\\ \label{H 1 Odd}H(1,\F) &=& \bigl\{x \in \Of
\mid \overline{x} \in \{\pm 1\}\bigr\}.\end{eqnarray} Suppose $\F$
is of even residual characteristic. Then
\begin{eqnarray} \label{D 1 even}D(1,\F) &=& \big\{x \in \Of \mid \overline{x} \notin
\{0,1\}\bigr\},\\ \label{H 1 even}H(1,\F) &=& \bigl\{x \in \Of
\mid \overline{x}=1 \bigr\},\end{eqnarray} for $1 \leq k  \leq e$:
\begin{equation} \label{D H 2k} D(2k,\F)=\varnothing, \, \,
H(2k,\F)=H(2k-1,\F)=1+\Pf^{k},\end{equation}  for all $1 \leq k
<e:$
\begin{equation} \label{D H 2k+1}D(2k+1,\F)=1+\Pf^{k}
\setminus 1+\Pf^{k+1},\end{equation} and
\begin{eqnarray} \label{D 2e+1}D(2e+1,\F) &=& \bigl\{1+b\pi^e \mid b \in \Of^*, \, \overline{b} \neq
\overline{w} \bigr\},\\ \label{H 2e+1}H(2e+1,\F) &=&
\bigl\{1+b\pi^e \mid b \in \Of, \, \overline{b} \in \{0,
\overline{w}\}  \bigr\}.\end{eqnarray}
\end{lem}
\begin{proof}
We start with $\F$ of odd residual characteristic. Suppose $x\in
\Of^*$. Then $\ab 1-x^2 \ab=1$ is equivalent to $1-\overline{x}^2
\neq 0$. This proves \eqref{D 1 Odd}. \eqref{H 1 Odd} follows
immediately since $H(1,F)=\Of^* \setminus D(1,F)$. Suppose now
that $\F$ is of even residual characteristic. \eqref{D 1 even} and
\eqref{H 1 even} are proven in the same way as \eqref{D 1 Odd} and
\eqref{H 1 Odd}. Note that in the case of even residual
characteristic $t=-t$ for all $t\in \rf$. Assume now $x=1+b\pi^m$,
$m\geq 1$, $b \in \Of^*$. We have:
$$\ab x^2-1\ab =\ab bw\pi^{m+e}+b^2\pi^{2m} \ab=\begin{cases} q^{-2m} &
 1\leq m <e \\ q^{-2e}\ab b+w \ab  & m=e \\ q^{m+e} & m>e
\end{cases}.$$
the rest of the assertions mentioned in this lemma follow at once.
\end{proof}
\begin{lem} \label{gauss} Assume $\chi$ is ramified and that $\psi$ is normalized. Fix $x\in
\Of^*$, and $n$, a non-negative integer.
\begin{equation} \label{gauss1}  \int_{\Of ^*}\psi(\pi^{-n}u)  \, du  =\begin{cases} 1-q^{-1} & n=0  \\
-q^{-1} & n=1 \\
0 & n>1
\end{cases}.\end{equation}
If $n>0$ then
\begin{equation}  \label{gauss2}
\int_{\Of ^*}\psi\bigl(\pi^{-m(\chi)}u(1-\pi^nx^2)\bigr) \chi(u)
\, du = \chi^{-1}(1-\pi^nx^2)G(\chi,\psi),\end{equation} and
\begin{equation} \label{gauss3} \int_{\Of ^*}\psi\bigl(\pi^{-n}u(1-x^2) \bigr)  \, du  =
\begin{cases} 1-q^{-1} & x\in H(n,\F)  \\
-q^{-1} & x \in D(n,\F) \\
0 & $otherwise$
\end{cases}.\end{equation}
If $n \neq m(\chi)$ then:
 \begin{equation} \label{gauss4}\int_{\Of ^*}\psi(\pi^{-n}u) \chi(u) \, du =0.
 \end{equation}
 \begin{equation} \label{gauss5}\int_{\Of ^*}\psi\bigl(\pi^{-n}u(1-x^2) \bigr) \chi(u) \, du  =
 \begin{cases} \chi^{-1}\bigl(\pi^{n-m(\chi)}(1-x^2)\bigr)G(\chi,\psi) & n \geq m(\chi) $ and$\\
 & x \in D(n-m(\chi)+1,\F) \\
0 & $otherwise$
\end{cases},\end{equation}
Suppose $2n \geq m(\chi)$. Then
\begin{equation} \label{gauss6}\int_{\Of^*} \chi(1+\pi^nu) \, du=\begin{cases}
1-q^{-1}& n\geq m(\chi) \\ -q^{-1} &  n=m(\chi)-1
\\ 0 & n \leq m(\chi)-2 \end{cases}. \end{equation}
Suppose in addition that $\chi'$ is another non-trivial character
and that $m(\chi)-n \neq m(\chi')$. Then
\begin{equation} \label{gauss7}\int_{\Of^*} \chi'(x) \chi(1+\pi^nu) \, du=0
.\end{equation}
\end{lem}
\begin{proof} To prove \eqref{gauss1} note that $$\int_{\Of ^*}\psi(\pi^{-n}u)  \, du  =
\int_{\Of}\psi(\pi^{-n}u)  \, du-\int_{\Pf}\psi(\pi^{-n}u)  \,
du,$$ and recall that for any compact group $G$ and a character
$\beta$ of $G$ we have $$\int_G \beta(g) \, dg=\begin{cases}
Vol(G)& \beta $ is trivial$ \\ 0 &  $otherwise$ \end{cases}.$$ The
proof of \eqref{gauss2} follows from the definition of
$G(\chi,\psi)$ by changing of $u \mapsto u(1-\pi^n x^{2})^{-1}$.
To prove \eqref{gauss3} note that the conductor of the additive
character $u \mapsto \psi\bigl(\pi^{-n}u(1-x^2) \bigr)$ is
$n-log_q \ab 1-x^2 \ab$ and repeat the proof of \eqref{gauss1}. We
now prove \eqref{gauss4}: If $n<m(\chi)$ then
$$\int_{\Of ^*}\psi(\pi^{-n}u) \chi(u) \, du= \sum_{t \in \Of^* /
1+\Pf^n} \psi(\pi^{-n} t)\chi(t)\int_{1+\Pf^n}\chi(u) \, du=0,$$
and if $n>m(\chi)$ then \begin{eqnarray*}\int_{\Of
^*}\psi(\pi^{-n}u) \chi(u) \, du &=& \sum_{t \in \Of^* /
1+\Pf^{m(\chi)}}\chi(t)\int_{1+\Pf^{m(\chi)}} \psi({\pi^{-n}ut})
\, du \\ &=& \sum_{t \in \Of^* /
1+\Pf^{m(\chi)}}\chi(t)\psi({\pi^{-n}t}) \int_{\Pf^{m(\chi)}}
\psi({\pi^{-n}k}) \, dk=0.\end{eqnarray*} To prove \eqref{gauss5}
write
$$\int_{\Of ^*}\psi\bigl(\pi^{-n}u(1-x^2) \bigr) \chi(u) \, du  =
\int_{\Of ^*}\psi\bigl(\pi^{-m(\chi)}(1-x^2)\pi^{m(\chi)-n}
u\bigr) \chi(u) \, du.$$ By \eqref{gauss4}, the integral vanishes
unless $(1-x^2)\pi^{m(\chi)-n} \in \Of^*$. This implies the second
case of \eqref{gauss5}. Changing  the integration variable
$u\mapsto u\bigl((1-x^2)\pi^{m(\chi)-n}\bigr)^{-1}$ proves the
first case. We move to \eqref{gauss6}: Since $2n \geq m(\chi)$ we
have $$\bigl(1+\pi^n(x+y)\bigr)=(1+\pi^n x)(1+\pi^n y) \,
\bigl(mod \, 1+\Pf^{m(\chi)}\bigr),$$ for all $x,y \in \Of$. This
implies that $u \mapsto \chi (1+\pi^n u)$ is an additive character
of $\Of$. It is trivial if $n\geq m(\chi)$ otherwise its conductor
is $m(\chi)-n$. The rest of the proof of \eqref{gauss6} is the
same as the proof of \eqref{gauss1}. The proof of \eqref{gauss7}
is now a repetition of the proof of \eqref{gauss4}.
\end{proof}
\subsubsection{Computation of $\int_{\F^*} \gamma_\psi^{-1}(u)
\chi_{(s)}(u) \psi(u) \, d^*u$} \label{int comp} In this
subsection we assume that $\psi$ is normalized. By Lemmas
\ref{reduction} and \ref{reduction psi}, the proof of Theorem
\ref{the formula} amounts to computing
$$\int_{\F^*} \gamma_\psi^{-1}(u) \chi_{(s)}(u) \psi(u) \, d^*u.$$
We write
\begin{equation} \label{decomp tate}\int_{\F^*} \gamma_\psi^{-1}(u)
\chi_{(s)}(u) \psi(u) d^*u \,
=A(\F,\psi,\chi,s)+B(\F,\psi,\chi,s),\end{equation} where
\begin{eqnarray} \label{A} A(\F,\psi,\chi,s) &=&\int_{\Of}
\gamma_\psi^{-1}(u) \chi_{(s)}(u) \, d^*u=\sum_{n=0}^\infty
\bigl(q^{-s}\chi(\pi)\bigr)^{n} \int_{\Of^*} \gamma_\psi^{-1}(\pi^n u)\chi(u) \, du \\
\nonumber &=& \frac {1} {1-\chi^2(\pi)q^{-2s}} \int_{\Of^*}
\gamma_\psi^{-1}(u)\chi(u)\, du+\frac {\chi(\pi)q^{-s}}
{1-\chi^2(\pi)q^{-2s}} \int_{\Of^*} \gamma_\psi^{-1}(\pi u)
\chi(u) \, du, \end{eqnarray} and  \begin{equation} \label{B}
B(\F,\psi,\chi,s)=\int_{\ab u \ab
> 1} \gamma_\psi^{-1}(u) \chi_{(s)}(u) \psi(u) \, d^*u=
\sum_{n=1}^\infty (\chi^{-1}(\pi)q^s)^{n}
J_n(\F,\psi,\chi),\end{equation} where
\begin{equation} \label{J def} J_n(\F,\psi,\chi)=\int_{\Of^*}
\gamma_\psi^{-1}(\pi^{-n} u) \psi(\pi^{-n} u)\chi(u) \, du.
\end{equation}
By \eqref{gamma-1 on pf-1}, if $k \in \N_{odd}$
 \begin{eqnarray}
 \label{j odd}
 \lefteqn{J_k(\F,\psi,\chi)}\\
\nonumber &=& \!
q^{-\half}c_\psi^{-1}(-1)\Biggl(\int_{\Of^*}\psi(\pi
^{-k}u)\chi(u) \, du +\sum_{n=1}^{e+1}q^n \int_{\Of^*}
\int_{\Of^*} \! \psi \bigl(u (\pi
^{-k}-\pi^{1-2n}x^2)\bigr)\chi(u) \, du \, dx
\Biggr),\end{eqnarray} while by \eqref{gamma-1 on of*}, if $k \in
\N_{even}$ then
\begin{eqnarray} \label{j even} \lefteqn{J_k(\F,\psi,\chi)} && \\ \nonumber
&& = c_\psi^{-1}(-1)\Biggl(\int_{\Of^*}\psi(\pi ^{-k}u)\chi(u) \,
du+\sum_{n=1}^e q^n \int_{\Of^*} \int_{\Of^*} \psi \bigl(u (\pi
^{-k}-\pi^{-2n}x^2)\bigr)\chi(u) \, du \, dx \Biggr).
\end{eqnarray}
 Due to Lemma \ref{reduction} and \eqref{decomp tate}, in order to
compute $C_\psi(\etac \otimes \gamma_\psi^{-1},s)$, it is
sufficient to compute
 $A(\F,\psi,\chi,s)$ and $B(\F,\psi,\chi,s)$.

\begin{lem} \label{A for chi^2}
$A(\F,\psi,\chi,s)=0$ unless $\chi^2$ is unramified, in which case
$$A(\F,\psi,\chi,s)=
c_\psi^{-1}(-1)\chi(-1)(1-q^{-1})q^{\frac {m(\chi)}
2}(1-\chi^2(\pi)q^{-2s})^{-1} G(\chi,\psi) \begin{cases}
\chi(\pi)q^{-s} & m(\chi)=2e+1 \\
1 &  m(\chi) \leq 2e \end{cases}.$$
\end{lem}
\begin{proof} First we show that \begin{eqnarray}
 \label{to be used} \lefteqn{\int_{\Of^*}
\gamma_\psi^{-1}(u)\chi(u)\, du} && \\ \nonumber&&
=c^{-1}_\psi(-1) \biggl( \int_{\Of^*}\chi(u) \, du+
\chi(-1)\int_{\Of^*} \chi(x^{-2}) \, dx \sum_{n=1}^e
q^n\int_{\Of^*} \psi(\pi^{-2n}u)\chi(u) \, du
\biggr),\end{eqnarray} and that
\begin{eqnarray} \label{to be usedp} \lefteqn{\int_{\Of^*}
\gamma_\psi^{-1}(\pi u)\chi(u)\, du} &&
\\ \nonumber&&
=q^{-\half}c^{-1}_\psi(-1) \biggl( \int_{\Of^*}\chi(u) \, du+
\chi(-1)\int_{\Of^*} \chi(x^{-2}) \, dx \sum_{n=1}^{e+1} q^n
\int_{\Of^*} \psi(\pi^{1-2n}u)\chi(u) \, du \biggr).\end{eqnarray}
Indeed, by Lemma \ref{gamma comp} \begin{eqnarray*}
\lefteqn{\int_{\Of^*} \gamma_\psi^{-1}(u)\chi(u)\, du} && \\ &=&
c^{-1}_\psi(-1) \int_{\Of^*}
 \biggl(1+ \sum_{n=1}^e q^n\int_{\Of^*} \psi^{-1}(\pi^{-2n}x^2u) dx
\biggr) \chi(u) \, du\\ &=&c^{-1}_\psi(-1)
 \biggl(\int_{\Of^*}\chi(u) \, du+ \sum_{n=1}^e q^n \int_{\Of^*}\int_{\Of^*}\psi(-\pi^{-2n}x^2u)\chi(u) \, du\, dx
\\ &=& c^{-1}_\psi(-1) \biggl( \int_{\Of^*}\chi(u) \, du+
\chi(-1)\int_{\Of^*} \chi(x^{-2}) \, dx \sum_{n=1}^e
q^n\int_{\Of^*} \psi(\pi^{-2n}u)\chi(u) \, du
\biggr).\end{eqnarray*} \eqref{to be usedp} follows in the same
way. If $\chi^2$ is ramified then \eqref{A}, \eqref{to be used}
and \eqref{to be usedp} implies that $A(\F,\psi,\chi,s)=0$. From
now on we assume $\chi^2$ is unramified. The fact, proven in Lemma
\ref{2e}, that $1+\Pf^{2e+1} \subset {\F^*}^2$ implies that
$m(\chi) \leq 2e+1$. Consider first the case $m(\chi) <2e+1$. By
Lemma \ref {usefull}:
\begin{eqnarray*}\int_{\Of^*} \gamma_\psi^{-1}(\pi u)\chi(u) \,
du&=&\gamma_\psi^{-1}(\pi )\int_{\Of^*}
\gamma_\psi^{-1}(u)\chip(u) \, du \\ &=&\gamma_\psi^{-1}(\pi
)\sum_{t \in {\Of^*}/{1+\Pf^{2e}}}\chi_\pi(t)\int_{1+\Pf^{2e}}
\gamma_\psi^{-1}(ut)\chi_\pi(u) \, du \\ &=&\gamma_\psi^{-1}(\pi
)\sum_{t \in {\Of^*}/{1+\Pf^{2e}}}\gamma_\psi^{-1}(t) \chi_\pi(t)
\int_{1+\Pf^{2e}}(u,\pi)_\F \, du=0.\end{eqnarray*} By \eqref{to
be used} and \eqref{gauss1}:
$$\int_{\Of^*} \gamma_\psi^{-1}(u) \, du=c^{-1}_\psi(-1)(1-q^{-1}).$$
 Recalling \eqref{A}, this lemma follows for unramified $\chi$. Assume
now that $1 \leq m(\chi) \leq 2e$. By Lemma \ref{m chi quad}
$m(\chi)$ is even. From \eqref{to be used}, \eqref{gauss2} and
\eqref{gauss4} it follows that
$$\int_{\Of^*} \gamma_\psi^{-1}(u)\chi(u) \, du=\chi(-1)c_\psi^{-1}(-1)(1-q^{-1})q^{\frac {m(\chi)} 2}
G(\chi,\psi).$$ By \eqref{A} we are done in this case also.
Finally, assume $m(\chi)=2e+1$. In this case, by Lemma
\ref{usefull}:
$$\int_{\Of^*} \gamma_\psi^{-1}(u)\chi(u) \, du=\sum_{t \in
{\Of^*}/{1+\Pf^{2e}}}\gamma_\psi^{-1}(t)\chi(t) \int_{1+\Pf^{2e}}
\chi(u) \, du=0.$$ Also, from \eqref{to be usedp}, \eqref{gauss2}
and \eqref{gauss4} it follows that $$\int_{\Of^*}
\gamma_\psi^{-1}(\pi u)\chi(u) \,
du=c_\psi^{-1}(-1)\chi(-1)(1-q^{-1})q^{\frac {m(\chi)} 2}
G(\chi,\psi).$$ With \eqref{A} we now finish.
\end{proof}

\begin{lem} \label{B poly}
$B(\F,\psi,\chi,s)=\sum_{n=1}^{max\bigl(2e+1,m(\chi)\bigr)}
(\chi^{-1}(\pi)q^s)^{n} J_n(\F,\psi,\chi)$
\end{lem}
\begin{proof} Put $k=max \bigl(2e+1,m(\chi)\bigr)$. By \eqref{B}, one should
prove that $J_n(\F,\psi,\chi)=0$ for all $n>k$. Indeed,
$$\int_{\Of^*} \gamma_\psi^{-1}(u) \psi(\pi^{-n}
u) (u,\pi^n)_\F \chi(u) \, du=\sum_{t \in
\Of^*/1+\Pf^k}\chi(t)\gamma_\psi^{-1}(t)
(t,\pi^n)_\F\int_{1+\Pf^k}\psi(\pi^{-n}tu) \, du=0.$$
\end{proof}
\begin{lem}  If $\chi$ is unramified then
$$B(\F,\psi,\chi,s)=c_\psi^{-1}(-1)\Bigl(\chi^{-2e-1}(\pi)q^{(2e+1)s-\half}+(1-q^{-1})\sum_{k=1}^{e}(\chi^{-1}(\pi)q^s)^{2k}\Bigr).$$

\end{lem}
\begin{proof} First note that by Lemma \ref{B poly} we have
$$B(\F,\psi,\chi,s)=\sum_{n=1}^{2e+1} (\chi^{-1}(\pi)q^s)^{n} J_n(\F,\psi,\chi).$$

Assume first that $\F$ of is of odd residual characteristic. In
this case we only want to show that
$J_1(\F,\psi,\chi)=q^{-\half}$. By \eqref{j odd} we have
$$J_1(\F,\psi,\chi)=q^{-\half} \Bigl(\int_{\Of^*}\psi(u\pi^{-1})
\, du+q\int_{\Of^*}
\int_{\Of^*}\psi^{-1}\bigl(u\pi^{-1}(x^2-1)\bigr) \, du \, dx
\Bigr).$$ Using \eqref{gauss1} and \eqref{gauss3} we now get
$$J_1(\F,\psi,\chi)=-q^{\frac {-3} {2}}+q^{\half} \Bigl(\mu\bigl (H(1,\F) \bigr
)(1-q^{-1})-q^{-1}\bigl(\mu(D(1,\F)\bigr) \Bigr).$$ Since from
Lemma \ref{H D comp} it follows that $\mu
\bigl(D(1,\F)\bigr)=1-3q^{-1}$ and $\mu\bigl
(H(1,\F)\bigr)=2q^{-1}$, we have completed the proof of this lemma
for $\F$ of odd residual characteristic. We now assume that the
residue characteristic of  $\F$ is even. Note that by Lemma
\ref{usefull}, if $m(\chi)<2e+1$ and $0 \leq k \leq e-1 $, then:
\begin{eqnarray} \label{j odd o}J_{2k+1}(\F,\psi,\chi)& =& \gamma_\psi^{-1}(\pi)\int_{\Of^*}
\gamma_\psi^{-1}(u) \chi_\pi(u) \psi(u\pi^{-2k-1})\\ \nonumber &=&
\gamma_\psi^{-1}(\pi )\sum_{t \in
{\Of^*}/{1+\Pf^{2e}}}\psi(t\pi^{-{2k+1}})\gamma_\psi^{-1}(t)\chi_\pi(t)\int_{1+\Pf^{2e}}(u,\pi)_\F
\, du=0.\end{eqnarray} Therefore,
\begin{eqnarray} \label{B mid even}
B(\F,\psi,\chi,s)&=&\sum_{k=1}^e (\chi^{-1}(\pi)q^s)^{2k}
\int_{\Of^*} \gamma_\psi^{-1}(u) \psi(u\pi^{-2k}) \, du \\
\nonumber &+& \chi^{-2e-1}(\pi)q^{(2e+1)s}\int_{\Of^*}
\gamma_\psi^{-1}(\pi u) \psi(u\pi^{-2e-1}) \, du.
\end{eqnarray}  Next, we use \eqref{j odd} and  \eqref{gauss1} and note that for $1
\leq k \leq e$:
$$J_{2k}(\F,\psi,\chi)=c_\psi^{-1}(-1) \sum_{n=1}^e q^n
\int_{\Of^*} \biggl(\int_{\Of^*} \psi\bigl((\pi^{-2k}-x^2
\pi^{-2n})u \bigr) \, du \biggr) \, dx.$$ Note that if $k \neq n$
then $\ab (\pi^{-2k}+x^2 \pi^{-2n}) \ab=q^{max(2k,2n)} \geq q^2$.
Therefore if $k \neq n$ the conductor of the character $u \mapsto
\psi\bigl((\pi^{-2k}-x^2 \pi^{-2n})u \bigr)$ is at least 2. Hence,
by  \eqref{gauss3} and by the same arguments used for
\eqref{gauss1} we get:
\begin{eqnarray} \label{I 2k} \int_{\Of^*} \gamma_\psi^{-1}(u)\psi(u\pi^{-2k}) \,
du &=& c_\psi^{-1}(-1)
 q^k \int_{\Of^*} \biggl(\int_{\Of^*}
\psi\bigl((1-x^2)\pi^{-2k}u \bigr) \, du \biggr)\\
\nonumber &=& c_\psi^{-1}(-1)q^k\Bigl(-q^{-1}\mu \bigl(D
(\F,2k)\bigr)+(1-q^{-1})\mu
\bigl(H(\F,2k)\bigr)\Bigr).\end{eqnarray} Similarly,
\begin{equation} \label {I 2e+1}
\int_{\Of^*} \gamma_\psi^{-1}(\pi u) \psi(u\pi^{-2e-1}) \,
du=c_\psi^{-1}(-1)q^{e+\frac 1 2} \Bigl(-q^{-1}\mu \bigl(D
(\F,2e+1)\bigr)+(1-q^{-1})\mu \bigl(H(\F,2e+1)\bigr)\Bigr).
\end{equation}
From Lemma \ref{H D comp} it follows that for $1 \leq k \leq e$
\begin{equation} \label{the hard comp 2k} \mu \bigl(D(2k,\F)\bigr)=0, \, \, \,
\mu \bigl(H(2k,\F)\bigr)=q^{-k},\end{equation} and that
\begin{equation} \label{the hard comp 2e+1}
\mu \bigl(D(2e+1,\F)\bigr)=\frac{q-2}{q^{e+1}}, \, \, \, \mu
\bigl(H(2e+1,\F)\bigr)=\frac{2}{q^{e+1}}. \end{equation}Combining
this with  \eqref{I 2k}, \eqref{I 2e+1} and \eqref{B mid even} the
lemma for the even residual characteristic case follows.

\end{proof}
The last lemma combined with the computation given in Lemma \ref{A
for chi^2} for unramified characters gives explicit formulas for
$A(\F,\psi,\chi,s)$ and $B(\F,\psi,\chi,s)$ for these characters.
A straight forward computation will now give the proof of Theorem
\ref{the formula} for unramified characters.
\begin{lem} \label{B chi quad}Suppose that $m(\chi)>0$ and that $\chi^2$ is unramified. Then:
\begin{eqnarray*}\lefteqn{B(\F,\psi,\chi,s)}  \\&=&c_\psi^{-1}(-1)G(\chi,\psi)\chi(-1)q^{\frac{m(\chi)}{2}-1}\Bigl(
(q-1)\sum_{k=1}^{e-\frac{m(\chi)}{2}}(\chi(\pi)^{-1}q^s)^{2k}-\chi(\pi)^{m(\chi)-2e-2}q^{\bigl(2e+2-m(\chi)\bigr)s}\Bigr),\end{eqnarray*}
where the sum in the right-hand side is to be understood as 0 if
$m(\chi) \geq 2e$.
\end{lem}
\begin{proof} We first assume that $\F$ is of odd residual characteristic.
Since $1+\Pf \subseteq {\Of^*}^2$, $m(\chi)=1$ ($\Of^* / 1+\Pf
\simeq \rf^*$, therefore $u\mapsto (u,\pi)$ is the only one
non-trivial quadratic character of $\Of^*$) . Due to Lemma \ref{B
poly} it is sufficient to show that
$$J_1(\F,\psi,\chi)=-q^{-\half}c_\psi^{-1}(-1)G(\chi,\psi)\chi(-1).$$
By \eqref{j odd} and by \eqref{gauss5} we get
\begin{eqnarray*} J_1(\F,\psi,\chi) &=&
q^{-\half}c_\psi^{-1}(-1)\Bigl(G(\chi,\psi)+q\int_{\Of^*}\int_{\Of^*}
\psi\bigl(u\pi^{-1}(1-x^2)\bigr)\chi(u)\,du \, dx \Bigr) \\ &=&
q^{-\half}c_\psi^{-1}(-1)G(\chi,\psi)\Bigl(1+q\int_{D(1,\F)}\chi^{-1}(1-x^2)
\, dx\Bigr).\end{eqnarray*} Fix $x\in D(1,\F)$. For any $u\in
1+\Pf$ we have $1-x^2u^2=(1-x^2)(1+\frac {x^2b}{1-x^2} \pi)$ for
some $b\in \Of$. This implies that $1-x^2u^2=1-x^2 \, (mod \,
1+\Pf)$. Therefore, $\chi^{-1}(1-x^2)=\chi^{-1}(1-x^2u^2)$. Thus,
\begin{eqnarray*} J_1(\F,\psi,\chi) &=& q^{-\half}c_\psi^{-1}(-1)G(\chi,\psi)\Bigl(1+\sum_{t
\in \Of^*, \, \overline{t} \neq \pm 1} \chi^{-1}(1-t^2)\Bigr)\\
&=& q^{-\half}c_\psi^{-1}(-1)G(\chi,\psi)\sum_{t \in \rf}
{\chi}'(1-t^2), \end{eqnarray*} Where $\chi'$ is the only
non-trivial quadratic character of $\rf^*$. Theorem 1 of Chapter 8
of \cite{I} combined with Exercise 3 of the same chapter implies
that $\sum_{t \in \rf} \chi'(1-t^2)=-\chi(-1)$.

From now on we assume that $\F$ is of even residual
characteristic. We start with the case $m(\chi) \leq 2e$. By Lemma
\ref{B poly} and by \eqref{j odd o} it is sufficient to prove that
$J_{2e+1}(\F,\psi,\chi)=0$ and that for all $1 \leq k \leq e$ we
have:
\begin{equation} \label{J(2k) chi quad}J_{2k}(\F,\psi,\chi)=
c_\psi^{-1}(-1)q^{\frac {m(\chi)} {2}
-1}G(\chi,\psi)\chi(-1)\begin{cases}q-1 &
 \frac {m(\chi)} {2}+k \leq e \\ -1 &  \frac {m(\chi)} {2}+k=e+1 \\ 0 &  \frac {m(\chi)} {2}+k \geq
 e+2
\end{cases}. \end{equation}
By \eqref{j odd} and by \eqref{gauss4} we have
$$J_{2e+1}(\F,\psi,\chi)=
c_\psi^{-1}(-1)\sum_{n=1}^{e+1}q^{n-\half} \Bigl(\int_{\Of^*}
\int_{\Of^*} \psi\bigl((\pi^{-2e-1}-\pi^{1-2n}x^2)u \bigr) \chi(u)
\, du \,dx \Bigr).$$ Since for $n<e+1, \, x \in \Of^*$ the
conductor of the character $u \mapsto
\psi\bigl((\pi^{-2e-1}-\pi^{1-2n}x^2)u \bigr)$ is $2e+1 \neq
m(\chi)$ it follows from \eqref{gauss4} that
$$J_{2e+1}(\F,\psi,\chi)=c_\psi^{-1}(-1)q^{e+\half} \Bigl(\int_{\Of^*}
\int_{\Of^*} \psi\bigl(\pi^{-2e-1}u(1-x^2) \bigr) \chi(u) \, du
\,dx \Bigr).$$ Using \eqref{gauss5} we obtain
$$J_{2e+1}(\F,\psi,\chi)=c_\psi^{-1}(-1)q^{e+\half}G(\chi,\psi)\int_{D(2e+2-m(\chi))}
\chi^{-1}\bigl(\pi^{m(\chi)-2e-1}(1-x^2)\bigr) \, dx.$$ Since by
Lemma \ref{m chi quad}, $m(\chi)$ is even, it follows from Lemma
\ref{H D comp} that ${D(2e+2-m(\chi))}=\varnothing$. We conclude
that $J_{2e+1}(\F,\psi,\chi)=0$ (note that we used only the facts
that $0<m(\chi)<2e+1$ and that $m(\chi)\in \N_{even}$). Suppose
now that $1 \leq k \leq e$ and that $2k \neq m(\chi)$. By \eqref{j
even} we
have:$$J_{2k}(\F,\psi,\chi)=c_\psi^{-1}(-1)\sum_{n=1}^{e}q^{n}
\Bigl(\int_{\Of^*} \int_{\Of^*}
\psi\bigl((\pi^{-2k}-\pi^{-2n}x^2)u \bigr) \chi(u) \, du \,dx
\Bigr).$$ For $n<k$ the conductor of $u \mapsto
\psi\bigl((\pi^{-2k}-\pi^{-2n}x^2)u \bigr)$ is $2k \neq m(\chi)$.
Hence, by  \eqref{gauss4} the first $k-1$ terms in the last
equation vanish and we have
\begin{equation} \label{for k not m}J_{2k}(\F,\psi,\chi)=c_\psi^{-1}(-1)\sum_{n=k}^{e}q^{n}
\Bigl(\int_{\Of^*} \int_{\Of^*}
\psi\bigl((\pi^{-2k}-\pi^{-2n}x^2)u \bigr) \chi(u) \, du \,dx
\Bigr).\end{equation} Since for $n>k$ the conductor of $u \mapsto
\psi\bigl((\pi^{-2k}-\pi^{-2n}x^2)u \bigr)$ is $2n$, it follows
from \eqref{gauss4} and \eqref{gauss5} that for $e\geq
k>\frac{m(\chi)}{2}$
$$J_{2k}(\F,\psi,\chi)=c_\psi^{-1}(-1)q^{k}G(\chi,\psi)\int_{D(2k-m(\chi)+1)}
\chi^{-1}\bigl(\pi^{m(\chi)-2k}(1-x^2)\bigr)\,dx.$$ By Lemma
\ref{H D comp}, in order to prove \eqref{J(2k) chi quad} for
$e\geq k>\frac{m(\chi)}{2}$, it is sufficient to show that
\begin{eqnarray} \label{for J(2k) k>m} \lefteqn{\int_{1+\Pf^{k-\frac{m(\chi)}{2}} \setminus
1+\Pf^{k-\frac{m(\chi)}{2}+1}}
\chi^{-1}\bigl(\pi^{m(\chi)-2k}(1-x^2)\bigr)\,dx} && \\
\nonumber && =q^{\frac{m(\chi)}{2}-k}\chi(-1)
\begin{cases}1-q^{-1}&\frac {m(\chi)} {2}+k \leq e \\ -q^{-1} &  \frac {m(\chi)} {2}+k=e+1 \\ 0
& \frac {m(\chi)} {2}+k \geq e+2 \end{cases}.\end{eqnarray} By
changing $x=1+u\pi^{k-\frac{m(\chi)}{2}}$ we have
$$\int_{D(2k-m(\chi)+1)}
\chi^{-1}\bigl(\pi^{m(\chi)-2k}(1-x^2)\bigr)\,dx=q^{\frac{m(\chi)}{2}-k}\int_{\Of^*}
\chi^{-1} \bigl((-u^2(1+u^{-1}w\pi^{e+\frac {m(\chi)}{2}-k})\bigr)
\, du.$$ (the term $q^{\frac{m(\chi)}{2}-k}$ comes form the change
of variables and from the fact that $dx$ is an additive measure).
Since $\chi^2$ is unramified, we get by changing $u^{-1}w \mapsto
u$:
$$\int_{D(2k-m(\chi)+1)}
\chi^{-1}\bigl(\pi^{m(\chi)-2k}(1-x^2)\bigr)\,dx=q^{\frac{m(\chi)}{2}-k}\chi(-1)\int_{\Of^*}
\chi \bigl(1+u\pi^{e+\frac {m(\chi)}{2}-k}) \, du.$$
\eqref{gauss6} implies now  \eqref{for J(2k) k>m}. We now prove
\eqref{J(2k) chi quad} for $1 \leq k<\frac{m(\chi)}{2}$. Equation
\eqref{for k not m} and arguments we have already used imply in
this case:
\begin{eqnarray*}J_{2k}(\F,\psi,\chi)&=&c_\psi^{-1}(-1)q^{\frac{m(\chi)}{2}}
\Bigl(\int_{\Of^*} \int_{\Of^*}
\psi\bigl(\pi^{-m(\chi)}u(\pi^{m(\chi)-2k}-x^2) \bigr) \chi(u) \,
du \,dx \Bigr)\\
&=& c_\psi^{-1}(-1)q^{\frac{m(\chi)}{2}}G(\chi,\psi)
\int_{\Of^*}\chi(\pi^{m(\chi)-2k}-x^2)\, dx \\
&=& c_\psi^{-1}(-1)q^{\frac{m(\chi)}{2}}G(\chi,\psi)
\chi(-1)\int_{\Of^*}\chi(x^2-\pi^{m(\chi)-2k})\,
dx.\end{eqnarray*} It is left to show that
$$\int_{\Of^*}\chi(x^2-\pi^{m(\chi)-2k})\,
dx=\begin{cases}1-q^{-1}&\frac {m(\chi)} {2}+k \leq e \\ -q^{-1} &
\frac {m(\chi)} {2}+k=e+1 \\ 0 & \frac {m(\chi)} {2}+k \geq e+2
\end{cases}.$$ This is done by changing $x=u+\pi^{\frac {m(\chi)}
{2}-k}$ and by using  \eqref{gauss6}. Finally, we prove
\eqref{J(2k) chi quad} for the case $2k=m(\chi)$: By \eqref{j
even}, \eqref{gauss5} and \eqref{gauss2} we have
\begin{eqnarray*} \lefteqn{J_{m(\chi)}(\F,\psi,\chi)} && \\ &=&
c_\psi^{-1}(-1)\Bigl(G(\chi,\psi)+\sum_{n=1}^{e}q^n
\int_{\Of^*}\int_{\Of^*}\psi
\bigl(u\pi^{-m(\chi)}(1-\pi^{m(\chi)-2n}x^2)\bigr)\chi(u) \, du \,
dx \Bigr) \\ &=& c_\psi^{-1}(-1)G(\chi,\psi)\Bigl(1+ \!
\sum_{n=1}^{\frac
{m(\chi)}{2}-1}q^n\int_{\Of^*}\chi(1-\pi^{m(\chi)-2n}x^2) \,
dx+q^{\frac {m(\chi)}{2}}\int_{D(1,\F)} \! \chi(1-x^2) \, dx
\Bigr).\end{eqnarray*} \eqref{J(2k) chi quad} will follow in this
case once we prove:
\begin{equation} \label{almost done}\int_{\Of^*}\chi(1-\pi^{m(\chi)-2n}x^2) \,
dx=\begin{cases}1-q^{-1}& n \leq e-\frac {m(\chi)}{2}
\\ -q^{-1} & n=e-\frac {m(\chi)}{2}+1
 \\ 0 & n>e-\frac {m(\chi)}{2}+1 \end{cases} \end{equation}
and \begin{equation} \label{done}\int_{\Of^* \setminus 1+\Pf
}\chi(1-x^2) \, dx =\begin{cases}\chi(-1)(1-2q^{-1})& m(\chi) \leq
e
\\ -q^{-1}\bigl(1+\chi(-1)) & m(\chi)=e+1
 \\ 0 &m(\chi) \geq e+2 \end{cases} .\end{equation} As for the first assertion:
Since $\chi^2$ is unramified:
$$\int_{\Of^*}\chi(1-\pi^{m(\chi)-2n}x^2) \, dx=
\int_{\Of^*}\chi(x^{-2}-\pi^{m(\chi)-2n}) \, dx.$$ We change $x
\mapsto x^{-1}$ and then $k=x-\pi^{\frac {m(\chi)}{2}-n}:$
$$\int_{\Of^*}\chi(1-\pi^{m(\chi)-2n}x^2) \,
dx=\int_{\Of^*}\chi\bigl(k^2(1+\pi^{\frac{m(\chi)}{2}-n+e} \omega
k^{-1})\bigr) \, dk.$$ We use once more the fact the $\chi^2$ is
unramified and we change $k=\omega u^{-1}$:
$$\int_{\Of^*}\chi(1-\pi^{m(\chi)-2n}x^2) \,
dx=\int_{\Of^*}\chi(1+\pi^{\frac{m(\chi)}{2}-n+e}u) \, du.$$
\eqref{gauss6} now implies \eqref{almost done}. As for
\eqref{done}, we change $k=1-x$ and then $u=-k^{-1}$ and get:
$$\int_{\Of^* \setminus 1+\Pf }\chi(1-x^2) \, dx =\chi(-1)
\int_{\Of^* \setminus 1+\Pf  }\chi(1+\omega\pi^{e}u) \, du.$$ It
is clear that if $m(\chi) \leq e$ then $$\int_{\Of^* \setminus
1+\Pf }\chi(1+\omega \pi^{e}u) \, du=\mu(\Of^* \setminus
1+\Pf)=1-\frac 2 q.
$$ If  $2e \geq m(\chi) > e$ we write
$$\int_{\Of^* \setminus 1+\Pf
}\chi(1- x^2) \, dx=\chi(-1) \Bigl( \int_{\Of^* }\chi(1+
\omega\pi^{e}u) \, du-\int_{1+\Pf} \chi(1+\omega\pi^{e}u) \, du
\Bigr).$$ By \eqref{gauss6} we have
$$\int_{\Of^* }\chi(1+\omega \pi^{e}u) \,
du=\begin{cases}-q^{-1}&m(\chi) =e+1 \\ 0 & m(\chi)  \geq e+2
\end{cases}.$$ We now compute $\int_{1+\Pf} \chi(1+\omega\pi^{e}u) \,
du$:  We change $u=1+\omega ^{-1}\pi k$. Since $$(1+ \omega \pi^e
u)=(1+\omega \pi^e)(1+ \pi^{e+1}k) \, \, \bigl(mod \,\Pf^
{2e+1}),$$ we have $$ \int_{1+\Pf} \chi(1+\omega \pi^{e}u)
\,du=\chi(3)\int_{\Of} \chi(1+\pi^{e+1}k)\, dk.$$As in the proof
of \eqref{gauss6}, since $m(\chi) \leq 2e$, we conclude that $k
\mapsto \chi(1+\pi^{e+1}k)$ is a character of $\Of$. It is trivial
if $m(\chi)=e+1$ and non-trivial otherwise. Also, since $-3 \in
1+\Pf^{2e}$, we have $\chi(3)=\chi(-1)$. Thus,$$ \int_{1+\Pf}
\chi(1+\omega\pi^{e}u) \,du=\begin{cases}\chi(-1)q^{-1}&
m(\chi)=e+1
\\ 0 &m(\chi) \geq e+2 \end{cases}.$$

We move to the case $m(\chi)=2e+1$. By Lemma \ref{B poly} it is
sufficient to show that $J_{n}(\F,\psi,\chi)=0$ for all $2 \leq n
\leq 2e+1$ and that
$$J_{1}(\F,\psi,\chi)=-q^{e-\half}c_\psi^{-1}(-1)G(\chi,\psi)\chi(-1).$$
From the fact that $m(\chi)=2e+1$ it follows that for all $1\leq k
\leq e$:
\begin{eqnarray} \label{j_2k m=2e+1} J_{2k}(\F,\psi,\chi)&=&\int_{\Of^*}\gamma_\psi^{-1}(u)
\chi(u)\psi(\pi^{-2k}u) \, du \\ \nonumber &=& \sum_{t \in
{\Of^*}/{1+\Pf^{2e}}}\gamma_\psi^{-1}(t)\chi(t)\psi^{-1}(\pi^{-2k}
t)\int_{1+\Pf^{2e}}\chi(u) \, du=0.\end{eqnarray} Next we show
that $J_{2e+1}(\F,\psi,\chi)=0$. By  \eqref{j odd},
\eqref{gauss4}, \eqref{gauss5} and \eqref{gauss2}:
\begin{eqnarray} \label{J_2k1+ m=2e+1} \lefteqn {J_{2e+1}(\F,\psi,\chi)} && \\
\nonumber &=& q^{-\half}c_\psi^{-1}(-1)\Bigl(
G(\chi,\psi)+\sum_{n=1}^{e+1}q^{n} \int_{\Of^*} \int_{\Of^*}
\psi\bigl(\pi^{-2e-1}u(1-\pi^{2(e+1-n)}x^2) \bigr)\chi(u) \, du
\,dx \Bigr) \\ \nonumber &=&
q^{-\half}c_\psi^{-1}(-1)G(\chi,\psi)\Bigl( 1 +
\sum_{n=1}^{e}q^{n} \int_{\Of^*} \chi(1-\pi^{2(e+1-n)}x^2) \,dx
\\ \nonumber
 && \hskip108pt+q^{e+1}\int_{D(1,\F)}\chi(1-x^2) \, dx \Bigr).\end{eqnarray} As
in the proof of \eqref{almost done} (using $u=x-\pi^{e+1-n}$
rather then $k=x-\pi^{\frac {m(\chi)}{2}-n}$) for all $1 \leq n
\leq e:$
$$\int_{\Of^*} \chi(1-\pi^{2(e+1-n)}x^2) \, dx=\int_{\Of^*}
\chi(1+\pi^{2e+1-n}u) \, du. $$ \eqref{gauss6} implies now that
$$\int_{\Of^*} \chi(1-\pi^{2(e+1-n)}x^2) \, dx=\begin{cases}-q^{-1}&n=1 \\ 0 &
n>1 \end{cases}.$$ It is clear now that if we prove that
\begin{equation} \label{chi(1-x^2)=0}\int_{D(1,\F)}\chi(1-x^2) \, dx =0,\end{equation}we will conclude
that $J_{2e+1}(\F,\psi,\chi)=0$. As in the proof of \eqref{done}
$$\int_{D(1,\F)}\chi(1-x^2) \, dx=\chi(-1) \Bigl( \int_{\Of^*  }\chi(1+w\pi^{e}u) \,
du-\int_{1+\Pf} \chi(1+w\pi^{e}u) \, du \Bigr),$$ and
$$\int_{1+\Pf} \chi(1+w\pi^{e}u) \, du=0.$$ We show that $$\int_{\Of^*  }\chi(1+w\pi^{e}u) \,du=0,$$ although
the map $u \mapsto \chi(1+w\pi^{e}u)$ is not a character of $\Of$:
$$\int_{\Of^*  }\chi(1+w\pi^{e}u) \,du=\sum_{t\in \Of^*
/ 1+\Pf} \int_{1+\Pf} \chi(1+w\pi^{e}tu) \,du.$$
 For any $t\in \Of^*, \, k\in \Of$: $$1+w\pi^{e}t(1+\pi k)=
 (1+w\pi^{e}t)(1+w\pi^{e+1}tk) \, (mod \, 1+\Pf^{2e+1}).$$ Hence, by
 changing $u=1+\pi k$ we get
$$\int_{\Of^*  }\chi(1+w\pi^{e}u) \,du=q^{-1}\sum_{t\in \Of^*
/ 1+\Pf}\chi(1+w\pi^{e}t) \int_{\Of} \chi(1+w\pi^{e+1}tk) \,dk.$$
Since for all $t \in \Of^*$, $k \mapsto \chi(1+w\pi^{e+1}tk)$ is a
non-trivial character of $\Of$, all the integrations in the right
wing vanish.

Finally, we compute $J_{2k+1}(\F,\psi,\chi)$ for $0 \leq k \leq
e-1$: By \eqref{j odd} and  \eqref{gauss4} we observe:
\begin{eqnarray} \label{J_2k+1
m_chi=2e+1} J_{2k+1}(\F,\psi,\chi) &=&
q^{e+\half}c_\psi^{-1}(-1)\int_{\Of^*}\int_{\Of^*}
\psi\bigl(\pi^{2e+1}u(\pi^{2(e-k)}-x^2) \bigr) \chi(u) \, du \, dx
\\ \nonumber &=& q^{e+\half}c_\psi^{-1}(-1)G(\chi,\psi)\int_{\Of^*}
\chi(\pi^{2(e-k)}-x^2) \,dx.\end{eqnarray} We finish by showing:
$$\int_{\Of^*}
\chi(\pi^{2(e-k)}-x^2) \,dx=\begin{cases}-\chi(-1)q^{-1}&k=0 \\ 0
& k>0 \end{cases}.$$ We note that $$\int_{\Of^*}
\chi(\pi^{2(e-k)}-x^2) \,dx=\chi(-1)\int_{\Of^*}
\chi(x^2-\pi^{2(e-k)}) \,dx.$$ We now change $u=x-\pi^{e-k}$ and
use  \eqref{gauss6}.
\end{proof}
The last lemma combined with the computation given in Lemma \ref{A
for chi^2} for these characters gives explicit formulas for
$A(\F,\psi,\chi,s)$ and $B(\F,\psi,\chi,s)$  for a character $\chi$ which staisfies
$m(\chi)>m(\chi^2)=0$. By a simple computation one can now finish the
proof of Theorem \ref{the formula} for a ramified $\chi$ provided that $\chi^2$ is unramified.
\begin{lem} If $\chi^2$ is ramified then $B(\F,\psi,\chi,s)$
is a non-zero monomial in $q^s$.
\end{lem}
\begin{proof} The fact that $B(\F,\psi,\chi,s)$
is non-zero follows from the fact, proven in Lemma \ref{A for
chi^2}, that  $A(\F,\psi,\chi,s)=0$ and from the fact, following
from Lemma \ref{reduction}, that $$\int_{F^*} \gamma_\psi^{-1}(u)
\chi_{(s)}(u) \psi(u) d^*u$$ can not be identically 0 as a
function of $s$. Assume first that $m(\chi)>2e+1$. By Lemma \ref{B
poly} it is enough to show that $J_n(\F,\psi,\chi)=0$ for all $1
\leq n<m(\chi)$: Put $k=max(n,2e+1)$. We have
$$J_n(\F,\psi,\chi)=\sum_{t \in
{\Of^*}/{1+\Pf^k}}\gamma_\psi^{-1}(\pi^{-n} t)\chi(t)\psi(\pi^{-n}
t)\int_{1+\Pf^k}\chi( u) \, du.$$ Since $k<m(\chi)$ that last
integral vanishes. This shows that
$$B(\F,\psi,\chi,s)=\chi^{-m(\chi)}(\pi)q^{m(\chi) s} J_{m(\chi)}(\F,\psi,\chi).$$
In order to complete the proof of this lemma for $\F$ of odd
residual characteristic it is left to consider the case
$m(\chi)=1$. In fact, there is nothing to prove here since by
Lemma \ref{B poly} we
get:$$B(\F,\psi,\chi,s)=\chi^{-1}(\pi)q^sJ_1(\F,\psi,\chi).$$ From
now on we assume that $\F$ is of even residual and that $m(\chi)
\leq 2e+1$. Assume first that $m(\chi)<2e+1$ and that $m(\chi)$ is
odd. By Lemma \ref{B poly}, once we prove $J_n(\F,\psi,\chi)=0$
for all $n<2e+1$ we will
conclude:$$B(\F,\psi,\chi,s)=\chi(\pi)^{-2e-1}q^{(2e+1)s}
J_{2e+1}(\F,\psi,\chi).$$ By \eqref{j odd o},
$J_n(\F,\psi,\chi)=0$ for all $0<n<2e+1$, $n\in \N_{odd}$. Let $1
\leq k \leq e$. Since $m(\chi)$ is odd, by \eqref{j even} and
\eqref{gauss4} we get:
$$J_{2k}(\F,\psi,\chi)=
c_\psi^{-1}(-1)\sum_{n=1}^{e}q^n \int_{\Of^*}\int_{\Of^*}\psi
\bigl((\pi^{-2k}-\pi^{-2n}x^2)u\bigr)\chi(u) \, du \, dx.$$ For $n
\neq k$, $x \in \Of^*$, the conductor of the character $u \mapsto
\psi \bigl( (\pi^{-2k}-\pi^{-2n}x^2)u\bigr)$ is $max(2k,2n)$.
Again, since $m(\chi)$ is odd, \eqref{gauss4} implies
$$J_{2k}(\F,\psi,\chi)=c_\psi^{-1}(-1)q^k \int_{\Of^*}\int_{\Of^*}\psi
\bigl(u\pi^{-2k}(1-x^2)\bigr)\chi(u) \, du \, dx.$$ Using
\eqref{gauss5} we conclude that if $2k<m(\chi)$ then
$J_{2k}(\F,\psi,\chi)=0$, and if $2k>m(\chi)$ then
$$J_{2k}(\F,\psi,\chi)=c_\psi^{-1}(-1)q^kG(\chi,\psi)\int_{D(2k-m(\chi)+1,\F)}
\chi^{-1}\bigl(\pi^{m(\chi)-2k}(1-x^2)\bigr) \,dx.$$ Lemma \ref{H
D comp} implies that $D(2k-m(\chi)+1,\F)=\varnothing$. Therefore,
we have shown: $J_{2k}(\F,\psi,\chi)=0$ for all $1 \leq k \leq e$.

We now assume that $2 \leq m(\chi) \leq 2e$ and that $m(\chi)$ is
even. Again, $J_n(\F,\psi,\chi)=0$ for all $n< 2e+1$, $n \in
\N_{odd}$. In this case, the proof of the fact that
$J_{2e+1}(\F,\psi,\chi)=0$ is a repetition of the proof that
$J_{2e+1}(\F,\psi,\chi)=0$ given in Lemma \ref{B chi quad} for the
case $ m(\chi) \leq 2e$. By Lemma \ref{B poly}, it is enough now
to prove that if $1 \leq k \leq e$, and if $2k \neq
2e+2m(\chi^2)-m(\chi)$ then $J_{2k}(\F,\psi,\chi)=0$ to conclude
that
$$B(\F,\psi,\chi,s)=\chi(\pi)^{-\bigl(2e+2m(\chi^2)-m(\chi)\bigr)}q^{\bigl(2e+2m(\chi^2)-m(\chi)\bigr)s}J_{\bigl(2e+2m(\chi^2)-m(\chi)\bigr)}(\F,\psi,\chi).$$
Assume that $2k<m(\chi)$. As in the proof of \eqref{J(2k) chi
quad} for the case $2k<m(\chi)$ we have:
$$J_{2k}(\F,\psi,\chi)=\chi(-1)c_\psi^{-1}(-1)q^{\frac{m(\chi)}{2}}G(\chi,\psi)
\int_{\Of^*}\chi^{-1}(x^2-\pi^{m(\chi)-2k})\, dx.$$ By changing
$x=u+\pi^{\frac {m(\chi)}{2}-k}$ and then $ wu^{-1} \mapsto u$ we
obtain
$$J_{2k}(\F,\psi,\chi)=\chi(w^{-2})\chi(-1)c_\psi^{-1}(-1)q^{\frac{m(\chi)}{2}}G(\chi,\psi)
\int_{\Of^*}\chi^{2}(u)\chi^{-1}(1+\pi^{e+\frac
{m(\chi)}{2}-k}u)\, du.$$ It follows from \eqref{gauss7} that if
$2k \neq 2e+2m(\chi^2)-m(\chi)$ then $J_{2k}(\F,\psi,\chi)=0$.
Assume now $2k>m(\chi)$.  As in the proof of \eqref{J(2k) chi
quad} for the case $2k>m(\chi)$ we have:
$$J_{2k}(\F,\psi,\chi)=c_\psi^{-1}(-1)q^{k}G(\chi,\psi)
\int_{1+\Pf^{k-\frac{m(\chi)}{2}} \setminus
1+\Pf^{k-\frac{m(\chi)}{2}+1}}
\chi^{-1}\bigl(\pi^{m(\chi)-2k}(1-x^2)\bigr)\,dx.$$ Changing
$x=1+u\pi^{k-\frac {m(\chi)}{2}}$ and then $ wu^{-1}  \mapsto u$
we reduce this computation to the previous case. As in the proof
of \eqref{J(2k) chi quad} for the case $2k=m(\chi)$ we have
\begin{eqnarray*} \lefteqn {J_{m(\chi)}(\F,\psi,\chi)} && \\ &&=
c_\psi^{-1}(-1)G(\chi,\psi)\Bigl(1+ \! \!  \sum_{n=1}^{\frac
{m(\chi)}{2}-1} \! q^n\int_{\Of^*}\chi^{-1}(1-\pi^{m(\chi)-2n}x^2)
\, dx \\ && \hskip102pt +q^{\frac {m(\chi)}{2}}\int_{D(1,\F)} \!
\chi^{-1}(1-x^2) \, dx \Bigr).\end{eqnarray*} It is sufficient now
to show that if $m(\chi) \neq 2e+2m(\chi^2)-m(\chi)$ then
\begin{equation} \label{sum=-1}\sum_{n=1}^{\frac
{m(\chi)}{2}-1}q^n\int_{\Of^*}\chi^{-1}(1-\pi^{m(\chi)-2n}x^2) \,
dx=-1,\end{equation} and \begin{equation} \label{after
sum=-1}\int_{\Of^* \setminus 1+\Pf}\chi^{-1}(1-x^2) \,
dx=0.\end{equation} We prove \eqref{sum=-1}. Let $1 \leq n
 \leq \frac {m(\chi)}{2}-1$. We change $x^{-1}=x'$, $u'=x'-\pi^{\frac
{m(\chi)}{2}-n}$. Note that $x'=u'(1+u'^{-1}\pi^{\frac
{m(\chi)}{2}-n})$ and that
$x'^2-\pi^{m(\chi)-2n}=u'^2(1+u'^{-1}\omega\pi^{e+\frac
{m(\chi)}{2}-n})$. Thus,
\begin{eqnarray*}\chi^{-1}(1-\pi^{m(\chi)-2n}x^2)&=& \chi^{-1}(x^2)\chi^{-1}
(x^{-2}-\pi^{m(\chi)-2n})=\chi^2(x')\chi^{-1}
(x'^2-\pi^{m(\chi)-2n}) \\ &=& \chi^2(u')
\chi^2(1+u'^{-1}\pi^{\frac
{m(\chi)}{2}-n})\chi^{-1}(u'^2)\chi^{-1}(1+u'^{-1}\omega
\pi^{e+\frac {m(\chi)}{2}-n}).\end{eqnarray*} Now setting
$u'=u^{-1}$ implies:
$$\int_{\Of^*}\chi^{-1}(1-\pi^{m(\chi)-2n}x^2) \,dx=
\int_{\Of^*}\chi^2(1+\pi^{\frac
{m(\chi)}{2}-n}u)\chi^{-1}(1+\pi^{e+\frac {m(\chi)}{2}-n}\omega
u)\, du.$$ The fact that $u\mapsto
\phi_1(u)=\chi^{-1}(1+\pi^{e+\frac {m(\chi)}{2}-n} \omega u)$ is
an additive character of $\Of$ follows from the same argument used
in the proof of \eqref{gauss6}. However, by the same argument, $u
\mapsto \phi_2(u)=\chi^2(1+\pi^{\frac {m(\chi)}{2}-n}u)$ is an
additive character of $\Of$ if and only if $m(\chi)-2n \geq
m(\chi^2)$. We first show that if $m(\chi)>2n> m(\chi)-m(\chi ^2)$
then

\begin{equation} \label{cor1} \int_{\Of^*} \phi_1(u)\phi_2(u) \, du=0.\end{equation} Define
$k=2n-m(\chi)+m(\chi^2)$. We have:
$$\int_{\Of^*}\phi_1(u)\phi_2(u) \, du= \sum_{t\in \Of^* / 1+\Pf^k}
\int_{1+\Pf^k}\chi^2(1+\pi^{\frac
{m(\chi)}{2}-n}tu)\chi^{-1}(1+\pi^{e+\frac {m(\chi)}{2}-n}\omega
tu)\, du.$$ In the right hand side we change $u=1+x\pi^k$ , $x \in
\Of$. We have
$$1+ut\pi^{\frac{m(\chi)}{2}-n}=(1+t\pi^{\frac{m(\chi)}{2}-n})(1+tx\pi^{\frac{m(\chi)}{2}-n+k})
\, (mod \, 1+\Pf^{m(\chi^2)}),$$ and
$$1+u \omega t\pi^{e+\frac{m(\chi)}{2}-n}=(1+ \omega t\pi^{e+\frac{m(\chi)}{2}-n})
(1+t \omega x\pi^{e+\frac{m(\chi)}{2}-n+k}) \, (mod \,
1+\Pf^{m(\chi)}).$$ Thus, we get
\begin{equation} \label{cor2} \int_{\Of^*}\phi_1(u)\phi_2(u) \, du=q^{-k} \sum_{t\in \Of^* / 1+\Pf^k}
\beta_1^t(1)\beta_2^t(1)\int_{\Of}\beta_1^t(x)\beta_2^t(x) \,
dx,\end{equation} where $\beta_1^t(x)=\chi^{-1}(1+t \omega
x\pi^{e+\frac{m(\chi)}{2}-n+k})$ and $ \beta_2^t(x)=\chi^{2}(1+t
 x\pi^{\frac{m(\chi)}{2}-n+k})$. Since
$$m(\chi)+2k-2n=2n-m(\chi)+2m(\chi^2)>m(\chi^2)$$ and
$$2e+m(\chi)+2k-2n>2e \geq m(\chi),$$ we conclude, using the proof of
\eqref{gauss6}, that for $t\in \Of^*$, $\beta_1^t$ and $\beta_2^t$
are two characters $\Of$. $\beta_2^t$ is a non-trivial character
and its conductor its $\frac {m(\chi)}{2}-n$. $\beta_1^t$ is
trivial if $\frac{3m(\chi)}{2}-n-e-m(\chi^2) \leq 0$, otherwise
its conductor is $\frac {3m(\chi)}{2}-n-e-m(\chi^2)$. Since we
assume $m(\chi)-m(\chi^2) \neq e$ we observe that $\beta_1^t$ and
$\beta_2^t$  have different conductors. In particular, $\beta_1^t
\beta_2^t$ is a non-trivial character of $\Of$. \eqref{cor2} now
implies \eqref{cor1}.

We now assume that $2n \leq m(\chi)-m(\chi^2)$. $\phi_1$ and
$\phi_2$ are now both characters of $\Of$. $\phi_1$ is trivial if
$n \leq e-\frac {m(\chi)}{2}$. Otherwise its conductor is $n+\frac
{m(\chi)}{2}-e$. $\phi_2$ is trivial if $n \leq \frac
{m(\chi)}{2}-m(\chi^2)$. Otherwise its conductor is $n-\frac
{m(\chi)}{2}+m(\chi^2)$. Recall that we assumed $m(\chi^2) \neq
m(\chi)-e$. Hence, if either $\phi_1$ or $\phi_2$ are non-trivial
then the conductor of $\phi_1 \phi_2$ is $n+max(m(\chi^2)-\frac
{m(\chi)}{2},\frac {m(\chi)}{2}-e)$. From \eqref{gauss1} it
follows now that if $m(\chi^2)
> m(\chi)-e $ then
$$\int_{\Of^*}\chi^2(1+\pi^{\frac
{m(\chi)}{2}-n}u)\chi^{-1}(1+\pi^{e+\frac {m(\chi)}{2}-n}u)\,
du=\begin{cases}1-q^{-1}& 1 \leq n \leq \frac
{m(\chi)}{2}-m(\chi^2)
\\ -q^{-1} & n=\frac
{m(\chi)}{2}-m(\chi^2)+1
 \\ 0 & n>\frac
{m(\chi)}{2}-m(\chi^2)+1 \end{cases},$$ and if $m(\chi^2) <
m(\chi)-e $ then
$$\int_{\Of^*}\chi^2(1+\pi^{\frac
{m(\chi)}{2}-n}u)\chi^{-1}(1+\pi^{e+\frac {m(\chi)}{2}-n}u)\,
du=\begin{cases}1-q^{-1}& 1 \leq n \leq e-\frac {m(\chi)}{2}
\\ -q^{-1} & n=e-\frac {m(\chi)}{2}+1
 \\ 0 & n>e-\frac {m(\chi)}{2}+1 \end{cases}.$$ Both cases imply
 \eqref{sum=-1}. We now prove \eqref{after sum=-1}. As in the proof of \eqref{done},
we change $u'=1-x$ and then $u=u'^{-1}$. Since
$1-x^2=-u'^2(1-\omega \pi^eu'^{-1})$, we get
$$\int_{\Of^* \setminus 1+\Pf} \! \chi^{-1}(1-x^2) \,
dx=\chi(-1)\Bigl(\int_{\Of^*} \chi^2(u)\chi^{-1}(1-u \omega \pi^e)
\, du-\int_{1+\Pf} \! \chi^2(u)\chi^{-1}(1-u \omega \pi^e) \, du
\Bigr).$$ Since $m(\chi)-e \neq m(\chi^2)$,  \eqref{gauss7}
implies
$$\int_{\Of^*} \chi^2(x)\chi^{-1}(1- \omega x\pi^e) \, dx=0.$$
If $m(\chi^2)>m(\chi)-e$ we have $$\int_{1+\Pf}
\chi^2(x)\chi^{-1}(1-x \omega \pi^e) \, dx=\sum_{t \in
1+\Pf^{m(\chi)-e}} \chi^2(t)\chi^{-1}(1-\omega
t\pi^e)\int_{1+\Pf^{m(\chi)-e}} \chi^2 (x) \, dx=0.$$ If
$m(\chi^2)<m(\chi)-e$ we have
$$\int_{1+\Pf} \chi^2(x)\chi^{-1}(1-\omega x\pi^e) \, dx=\sum_{t \in
1+\Pf^{m(\chi^2)}} \chi^2(t)\chi^{-1}(1- \omega
t\pi^e)\int_{1+\Pf^{m(\chi^2)}} \chi^{-1} (1- \omega xt\pi^e) \,
dx.$$ The proof of \eqref{after sum=-1} will be finished once we
show that for all $t\in \Of^*$:
$$\int_{1+\Pf^{m(\chi^2)}} \chi^{-1} (1- \omega xt\pi^e)=0.$$ We
change $x=1+u\pi^{m(\chi^2)}$. Since
$$1-\omega xt\pi^e=(1-\omega t\pi^e)(1- \omega tu\pi^{e+m(\chi^2)})\, (mod \,
1+\Pf^{2e}),$$ we get
$$\int_{1+\Pf^{m(\chi^2)}} \chi^{-1} (1- \omega xt\pi^e)=
q^{-m(\chi^2)}\chi^{-1}(1+t\pi^e)\int_{\Of}  \chi^{-1} (1-\omega
ut\pi^{e+m(\chi^2)}) \, du.$$ The integral on the right hand side
vanishes since $u \mapsto \chi^{-1} (1-\omega
ut\pi^{e+m(\chi^2)})$ is a non-trivial character of $\Of$.

It is left to prove the lemma for the case $m(\chi)=2e+1$. By
Lemma \ref{B poly} and by \eqref{j_2k m=2e+1}, we only have to
show that for all $0 \leq k \leq e$, unless $k+1=m(\chi^2)$,
 $J_{2k+1}(\F,\psi,\chi)=0$, to conclude that
$$B(\F,\psi,\chi,s)=\chi(\pi)^{-(2m(\chi^2)-1)s}q^{(2m(\chi^2)-1)s} J_{2m(\chi^2)-1}(\F,\psi,\chi).$$
Assume first that $0 \leq k<e$. The same argument we used for
\eqref{J_2k+1 m_chi=2e+1} shows that $$
J_{2k+1}(\F,\psi,\chi)=q^{e+\half}c_\psi^{-1}(-1)G(\chi,\psi)\int_{\Of^*}
\chi^{-1}(\pi^{2(e-k)}-x^2) \,dx.$$ We change $x=\pi^{e-k}-u'$. We
have $\pi^{2(e-k)}-x^2=-u'^2(1-u'^{-1}\pi^{2e-k}\omega)$. Thus,
$$\int_{\Of^*} \chi^{-1}(\pi^{2(e-k)}-x^2) \,dx=
\chi^{-1}(-1)\int_{\Of^*}\chi^{-1}(u'^2)\chi^{-1}(1-u'^{-1}\pi^{2e-k}
\omega) \, du'.$$ Next we change $u=-u'^{-1}\omega$ and obtain:
$$\int_{\Of^*} \chi^{-1}(\pi^{2(e-k)}-x^2)
\,dx=\chi(-\omega^{-2})\int_{\Of^*}
\chi^2(u)\chi^{-1}(1+u\pi^{2e-k}) \, du.$$ \eqref{gauss7} implies
now that if $k+1 \neq m(\chi^2)$ then
$$\int_{\Of^*} \chi^{-1}(\pi^{2(e-k)}-x^2) \,dx=0.$$
Since $m(\chi)=2e+1$ implies $1 \leq m(\chi^2) \leq e+1$ it is
left to show that if  $m(\chi^2) < e+1$ then $J_{2e+1}=0$. As in
the proof of \eqref{J_2k1+ m=2e+1}
\begin{eqnarray*} J_{2e+1}(\F,\psi,\chi)= q^{-\half}c_\psi^{-1}(-1)G(\chi,\psi)\Bigl(
1&+&\sum_{n=1}^{e}q^{n} \int_{\Of^*}
\chi^{-1}(1-\pi^{2(e+1-n)}x^2) \,dx \\
&+&q^{e+1}\int_{D(1,\F)}\chi^{-1}(1-x^2) \, dx
\Bigr).\end{eqnarray*} As we have seen before:
\begin{eqnarray*} \lefteqn{\int_{D(1,\F)}\chi^{-1}(1-x^2) \,
dx} && \\ &&
=\chi(-\omega^{-2})\Bigl(\int_{\Of^*}\chi^2(u)\chi^{-1}(1+w\pi^eu)
\, du-\int_{1+\Pf}\chi^2(u)\chi^{-1}(1+w\pi^eu) \, du
\Bigr).\end{eqnarray*} We show that the last two integrals vanish.
Similar steps to those used in the proof of \eqref{chi(1-x^2)=0}
shows that
\begin{eqnarray*} \lefteqn{\int_{\Of^*}\chi^2(u)\chi^{-1}(1+\omega\pi^eu) \,
du} && \\ &=& \sum_{t\in \Of^* / 1+\Pf^{m(\chi^2)}}\chi^2(t)
\int_{1+\Pf^{m(\chi^2)}}\chi^{-1}(1+w\pi^etu) \, du \\
&=& \sum_{t\in \Of^* / 1+\Pf^{m(\chi^2)}}\chi^2(t)
\chi^{-1}(1+w\pi^et) \int_{\Of}\chi^{-1}(1+w\pi^{e+m(\chi^2)}k) \,
dk.\end{eqnarray*} Since $\frac {m(\chi)}{2} \leq
e+m(\chi^2)<m(\chi)$, \eqref{gauss6} implies that all the
integrals in the last sum vanish. The fact that
$$\int_{1+\Pf}\chi^2(u)\chi^{-1}(1+w\pi^eu) \, du = 0$$ is proven
by a similar argument. Last, we have that to show that $1 \leq
m(\chi^2) \leq e$ implies that
$$\sum_{n=1}^{e}q^{n} \int_{\Of^*} \chi^{-1}(1-\pi^{2(e+1-n)}x^2)
\,dx =-1.$$ This is done by a similar argument to the one we use
in the proof of \eqref{sum=-1}.
\end{proof}
In Lemma \ref{A for chi^2} we proved that if $\chi^2$ is ramified
then $A(\F,\psi,\chi,s)=0$. Thus, the last lemma completes the
proof of Theorem \ref{the formula} for these characters.
\subsection{Real case} \label{real case}
\subsubsection{Notations and main result} \label{real note} Any
non-trivial character of $\R$ has the form $\psi_b(x)=e^{ibx}$ for
some $b$ in $\R^*$. Any character of $\R^*$ has the form
$\chi=\chi_{x,n}(a)=\bigl(sign(a)\bigr)^n \ab a \ab^x$ for some $x
\in \C, \, n \in \{0,1\}$. We may assume that $x=0$, i.e., that
$\chi$ is either the trivial character of the sign character. Thus,
$\chi_{(s)}$ will denote the character $a \mapsto
\bigl(sign(a)\bigr)^n \ab a \ab^s$ where $ n \in \{0,1\}$.
\begin{lem} \label{the real comp}
$$C_{\psi_{b}}(\etac \otimes \gamma_{\psi_a}^{-1},s)=\frac
{e^{-\frac {i \pi \chi(-1) sign(a)} 4}} {2\pi} \frac {\Gamma(\frac
{1-s} 2+ \frac {sign(ab)\chi(-1)} 4)\Gamma((\frac {1+s} 2- \frac
{sign(ab)\chi(-1)} 4)} {\Gamma(s)}.$$
\end{lem}
We shall prove Lemma \ref{the real comp} in the next subsection.
An immediate corollary of this lemma and of the classical
duplication formula $${\Gamma(-\frac s 2+ \frac 1 4)}\Gamma(-\frac
s 2+\frac 3 4)=2^{s+\half}\sqrt{\pi} {\Gamma(-s +\half)}$$ is the
following:
\begin{thm} \label{real thm}
$$C_{\psi_a}(\chi \otimes \gamma_{\psi_a}^{-1},s)=\frac {e^{-\frac {i \pi
\chi(-1)sign(a)} 4}2^{s-\half}} {\pi^s} \frac{L(\chi,s+\half)}
{L(\chi,-s+\half)} \frac{L(\chi^2,-2s+1)} {L(\chi^2,2s)}.$$
\end{thm}
Recall that the local L-function for $\R$ is defined by
$$L_{\R}(\chi_{0,n},s)=\begin{cases} \pi^{\frac {-s}2}\Gamma(\frac s 2) &
n=0 \\ \pi^{\frac {-(s+1)}2}\Gamma(\frac {s+1} 2) & n=1
\end{cases}.$$
\subsubsection{Proof of Lemma
\ref{the real comp}} We shall see that the computation of the
local coefficient for $\overline{SL_2(\R)}$ is done by the same
methods as the computation for $SL_2(\R)$. Namely, we shall use
the Iwasawa decomposition,
$\overline{SL_2(\R)}=\overline{B(\R)}\overline{SO_2(\R)}$, and the
fact that, as an inverse image of a commutative group,
$\overline{SO_2(\R)}$ is commutative; see Section \ref{kubota}.

All the characters of $\overline{SO_2(\R)}$ are given by
$$\vartheta_n \bigl(k(t),\epsilon
\bigr)=e^{in\phi^{-1}\bigl(k(t),\epsilon \bigr)},$$ where $2n \in
\Z$ and $\phi:\R / 4\pi \Z \rightarrow \overline{SO_2(\R)}$ is the
isomorphism from Lemma \ref{so2}. Denote by $I(\chi \otimes
\gamma_{\psi_a}^{-1},s)_n$ the subspace of $I(\chi \otimes
\gamma_{\psi_a}^{-1},s)$ on which $\overline{SO_2(\R)}$ acts by
$\vartheta_n$. From the Iwasawa decomposition it follows that
$$\dim \, \,I(\chi \otimes \gamma_{\psi_a}^{-1},s)_n \leq 1.$$
Furthermore, since $\overline{B(\R)} \cap\overline{SO_2(\R)}$ is a
cyclic group of order 4 which is generated by $(-I_2,1)$, it
follows that $I(\chi \otimes \gamma_{\psi_a}^{-1},s)_n \neq \{0\}$
if and only if
$e^{in\phi^{-1}(-I_2,1)}=\chi(-1)\gamma_{\psi_a}^{-1}(-1)$. Recall
that
$$\gamma_{\psi_a}(y)=\begin{cases} 1 & y>0
 \\ -sign(a)i  & y<0  \end{cases}.$$ Therefore, the last condition
is equivalent to $e^{in\pi}=i\chi(-1)sign(a)$. Thus, we proved:
\begin{lem} \label{n}$\dim I(\chi \otimes \gamma_{\psi_a}^{-1},s)_n =1$ if and
only if $n \in \frac{\chi(-1)sign(a)}{2}+2\Z$.
\end{lem}
Let $n$ be as in Lemma \ref{n}. Define $f_{s,\chi,a,n}$ to be the
unique element of $I(\chi_{(s)} \otimes \gamma_{\psi_a}^{-1})_n$
which satisfies $f_{s,\chi,a,n}(I_2,1)=1$. Since $A(s)I(\chi_{(s)}
\otimes \gamma_{\psi_a}^{-1})_n \subset I(\eta(\chi,-s) \otimes
\gamma_{\psi_a}^{-1})_n$, it follows that $$\frac
{A_s(f_{s,\chi,a,n})}{A_s(f_{s,\chi,a,n})(I_2,1)}=f_{-s,\chi,a,n}.$$
Therefore: \begin{equation}\label{lcn}C_{\psi_b}(\etac \otimes
\gamma_{\psi_a}^{-1},s)=\frac{\lambda(\psi_b,\chi_{(s)} \otimes
\gamma_{\psi_a}^{-1})(f_{s,\chi,a,n})}{\lambda(\psi_b,\eta(\chi,-s)
\otimes \gamma_{\psi_a}^{-1})(f_{-s,\chi,a,n})} \bigl(
A_s(f_{s,\chi,a,n})(I_2,1)\bigr)^{-1}.\end{equation} We have:
$$A_s(f_{s,\chi,a,n})(I_2,1)=\int_{\R}f_{s,\chi,a,n}\Biggl(\begin{pmatrix}
_{0} & _{1}\\_{-1} & _{0}\end{pmatrix}\begin{pmatrix} _{1} &
_{x}\\_{0} & _{1}\end{pmatrix},1 \Biggr)\,
dx=\int_{\R}f_{s,\chi,a,n}(u_x k_x,1)\, dx,$$ where
$$u_x=\begin{pmatrix}
{\frac 1{\sqrt{1+x^2}}} &{\frac {-x}{\sqrt{1+x^2}}}\\ _{0} &
_{\sqrt{1+x^2}}\end{pmatrix}, \, \, \, \, k_x=\begin{pmatrix}
{\frac {-x}{\sqrt{1+x^2}}} &{\frac 1{\sqrt{1+x^2}}}\\ {\frac
{-1}{\sqrt{1+x^2}}} & {\frac {-x}{\sqrt{1+x^2}}}\end{pmatrix}.
$$
Note that $c(u_x,k_x)=1$, and that since $\frac
1{\sqrt{1+x^2}}>0$, $\phi^{-1}(k_x,1)=t$, where $0<t<\pi$ is the
unique number that satisfies $e^{it}=\frac{i-x}{\sqrt{1+x^2}}$.
Therefore, \begin{eqnarray*}
A_s(f_{s,\chi,a,n})(I_2,1)&=&\int_{\R}
\Bigl(1+x^2\Bigr)^{-\frac{s+1}2}\Bigl(\frac{i-x}{\sqrt{1+x^2}}\Bigr)^n\,
dx \\ &=& e^{\frac{i\pi n}2}\int_{\R}(1+ix)^{\frac
{-s-1+n}2}(1-ix)^{\frac {-s-1-n}2}\, dx.\end{eqnarray*} By Lemma
53 of \cite{Wa} we conclude:
\begin{equation} \label{IO f_n}A_s(f_{s,\chi,a,n})(I_2,1)=e^{\frac{i\pi n}2}\pi2^{1-s}\frac
{\Gamma(s)}{\Gamma(\frac {s+1+n}2) \Gamma(\frac
{s+1-n}2)}.\end{equation} For a similar computation, see
\cite{Gaweb}. In the rest of the proof we assume that $b>0$. We
compute $C_{\psi_b}(\etac \otimes \gamma_{\psi_a}^{-1},s)$:
\begin{eqnarray*}\lambda(s,\chi,\psi_b)(f_{s,\chi,a,n})&=& \int_{\R} f_{s,\chi,a,n}
\Biggl(\begin{pmatrix} _{0} & _{1}\\_{-1} &
_{0}\end{pmatrix}\begin{pmatrix} _{1} & _{x}\\_{0} &
_{1}\end{pmatrix},1 \Biggr) \psi_b^{-1}(x) dx \\ &=& \int_{\R}
\Bigl(1+x^2\Bigr)^{-\frac{s+1}2}\Bigl(\frac{i-x}{\sqrt{1+x^2}}\Bigr)^n
e^{-ibx} \, dx\\ &=& e^{i\pi n}\int_{\R} \ab x+i
\ab^{n-(s+1)}(x+i)^{-n} e^{-ibx} \, dx \\ &=& e^{i\pi
n}b^s\int_{\R} \ab x+bi \ab^{n-(s+1)}(x+bi)^{-n} e^{-ix} \,
dx.\end{eqnarray*} The rest of the computation goes word for word
as the computations in pages 283-285 of \cite{J}:
$$\lambda(s,\chi,\psi_b)(f_{s,\chi,a,n})=2^{\frac{1-s}2}e^{\frac 5 2
i \pi n}\Gamma^{-1}(\frac {s+1+n}{2})W_{\frac n 2, -\frac s
2}(2b),$$ where
$$W_{k,r}(u)=e^{-\half u}u^k \frac 1 {\Gamma(\half+r-k)}
\int_0^\infty v^{r-\half-k}(1+\frac v u)^{r-\half-k} e^{-v} \,
dv$$ is the Whittaker function defined in Chapter 16 of
\cite{Whi}. This function satisfies $$W_{k,r}=W_{k,-r}.$$ It
follows that
\begin{equation} \label{ratio}
\frac{\lambda(s,\chi,\psi_b)(f_{s,\chi,a,n})}{\lambda((-s,\chi,\psi_b)(f_{-s,\chi,a,n})}=2^{-s}\frac
{\Gamma(\frac {-s+1+n}{2})}{\Gamma(\frac {s+1+n}{2})}
.\end{equation} Plugging \eqref{ratio} and \eqref{IO f_n} into
\eqref{lcn} and using Lemma \ref{n} we conclude that
\begin{equation} \label{lc b>0} C_{\psi_b}(\etac \otimes
\gamma_{\psi_a}^{-1},s)=\frac {e^{-\frac {i \pi n} 2}} {2\pi}
\frac {\Gamma(\frac {-s+1+n} 2)\Gamma(\frac {s+1-n}2)}
{\Gamma(s)},\end{equation} where $n \in
\frac{\chi(-1)sign(a)}{2}+2\Z$. We now compute
$C_{\psi_{-b}}(\etac \otimes \gamma_{\psi_a}^{-1},s)$ and finish
the prove of Lemma \ref{the real comp}:
\begin{eqnarray*}\lambda(s,\chi,\psi_{-b})(f_{s,\chi,a,n})&=&e^{i\pi n}b^s\int_{\R} \ab
x+bi \ab^{n-(s+1)}(x+bi)^{-n} e^{ix} \, dx\\ &=& b^s\int_{\R} \ab
x+bi \ab^{n-(s+1)}(x-bi)^{-n} e^{-ix} \, dx.\end{eqnarray*} We
note that $(x-bi)^{-n}=\ab x+bi \ab^{-2n} (x+bi)^n$. Thus,
repeating the previous computation we get:

$$\frac{\lambda(s,\chi,\psi_{-b})(f_{s,\chi,a,n})}{\lambda(-s,\chi,\psi_{-b})(f_{-s,\chi,a,n})}=2^{-s}\frac
{\Gamma(\frac {-s+1-n}{2})}{\Gamma(\frac {s+1-n}{2})} ,$$ which
implies \begin{equation} \label{lc b<0}C_{\psi_{-b}}(\etac \otimes
\gamma_{\psi_a}^{-1},s)=\frac {e^{-\frac {i \pi n} 2}} {2\pi}
\frac {\Gamma(\frac {-s+1-n} 2)\Gamma(\frac {s+1+n}2)}
{\Gamma(s)},\end{equation} where $n \in
\frac{\chi(-1)sign(a)}{2}+2\Z$.

{\bf Remark}: When one computes
$$C_{\psi_a}^{SL_2(\R)}\bigl(B_{SL_2}(\R),s,\chi,(\begin{smallmatrix}
& {0}& {-1} \\
& {1} & {0}\end{smallmatrix})\bigr)$$ one obtains the same
expressions on the right-hand sides of \eqref{lc b>0} and
\eqref{lc b<0}, only that one should assign $n \in 2\Z$ if
$\chi=\chi_{0,0}$ and $n \in 2\Z+1$ if $\chi=\chi_{0,1}$.

\subsection{Complex case} \label{easy complex sec}
Since $\overline{SL_2(\C)}=SL_2(\C) \times \{\pm 1 \}$ and since
$\gamma_\psi(\C^*)=1$ it follows that
\begin{equation} \label{msl is sl}C_\psi^{SL_2(\C)}\bigl(B_{SL_2}(\C),s,\chi,(\begin{smallmatrix}
& {0}& {1} \\
& {-1} &
{0}\end{smallmatrix})\bigr)=C_\psi^{\overline{SL_2(\C)}}\bigl(\overline{B_{SL_2}(\C)},s,\chi,(\begin{smallmatrix}
& {0}& {1} \\
& {-1} & {0}\end{smallmatrix})\bigr).\end{equation} Theorem 3.13
of \cite{Sha85} states that
\begin{equation} \label{shahidi complex}C_\psi^{SL_2(\C)}\bigl(B_{SL_2}(\C),s,\chi,(\begin{smallmatrix}
& {0}& {1} \\
& {-1} & {0}\end{smallmatrix})\bigr) =c'(s) \frac {L(
\chi^{-1},1-s)} {L( \chi,s)}, \end{equation} where $c'(s)$ is an
exponential factor. Recall that any character of $\C^*$ has the
form
$$\chi(re^{i\theta})=\chi_{n,s_0}(re^{i\theta})=r^{s_0}e^{in\theta},$$
for some $s_0\in \C$, $n \in \Z$. We may assume that $s_0=0$ or
equivalently that $\chi=\chi_{n,0}$. The corresponding local
L-function is defined by
$$L_\C(\chi_{n,0},s)=(2\pi)^{-(s+\frac {\ab n \ab}2)} \Gamma(s+\frac {\ab n
\ab}2).$$
\begin{lem} \label{easy complex lem} There exists an exponential factor $c(s)$
such that

\begin{equation} \label{c to show}C_\psi^{SL_2(\C)}\bigl(B_{SL_2}(\C),s,\chi,(\begin{smallmatrix}
& {0}& {1} \\
& {-1} &
{0}\end{smallmatrix})\bigr)=c(s)\frac{\gamma(\chi^2,2s,\psi)}{\gamma(\chi,s+\half,\psi)}.
\end{equation}

\end{lem}
\begin{proof}
Due to \eqref{msl is sl} and \eqref{shahidi complex} we only have
to show that
$$\frac{\Gamma(1+\frac {\ab n\ab} 2 -s)}{\Gamma (\frac{\ab n\ab} 2+s)}=2\frac{\Gamma(1+n-2s)\Gamma(\half+\frac{\ab n\ab} 2+s)}
{\Gamma(n+2s) \Gamma(\half+\frac{\ab n\ab} 2-s)}$$ This fact follows from the classical
duplication formula
$$\Gamma(z)\Gamma(z+\half)=2^{1-2z} \sqrt{\pi} \Gamma(2z).$$
\end{proof}

\subsection{A remark about the chosen parameterization} \label{par
remark} From Theorem \ref{the formula}, Lemma \ref{a to 1} and
Lemma \ref{the real comp} it follows that if $\chi^2$ is
unramified $$C_{\psi}(\etac \otimes \gamma_\psi^{-1},s)$$ and $$
C_{\psi_a}(\etac \otimes \gamma_\psi^{-1},s)$$ might not have the
same of zeros and poles. This phenomenon has no analog in
connected reductive algebraic groups over a local field. Many
authors consider $\psi^{-1}$-Whittaker functionals rather than
$\psi$-Whittaker functionals, see \cite{BFH} for example. It
follows from Theorem \ref{the formula}, Lemma \ref{a to 1} and
Lemma \ref{the real comp} that $C_{\psi^{-1}}(\etac \otimes
\gamma_\psi^{-1},s)$ will contain twists of $\chi$ by the
quadratic character $x \mapsto (x,-1)_{\F}$ whose properties
depend on $\F$. The situation is completely clear in the real
case. The case of p-adic fields of odd residual characteristic is
also easy since a character of $\F^*$, $x \mapsto (x,-1)_{\F}$ is
trivial if $-1 \in {\F^*}^2,$ otherwise it is non-trivial but
unramified. Thus, Theorem \ref{the formula} and Lemma \ref{a to 1}
give a formula for $C_{\psi^{-1}}(\etac \otimes
\gamma_\psi^{-1},s)$.

However, in the case of p-adic fields of even residual
characteristic, the properties of $x \mapsto (x,-1)_{\F}$ depend
so heavily on the field, that it is not possible to give a general
formula for $C_{\psi^{-1}}(\etac \otimes \gamma_\psi^{-1},s)$. If
$\chi^2$ is ramified it follows from Theorem \ref{the formula} and
from Lemma \ref{a to 1} that $C_{\psi^{-1}}(\etac \otimes
\gamma_\psi^{-1},s)$ is a monomial in $q^{-s}$. But, if $\chi^2$
is unramified, things are much more complicated. For example,
assume that $\psi$ is normalized and that $\chi$ is the trivial
character of $\F^*$. By the same techniques we used in Section
\ref{padic comp} the following can be proved: If $e$ is odd then
$$C_{\psi^{-1}}(\chi \otimes \gamma_\psi^{-1},s)=c_\psi(-1)q^{(1-e)s}\frac{L(\chi^{-1},-2s+1)}
{L(\chi,2s)}.$$ If $e$ is even, we have:
\begin{equation} \label{e even} C_{\psi^{-1}}(\chi \otimes \gamma_\psi^{-1},s)= \frac
{q-1} q\biggl((1-q^{-2s})^{-1}+\sum_{1 \leq k \leq \frac {e}
{2}}q^{2ks} \biggr)+P_\F(s), \end{equation} where $P_\F(s)$  is
defined to be
$$\sum_{\frac {e} {2} < k \leq e} q^{2ks}\biggl(
\bigl(c_1(\F,2k)(-q^{-1})+c_2(\F,2k)(1-q^{-1}) \bigr)q^k
\biggr)+$$ $$ q^{s(2e+1)}\biggl(
\bigl(c_1(\F,2e+1)(-q^{-1})+c_2(\F,2e+1)(1-q^{-1})\bigl)q^{e+\frac
1 2} \biggr),$$ where \begin{equation}  c_1(\F,n)=\mu \bigl(
\{x\in \Of^* \mid \ab x^2+1 \ab=q^{1-n} \} \bigr),
\end{equation} \begin{equation}  c_2(\F,n)=\mu \bigl(
\{x\in \Of^* \mid \ab x^2+1 \ab \leq q^{-n} \} \bigr).
\end{equation}
For the 6 ramified quadratic extensions of $\Q_2$, \eqref{e even}
gives the following: If $\F$ is either $\Q_2(\sqrt{2}), \,
\Q_2(\sqrt{-2}), \, \Q_2(\sqrt{6}), \, \Q_2(\sqrt{-6})$ we have:
$$C_{\psi^{-1}}(\chi \otimes \gamma_\psi^{-1},s)=
c_\psi(-1)q^{-s}\frac{L(\chi^{-1},-2s+1)} {L(\chi,2s)},$$ if $\F=
\Q_2(\sqrt{-1})$ we have
$$C_{\psi^{-1}}(\chi \otimes \gamma_\psi^{-1},s)=C_{\psi}(\chi \otimes
\gamma_\psi^{-1},s)=q^{-4s} \frac{L_{\F}(\chi,s+\half)}
{L_{\F}(\chi,-s+\half)} \frac{L_{\F}(\chi^{-1},-2s+1)}
{L_{\F}(\chi,2s)}$$ and if $\F=(\sqrt{3})$ we have
$$C_{\psi^{-1}}(\chi \otimes \gamma_\psi^{-1},s)=c_\psi(-1)q^{-4s} \frac{L_{\F}(\chi^{-1},-s+\half)}
{L_{\F}(\chi,s+\half)} \frac{L_{\F}(\chi^{-1},-2s+1)}
{L_{\F}(\chi,2s)}.$$ These observations, i.e. the appearance of
twists and the incomplete picture in p-adic fields of even
residual characteristic, justifies our parameterization, that is ,
our choice to calculate $C_{\psi}(\etac \otimes
\gamma_\psi^{-1},s)$ rather then $C_{\psi^{-1}}(\etac \otimes
\gamma_\psi^{-1},s)$. We do note that analytic properties of
$$C _{\psi_a}(\etac \otimes \gamma_\psi^{-1},s)C _{\psi_a}(\chi^{-1}\otimes \gamma_\psi^{-1},-s)$$ do not depend on $a$. It
is the last expression which is important for irreducibility; see
Chapter \ref{Irreducibility theorems}.
\subsection{A remark about the similarity between $C_{\psi}(\etac \otimes \gamma_\psi^{-1},s)$ and the Tate $\gamma$-factor}
\label{tate remark} Let $\F$ be a p-adic field and let $\chi$ and
$\psi$ be as before. For simplicity we assume that $\psi$ is
normalized and that $\chi(\pi)=1$ (the latter can be achieved by
shifting the complex parameter $s$). Let $S(\F)$ be the space of
Schwartz-Bruhat functions on $\F$. For $\phi \in S(\F)$ define its
Fourier transform, $\widehat{\phi} \in S(\F)$, by
$$\widehat{\phi}(x)=\int_{\F} \phi(y)\psi(xy) \, dy.$$ It is well
known that for $Re (s)>0$,

$$\int_{\F^*} \phi(x) \chi_{(s)}(x) \, d^*x $$ converges absolutely
for any $\phi \in S(\F)$ and that it has a meromorphic
continuation to the whole complex plane. This continuation is the
Mellin transform, $\zeta(\phi,\chi,s)$. It is a rational function
in $q^{-s}$. Let $\gamma(\chi,\psi,s)$ be the Tate $\gamma$-factor,
see \cite{T}, defined by the functional equation
\begin{equation} \label{tate}
\gamma(\chi,\psi^{-1},s)\zeta(\phi,\chi,s)=\zeta(\widehat{\phi},\chi^{-1},1-s).\end{equation}
where $\phi$ is any function in $S(\F)$. It is clearly a rational
function in $q^{-s}$. It as well known that for $Re(s)>0$
\begin{equation} \label{tate integral}\gamma(\chi^{-1},\psi^{-1},1-s)=\int_{0 < \ab u \ab \leq q^m}  \chi_{(s)}(u)
\psi(u) \,  d^*u,\end{equation} where $m \geq m(\chi)$. We have
seen that for $k \geq {max \bigl (m(\chi),2e+1 \bigr)}$
\begin{equation} \label{my tate integral}C_{\psi}(\etac \otimes
\gamma_\psi^{-1},s)^{-1}= \int_{0 < \ab u \ab \leq q^k}
\gamma_\psi^{-1}(u) \chi_{(s)}(u) \psi(u) \, d^*u \end{equation}
whenever the integral on the right-hand side converges, that is,
for $Re(s)>0$.

The similarity between the integrals in the right-hand sides of
\eqref{tate integral} and \eqref{my tate integral} and the fact
that
\begin{equation} \label{last correction}{C_{\psi}^{SL_2(\F)}\bigl(B_{SL_2}(\F),s,\chi,(\begin{smallmatrix}
& {0}& {-1} \\
& {1} &
{0}\end{smallmatrix})\bigr)}^{-1}=\gamma(\chi^{-1},\psi^{-1},1-s)
\end{equation} raise the question whether one can define a meromorphic
function $\widetilde{\gamma}(\chi,\psi,s)$ by a similar method to
the one used by Tate, replacing to role of the Fourier transform
with $\phi \mapsto \widetilde{\phi}$ defined on $S(\F)$ by
$$\widetilde{\phi}(x)=\int_{\F}
\phi(y)\psi(xy)\gamma_{\psi}^{-1}(xy) \, dy$$ (we define
$\gamma_{\psi}^{-1}(0)=0$). We shall show that such a definition
is possible and that there exits a metaplectic analog to
\eqref{last correction}, namely that
\begin{equation} \label{meta lc is gamtil}C_\psi^{\overline{SL_2(\F)}}\bigl(\overline{B_{SL_2(\F)}},s,\chi,(\begin{smallmatrix}
& {0}& {-1} \\
& {1} &
{0}\end{smallmatrix})\bigr)^{-1}=\widetilde{\gamma}(\chi^{-1},\psi^{-1},1-s).\end{equation}
\begin{lem} \label{new mellin} For any $\phi \in S(\F)$
$$\int_{\F^*} \widetilde{\phi}(x) \chi_{(s)}(x) \, d^*x $$converges absolutely
for $0<Re(s)<1$. This integral extends to a rational function in
$q^{-s}$.
\end{lem}
\begin{proof}
 For $A \subseteq \F$, denote by $ 1 _A$ the characteristic function of $A$.
Regarding the convergence of $$\int_{F^*} \widetilde{\phi}(x)
\chi_{(s)}(x)d^* x, $$ we may assume that $\phi=1_{\Pf^{n}}$ or
that $\phi=1_{a+\Pf^{n}}$, where $\ab a \ab  \geq q^{2e+1-n}$
since these functions span $S(\F)$.

Fix $n \in \Z$, $a \in \F^*$ such that $\ab a \ab \geq
q^{2e+1-n}$.

\begin{eqnarray*} \widetilde{1_{a+
\Pf^{n}}}(y) &=& \int_{a+\Pf^{n}}\psi(xy)\gamma_{\psi}^{-1}(xy) \,
dx = \psi(ay)\int_{\Pf^{n}}\psi(xy)\gamma_{\psi}^{-1}\bigl((a+x)y
\bigr) \, dx\\
&=& \ab y \ab^{-1}\psi(ay)\int_{\ab x \ab \leq \ab y \ab
q^{-n}}\psi(x)\gamma_{\psi}^{-1}\bigl(ay(1+xy^{-1}a^{-1})\bigr) \,
dx. \end{eqnarray*} Note that for $x,y$ and $a$ in the last
integral we have
$$\ab xy^{-1}a^{-1} \ab \leq  \ab q^{-n}a^{-1} \ab \leq q^{-2e-1}.$$
Therefore,
$$\widetilde{1_{a+\Pf^{n}}}(y)=\ab y \ab^{-1}\psi(ay)\gamma_\psi^{-1}(ay)\int_{\ab x \ab \leq \ab y \ab q^{-n}}\psi(x) \, dx.$$
This implies that for $y \in \F^*$
\begin{equation}  \label{widetilde a+Pfn} \widetilde{1_{a+\Pf^{n}}}(y)=\begin{cases} 0 & \ab y \ab>q^n  \\ \psi(ay)\gamma_\psi^{-1}(ay)q^{-n} & \ab y \ab \leq q^n
\end{cases}.\end{equation}
Thus,
$$\int_{\F^*}\ab \widetilde{1_{a+\Pf^{n}}}(x) \chi_{(s)}(x) \ab d^*x =q^{-n}
\int_{0<\ab x \ab \leq q^n}\ab x \ab^{Re(s)} \,
d^*x=q^{-n}\sum_{k=-n}^{\infty}q^{-kRe(s)}.$$ This shows that
$$\int_{\F^*} \widetilde{1_{a+\Pf^{n}}}(x) \chi_{(s)}(x)  \,
d^*x$$ converges absolutely for $Re(s)>0$. Furthermore, for
$Re(s)>0$,
\begin{eqnarray*}\int_{\F^*} \widetilde{1_{a+\Pf^{n}}}(x) \chi_{(s)}(x)  d^*x
&=& q^{-n}\int_{0<\ab x \ab \leq q^n} \psi(ax)\gamma_\psi^{-1}(ax)
\chi_{(s)}(x) \, d^*x \\
&=& \chi_{(s)}(a^{-1})q^{-n}\int_{0<\ab x \ab  \leq \ab a \ab q^n}
\psi(x)\gamma_\psi^{-1}(x) \chi_{(s)}(x) \, d^*x.\end{eqnarray*}
It was shown in Section \ref{int comp} that the last integral is a
rational function in $q^{-s}$.

We now move to $\phi=1_{\Pf^{n}}$. By an argument that we used
already, for $y\in \F^*$ we have

$$\widetilde{1_{\Pf^{n}}}(y)=\ab y \ab^{-1}=\int_{\ab x \ab \leq  \ab y \ab q^{-n}} \psi(x)\gamma_\psi^{-1}(x) \, dx.$$
Put $\ab y \ab= q^{m}.$ Assume first that $m \leq n$:
$$\widetilde{1_{\Pf^{n}}}(y)=q^{-m}\int_{\ab x \ab \leq  q^{m-n}} \gamma_\psi^{-1}(x) \, dx=
q^{-m}\sum_{k=n-m}^{\infty} q^{-k}\int_{\Of^*}
\gamma_\psi^{-1}(u\pi^k) \, du.$$ It was shown in lemma \ref{A for
chi^2} that
$$\int_{\Of^*} \gamma_\psi^{-1}(u\pi^k) \, du=
\begin{cases} 0 & k \in \N_{odd}  \\ c_\psi^{-1}(-1)(1-q^{-1})&  k \in \N_{even}
\end{cases}.$$
This implies
\begin{equation} \label{widetilde Pfn small} \widetilde{1_{\Pf^{n}}}(y)=q^{-m}c_\psi^{-1}(-1)(1-q^{-1})\sum_{k \in \N_{even}, \, k \geq n-m} q^{-k}=
\frac 1 {c_\psi(-1)(1+q^{-1})q^n} \cdot \begin{cases} 1 & n-m \in
\N_{even}  \\ q^{-1}&  n-m \in \N_{odd}
\end{cases}.\end{equation}
If $m>n$ then $$\widetilde{1_{\Pf^{n}}}(y)=q^{-m}\int_{\ab x \ab
\leq q^{m-n}} \psi(x) \gamma_\psi^{-1}(x) \, dx.$$ Put
$$c_{m}=\int_{\ab x \ab \leq  q^{m-n}} \psi(x) \gamma_\psi^{-1}(x) \,
dx.$$ with this notation, if $\ab y \ab= q^{m}>q^n$ then
\begin{equation} \label{widetilde Pfn big}\widetilde{1_{\Pf^{n}}}(y)=c_{m}\ab y \ab^{-1}.\end{equation}
It was shown in Lemma \ref{B poly} that $c_{m}$ stabilizes, more
accurately, that there exists $c \in \C$ such that $c_m=c$ for all
$m>n+2e+1$.

Since $ \widetilde{1_{\Pf^{n}}}(x)$ is bounded, it now follows
that there exists two positive constants $c_1$ and $c_2$ such that
\begin{eqnarray*} \lefteqn{\int_{\F^*}\ab
\widetilde{1_{\Pf^{n}}}(x) \chi_{(s)}(x) \ab d^*x} \\  &\leq& c_1
\int_{0<\ab x \ab \leq q^n}\ab x \ab^{Re(s)} \, d^*x+c_2\int_{\ab
x \ab > q^n}\ab x \ab^{Re(s)-1} \, d^*x\\ &=&
c_1\sum_{k=-n}^{\infty}q^{-nRe(s)}+c_2\sum_{k=n}^{\infty}q^{n\bigl(Re(s)-1
\bigr)}.\end{eqnarray*} This shows that  $$\int_{\F^*}
\widetilde{1_{\Pf^{n}}}(x) \chi_{(s)}(x)  \, d^*x$$ converges
absolutely for $1>Re(s)>0$. From \eqref{widetilde Pfn small} and
\eqref{widetilde Pfn big} it follows that there exist complex
constants such that for $1>Re(s)>0$ we have
\begin{eqnarray*} \lefteqn{\int_{\F^*} \widetilde{1_{\Pf^{n}}}(x)
\chi_{(s)}(x) \, d^*x}\\  &=& k_1  \sum_{k \leq m, \, k \in
Z_{even}}\int_{\ab x \ab=q^k} \chi_{(s)}(x) \, d^*x + k_2\sum_{k
\leq m, \, k \in Z_{odd}}\int_{\ab x
\ab=q^k} \chi_{(s)}(x) \, d^*x\\
&+&  \sum_{k=n}^{n+2e+1}c_k\int_{\ab x \ab=q^k} \chi_{(s-1)}(x)
 \, d^*x + c \sum_{k>n+2e+1}\int_{\ab x \ab=q^k}  \chi_{(s-1)}(x)
 \, d^*x.\end{eqnarray*} It is clear that the right hand side of the
last equation is rational in $q^{-s}$.
\end{proof}
From the lemma just proven it follows that it makes sense to
consider $\zeta(\widetilde{\phi},\chi,s)$, although
$\widetilde{\phi}$ does not necessarily lie in $S(\F)$.
\begin{thm} \label{tate analog}
There exists a meromorphic function
$\widetilde{\gamma}(\chi,\psi,s)$ such that
\begin{equation} \label{my tate}
\widetilde{\gamma}(\chi,\psi^{-1},s)\zeta(\phi,\chi,s)=\zeta(\widetilde{\phi},\chi^{-1},1-s)\end{equation}
For every $\phi \in S(\F)$. In fact, it is a rational function in
$q^s$. Furthermore, \eqref{meta lc is gamtil} holds.
\end{thm}
\begin{proof}
It is sufficient to prove \eqref{my tate} for $0<Re(s)<1$. By
lemma \ref{new mellin}, in this domain both $\zeta(\phi,\chi,s)$
and $\zeta(\widetilde{\phi},\chi^{-1},1-s)$ are given by
absolutely convergent integrals. The proof goes word for word as
the well known proof of \eqref{tate}; see Theorem 7.2 of \cite{RV}
for example, replacing the role of the Fourier transform with
$\phi \mapsto \widetilde{\phi}$. We now prove \eqref{meta lc is
gamtil}. Due to Lemmas \ref{reduction} and \ref{B poly} it is
sufficient to prove that
\begin{equation} \label{my tate integral
again}\widetilde{\gamma}(\chi^{-1},\psi^{-1},1-s)= \int_{0 < \ab x \ab
\leq q^{max(m(\chi),2e+1)}} \gamma_\psi^{-1}(u) \chi_{(s)}(u)
\psi(u) \, d^*u .\end{equation} To prove \eqref{my tate integral
again} we use $\phi=1_{1+\Pf^k}$, where $k \geq max
\bigl(m(\chi),2e+1 \bigr)$. An easy computation shows that
$$\zeta(1_{1+\Pf^k},\chi^{-1},1-s)=q^{-k}.$$ Thus, from \eqref{my
tate} it follows that
$$\widetilde{\gamma}(\chi^{-1},\psi^{-1},1-s)=q^k \zeta(\widetilde{1_{1+\Pf^k}},\chi,s).$$
\eqref{my tate integral again} follows now from \eqref{widetilde
a+Pfn}.
\end{proof}
Many properties of the Tate $\gamma$-factor follow from the
properties of the Fourier transform. What properties of
$\widetilde{\gamma}$ follow from the properties of the transform
$\phi \mapsto \widetilde{\phi}$ is left for a future research.

\subsection{Computation of $C_{\psi}^{\mspm}\bigl(\overline{P_{m;0}(\F)},s,\tau,
\omega_m'^{-1} \bigr)$ for principal series representations}
\label{prin comp with soo}In this subsection we assume that $\F$
is either $\R$, $\C$ or a p-adic field. Let
$\alpha_1,\alpha_2,\ldots,\alpha_m$ be $m$ characters of $\F^*$
and let $\mu$ be the character of $\tglm$ defined by
$$Diag(t_1,t_2,\ldots,t_m) \mapsto \prod_{i=1}^{m}\alpha_i(t_i).$$
We also regard  $\mu$ as a character of $\bglm$. Define $\tau$ to
be the corresponding principal series representation:
$$\tau=I(\mu)=Ind^{\glm}_{B_{\glm}} \mu.$$

\begin{lem}\label{lc zigel prin}  There exists $d \in \{ \pm 1 \}$ such
that
\begin{eqnarray} \label{for soudty said}&& \lefteqn{{C_{\psi}^{\mspm}(\overline{P_{m;0}(\F)},s,\tau,
\omega_m'^{-1} )}} \\ \nonumber &&= d\prod_{i=1}^m
C_\psi^{\overline{SL_2(\F)}}\Bigl(\overline{B_{SL_2}(\F)},s,\alpha_i,(\begin{smallmatrix}
& {0}& {1} \\
& {-1} & {0}\end{smallmatrix})\Bigr) \prod_{i=1}^{m-1}
C_{\psi}^{GL_{m+1-i}(\F)}\bigl(P_{1,m-i}^0(\F),(s,-s),\alpha_i
\otimes \widehat{\tau_i}, \varpi_{1,m-i}^{-1}\bigr) \end{eqnarray}
 where for $1 \leq i
\leq m-2$, $\tau_i=Ind^{GL_{m-i}(\F)}_{B_{GL_{m-i}}(\F)}
\otimes_{j=i+1}^m \alpha_j$ and $\tau_{m-1}=\alpha_m$.
Furthermore, $d=1$ if $\F$ is a p-adic field of odd residual
characteristic and $\tau$ is unramified.
\end{lem}

\begin{proof} We prove this lemma by induction. For $m=1$ there is nothing to prove.
Suppose now that the theorem is true for $m-1$. With our
enumeration this means that there exists $d' \in \{ \pm 1 \}$ such
that
\begin{eqnarray*} && \lefteqn{{C_{\psi}^{\widehat{Sp_{2(m-1)}(\F)}}(\overline{P_{m-1;0}(\F)},s,\tau_1, \omega_{m-1}'^{-1}
)}} \\ &&=d'\prod_{i=2}^m
C_\psi^{\overline{SL_2(\F)}}\bigl(\overline{B_{SL_2(\F)}},s,\alpha_i,(\begin{smallmatrix}
& {0}& {1} \\
& {-1} & {0}\end{smallmatrix})\bigr) \prod_{i=2}^{m-1}
C_{\psi}^{GL_{m-i}(\F)}\bigl(P_{1,m-i}(\F)^0,(s,-s),\alpha_i
\otimes \widehat{\tau_i},
\varpi_{1,m-1-i}^{-1}\bigr)\end{eqnarray*} and that $d'=1$ if $\F$
is a p-adic field of odd residual characteristic and $\tau$ is
unramified. The proof is done now once we observe that since $\tau
\simeq Ind^{\glm}_{P^0_{1,m-1}} \alpha_1 \otimes \tau_1$ it
follows from Lemmas \ref{ind by parts}, \ref{lem heart} and
\ref{id with c} that \begin{eqnarray*} &&
{C_{\psi}^{\mspm}(\overline{P_{m;0}(\F)},s,\tau, \omega_m'^{-1}
)}=
d''C_\psi^{\overline{SL_2(\F)}}\bigl(\overline{B_{SL_2}(\F)},s,\alpha_1,(\begin{smallmatrix}
& {0}& {1} \\
& {-1} & {0}\end{smallmatrix})\bigr) \\ &&
C_{\psi}^{\overline{Sp_{2(m-1)}(\F)}}\bigl(\overline{P_{m-1;0}(\F)},s,
\tau_1  ,
\omega_{m-1}'^{-1}\bigr)C_{\psi}^{GL_{m}(\F)}\bigl(P_{1,m-1}^0(\F),(s,-s),\alpha_1
\otimes \widehat{\tau_1},
\varpi_{1,m-1}^{-1}\bigr).\end{eqnarray*} for some  $d'' \in \{
\pm 1 \}$. If $\F$ is a p-adic field of odd residual
characteristic and $\tau$ is unramified then $d''=1$.
\end{proof}

\begin{thm} \label{soudry said it}
There exists an exponential function $c(s)$ such that

\begin{equation} \label{soudry said}
C_{\psi}^{\mspm}\bigl(\overline{P_{m;0}(\F)},s,\tau,
\omega_m'^{-1} \bigr)=c(s)\frac
{\gamma(\tau,sym^2,2s,\psi)}{\gamma(\tau,s+\half,\psi)}\end{equation}
If $\F$ is a p-adic field of odd residual characteristic, $\psi$
is normalized and $\tau$ is unramified then $c(s)=1$.
\end{thm}
In Section \ref{true for cusp} we shall show that \eqref{soudry
said} holds for every irreducible admissible generic representation $\tau$ of
$\glm$; see Theorem \ref{soudry said always}.
\begin{proof}
We keep the notations of Lemma \ref{lc zigel prin}. From
\eqref{gama def gl} and from the known properties of
$\gamma(\tau,s,\psi)$ (see \cite{Sha84} or \cite{JPS}) it
follows that for every $1 \leq i \leq m-1$ there exists $d_i \in
\{ \pm 1 \}$ such that
\begin{equation} \label{product 1}C_{\psi}^{GL_{m+1-i}(\F)}\bigl(P_{1,m-i}^0(\F),(s,-s),\alpha_i
\otimes \widehat{\tau_i}, \varpi_{1,m-i}^{-1}\bigr)=d_i
\prod_{j=i+1}^m \gamma(\alpha_i \alpha_j,2s,\psi).\end{equation}
and that $d_i=1$ provided that $\F$ is a p-adic field and $\tau$
is unramified. From \eqref{res g} it follows that \begin{equation}
\label{product 2}\prod_{i=1}^m
C_\psi^{\overline{SL_2(\F)}}\Bigl(\overline{B_{SL_2}(\F)},s,\alpha_i,(\begin{smallmatrix}
& {0}& {1} \\
& {-1} & {0}\end{smallmatrix})\Bigr)=
c'(s)\prod_{i=1}^m\frac{\gamma(\alpha_i^2,2s,\psi)}{\gamma(\alpha_i,s+\half,\psi)},\end{equation}
where $c_\F'(s)$ is an exponential function that equals 1 if $\F$
is a p-adic field of odd residual characteristic, $\psi$ is
normalized and $\tau$ is unramified. Plugging \eqref{product 1}
and \eqref{product 2} into \eqref{soudry said} we get
$$C_{\psi}^{\mspm}\bigl(\overline{P_{m;0}(\F)},s,\tau,
\omega_m'^{-1}
\bigr)=c(s)\prod_{i=1}^m\Bigl(\frac{\gamma(\alpha_i^2,2s,\psi)}{\gamma(\alpha_i,s+\half,\psi)}
\prod_{j=i+1}^m \gamma(\alpha_i \alpha_j,2s,\psi)\Bigr).$$ By
definition, \eqref{soudry said} now follows.
\end{proof}
\newpage

\section{An analysis of Whittaker coefficients of an Eisenstein series} \label{global sec}  In this
chapter we prove a global functional equation satisfied by
$\gamma(\msig \times \tau,s,\psi)$; see Theorem \ref {crude} of
Section \ref{crude equ}. It is a metaplectic analog to Theorem 4.1
\cite{Sha 1} (see also Part 4 of Theorem 3.15 of \cite{Sha 3}).
The argument presented here depends largely on the theory of
Eisenstein series developed by Langlands, \cite{L}, for reductive
groups. Moeglin and Wladspurger extended this theory to coverings
groups; see \cite{MW}. This chapter is organized as follows. In
Section \ref{anunramcomp} we introduce some unramified
computations which follow from \cite{BFH} and \cite{CS}. These
computations will be used in Section \ref{crude equ} where we
prove a global functional equation by analyzing the
$\psi$-Whittaker coefficient of a certain Eisenstein series. As a
consequence of a particular case of the global functional equation
we show in Theorem \ref{soudry said always} of Section \ref{true
for cusp} that if $\F$ is a p-adic field and $\tau$ is an
irreducible admissible generic representation of $\glm$ then there
exists an exponential function $c(s)$ such that
$$ C_{\psi}^{\mspm}\bigl(\overline{P_{m;0}(\F)},s,\tau,
\omega_m'^{-1} \bigr)=c(s)\frac
{\gamma(\tau,sym^2,2s,\psi)}{\gamma(\tau,s+\half,\psi)}.$$

\subsection{Unramified computations} \label{anunramcomp}
We keep all the notations we used in Lemma \ref{princcomp} but we
add the following restrictions; we assume that $\F$ is a p-adic
field of odd residual characteristic, that $\psi$ is normalized
and that $\chi$ and $\mu$ are unramified. We define the local
unramified $L$-function of $\msig \times \tau$ with respect to
$\psi$:
\begin{equation} \label{l def} L_\psi(\msig \times \tau,s)=
\prod_{1 \leq i \leq k} \prod_{1 \leq j \leq m} L(\eta_i
\alpha_j,s)L(\eta_i^{-1} \alpha_j,s).
\end{equation} The subscript
$\psi$ is in place due to the dependence on $\gamapsi$ in the
definition of $\msig$.

 Similar to the algebraic case,
$I(\chi_{(s)})$ has a one dimensional $\kappa_{2k}\bigl(\Kk
\bigr)$ invariant subspace. Let $f^0_{\chi_{(s)}}$ be the
normalized spherical vector of $I(\chi_{(s)})$, i.e., the unique
$\kappa_{2k}\bigl(\Kk \bigr)$ invariant vector with the property
$f^0_{\chi_{(s)}}(I_{2k},1)=1$. For $f \in I(\chi_{(s)})$ the
corresponding Whittaker function is defined by $$ W_f(g)= \frac
{1}{C_{\chi_{(s)}}}\int_{\zspk} f \bigl((J_{2k}u,1)g \bigr)
\psi^{-1}(u) \, du, $$ where $$ \label{bfh}
C_{\chi_{(s)}}=\prod_{i=1}^k \bigl (1+ \eta_i(\pi) q^{-(s+\half)}
\big) \prod_{1 \leq i<j \leq k} \bigl( (1-q^{-1}\eta_i(\pi)
\eta_j(\pi)^{-1})(1-q^{-1}\eta_i(\pi) \eta_j(\pi)q^{-2s})
\bigr).$$

With the normalization above, Theorem 1.2 of \cite{BFH} states
that $W_{f^0_{\chi_{(s)}}}=W^0_{\chi_{(s)}}$, where
$W^0_{\chi_{(s)}}$ is the normalized spherical function in
$W(I(\chi_{(s)}),\psi)$. To be exact we note that in \cite{BFH},
the $\psi^{-1}$-Whittaker functional is computed. This difference
manifests itself only in the $\msl$ computation presented in page
387 of \cite{BFH}. Consequently the left product defining
$C_{\chi_{(s)}}$ presented in \cite{BFH} differs slightly from the
one given here. It is the same difference discussed in Lemma
\ref{a to 1} and Section \ref{par remark}.

Let $f_\mu^0$ be the normalized spherical vector of $I(\mu)$.
Define $$ W_f(g)= \frac {1}{D_\mu}\int_{\zglm} f (\omega_mug)
\psi^{-1}(u)\, du,  $$ where $$ D_\mu=\prod_{1 \leq i<j \leq
m}(1-q^{-1}\alpha_i \alpha_j^{-1}).$$

Denote by $W^0_\mu$ the normalized spherical function of
$W(I(\mu),\psi)$. Theorem 5.4 of \cite{CS} states that
$W_{f^0_\mu}=W^0_{\mu}$. Let $\lambda_{\tau,\psi}$ and
$\lambda_{\msig,\psi}$ be Whittaker functionals on $I(\mu)$ and
$I(\chi)$ respectively. Note that
\begin{equation} \label{relation}
\lambda_{\tau,\psi}\bigl(\tau(g)f \bigr)=W_f(g), \quad \quad
\lambda_{\msig,\psi}\bigl(\msig(s)f\bigr)=W_f(s). \end{equation}
Similar to Chapter \ref{cha gam}, we realize
$$I_1=Ind^{\mspn}_{\overline{P_{m;k}(\F)}} (\gamapsi \otimes  \tau_{(s)}) \otimes \msig $$ as a space of complex
functions on $\mspn \times \glm \times \mspk$ which are smooth
from the right in each argument and which satisfy
\begin{eqnarray*} \lefteqn{f
\bigl((j_{m,n}(\widehat{g}),1)i_{k,n}(y)nh,bg_0,(b',\epsilon)y_0\bigr)}
\\ && = \epsilon \gamapsi \bigl (\det(g)\det(b') \bigr) \ab det(g)
\ab^{s+\frac {n+k+1}{2}} \delta_{B_{\glm}}(b)\delta_{\bspn}(b')
\mu(b) \chi(b)f(h,g_0g,y_0y),\end{eqnarray*} For all $g,g_0 \in
\glm, \,y,y_0 \in \mspk, \,n \in (N_{m;k},1), \, h \in \mspn, \, b
\in B_{\glm}, \, (b',\epsilon) \in \overline{B_{\spk}}.$ We
realize
$$I_1'=Ind^{\mspn}_{\overline{P_{m;k}(\F)}} \bigl(\gamapsi \otimes  W_{(s)}(\tau,\psi)\bigr) \otimes W(\msig,\psi)$$ as we
did in Lemma \ref{ind by parts}. An isomorphism $T_1:I_1
\rightarrow I_1'$ is given by
$$(T_1f)(h,g,y)= \frac {1}{C_\chi D_\mu} \int_{n_1 \in Z_{\glm}} \int_{n_2 \in
Z_{\spk}}f \bigl(s,\omega_mn_1 g,(J_{2k}n_2,1)y \bigr)
\psi^{-1}(n_1) \psi^{-1}(n_2) .$$ Let $f^0_{I_1} \in I_1$ be the
unique function such that
$$f^0_{I_1}\bigl((I_{2n},1),I_m,(I_{2k},1) \bigr)=1$$ and such that
for all $o \in \kappa_{2n} \bigl(\Kn \bigr), \, g \in \glm, \, y
\in \mspk$ we have $$f^0_{I_1}(o,g,y)=f^0_\mu(g) \cdot
f^0_\chi(y).$$ Let $f^0_{I_1'} \in I_1'$ be the unique function
such that
$$f^0_{I_1'}\bigl((I_{2n},1),I_m,(I_{2k},1) \bigr)=1$$ and such that
for all $o \in \kappa_{2n}\bigl(\Kn \bigr), \, g \in \glm, \, y
\in \mspk$ we have $$f^0_{I_1'}(o,g,y)=W^0_\mu(g) \cdot
W^0_\chi(y).$$ According \eqref{relation}:
$T_1(f^0_{I_1})=f^0_{I_1'}$. We denote by $\lambda'(s,\tau \otimes
\msig,\psi \bigr)$ the Whittaker functional on $I_1'$ constructed
in the usual way.
\begin{lem}
\begin{equation} \label{proven up} \lambda'(s,\tau \otimes \msig,\psi
\bigr)f^0_{I_1'}=\frac{L(\tau,s+\half)}{L(\tau,sym^2,2s+1)L_\psi(\msig
\otimes \tau,s+1)}.\end{equation}
\end{lem}
\begin{proof}
For $f \in I_1$ we have:
\begin{eqnarray} \label {unmram lambda 1}  &&\lefteqn{\lambda'(s,\tau \otimes \msig,\psi
\bigr)(T_1(f))}  \\ \nonumber &&=\int_{N_{m;k}(\F)}  (T_1f)
\bigl((j_{m,n}(\omega_m') u,1),I_m,(I_{2k},1)  \bigr)
\psi^{-1}(u)du  \\ \nonumber &&= \frac {1}{C_\chi D_\mu}
\int_{N_{\spn}} f \bigl( (J_{2n}u,1),I_m,(I_{2k},1) \bigr)
\psi^{-1}(u)du.\end{eqnarray} In particular $$ \lambda'(s,\tau
\otimes \msig,\psi \bigr)f^0_{I_1'}= \frac {1}{C_\chi D_\mu}
\int_{N_{\spn}} f^0_{I_1} \bigl( (J_{2n}u,1),I_m,(I_{2k},1) \bigr)
\psi^{-1}(u)du.$$ Let $\mu_{(s)} \otimes \chi$ be the character of
$\tspn$ defined by
$$\bigl(j_{m,n}(\widehat{t})i_{k,n}(t') \bigr) \mapsto \ab \det(t)
\ab ^s \mu(t) \chi(t'),$$ where $t \in \tglm, \, t' \in \tspk$. We
realize $I(\mu_{(s)} \otimes \chi)=Ind^{\mspn}_{\overline{\bspn}}
\gamapsi \otimes \mu_{(s)} \otimes \chi$,  in the obvious way. For
the Whittaker functional defined on this representation space,
$$\lambda(s,\chi \otimes \mu)(f)= \int_{N_{\spn}} f(J_{2n}u,1)
\psi^{-1}(u)du, $$ we have
\begin{equation} \label {unmram lambda 2 prop}\lambda(s, \mu \otimes \chi )(f^0_{I(\mu _s \otimes
\chi)})=C_{\mu_{(s)} \otimes \chi}.\end{equation} The isomorphism
$T_2: I_1 \rightarrow I(\mu_{(s)} \otimes \chi)$ defined by
$$(T_2f)(h)=f\bigl(h,I_m,(I_{2k},1) \bigr),$$ whose inverse is given
by $$(T_2^{-1}f)(h,g,y)=\gamma_{\psi}\bigl(\det (g) \bigr) \ab
det(g) \ab^{\frac{n+k+1}{-2}}f
\bigl((j_{m,n}(\widehat{g}),1)i_{k,n}(y)h \bigl),$$ has the
property:
\begin{equation} \label{T property}T_2(f^0_{I_1})=f^0_{I(\mu _s \otimes
\chi)}. \end{equation} Using \eqref{unmram lambda 1},
\eqref{unmram lambda 2 prop} and \eqref{T property} we observe
that

\begin{equation} \label{preproven up}C_{\mu_{(s)} \otimes \chi}=
\lambda(s,\chi \otimes \mu)(f^0_{I(\chi _s \otimes \mu)})=
\int_{N_{\spn}} \bigl(T_2(f^0_{I_1})\bigr)
\bigl(J_{2n}u,1)\psi^{-1}(u)du=$$
$$ \int_{N_{\spn}}
f^0_{I_1} \bigl( (J_{2n}u,1),I_m,(I_{2k},1) \bigr)
\psi^{-1}(u)du=C_\chi D_\mu \lambda'(s,\tau \otimes \msig,\psi
\bigr)f^0_{I_1'}.
\end{equation} Since
$$ \frac {C_{\mu_{(s)} \otimes \chi}}{C_\chi
D_\mu}=\frac{L(\tau,s+\half)}{L(\tau,sym^2,2s+1)L_\psi(\msig
\otimes \tau,s+1)},$$ the lemma is proved.
\end{proof}
{\bf Remark:} In the case $k=0$ \eqref{proven up} reduces to
$$\lambda'(s,\tau \otimes \msig,\psi
\bigr)f^0_{I_1'}=\frac{L(\tau,s+\half)}{L(\tau,sym^2,2s+1)}.$$
This case appears in the introduction of \cite{BFH}.

Let $A_{j_{m,n}(\omega_m'^{-1})}^{\mu_{(s)} \otimes \chi},
A_{j_{m,n}(\omega_m'^{-1})}, \, A_{j_{m,n}(\omega_m'^{-1})}'$ be
the intertwining operators defined on $I(\mu_{(s)} \otimes
\chi),\, I_1, \, I_1'$ respectively.
\begin{lem}\begin{eqnarray}
\label{proven down}&& \lefteqn{\lambda'(-s,\widehat{\tau} \otimes
\msig,\psi \bigr)
\bigl(A_{j_{m,n}(\omega_m'^{-1})}'(f^0_{I_1'})\bigr)} \\
\nonumber && = \frac{L(\widehat{\tau},-s+\half)L_\psi(\msig
\otimes \tau,s)L(\tau,sym^2,2s)}
{L(\widehat{\tau},sym^2,-2s+1)L_\psi(\msig \otimes
\widehat{\tau},-s+1)L_\psi(\msig \otimes
\tau,s+1)L(\tau,sym^2,2s+1)}.\end{eqnarray} \end{lem}

\begin{proof}
Application of Lemma 3.4 of \cite{BFH} to the relevant Weyl
element proves that $$ A_{j_{m,n}(\omega_m'^{-1})}^{\mu_{(s)}
\otimes \chi}(f^0_{\mu_{(s)} \otimes \chi})=K_{\mu_{(s)} \otimes
\chi}f^0_{I(\mu^{-1}_{(-s)} \otimes \chi)}, $$ where
$$K_{\mu_{(s)} \otimes \chi}=\frac{L(\msig \otimes \tau,s)}{L(\msig
\otimes \tau,s+1)}\frac{L(\tau,sym^2,2s)}{L(\tau,sym^2,2s+1)}.$$
We define $$\widetilde{T_1}:Ind^{\mspn}_{\overline{P_{m;k}(\F)}}
(\gamapsi \otimes  \tau^{(-1)}_{(-s)}) \otimes \msig \rightarrow
Ind^{\mspn}_{\overline{P_{m;k}(\F)}} \bigl(\gamapsi \otimes
W_{(-s)}(\tau^{(-1)},\psi)\bigr) \otimes W(\msig,\psi)$$ and
$$\widetilde{T_2}:Ind^{\mspn}_{\overline{P_{m;k}(\F)}}
(\gamapsi \otimes  \tau^{(-1)}_{(-s)}) \otimes \msig \rightarrow
Ind^{\mspn}_{\overline{\bspn}} \gamapsi \otimes \mu^{-1}_{(-s)}
\otimes \chi $$ by analogy with $T_1$ and $T_2$. Note that
$\widetilde{T_1}$ commutes with $A_{j_{m,n}}$ and that
$\widetilde{T_2}$ commutes with $A_{j_{m,n}(\omega_m'^{-1})}$.
Therefore,
$$\widetilde{T_1} \bigl(A_{j_{m,n}(\omega_m'^{-1})}(I_1) \bigr)
\subseteq A_{j_{m,n}(\omega_m'^{-1})}'(I_1')$$ and
$$\widetilde{T_2}: \bigl(A_{j_{m,n}(\omega_m'^{-1})}(I_1)\bigr) \subseteq
A_{j_{m,n}(\omega_m'^{-1})}^{\mu_{(s)} \otimes \chi}
\bigl(I(\mu_{(s)} \otimes \chi) \bigr).$$ We denote by
$\widetilde{f^0_{I_1}}$ and $\widetilde{f^0_{I_1'}}$ the spherical
functions of $A_{j_{m,n}(\omega_m'^{-1})}(I_1)$ and
$A_{j_{m,n}(\omega_m'^{-1})}'(I_1')$ respectively. Since $T_1, \,
T_2, \, ,\widetilde{T_1}, \,\widetilde{T_2}$ map a normalized
spherical function to a normalized spherical function and since a
straightforward computation shows that $$\widetilde{T_1}
\widetilde{T_2}^{-1}A_{j_{m,n}(\omega_m'^{-1})}^{\mu_{(s)} \otimes
\chi}= A_{j_{m,n}(\omega_m'^{-1})}' T_1 T_2^{-1},$$ we have
$$A_{j_{m,n}(\omega_m'^{-1})}'(f^0_{I_1'})=\widetilde{T_1}
\widetilde{T_2}^{-1}A_{j_{m,n}(\omega_m'^{-1})}^{\mu_{(s)} \otimes
\chi} T_2 T_1^{-1} (\widetilde{f^0_{I_1'}})= K_{\mu_{(s)} \otimes
\chi}\widetilde{f^0_{I_1'}}.$$ From this and from \eqref{preproven
up} we conclude that
$$\lambda'(-s,\widehat{\tau} \otimes \msig,\psi \bigr) \bigl(A_{j_{m,n}(\omega_m'^{-1})}'(f^0_{I_1'})\bigr)
=K_{\mu_{(s)} \otimes \chi}\lambda'(-s,\widehat{\tau} \otimes
\msig,\psi \bigr)(f^0_{I_1'})=K_{\mu_{(s)} \otimes \chi} \frac
{C_{\widehat {\mu}_{-s} \otimes \chi}        }{C_\chi
D_{\breve{\mu}}}.$$ This finishes the proof of this lemma.
\end{proof}

{\bf Remark}: Combining \eqref{proven up} and \eqref{proven down}
we get
\begin{eqnarray} \label{uram lc} && C_{\psi}^{\mspn}(\overline{P_{m;k}(\F)},s,\tau\otimes \msig,
j_{m,n}(\omega_m'^{-1}))=\frac {\lambda '(s,\tau \otimes
\msig,\psi \bigr)(f^0_{I_1'})}{\lambda' (-s,\widehat{\tau} \otimes
\msig,\psi \bigr)(A'_{j_{m,n}(\omega_m'^{-1})} (f^0_{I_1'}))}\\
\nonumber && =
 \frac{L(\tau,s+\half)L(\widehat{\tau},sym^2,-2s+1)L_\psi(\msig \otimes
\widehat{\tau},1-s)} {L(\widehat{\tau},-s+\half)L(\tau,sym^2,2s )
L_\psi(\msig \otimes \tau,s)}.\end{eqnarray} In particular,
\begin{equation} \label{uram lc par}
C_{\psi}^{\mspm}(\overline{P_{m;0}(\F)},s,\tau, \omega_m'^{-1}
)=\frac{L(\tau,s+\half)L(\widehat{\tau},sym^2,-2s+1)}
{L(\widehat{\tau},-s+\half) L(\tau,sym^2,2s)}.
\end{equation}
Dividing \eqref{uram lc} by \eqref{uram lc par} we get
\begin{equation} \label{unramgamma} \gamma(\overline{\sigma} \times
\tau,s,\psi)=\frac{L_\psi(\msig \times
\widehat{\tau},1-s)}{L_\psi(\msig \times \tau,s)}. \end{equation}
Thus, the computations presented in this subsection provide an
independent proof of Corollary \ref{princcomp} in the case where
$\F$ is a p-adic field of odd residual characteristic and $\msig$
and $\tau$ are unramified.

\subsection{Crude functional equation.} \label{crude equ}
Throughout this section, $\F$ will denote a number field. For
every place $\nu$ of $\F$, denote by $\F_\nu$ the completion of
$\F$ at $\nu$. Let $\A$ be the adele
ring of $\F$. We fix a non-trivial character $\psi$ of $\F
\backslash \A$. We write $\psi(x)=\prod_\nu \psi_\nu(x_\nu)$,
where for almost all finite $\nu$, $\psi_\nu$ is normalized. As in
the local case, $\psi$ will also denote a character of
$Z_{GL_m}(\A)$, $\overline{Z_{Sp_{2n}}(\A)}$ and of their
subgroups.

Let $\tau$ and $\msig$ be a pair of irreducible automorphic
cuspidal representations of $GL_m(\A)$ and
$\overline{Sp_{2k}(\A)}$ respectively. Let $\tau$ and $\msig$ act
in the spaces $V_\tau$ and $V_\msig$ respectively. We assume that
$\msig$ is genuine and globally $\psi$-generic, i.e., that

\begin{equation} \label{our assu} \int_{n \in \zspk \backslash Z_{Sp_{2k}}(\A)}
\phi_\msig (n,1) \psi^{-1}(n) \, dn \neq 0 \end{equation} for some
$\phi_\msig \in V_\msig$. Fix isomorphisms $T_1:\otimes'_\nu \tnu
\rightarrow \tau$ and $T_2:\otimes'_\nu \msignu \rightarrow
\msig$. Here, for each place $\nu$ of $\F$, $\tnu$ and $\msignu$
are the local components. Outside a finite set of places $S$,
containing the even places and those at infinity, $\tnu$ and
$\msignu$ come together with a chosen spherical vectors
$\alpha_\tnu^0$ and $\beta_\msignu^0$ respectively. We may assume,
and in fact do, that $\psi_\nu$ is normalized for all $\nu \notin
S$.

Let $T=T_1 \otimes T_2$. We identify $(\otimes_\nu'\tnu)\otimes
(\otimes_\nu' \msignu)$ with $\otimes_\nu'(\tnu \otimes \msignu)$
in the obvious way. We also identify the image of $T$ with the
space of cusp forms on $GL_m(\A) \times \overline{Sp_{2k}(\A)}$
generated by the functions $(g,h) \mapsto
\phi_\tau(g)\phi_\msig(h)$, here $g \in GL_m(\A)$, $h \in
\overline{Sp_{2k}(\A)}$, $\phi_\tau \in V_\tau$ and $\phi_\msig
\in V_\msig$. $T$ then is an isomorphism $T:\otimes_\nu'(\tnu
\otimes \msignu) \rightarrow \tau \otimes \msig$. Denote for $\phi
\in V_{ \tau \otimes \msig}$.

$$W_{\phi}(g,h)=\int_{n_1 \in \zglm \backslash Z_{GL_m}(\A)} \int_{n_2 \in \zspk \backslash Z_{Sp_{2k}}(\A)}
\phi \bigl(n_1g,(n_2,1)h \bigr) \psi^{-1}(n_1) \psi^{-1}(n_2) dn_2
\, dn_1.$$ By our assumption \eqref{our assu}, there exists $\phi \in
V_{ \tau \otimes \msig}$ such that $W_{\phi} \neq 0$ is not the zero function. Note that the linear functional
$$\lambda_{\tau \otimes \msig,\psi}(\phi)=
W_{\phi}\bigl(I_m,(I_{2n},1)\bigr)$$ is a non-trivial (global)
$\psi$-Whittaker functional on $V_{\tau \otimes \msig}$, i.e,
$$\lambda_{\tau \otimes \msig,\psi}\bigl(\tau \otimes \msig
(n_1,n_2)\phi \bigr)=\psi(n_1)\psi(n_2)\lambda_{\tau \otimes
\msig,\psi}(\phi),$$ for all $(n_1,n_2) \in GL_m(\A) \times
\overline{Sp_{2k}(\A)}$. The last fact and the local uniqueness of
Whittaker functional imply that
\begin{lem} \label{whi1}There exists a unique, up to scalar, global $\psi$-Whittaker functional
on $\tau \otimes \msig$:
$$\phi \mapsto \lambda_{\tau \otimes \msig,\psi}(\phi)=
\int_{\zglm \backslash Z_{GL_m}(\A)} \int_{\zspk \backslash
Z_{Sp_{2k}}(\A)} \phi \bigl(n_1,(n_2,1) \bigr) \psi^{-1}(n_1)
\psi^{-1}(n_2) dn_2 \, dn_1 .$$ For each $\nu$ let us fix a
non-trivial $\psinu$ Whittaker functional $\lambda_{\tnu \otimes
\msignu, \psi_\nu}$ on $V_{\tnu \otimes \msignu}$ at each place
$\nu$, such that if $\tau_\nu \otimes \msignu$ is unramified then
$$\lambda_{\tnu \otimes \msignu, \psi_\nu}(\alpha_\tnu^0 \otimes
\beta_\msignu^0)=1.$$ Then, by normalizing $\lambda_{\tnu \otimes
\msignu, \psi_\nu}$ at one ramified place, we have
$$\lambda_{\tau \otimes \msig,\psi}
(\phi)=\prod_\nu \lambda_{\tnu \otimes \msignu, \psi_\nu}(v_\tnu
\otimes v_\msignu),$$ where $\phi=T \bigl(\bigotimes_\nu (v_\tnu
\otimes v_\msignu) \bigr)$, i.e., $\phi$ corresponds to a pure
tensor.

\end{lem}
We shall realize each local representation
$$I_\nu({\tnu}_{(s)},\msignu)=Ind^{\mspnnu}_{\overline{P_{m;k}(\fnu)}}(\gamapsinu
\otimes  {\tnu}_{(s)})\otimes \msignu$$ as the space of smooth
from the right functions $$f: \mspnnu \rightarrow V_\tnu \otimes
V_\msignu$$ satisfying
$$f\bigl((j_{m,n}(\widehat{g}),1)i_{k,n}(y)nh \bigl)=
\gamapsinu(g)\ab \det(g) \ab_\nu ^{s+\frac{n+k+1}{2}} \tnu(g)
\otimes\msignu(y) f(h)$$ for all $g \in GL_m(\fnu), \, y \in
\overline{Sp_{2k}(\fnu)}, \,n \in  \bigl(N_{m,k}(\fnu),1 \bigr),
\, h \in \overline{Sp_{2n}(\fnu)}.$ For each place where $\tnu$
and $\msignu$ are unramified we define $f_\nu^{0,s} \in
I_\nu({\tnu}_{(s)},\msignu)$ to be the normalized spherical
function, namely, $f_\nu^{0,0}(I_{2n},1)=\alpha_\tnu^0 \otimes
\beta_\msignu^0$. We shall realize the global representation
$$I({\tau}_{(s)},\msig)=Ind^{\overline{Sp_{2n}(\A)}}_{\overline{P_{m;k}(\A)}}(\gamapsi
\otimes \tau_{(s)})\otimes \msig$$ as a space of functions $$f:
\overline{Sp_{2n}(\A)} \times GL_m(\A) \times
\overline{Sp_{2k}(\A)} \rightarrow V_{\tau \otimes \msig}$$ smooth
from the right in the first variable such that
$$f \bigl((j_{m,n}(\widehat{g}),1)i_{k,n}(y)nh ,g_0,y_0 \bigr)=
\gamapsi(g)\ab \det(g) \ab ^{s+\frac{n+k+1}{2}} \tau(g)
\otimes\msig(y) f(h,g_0g,y_0y),$$ for all $g,g_0 \in GL_m(\A), \,
y,y_0 \in \overline{Sp_{2k}(\A)}, \,n \in  \bigl(N_{m,k}(\A),1
\bigr), \, h \in \overline{Sp_{2n}(\A)},$ and such that for all $h
\in \overline{Sp_{2n}(\A)}$ the map $(g,y) \mapsto f(h,g,y)$ lies
in $V_{\tau \otimes \msig}$.

$I({\tau}_{(s)},\msig)$ is spanned by functions of the form
$f(g)=T\bigl((\otimes_\nu f_\nu(g_\nu)\bigr)$, where $f_\nu \in
I_\nu({\tnu}_{(s)},\msignu)$ and for almost all $\nu$:
$f_\nu=f^0_\nu$ (for a fixed $g$, $f(g)$ is a cuspidal automorphic
form corresponding to a pure tensor).

We note that for $(p,1) \in \bigl(P_{m;k}(\F),1 \bigr)$ we have $f
\bigl((p,1)g\bigr)=f (g)$. This follows from the fact that
$\prod_\nu \gamapsinu(a)=1$ for all $a \in \F^*$. Hence, it makes
sense to consider Eisenstein series: For a holomorphic smooth
section  $f_s \in I({\tau}_{(s)},\msig)$ define
$$E(f_s,g)=\sum_{\gamma \in P_{m;k}(\F) \backslash
\spn}f_s \bigl((\gamma,1) g,I_m,(I_{2k},1) \bigr)$$ It is known
that the series in the right-hand side converges absolutely for
$Re(s)>>0$, see Section II.1.5 of \cite{MW} and that it has a
meromorphic continuation to the whole complex plane, see Section
IV.1.8 of \cite{MW}. We continue to denote this continuation by
$E(f_s,g)$.

We introduce the $\psi$- Whittaker coefficient
$$E_\psi(f_s,g)=\int_{u \in Z_{Sp_{2n}(\F)} \backslash Z_{Sp_{2n}}(\A)}
 E \bigl(f_s,(u,1)g \bigr) \psi^{-1}(u) \, du. $$
Note that no question of convergence arises here since
$Z_{Sp_{2n}}(\F) \backslash Z_{Sp_{2n}}(\A)$ is compact. It is
also clear that $E_\psi(f_s,g)$ is meromorphic in the whole
complex plane.
\begin{lem}
$\prod_{\nu \notin S}L_{\psi_\nu}(\msignu \otimes \tnu,s)$
converges absolutely for $Re(s)>>0$. This product has a
meromorphic continuation on $\C$. We shall denote this
continuation by $L_\psi^S(\msig \otimes \tau,s)$. We have:
\begin{equation} \label{E psi formula}E_\psi
\bigl(f_s,(I_{2n},1)\bigr)=\frac{L^S(\tau,s+\half)}{L^S(\tau,sym^2,2s+1)L_\psi^S(\msig
\otimes \tau,s+1)}\prod_{\nu \in S} \lambda \bigl(s,\tnu \otimes
\msignu,\psi \bigr)(f_\nu). \end{equation}
\end{lem}
Recall that $\prod_{\nu \notin S} L(\tnu,s)$ and $\prod_{\nu
\notin S} L(\tnu,sym^2,s)$ converge absolutely for $Re(s)>>0$ and
that these products have a meromorphic continuation on $\C$; see
\cite{L71}. These continuations are denoted by $L^S(\tau,s)$ and
$L^S(\tau,sym^2,s)$ respectively.
\begin{proof}
Recall that in Section \ref{The symplectic group} we have denoted
by  $W_{Sp_{2n}}$ the Weyl group of $Sp_{2n}(\F)$. We now and
denote by $W_{M_{m;k}}$ the Weyl group of $M_{m;k}$. We fix
$\Omega$, a complete set of representatives of $W_{M_{m;k}}
\backslash W_{Sp_{2n}}$. Recall the Bruhat decomposition
$$Sp_{2n}(\F)=\dotunion_{w \in \Omega} P_{m,k}(\F)wB_{Sp_{2n}}(\F).$$ Clearly
for $w \in \Omega$: $$P_{m,k}(\F)wB_{Sp_{2n}}(\F)=
P_{m,k}(\F)wZ_{Sp_{2n}}(\F).$$ Also, for $w \in \Omega,\, p_1,p_2
\in P_{m;k}(\F), \, u_1,u_2 \in Z_{Sp_{2n}}(\F)$ we have: If
$p_1wu_1=p_2wu_2$ then $$u_2 u_1^{-1} \in Z_w(\F)=Z_{Sp_{2n}}(F)
\cap w^{-1}Z_{Sp_{2n}}(\F)w.$$ Thus, every element $\gamma$ of
$Sp_{2n}(\F)$ can be expressed as $g=pwu$, where $p \in
P_{m,k}(\F)$ and $w \in \Omega$ are determined uniquely and $u \in
Z_{Sp_{2n}}(\F)$ is determined uniquely modulo $Z_w(\F)$ from the
left (note that if $W_{M_{m;k}}w_1=W_{M_{m;k}}w_2$ it does not
follow that $Z_{w_1}(\F)=Z_{w_1}(\F)$. This is why we started from
fixing $\Omega$). Thus, for $f_s \in I({\tau}_{(s)},\msig)$ and
$Re(s)>>0$ we have:
\begin{eqnarray} \label{eis 1}
&& \lefteqn{E_\psi \bigl(f_s,(I_{2n},1) \bigr)} \\
\nonumber &&=\int_{u \in \zspn \backslash Z_{Sp_{2n}}(\A)}
\sum_{\gamma \in P_{m;k}(\F) \backslash \spn}f_s \bigl((\gamma
u,1),I_m,(I_{2k},1) \bigr) \psi^{-1}(u)du \\ \nonumber && =\int_{u
\in \zspn \backslash Z_{Sp_{2n}}(\A)} \psi^{-1}(u) \sum_{w \in
\Omega} \, \sum_{n \in Z_w(\F) \backslash \zspn} f_s
\bigl((wnu,1),I_m,(I_{2k},1) \bigr) du \\ \nonumber && = \sum_{w
\in \Omega} \int_{u \in \zspn \backslash Z_{Sp_{2n}}(\A)}
\psi^{-1}(u) \sum_{n \in Z_w(\F) \backslash \zspn} f_s
\bigl((wnu,1),I_m,(I_{2k},1) \bigr) du \\ \nonumber &&=
 \sum_{w \in \Omega} \int_{u \in Z_w(\F) \backslash
Z_{Sp_{2n}}(\A)} \psi^{-1}(u) f_s\bigl(wu,I_m,(I_{2k},1) \bigr)
du.
 \end{eqnarray}
We now choose $w_0= w_l'(m;k)$ (see \eqref{good long rep}) as the
representative of $J_{2n}$ in $\Omega$. We note that
$Z_{w_0}=Z_{Sp_{2n}} \cap M_{n,k}$. By the same argument used in
page 182 of \cite{Sha78} one finds that
\begin{equation} \label{only wo} \int_{u \in Z_w(\F) \backslash
Z_{Sp_{2n}}(\A)} \psi^{-1}(u) f_s\bigl(wu,I_m,(I_{2k},1) \bigr)=0,
\end{equation} for all  $w \in \Omega, \, w \neq w_0$. Thus, from
\eqref{eis 1} we have:
\begin{eqnarray} \label{eis 2} && \lefteqn{E_\psi
\bigl(f_s,(I_{2n},1)\bigr)} \\ \nonumber && =\int_{u \in \zspn
\cap M_{n,k}(\F) \backslash Z_{Sp_{2n}}(\A)} \psi^{-1}(u)
f_s\bigl(w_0u,I_m,(I_{2k},1) \bigr) du\\ \nonumber && =\int_{u \in
\zspn \cap M_{n,k}(\F) \backslash \bigl(Z_{Sp_{2n}}(\A) \cap
M_{n,k}(\A) \bigr)N_{m;k}(\A)} \psi^{-1}(u)
f_s\bigl(w_0u,I_m,(I_{2k},1) \bigr) du\\ \nonumber && = \int_{n
\in N_{m;k}(\A)}\psi^{-1}(n) \int_{u \in  \zspn \cap M_{n,k}(\F)
\backslash \bigl(Z_{Sp_{2n}}(\A) \cap M_{n,k}(\A) \bigr)}
\psi^{-1}(u) f_s\bigl(w_0un,I_m,(I_{2k},1) \bigr) du \, dn\\
\nonumber && = \int_{n \in N_{m;k}(\A)} \! \psi^{-1}(n) \!  \! \!
\int_{n_1 \in \zglm \backslash \zglma} \int_{n_2 \in \zspk
\backslash \zspka} \! \! \! \! \! \psi^{-1}(n_2)\psi^{-1}(n_1)
f_s\bigl(w_0n,n_1,(n_2,1) \bigr) du \, dn.\end{eqnarray} Recall
that $S$ is a finite set of places of $\F$, such that for all $\nu
\notin S$, $\nu$ is finite and odd, $\tnu \otimes \msignu$ is
unramified and $\psinu$ is normalized. Assume now that $f_s$
corresponds to the following pure tensor of holomorphic smooth
sections, $f_s(g)=T\bigl(\otimes_\nu f_{s,\nu}(g_\nu)\bigr)$,
where $f_{s,\nu} \in V_{I_\nu({\tau_{\nu}}_{(s)},\msig)}$ and for $\nu
\notin S$: $f_{s,\nu}=f_\nu^{0,s}$. By Lemma \ref{whi1}, we have
\begin{equation} \label{eis 3} E_\psi \bigl(f_s,(I_{2n},1)\bigr)=
\prod_\nu \int_{n \in N_{m;k}(\fnu)} \lambda_{\tnu \otimes
\msignu, \psi_\nu} \Bigl(\bigr( \rho(w_0n){f_s}_\nu\bigl)\Bigr)
\psi^{-1}(n) \, dn=\prod_\nu \lambda \bigl(s,\tnu \otimes
\msignu,\psi \bigr)({f_s}_\nu)
\end{equation}
(see Section 3.3 of \cite{T} for the general arguments about
Eulerian integrals). The last equation should be understood as a
global metaplectic analog to Rodier`s local algebraic heredity.
\eqref{eis 3} and \eqref{proven up} imply that for $Re(s)>>0$
\begin{equation} \label{almost E psi formula} E_\psi
\bigl(f_s,(I_{2n},1)\bigr)=\frac{L^S(\tau,s+\half)}{L^S(\tau,sym^2,2s+1)\prod_{\nu
\notin S}L_\psi(\msignu \otimes \tnu,s+1)}\prod_{\nu \in S}
\lambda \bigl(s,\tnu \otimes \msignu,\psi
\bigr)(f_{s,\nu}).\end{equation}

We claim that we may choose $f_s$ as above such that for all $\nu \in S$
\begin{equation} \label{last eq}\lambda \bigl(s,\tnu \otimes \msignu,\psi \bigr)(f_{s,\nu})=1\end{equation}
for all $s\in \C$. Indeed, we choose $f_{s,\nu}$ which is supported on the open Bruhat cell
$$\overline{P_{m;k}(\fnu)}\bigl(w_l'(m;k)Z_{m;k}(\fnu),1 \bigr)$$ which satisfies
$$f_{s,\nu}(w_l'(m;k)z,1)(g,y)=\phi(z)W_{\tau_\nu}(g)W_{\msignu}(y).$$
Here $z\in Z_{m;k}(\fnu)$, $g\in \glmnu$, $y \in \mspknu$, $\phi$ is a properly chosen smooth compactly supported function on
$Z_{m;k}(\fnu)$ and $W_\tnu(I_m)=W\msignu(I_k,1)=1$. We have
$$\lambda \bigl(s,\tnu \otimes \msignu,\psi \bigr)(f_\nu)=\int_{Z_{m;k}(\fnu)}\phi(z)\psi^{-1}(z_{n,n}) \, dn.$$ We now choose
$\phi$ such that \eqref{last eq} holds. For such a choice we have
$$E_\psi
\bigl(f_s,(I_{2n},1)\bigr)=\frac{L^S(\tau,s+\half)}{L^S(\tau,sym^2,2s+1)\prod_{\nu
\notin S}L_\psi(\msignu \otimes \tnu,s+1)}.$$ The absolute
convergence of $\prod_{\nu \notin S}L_\psi(\msignu \otimes
\tnu,s)$ for $Re(s)>>0$ is clear now. Furthermore, the fact that
this product has a meromorphic continuation to $\C$ follows from
the meromorphic continuations of $E_\psi
\bigl(f_s,(I_{2n},1)\bigr)$, $L^S(\tau,s)$ and
$L^S(\tau,sym^2,s)$. Finally, the validity of \eqref{E psi
formula} for all $s$ follows from \eqref{almost E psi formula}.

\end{proof}
\begin{thm} \label{crude} \begin{equation} \label{unramgamm} \prod_{\nu \in
S} \gamma(\msignu \times \tnu,s,\psi_\nu)=\frac{L_\psi^S(\msig
\times \tau,s)}{L_\psi^S(\msig \times \widehat{\tau},1-s)}.
\end{equation}
\end{thm}
\begin{proof}
The global functional equation for the Eisenstein series states
that
$$E(f_s,g)=E\bigl(A(f_s,g)\bigr),$$
where $A$ is the global intertwining operator; see Section IV.1.10
of \cite{MW}. We compute the $\psi$-Whittaker coefficient of both sides of the last equation.
By \eqref{proven down} we have
 \begin{eqnarray*} && \frac{L^S(\tau,s+\half)}{L^S(\tau,sym^2,2s+1)L_\psi^S(\msig \otimes
\tau,s+1)}\prod_{\nu \in S} \lambda \bigl(s,\tnu \otimes
\msignu,\psi \bigr)({f_s}_\nu) \\ &&=
\frac{L^S(\widehat{\tau},-s+\half)L_\psi^S(\msig \otimes
\tau,s)L^S(\tau,sym^2,2s)}
{L^S(\widehat{\tau},sym^2,-2s+1)L_\psi^S(\msig \otimes
\widehat{\tau},-s+1)L_\psi^S(\msig \otimes
\tau,s+1)L^S(\tau,sym^2,2s+1)}\\ && \prod_{\nu \in S}\lambda
\bigl(-s,\widehat{\tnu} \otimes \msignu,\psi\bigr)
\bigl(A_{j_{m,n}(\omega_m'^{-1})}(f^0_{I_1'})\bigr),\end{eqnarray*}
Or equivalently, by the definition of the local coefficients
\begin{equation} \label{crude1}
\prod_{\nu \in
S}C_{\psi_\nu}^{\overline{Sp_{2n}(\fnu)}}\bigl(\overline{{P_{m;k}(\fnu)}},s,\tnu\otimes
\msignu, j_{m,n}(\omega_m'^{-1})\bigr)= \frac
{L^S(\widehat{\tau},-s+\half)L^S(\tau,sym^2,2s )L_\psi^S(\msig
\otimes \tau,s)}
{L^S(\tau,s+\half)L^S(\widehat{\tau},sym^2,-2s+1)L_\psi^S(\msig
\otimes \widehat{\tau},1-s)}.
\end{equation}
In particular, for $k=0$
\begin{equation} \label{crude1 k=0}
\prod_{\nu \in
S}C_{\psi_\nu}^{\overline{Sp_{2m}(\fnu)}}\bigl(\overline{{P_{m;0}(\fnu)}},s,\tnu,
j_{m,n}(\omega_m'^{-1})\bigr)= \frac
{L^S(\widehat{\tau},-s+\half)L^S(\tau,sym^2,2s)}{L^S(\tau,s+\half)L^S(\widehat{\tau},sym^2,-2s+1)}.
\end{equation}

Dividing \eqref{crude1} and \eqref{crude1 k=0} we get
\eqref{unramgamm}.
\end{proof}

\subsection{Computation of $C_{\psi}^{\mspm}\bigl(\overline{P_{m;0}(\F)},s,\tau,
\omega_m'^{-1} \bigr)$ for generic representations}
\label{true for cusp}

\begin{thm} \label{soudry said cusp}Let $\F$ be a p-adic field and let
$\tau$ be an irreducible admissible supercuspidal representation
of $\glm$. There exists an exponential function $c_\F(s)$ such
that
$$ C_{\psi}^{\mspm}\bigl(\overline{P_{m;0}(\F)},s,\tau,
\omega_m'^{-1} \bigr)=c_\F(s)\frac
{\gamma(\tau,sym^2,2s,\psi)}{\gamma(\tau,s+\half,\psi)}.$$
\end{thm}

\begin{proof}
Since $\tau$ is supercuspidal it is also generic. Proposition 5.1
of \cite{Sha 3} implies now that there exists a number field $\K$,
a non-degenerate character $\widetilde{\psi}$ of
$Z_{GL_n(\K)}\backslash Z_{GL_n(\A)}$ and an irreducible cuspidal representation  $\pi \simeq\otimes_\nu
\pi_\nu$ of $GL_n(\A)$ such that\\
1. $\K_{\nu_0}=\F$ for some place $\nu_0$ of $\K$.\\
2. $\widetilde{\psi}_{\nu_0}=\psi$.\\
3. $\pi_{\nu_0}=\tau$.\\
4. For any finite place $\nu \neq \nu_0$ of $\K$, $\pi_\nu$ is
unramified.\\
Define $S$ to be the finite set of places of $\K$ which consists
of $\nu_0$, of the infinite and even places and of the finite
places where $\widetilde{\psi}_{\nu}$ is not normalized. From the
fourth part of Theorem 3.5 of \cite{Sha 2} it follows that
$$\prod_{\nu \in
S}\gamma_{\fnu}(\pi_\nu,s,\psi_\nu)=\frac{L^S(\pi,s)}{L^S(\widehat{\pi},1-s)}$$
and that
$$\prod_{\nu \in
S}\gamma_{\fnu}(\pi_\nu,s,sym^2,\psi_\nu)=\frac{L^S(\pi,sym^2,s)}{L^S(\widehat{\pi},sym^2,1-s)}.$$
Therefore, \eqref{crude1 k=0} can be written as

$$
\prod_{\nu \in
S}C_{\psi_\nu}^{\overline{Sp_{2m}(\fnu)}}\bigl(\overline{{P_{m;0}(\fnu)}},s,\pi_\nu,
j_{m,n}(\omega_m'^{-1})\bigr)=\prod_{\nu \in S} \frac
{\gamma(\pi_\nu,sym^2,2s,\psi)}{\gamma(\pi_\nu,s+\half,\psi)}.
$$

This implies that this theorem will be proven once we show that
for all $\nu \in S$, $\nu \neq \nu_0$, there exists an exponential
 function $c_\nu(s)$ such that
$$C_{\psi_\nu}^{\overline{Sp_{2m}(\fnu)}}\bigl(\overline{{P_{m;0}(\fnu)}},s,\pi_\nu,
j_{m,n}(\omega_m'^{-1})\bigr)= c_\nu(s)\frac
{\gamma_\fnu(\pi_\nu,sym^2,2s,\psi)}{\gamma_\fnu(\pi_\nu,s+\half,\psi)}.
$$
Since for all for all $\nu \in S$, $\nu \neq \nu_0$, $\pnu$ is the
generic constituent of a principal series representation series,
this follows from Theorem \ref{soudry said it}.
\end{proof}
\begin{thm} \label{soudry said always} Let $\F$ be a p-adic field and let
$\tau$ be an irreducible admissible generic representation of
$\glm$. There exists an exponential function $c_\F(s)$ such that
\begin{equation} \label{always always} C_{\psi}^{\mspm}\bigl(\overline{P_{m;0}(\F)},s,\tau,
\omega_m'^{-1} \bigr)=c_\F(s)\frac
{\gamma(\tau,sym^2,2s,\psi)}{\gamma(\tau,s+\half,\psi)}.\end{equation}
\end{thm}
\begin{proof} By Chapter II of \cite{BZ}, $\tau$ may be realized
as a sub-representation of $$\tau'=Ind^\glm_{Q(\F)}
(\otimes_{i=1}^r \tau_i),$$ where $Q(\F)$ is a standard parabolic
subgroup of $\glm$ whose Levi part, $M(\F)$, is isomorphic to
$$GL_{n_1}(\F) \times GL_{n_2}(\F) \ldots \times GL_{n_r}(\F)$$
and where for all $1 \leq i \leq r$, $\tau_i$ is  an irreducible
admissible supercuspidal representation of $GL_{n_i}(\F)$. Since
for all $1 \leq i \leq r$, $\tau_i$ has a unique Whittaker model
it follows from the heredity property of the Whittaker model that
$\tau'$ has a unique Whittaker model; see \cite{Rod}. This implies
that $\tau$ is the generic constituent of $\tau'$. Hence,
$$ C_{\psi}^{\mspm}\bigl(\overline{P_{m;0}(\F)},s,\tau,
\omega_m'^{-1} \bigr)=
C_{\psi}^{\mspm}\bigl(\overline{P_{m;0}(\F)},s,\tau',
\omega_m'^{-1} \bigr).$$ Thus, it is sufficient to prove
\eqref{always always} replacing $\tau$ with $\tau'$. By similar
arguments to those used in Lemma \ref{lc zigel prin} and Theorem
\ref{soudry said it} of Section \ref{prin comp with soo} one shows
that there exists $d \in \{ \pm 1 \}$ such that
$$C_{\psi}^{\mspm}\bigl(\overline{P_{m;0}(\F)},s,\tau',
\omega_m'^{-1} \bigr)=d
\prod_{i=1}^r\Bigl(C_{\psi}^{\overline{Sp_{n_i}(\F)}}\bigl(\overline{P_{n_i;0}(\F)},s,\tau_i,
\omega_{n_i}'^{-1} \bigr) \prod_{j=i+1}^r \gamma(\tau_i \times
\tau_j,2s,\psi)\Bigr).$$ Since for $1 \leq i \leq r$, $\tau_i$ are
irreducible admissible supercuspidal representations it follows
from Theorem \ref{soudry said cusp} that there exits an exponetial
factor, $c_\F(s)$, such that
$$C_{\psi}^{\mspm}\bigl(\overline{P_{m;0}(\F)},s,\tau',
\omega_m'^{-1} \bigr)=c_\F(s) \prod_{i=1}^r\Bigl(\frac
{\gamma(\tau_i,sym^2,2s,\psi)}{\gamma(\tau_i,s+\half,\psi)}
\prod_{j=i+1}^r \gamma(\tau_i \times \tau_j,2s,\psi)\Bigr).$$
Using the known multiplicativity of the symmetric square
$\gamma$-factor (see Part 3 of Theorem 3.5 of \cite{Sha 2}) we
finish.

\end{proof}

\newpage
\section{Irreducibility theorems} \label{Irreducibility theorems}
Let $\F$ is a p-adic field. In this chapter we prove several
irreducibility theorems for parabolic induction on the metaplectic
group. Through this chapter we shall assume that the inducing
representations are smooth irreducible admissible cuspidal unitary
genuine generic.

In Section \ref{irr prin uni sec} we prove the irreducibility of
principal series representations (induced from a unitary
characters); see Theorem \ref{irr unitary prinipal}. A
generalization of the argument given in this theorem is used in
Section \ref{General parabolic case} where we prove a general
criterion for irreducibility of parabolic induction; see Theorem
\ref{big thm for max}. As a corollary we reduce the question of
irreducibility of parabolic induction on $\mspn$ to irreducibility
of representations of metaplectic groups of smaller dimension
induced from representations of maximal parabolic groups; see
Corollary \ref{irr from irr}. In Section \ref{irr max sec} we
reduce the question of irreducibility of representations of
$\mspn$ induced from the Siegel parabolic subgroup to the question
of irreducibility of a parabolic induction on $\soo$; see Theorem
\ref{metaplectic and so}. We then use \cite{Sha 4} to give a few
corollaries.

Our method here is an application of Theorem \ref{silberger} to
the computations and results presented in Chapters \ref{cha gam},
\ref{second paper} and \ref{global sec}. The link between Theorem
\ref{silberger} and the local coefficients is explained in Chapter
\ref{chapter def lc}; see \eqref{beta is lc}. Throughout this
chapter we shall use various definitions and notation given in
previous chapters. Among them are $\pi^w$ and the notion of a
regular and of a singular representation (see Section \ref{app
bruhat theory}), $W_\pt$ (see Section \ref{The symplectic group}),
$W(\pi)$ (see \eqref{wpi def}) and $\Sigma_{\pt}$ (see Section
\ref{knapp sec}).

\subsection{Irreducibility of principal series representations of
$\mspn$ induced from unitary characters.} \label{irr prin uni sec}
\begin{lem} Let $\beta_1$ and $\beta_2$ be two characters of $\F^*$.
Denote $\beta=\beta_1\beta_2^{-1}$. If $\F$ is a p-adic field then
\begin{equation}\label{lc gl2}C_{\psi}^{GL_2(\F)}\bigl(B_{GL_2}(\F),(s_1,s_2),\beta_1 \otimes
\beta_2, (\begin{smallmatrix}
& {0}& {1} \\
& {1} & {0}\end{smallmatrix})
\bigr)=\frac{\beta(\pi^{m(\beta)-n})q^{(n-m(\beta)
n)(s_1-s_2)}}{G(\beta,\psi^{-1})} \frac{L \bigl
(\beta^{-1},1-(s_1-s_2)\bigr)}{L \bigl(\beta,(s_1-s_2) \bigr)},
\end{equation} where $n$ is the
conductor of $\psi$. If $\F=\R$ then

\begin{equation}\label{lc gl2 r}C_{\psi}^{GL_2(\R)}\bigl(B_{GL_2}(\R),(s_1,s_2),\beta_1 \otimes
\beta_2, (\begin{smallmatrix}
& {0}& {1} \\
& {1} & {0}\end{smallmatrix}) \bigr)=\beta_2(-1)e^{\frac
{-i\pi\beta(-1)}{2}} (2\pi)^{-(s_1-s_2)}\frac{L \bigl
(\beta^{-1},1-(s_1-s_2)\bigr)}{L \bigl(\beta,(s_1-s_2) \bigr)}.
\end{equation}
\end{lem}

This lemma can be proved by similar computations to those
presented in Sections \ref{padic comp} and \ref{real case}. The
p-adic computations here are much easier than those presented in
this dissertation. The lemma also follows from more general known
results; see Lemma 2.1 of \cite{Sha83} for the p-adic case and see
Theorem 3.1 of \cite{Sha85} for the real case.

\begin{thm} \label{irr unitary prinipal}
Let $\alpha_1,\ldots,\alpha_n$ be $n$ unitary characters of
$\F^*$. Let $\alpha$ be the character of $\overline{\tspn}$
defined by
$$\bigl(diag(a_1,\ldots,a_n,a^{-1},\ldots,a_n^{-1}),\epsilon \bigr)
\mapsto \epsilon \gamma_\psi^{-1}(\prod_{i=1}^na_i)\prod_{i=1}^n
\alpha_i(a_i).$$ Then $I(\alpha)$ is irreducible.
\end{thm}
\begin{proof}
Since $\alpha$ is unitary, $I(\alpha)$ is also unitary. Therefore,
the irreducibility of $I(\alpha)$ will follow once we show that
$$Hom_{\mspn} \bigl(I(\alpha),I(\alpha)\bigr) \simeq \C.$$ For
$1\leq i<j \leq n$ define
$$w_{(i,j)}=\widehat{
\begin{pmatrix} _{I_{i-1}} & _{ } & _{ } & _{ } &
_{ } \\ _{ } & _{ }  & _{ } & _{1}  & _{ }
\\ _{ } & _{ } & _{I_{j-i-2}} & _{ } & _{ } \\
_{ } & _{1}& _{ }  & _{ } & _{ } \\
_{ } & _{ }  & _{ } & _{ } & _{I_{n-j+1}}
\end{pmatrix}}$$
and
$$w'_{(i,j)}=\tau_{\{i,j\}}w_{(i,j)}.$$
A routine exercise shows that \begin{equation} \label{sigma borel}
\Sigma_{\bspn}=\{w_{(i,j)} \mid 1 \leq i < j \leq n \} \cup
 \{ \tau_{\{r\}} \mid 1 \leq  r \leq n \} \cup \{w'_{(i,j)} \mid 1
\leq<i < j \leq n \}.\end{equation} Note that,
$$w_{(i,j)}\widehat{diag(a_1,\ldots,a_n)}w_{(i,j)}^{-1}=
\widehat{diag(a_1,\ldots,a_{i-1},a_j,a_{i+1},\ldots,a_{j-1},a_i,a+{j+1},\ldots,a_n)},$$

$$w_{\{r\}}\widehat{diag(a_1,\ldots,a_n)}w_{\{r\}}^{-1}=
\widehat{diag(a_1,\ldots,a_{r-1},a_r^{-1},a_{r+1},\ldots,a_n)}$$
and that

$$w'_{(i,j)}\widehat{diag(a_1,\ldots,a_n)}w_{(i,j)}'^{-1}=
\widehat{diag(a_1,\ldots,a_{i-1},a_j^{-1},a_{i+1},\ldots,a_{j-1},a_i^{-1},a+{j+1},\ldots,a_n)}.$$
Therefore $w_{(i,j)} \in W(\alpha)$ $\Leftrightarrow$
$\alpha_i=\alpha_j$, $\tau_{\{r\}} \in W(\alpha)$
$\Leftrightarrow$ $\alpha_r$ is quadratic and $w_{(i,j)}' \in
W(\alpha)$ $\Leftrightarrow$ $\alpha_i=\alpha_j^{-1}$.
Furthermore, $W(\alpha)$ is generated by $\Sigma_{B_\spn} \cap
W(\alpha)$. Thus, using Theorem \ref{silberger} and \eqref{beta is
lc}, the proof of this theorem amounts to showing that
\begin{equation} \label{what to show}C_{\psi}^{\mspn}(\mpt,\overrightarrow{s},(\otimes_{i=1}^n
\alpha_i)\otimes,w \bigl)
C_{\psi}^{\mspn}(\mpt,\overrightarrow{s^w},(\otimes_{i=1}^n
\alpha_i)^{w},w^{-1} \bigl)=0 \end{equation} for all $w \in
\Sigma_{\bspn} \cap W(\alpha)$. We prove it for each of the three
types in the right-hand side of \eqref{sigma borel}.

Suppose that $w_{(i,j)} \in W(\alpha)$. We write:
\begin{equation} \label{decompose wij} w_{(i,j)}=w_{(i,i+1)}w_{(i+1,i+2)}\ldots
w_{(j-1,j)}w_{(j-2,j-1)}\ldots w_{(i+1,i)}.\end{equation} We claim
that the expression in the right-hand side of \eqref{decompose
wij} is reduced. Indeed,

$$\Sigma'_{\bspn}=\{w_{(i,+1)} \mid 1
\leq<i <n \} \cup
 \{ \tau_{\{1\}}\} \subset \Sigma_{B_\spn}$$
is the subset of reflections corresponding to simple roots and the
length of $w_{(i,j)}$ is $2(j-i)-1$ (any claim about the length of
a given Weyl element $w$ may be verified by counting the number of
positive root subgroups mapped by $w$ to negative root subgroups).
Thus, we may use the same argument as in Lemma \ref{lem heart} and
conclude that there exists $c \in \{ \pm 1 \}$ such that
\begin{eqnarray} \label{first type factor}&&C_{\psi}^{\mspn}\bigl(\bspn,\overrightarrow{s},(\otimes_{i=1}^n
\alpha_i)),w_{(i,j)} \bigl)= \\ \nonumber && c
\Bigl(\prod_{k=i}^{j-2}f_k(s)f_k'(s)\Bigr)
C_{\psi}^{\mspn}\bigl(\bspn,\overrightarrow{s}^{\widetilde{w}},(\otimes_{i=1}^n
\alpha_i)^{\widetilde{w}} ,w_{(j-1,j)}\bigr),
\end{eqnarray}
where \begin{eqnarray*} &&
f_k(\overrightarrow{s})=C_{\psi}^{\mspn}\bigl(\bspn,\overrightarrow{s}^{w_{(k)}},(\otimes_{i=1}^n
\alpha_i)^{w_{(k)}},w_{(k,k+1)} \bigl) \\ &&
f_k'(\overrightarrow{s})=C_{\psi}^{\mspn}\bigl(\bspn,\overrightarrow{s}^{w'_{(k)}},(\otimes_{i=1}^n
\alpha_i)^{w'_{(k)}},w_{(k,k+1)} \bigl), \end{eqnarray*} where
\begin{eqnarray*} w_{(k)}&=&w_{(i,i+1)}w_{(i+1,i+2)}\ldots w_{(k-1,k)},\\
\widetilde{w}&=&w_{(i,i+1)}w_{(i+1,i+2)}\ldots
w_{(j-1,j)}w_{(j-2,j-1)},\\
w'_{(k)}&=&w_{(i,i+1)}w_{(i+1,i+2)}\ldots
w_{(j-1,j)}w_{(j-2,j-1)}w_{(j-1,j)}w_{(j-2,j-1)}\ldots
w_{(k+1,k+2)}.
\end{eqnarray*}
Since all the local coefficients in the right-hand side of
\eqref{first type factor} correspond to simple reflections we may
use the same argument as in Lemma \ref{id with c} and conclude
that

\begin{eqnarray} \label{first type id} \\ \nonumber && \! \! \! \! \! \! \! C_{\psi}^{\mspn}\bigl(\bspn,\overrightarrow{s}^{w_{(k)}},(\otimes_{i=1}^n
\alpha_i)^{w_{(k)}},w_{(k,k+1)} \bigl)=
C_\psi^{GL_2(\F)}\bigr(B_{GL_2}(\F),(s_i,s_{i+k}),\alpha_i\otimes
\alpha_{i+k},{\omega_2},(\begin{smallmatrix} _{0} & _{1}\\_{1} &
_{0}\end{smallmatrix})  \bigl),\\ \nonumber && \! \! \! \! \! \!
\!
C_{\psi}^{\mspn}\bigl(\bspn,\overrightarrow{s}^{\widetilde{w}},(\otimes_{i=1}^n
\alpha_i)^{\widetilde{w}} ,w_{(j-1,j)}\bigr)=
C_\psi^{GL_2(\F)}\bigr(B_{GL_2}(\F),(s_i,s_j),\alpha_i\otimes
\alpha_j,{\omega_2},(\begin{smallmatrix} _{0} & _{1}\\_{1} &
_{0}\end{smallmatrix})  \bigl),\\ \nonumber && \! \! \! \! \! \!
\!
C_{\psi}^{\mspn}\bigl(\bspn,\overrightarrow{s}^{w'_{(k)}},(\otimes_{i=1}^n
\alpha_i)^{w'_{(k)}},w_{(k,k+1)} \bigl) =
C_\psi^{GL_2(\F)}\bigr(B_{GL_2}(\F),(s_{i+k},s_j),\alpha_{i+k}\otimes
\alpha_{j},{\omega_2},(\begin{smallmatrix} _{0} & _{1}\\_{1} &
_{0}\end{smallmatrix})  \bigl). \end{eqnarray} Since
$\alpha_1,\ldots,\alpha_k$ are unitary, \eqref{first type id} and
\eqref{lc gl2} implies that for $i \leq k \leq j-2$,
$$C_{\psi}^{\mspn}\bigl(\bspn,\overrightarrow{s}^{w_{(k)}},(\otimes_{i=1}^n
\alpha_i)^{w_{(k)}},w_{(k,k+1)} \bigl)$$ and
$$C_{\psi}^{\mspn}\bigl(\bspn,\overrightarrow{s}^{w'_{(k)}},(\otimes_{i=1}^n
\alpha_i)^{w'_{(k)}},w_{(k,k+1)} \bigl)$$ are holomorphic at
$\overrightarrow{s}=0$. Also, since $w_{(i,j)} \in W(\alpha)$
implies that $\alpha_i=\alpha_j$, \eqref{first type id} and
\eqref{lc gl2} imply that
$$C_{\psi}^{\mspn}\bigl(\bspn,\overrightarrow{s}^{\widetilde{w}},(\otimes_{i=1}^n
\alpha_i)^{\widetilde{w}} ,w_{(j-1,j)}\bigr)$$ vanishes for
$\overrightarrow{s}=0$. Recalling \eqref{first type factor} we now
conclude that if $w=w_{(i,j)} \in W(\alpha)$ then \eqref{what to
show} holds.

Suppose now that $\tau_{\{r\}} \in W(\alpha)$. We write
\begin{equation} \label{decompose taur} \tau_{\{r\}}=w_{(r,r+1)}w_{(r+1,r+2)}\ldots
w_{(n-1,n)}\tau_{\{n\}}w_{(n-1,n)}w_{(n-2,n-1)}\ldots
w_{(r+1,r)}.\end{equation} The reader may check that the
expression in the right-hand side of \eqref{decompose taur} is
reduced. We now use the same arguments we used for $w=w_{(i,j)}$:
We decompose
$$C_{\psi}^{\mspn}\bigl(\bspn,\overrightarrow{s},(\otimes_{i=1}^n
\alpha_i),\tau_{\{r\}} \bigr)$$ into $1+2(n-i)$ local
coefficients. $2(n-i)$ of them are of the form \eqref{lc gl2}.
These factors are holomorphic at $\overrightarrow{s}=0$. The
additional local coefficient, the one corresponding to
$\tau_{\{n\}}$ is
$$C_\psi^{\msl}\bigr(\overline{B_{SL_2}(\F)},s_r,\alpha_r,(\begin{smallmatrix} _{0} & _{1}\\_{-1} &
_{0}\end{smallmatrix}) \bigl).$$ Theorem \ref{the formula} implies
that there exists $c \in \C^*$ such that
 \begin{eqnarray*} &&C_\psi^{\msl}(\overline{B_{SL_2}(\F)},s,\chi,(\begin{smallmatrix} _{0} & _{1}\\_{-1} &
_{0}\end{smallmatrix}) \bigl)
C_\psi^{\msl}(\overline{B_{SL_2}(\F)},-s,\chi^{-1},(\begin{smallmatrix}
_{0} & _{1}\\_{-1} & _{0}\end{smallmatrix}) \bigl) \\ &&=c
\frac{L_{\F}(\chi^{2},-2s+1)}
{L_{\F}(\chi^{2},2s)}\frac{L_{\F}(\chi^{2},2s+1)}
{L_{\F}(\chi^{2},-2s)}.\end{eqnarray*} Since $\tau_{\{r\}} \in
W(\alpha)$ implies that $\alpha_r$ is quadratic, we now conclude
that \eqref{what to show} holds for  $w=\tau_{\{r\}}$.

Finally, assume that $w'_{(i,j)} \in W(\alpha)$. We write it as a
reduced product of simple reflections:
\begin{eqnarray*} w'_{(i,j)}= &&
w_{(j,j+1)}w_{(j+1,j+2)}\ldots
w_{(n-1,n)}\tau_{\{n\}}w_{(n-1,n)}w_{(n-2,n-1)}\ldots w_{(j+1,i)}
\\ && w_{(i,i+1)}w_{(i+1,i+2)}\ldots
w_{(j-1,j)}w_{(j-2,j-1)}\ldots w_{(i+1,i)} \\
&&w_{(j,j+1)}w_{(j+1,j+2)}\ldots
w_{(n-1,n)}\tau_{\{n\}}w_{(n-1,n)}w_{(n-2,n-1)}\ldots w_{(j+1,i)}.
\end{eqnarray*}
We then decompose
$$C_{\psi}^{\mspn}\bigl(\bspn,\overrightarrow{s},(\otimes_{i=1}^n
\alpha_i),w'_{(i,j)} \bigr)$$ into $1+2(n-i)$ local coefficients
coming either from $GL_2(\F)$ or from $\msl$. All these local
coefficients are holomorphic at $\overrightarrow{s}=0$. The factor
corresponding to $w_{(j,j-1)}$ equals
$$C_\psi^{GL_2(\F)}\bigr(B_{GL_2}(\F),(s_i,-s_j),\alpha_i\otimes
\alpha_j^{-1},(\begin{smallmatrix} _{0} & _{1}\\_{1} &
_{0}\end{smallmatrix})  \bigl).$$ Since $w'_{(i,j)} \in W(\alpha)$
implies that $ \alpha_i=\alpha_j^{-1}$ we conclude, using
\eqref{lc gl2}, that \eqref{what to show} holds for
$w=w'_{(i,j)}$, provided that $w'_{(i,j)} \in W(\alpha)$.
\end{proof}{\bf Remarks}:\\
1. Assume that $\F$ is a p-adic field of odd residual
characteristic. For the irreducibility of principal series
representations of $\msl$, induced from unitary characters see
\cite{Mo88}. For the irreducibility of principal series
representations induced from unitary characters to the $\C^1$
cover of $Sp_4(\F)$  see \cite{Zo09}. During the final
preparations of this manuscript the author of these lines
encountered a preprint which proves Theorem \ref{irr unitary
prinipal} using the theta
correspondence; see \cite{HMa}.\\
2. One can show that Theorem \ref{silberger} applies also to the
field of real numbers in the case of a parabolic induction from
unitary characters of $\overline{B_{Sp_{2n}}(\R)}$. Thus,
repeating the same argument used in this section, replacing
Theorem \ref{the formula} with Theorem \ref{real thm} and
\eqref{lc gl2} with \eqref{lc gl2 r}, one concludes that Theorem
\ref{irr unitary prinipal} applies for the real case as well.
\subsection{Irreducibility criteria for parabolic induction} \label{General parabolic case}
\begin{thm} \label{big thm for max}Let $\overrightarrow{t}=(n_1,n_2,\ldots,n_r;k)$ where
$n_1,n_2,\ldots,n_r,k$ are $r+1$ non-negative integers whose sum
is $n$. For $1 \leq i \leq r$ let $\tau_i$ be an irreducible
admissible supercuspidal unitary  representation of $GL_{n_i}(\F)$
and let $\msig$ be an an irreducible admissible supercuspidal
$\psi$-generic genuine representation of $\mspk$. Denote
$\pi=\bigl(\otimes_{i=1}^r (\gamma_{\psi}^{-1}\otimes
{\tau_i})\bigr) \otimes \msig$. $I(\pi)$ is reducible if and only
if there exists $1 \leq i \leq r$ such that $\tau_i$ is self dual
and
\begin{equation} \label{the condition}C_{\psi}^{\overline{Sp_{2(k+n_i)(\F)}}}\bigl(\overline{P_{n_i;k}(\F)},0,\tau_i\otimes
\msig, j_{n_i,k+n_i}(\omega_{n_i}'^{-1})\bigr) \neq 0
\end{equation}
\end{thm}
\begin{proof}
Since $j_{n_i,n}(\omega_{n_i}'^{-1})$ is of order two as a Weyl
element it follows that if $\tau_i$ is self dual then
$$C_{\psi}^{\overline{Sp_{2(k+n_i)(\F)}}}\bigl(\overline{P_{n_i;k}(\F)},0,\tau_i\otimes
\msig, j_{n_i,k+n_i}(\omega_{n_i}'^{-1})\bigr)=0$$ if and only if
\begin{equation} \label{vanish tau self sual}C_{\psi}^{\overline{Sp_{2(k+n_i)(\F)}}}\bigl(\overline{P_{n_i;k}(\F)},s,\tau_i\otimes
\msig, j_{n_i,k+n_i}(\omega_{n_i}'^{-1})\bigr)
C_{\psi}^{\overline{SP_{k+2n_i}(\F)}}\bigl(\overline{P_{n_i;k}(\F)},-s,\widehat{\tau_i}\otimes
\msig, j_{n_i,k+n_i}(\omega_{n_i}'^{-1})\bigr)\end{equation}
vanishes at $s=0$. Thus, since $I(\pi)$ is unitary, we only have
to show that
\begin{equation} \label {dim hom more 1}dim \bigl(Hom_{\mspn} \bigl(I(\pi),I(\pi)\bigr)\bigr) >1
\end{equation} if and only if there exits $1 \leq i \leq r$ such that
$\tau_i$ is self dual and \eqref{vanish tau self sual} does not
vanish at $s=0$.

Suppose first that there exits $1 \leq i \leq r$ such that
$\tau_i$ is self dual and \eqref{vanish tau self sual} does not
vanish at $s=0$. Since for any $w \in W_\pt$, $I(\pi)$ and
$I(\pi^w)$ have the same Jordan Holder series we may assume that
$i=r$. It follows from \eqref{w action} that
$w_0=j_{n_r,n}(\omega_{n_r}'^{-1}) \in \sigma_\pt \cap W(\pi)$.
Since $w_0$ is a simple reflection we may use a similar argument
to the one used in Lemma \ref{id with c} and conclude that
$$
C_{\psi}^{\mspn}\bigl(\mpt,\overrightarrow{s},\bigl(\otimes_{i=1}^r
\tau_i \bigr)\otimes \msig,w_0\bigr)=
C_{\psi}^{\overline{Sp_{2(k+n_i)}(\F)}}\bigl(\overline{P_{n_i;k}(\F)},s_r,\tau_i\otimes
\msig, w_0 \bigr),$$where
$\overrightarrow{s}=(s_1,s_2,\ldots,s_r)$. Thus, our assumption
implies that
\begin{equation} \label{vanish or not}C_{\psi}^{\mspn}(\mpt,\overrightarrow{s},(\otimes_{i=1}^r
\tau_i)\otimes \overline{\sigma},w_0 \bigl)
C_{\psi}^{\mspn}(\mpt,\overrightarrow{s^{w_0}},((\otimes_{i=1}^r
\tau_i)\otimes \overline{\sigma})^{w_0},w_0 \bigl)\end{equation}
does not  vanish at $s=0$. Theorem \ref{silberger} and \eqref{beta
is lc} imply now that \eqref{dim hom more 1} holds.

We now assume that for any $1 \leq i \leq r$, if $\tau_i$ is self
dual then \eqref{vanish tau self sual} vanishes at $s=0$. Again,
by theorem \ref{silberger} and \eqref{beta is lc} we only have to
show that \eqref{vanish or not} vanishes at $\overrightarrow{s}=0$
for any $w_0 \in \Sigma_{\pt} \cap W(\pi)$. Similar to the proof
of Theorem \ref{irr unitary prinipal}, there are three possible
types of $w_0 \in
\Sigma_{\pt} \cap W(\pi)$:\\
Type 1. $\tau_i \simeq \tau_j$ for some $1 \leq i<j \leq r$ then
the if and only if the Weyl element that inter change the
$GL_{n_i}(\F)$ and the $GL_{n_j}(\F)$ blocks lies in
$\Sigma_{\pt} \cap W(\pi)$.\\
Type 2. $\tau_i \simeq \widehat{\tau_j}$ for some $1 \leq i<j \leq
r$  if and only if the Weyl element that inter change the
$GL_{n_i}(\F)$ with the "dual" $GL_{n_j}(\F)$ blocks lies in $\Sigma_{\pt} \cap W(\pi)$.\\
Type 3. $\tau_i$ is self dual for some $1 \leq i \leq r$ if and
only if the Weyl element that inter change the $GL_{n_i}(\F)$ with
its "dual" block lies in $\Sigma_{\pt} \cap W(\pi)$.\\

In fact, by switching from $\pi$ to $\pi^w$ for some $w \in
W_\pt$, we may assume that there are no elements in $\Sigma_{\pt}
\cap W(\pi)$ of type 2. Indeed, Let $I \subseteq \{1,2,\ldots,r\}$
such that
$$\{1,2,\ldots,r\}=\dotunion_{i \in I} A_i,$$ where $A_i$ are the
equivalence classes $$A_i=\{1 \leq j \leq r \mid \tau_i \simeq
\tau_j \, \, or \, \, \tau_i \simeq \widehat{\tau_j}\}.$$ By
choosing $w \in W_\pt$ properly we may assume that $\tau_i \simeq
\tau_j$ for all $j \in A_i$. Thus, we only prove that
\eqref{vanish or not} vanishes at $\overrightarrow{s}=0$ for any
$w_0 \in \Sigma_{\pt} \cap W(\pi)$ of type 1 or 3.

Assume that $\tau_i \simeq \tau_j$. Let $w_0\in \Sigma_{\pt} \cap
W(\pi)$ be the corresponding Weyl element. We decompose
\eqref{vanish or not} into a product of local coefficients
corresponding to simple reflections which may be shown to be equal
to local coefficients of the form
\begin{equation} \label{gl lc}C_{\psi}^{GL_{n_p+n_q}(\F)}
\bigr(P^0_{n_p,n_q}(\F),(s_p,s_q),\tau_p \otimes
\tau_q,\varpi_{q,p} \bigl) C_{\psi}^{GL_{n_q+n_p}(\F)}
\bigr(P^0_{n_q,n_p}(\F),(s_q,s_p),\tau_q \otimes
\tau_p,\varpi_{p,q} \bigl). \end{equation}

All these factors are analytic at $(0,0)$;  see Theorem 5.3.5.2 of
\cite{Sil book}. One of these factors corresponds to
$(p,q)=(i,j)$. Since by assumption $\tau_i \simeq \tau_j$, the
well-known reducibility theorems for parabolically induced
representation of $\gln$ (see the first remark on page 1119 of
\cite{Gol94}, for example) implies that the factor that
corresponds to $p=i$, $q=j$ vanishes at $(0,0)$. This shows that
\eqref{vanish or not} vanishes at $\overrightarrow{s}=0$ for any
$w_0\in \Sigma_{\pt} \cap W(\pi)$ of type 1.

Assume that $\tau_i$ is self dual. Let $w_0\in \Sigma_{\pt} \cap
W(\pi)$ be the corresponding Weyl element. We decompose
\eqref{vanish or not} into a product which consist of elements of
the form \eqref{gl lc} and of factor of the form \eqref{vanish tau
self sual}. All the factors of the form \eqref{gl lc} are analytic
and $(0,0)$. Since $\tau_i$ is self dual, by our assumption the
other factor vanishes at $s=0$. This shows that \eqref{vanish or
not} vanishes at $\overrightarrow{s}=0$ for any $w_0\in
\Sigma_{\pt} \cap W(\pi)$ of type 3. \
\end{proof}

\begin{cor} \label{big thm for max again} We keep the notations and
assumptions of Theorem \ref{big thm for max}. $I(\pi)$ is
reducible if and only if there exists $1 \leq i \leq r$ such that
$\tau_i$ is self dual and
\begin{equation} \label{the condition again}\gamma(\msig \times \tau_i,0,
\psi) \gamma(\tau_i,sym^2,0,\psi) \neq 0
\end{equation}
\end{cor}

{\begin{proof} Let $\tau$ be an irreducible admissible generic
representation of $\glm$. From the definition of $\gamma(\msig
\times \tau,s, \psi)$, \eqref{gama def}, and from Theorem
\ref{soudry said always} it follows that
$$C_{\psi}^{\mspn}
\bigl(\overline{P_{m;k}(\F)},s,\tau\otimes \msig,
j_{m,n}(\omega_m'^{-1}) \bigr)=c(s)\gamma(\msig \times \tau,s,
\psi)\frac
{\gamma(\tau,sym^2,2s,\psi)}{\gamma(\tau,s+\half,\psi)}$$ for some
exponential factor $c(s)$. By (6.1.4) in  page 108 of \cite{Sha84}
we have
\begin{equation} \label{gama tau tau hat}
\gamma(\widehat{\tau},1-s,\psi)\gamma(\tau,s,\psi)=\tau(-I_m) \in
\{\pm 1\}.\end{equation}
 Therefore, if we
assume in addition that $\tau$ is self dual we know that
${\gamma(\tau, \half,\psi)}\in \{\pm 1\}$. This implies that
\eqref{the condition} may be replaced with \eqref{the condition
again} \end{proof} The following two corollaries follow
immediately from Theorem \ref{big thm for max}.
\begin{cor} \label{irr from irr}
With the notations and assumptions of Theorem \ref{big thm for
max}, $I(\pi)$ is irreducible if and only if $I(\tau_i, \msig)$ is
irreducible for every $1 \leq i \leq r$.
\end{cor}
\begin{cor} \label{nice for max} Let $\overrightarrow{t}=(n_1,n_2,\ldots,n_r;0)$ where
$n_1,n_2,\ldots,n_r$ are $r$ non-negative integers whose sum is
$n$. For $1 \leq i \leq r$ let $\tau_i$ be an irreducible
admissible supercuspidal unitary representation of $GL_{n_i}(\F)$.
Denote $\pi=\otimes_{i=1}^r (\gamma_{\psi}^{-1}\otimes {\tau_i})
$. $I(\pi)$ is reducible if and only if there exits $1 \leq i \leq
r$ such that $\tau_i$ is self dual and $\gamma(\tau_i,sym^2,0)
\neq 0$.
\end{cor}

\subsection{A comparison with $\soo$} \label{irr max sec}
Let $\soo$ be the special orthogonal group:
$$\soo=\{g \in GL_{2n+1}(\F) \mid
g J'_{2n+1}g^t=J'_{2n+1}, \, \det (g)=1 \},$$ where
$J'_{n}=\left(\begin{smallmatrix}
_{ } & _{ } & _{ } & _{1} \\
_{ } & _{ } & _{1} & _{ } \\
_{ } & _{\upddots} & _{ } & _{ } \\
_{1} & _{ } & _{ } & _{ }
\end{smallmatrix}\right).$
Denote by $B_{SO_{2n+1}}(\F)$, $N_{SO_{2n+1}}(\F)$ the standard
Borel subgroup and its unipotent radical respectively (see page 2
of \cite{So93} for example). Let $\psi$ be a non-trivial character
of $\F$. We continue to denote by $\psi$ the character of
$N_{SO_{2n+1}}(\F)$ defined by

$$\psi(u)=\psi\bigl(\sum_{k=1}^{n} u_{k,k+1}\bigr).$$
We also view $\psi$ a character of any subgroup of
$N_{SO_{2n+1}}(\F)$. Let $P_{SO_{2n+1}}(\F)$ be the standard
parabolic subgroup of $\soo$ whose Levi part and unipotent radical
are
$$M_{SO_{2n+1}}(\F)=\{ \begin{pmatrix} &_{g} & _{} & _{}\\ &_{}& _{1} &_{}
\\ &_{}& _{} &_{g*}\end{pmatrix} \mid g \in \gln \} \simeq \gln,$$
$$U_{SO_{2n+1}}(\F)=\{ \begin{pmatrix} &_{I_n} & _{x} & _{z}\\ &_{}& _{1} &_{x'}
\\ &_{}& _{} &_{I_n}\end{pmatrix} \in \soo \},$$
where $g^*=J'_{2n+1}{^tg^{-1}}, x'=-{^tx}J'_n$. Define
$$\omega''_n=\begin{pmatrix}
& _{ }& _{ }& _{\omega_n } \\
& _{ } & _{(-1)^n} & _{ } \\ & _{\omega_n} & _{ }\end{pmatrix} \in
SO_{2n+1},$$ and let $\tau$ be a generic representation of $\gln$
identified with $M_{SO_{2n+1}}(\F)$. The local coefficient
$${C_{\psi}^{\soo}(P_{SO_{2n+1}}(\F),s,\tau, \omega_n''^{-1}
)}$$ is defined in the same way as in Chapter \ref{chapter def lc}
via Shahidi`s general construction; see Theorem 3.1 of \cite{Sha
1}. From the second part of Theorem 3.5 of \cite{Sha 2} it follows
there exists $c \in \C^*$ such that
$$C_{\psi}^{\soo}(P_{SO_{2n+1}}(\F),s,\tau, \omega_n''^{-1}
)=c\gamma(\tau,sym^2,2s,\psi).$$ Furthermore, if $\tau$ is
unramified then $c=1$. In Theorem \ref{soudry said always} we have
proven that
$$C_{\psi}^{\mspm}\bigl(\overline{P_{m;0}(\F)},s,\tau,
\omega_m'^{-1} \bigr)=c_\F(s)\frac
{\gamma(\tau,sym^2,2s,\psi)}{\gamma(\tau,s+\half,\psi)},$$ where
$c_\F(s)$ is an exponential factor which equals 1 if  $\F$ is a
p-adic field of odd residual characteristic,
 $\psi$ is normalized and $\tau$ is unramified. Recalling \eqref{gama
tau tau hat} we have proved the following.
\begin{lem} \label{lc soo lc meta}Let $\tau$ be an irreducible admissible
generic representation of $\gln$. There exits an exponential
function $c(s)$ such that
\begin{eqnarray} \label{s0=meta cusp}
&&{C_{\psi}^{\soo}(P_{SO_{2n+1}}(\F),s,\tau, \omega_n''^{-1}
)}{C_{\psi}^{\soo}(P_{SO_{2n+1}}(\F),-s,\widehat{\tau},
\omega_n''^{-1} )}\\ \nonumber &&=
c(s){C_{\psi}^{\mspn}(\overline{P_{n;0}(\F)},s,\tau,
\omega_n'^{-1}
)}{C_{\psi}^{\mspn}(\overline{P_{n;0}(\F)},-s,\widehat{\tau},
\omega_n'^{-1} )}.\end{eqnarray} $c(s)=1$ provided that $\F$ is a
p-adic field of odd residual characteristic,
 $\psi$ is normalized and $\tau$ is unramified.
\end{lem}
\begin{thm} \label{metaplectic and so}Let $\tau$ be an irreducible admissible self dual
supercuspidal representation of $\gln$. Then,
$$I(\tau)=Ind^{\mspn}_{\overline{P_{n;0}(\F)}}
\bigl((\gamma_{\psi}^{-1} \! \! \circ \! \det) \otimes \tau
\bigr)$$ is irreducible if and only if
$$I'(\tau)=Ind^{\soo}_{P_{SO_{2n+1}}(\F)} \tau $$ is irreducible.
\end{thm}
\begin{proof} In both cases we are dealing with a representation induced
from a singular representation of a maximal parabolic subgroup.
Therefore, applying Theorem \ref{silberger} and \eqref{beta is lc}
to these representations, the theorem follows from Lemma \ref{lc
soo lc meta}.
\end{proof}

{\bf Remarks:}

1. One can replace the assumption that $\tau$ is self dual and
replace it with the assumption that $\tau$ is unitary, since by
Theorem \ref{bruhat dimension} the commuting algebras of these
representations are one dimensional if $\tau$ is not self dual.

2. Theorem \ref{metaplectic and so} may be proved without a direct
use of Lemma \ref{lc soo lc meta}. One just has to recall
Corollary \ref{nice for max} and the well known fact that
$I'(\tau)$ is irreducible if and only if $\gamma(\tau,sym^2,0)\neq
0$;
 see \cite{Sha 4}. However, the last fact follows also from the Knapp-Stein dimension theory
and from the theory of local coefficients. In fact, Lemma \ref{lc
soo lc meta} gives more information than Theorem \ref{metaplectic
and so}. This Lemma implies that $\beta(s,\tau,\omega_n'^{-1})$
has the same analytic properties as the Plancherel measure
attached to $\soo$, $P_{SO_{2n+1}}(\F)$ and $\tau$.

3. In general, we expect the same connection between the the
parabolic inductions $Ind^{\mspn}_{\overline{Q(\F)}} \bigl((\tau
\otimes \gamma_\psi^{-1} \circ \det)   \otimes \msig \bigr)$ and
$Ind^{SO_{2n+1}(\F)}_{Q'(\F)} \bigl(\tau \otimes
\theta_\psi(\msig) \bigr )$ , where $\F$ is a p-adic field,
$Q(\F)$ is the standard parabolic subgroup of $\spn$ which has
$GL_m(\F) \times Sp_{2k}(\F)$ as its Levi part, $Q'(\F)$ is the
standard parabolic subgroup of $SO_{2n+1}(\F)$ which has $GL_m(\F)
\times SO_{2k+1}(\F)$ as its Levi part ($r+m=n$), $\tau$ is an
irreducible supercuspidal generic representation of $GL_m(\F)$ and
$\msig$ is an irreducible genuine supercuspidal generic
representation of $\overline{Sp_{2k}(\F)}$. Here
$\theta_\psi(\msig)$ is the generic representation of
$SO_{2k+1}(\F)$ obtained from $\msig$ by the local theta
correspondence. See \cite{JS} and \cite{F} for more details on the
theta correspondence between generic representations of
$\overline{Sp_{2k}(\F)}$ and $SO_{2k+1}(\F)$. In a recent work of
Hanzer and Muic, a progress in this direction was made; see
\cite{HM}.\\

\begin{cor} \label{irr with gama}Let  $\tau$ be an irreducible admissible self dual supercuspidal
representation of $\glm$. Let $\msig$ be a generic genuine
irreducible admissible supercuspidal representation of $\mspk$. If
$I(\tau)$ is irreducible then $I(\tau,\msig)$ is irreducible if
and only if $\gamma(\overline{\sigma} \times \tau,0,\psi)=0$
\end{cor}

\begin{proof} Recalling Theorem \ref{big thm for max}, we only have to
show that  $\gamma(\overline{\sigma} \times \tau,0,\psi)=0$ if and
only if
$${C_{\psi}^{\mspn}(\overline{P_{m;k}(\F)},0,\tau\otimes \msig,
j_{m,n}(\omega_m'^{-1}))}=0.$$ Therefore, from \eqref{gama def},
the definition of $\gamma(\overline{\sigma} \times \tau,0,\psi)$,
the proof is done once we show that
$$C_{\psi}^{\mspm}(\overline{P_{m;0}(\F)},s,\tau, \omega_m'^{-1})$$
is analytic and non-zero in $s=0$. The analyticity of this local
coefficient at $s=0$ follows since by \eqref{s0=meta cusp}
$$C_{\psi}^{\mspm}(\overline{P_{m;0}(\F)},s,\tau, \omega_m'^{-1})
C_{\psi}^{\mspm}(\overline{P_{m;0}(\F)},-s,\tau, \omega_m'^{-1})$$
has the same analytic properties as
$${C_{\psi}^{SO_{2m+1}(\F)}(P_{SO_{2m+1}}(\F),s,\tau, \omega_m''^{-1}
)}{C_{\psi}^{SO_{2m+1}(\F)}(P_{SO_{2m+1}}(\F),-s,\widehat{\tau},
\omega_m''^{-1} )}$$ which is known to be analytic in $s=0$; see
Theorem 5.3.5.2 of \cite{Sil book} (note that that last assertion
does not relay on the fact that $\tau$ is self dual). The fact
that
$$C_{\psi}^{\mspm}(\overline{P_{m;0}(\F)},s,\tau, \omega_m'^{-1}) \neq
0$$ follows from Theorem \ref{silberger} and the assumption that
$I(\tau)$ is irreducible.
\end{proof}

The corollaries below follow from \cite{Sha 4}:

\begin{cor} \label{cor 1}Let $\tau$ be as in Theorem \ref{metaplectic and so}.
Assume that $n \geq 2$. Then $I(\tau)$ is irreducible if and only
if
$$I''(\tau)=Ind^{\spn}_{P_{n;0}(\F)} \tau$$ is reducible.
\end{cor}

\begin{proof}
Theorem 1.2 of \cite{Sha 4} states that $I''(\tau)$ is irreducible
if and only if $I'(\tau)$ is reducible.
\end{proof}

\begin{cor} \label{cor 2} Let $\tau$ be as in Theorem \ref{metaplectic and
so}. If $n$ is odd then $I(\tau)$ is irreducible.
\end{cor}

\begin{proof} Corollary 9.2 of \cite{Sha 4} states that under the
conditions in discussion $I''(\tau)$ is reducible.
\end{proof}


\newpage

\begin{thebibliography}{}

\bibitem{BanPhD} Banks W.D., {\it Exceptional representations on
the metaplectic group}, Dissertation, Stanford University, (1994).

\bibitem{Ban}  Banks W.D., {\it Heredity of Whittaker models on the metaplectic
group.} Pacific J. Math. { 185} (1998), no. 1, pp. 89-96.

\bibitem{Bar} Barthel L., {\it Local Howe correspondence for groups of
similitudes.} J. Reine Angew. Math. { 414}, (1991), pp. 207-220.

\bibitem{Bru61} Bruhat F., {\it Distributions sur un groupe
 localement compact et applications  l`etude des representations des groupes $p$-adiques.}
 Bull. Soc. Math. France 89, (1961), pp. 43-75.

\bibitem{Bud} Budden M., {\it Local coefficient matrices of metaplectic groups.}
J. Lie Theory 16, no. 2, (2006), pp. 239-249.

\bibitem {B} Bump D., {\it Automorphic Forms and Representations}, Cambridge
U. Press, 1988.

\bibitem{BFH} Bump D., Friedberg S., Hoffstein J., {\it p-adic
Whittaker functions on the metaplectic group}, Duke Math. Journal,
Vol { 31}, No. { 2}, (1991), pp. 379-397.

\bibitem {BZ}  Bernstein I. N., Zelevinsky A. V., {\it Representations of the group
$GL(n,F)$ where $F$ is a non-archimedean local field.} Russian
Math. Surveys { 31}, (1976), pp. 1-68.

\bibitem{BZ77} Bernstein I. N., Zelevinsky A. V., {\it
Induced representations of reductive $p$-adic groups},  I. Ann.
Sci. ecole Norm. Sup. (4) 10, no. 4, (1977), 441--472. 22E50

\bibitem{CS}  Casselman W., Shalika J. {\it The unramified principal
series of $p$-adic groups. II. The Whittaker function.} Compositio
Math. 41 (1980), no. 2, pp. 207-231.

\bibitem{F86}  Flicker Y., Kazhdan D.A., {\it Metaplectic
correspondence.}
 Inst. Hautes \'{E}tudes Sci. Publ. Math. No. 64 (1986), pp. 53-110.

\bibitem{F} Furusawa M., {\it On the theta lift from
${\rm SO}\sb {2n+1}$ to $\widetilde{\rm Sp}\sb n$.} J. Reine
Angew. Math. 466 (1995), pp. 87--110.

\bibitem{Gaweb} Garrett P., {\it irreducibles as kernels of
intertwinings among principal series.} Avaliable at
http://www.math.umn.edu/~garrett/m/v/ios.pdf

\bibitem{G} Gelbart S. S., {\it Weil's representation and
the spectrum of the metaplectic group.} Springer-Verlag (1976).

\bibitem{GHP} Gelbart S., Howe R., Piatetski-Shapiro I., {\it
Uniqueness  and existence of Whittaker models for the metaplectic
group.} Israel Journal of Mathematics, Vol { 34}, nos 1-2, (1979),
pp. 21-37.

\bibitem{GS}  Gelbart S., Shahidi F., {\ it Analytic properties of automorphic $L$-functions.}
Perspectives in Mathematics, 6. Academic Press, Inc., Boston, MA,
(1988).

\bibitem{GK} Gelfand I.M., Kajdan D.A., {\it Representations of
the group $GL(n,K)$, where $K$ is a local field.} In {\it Lie
groups and their representations.} Budapest (1971), pp 95-117.

\bibitem{Gol94} Goldberg, D. {\it Reducibility of induced
representations for $\spn$ and $\soo)$.} American Journal of
Mathematics, vol. 116 , no. 5, (1994) pp. 1101-1151.

\bibitem {HMa} Hanzer M., Matic I., {\it Irreducibility of the unitary principal
series of $p$-adic $\widetilde{Sp(n)}$.} Preprint. Available at
http://www.mathos.hr/~imatic

\bibitem{HM} Hanzer M., Muic G., {\it Parabolic induction and
Jacquet functors for the metaplectic groups.} Preprint.

\bibitem{Har73}  Harish-Chandra, {\it  Harmonic analysis on
reductive $p$-adic groups} in  {\it Harmonic analysis on
homogeneous spaces} ( Proceeding of symposia in Pure Mathematics,
Vol. XXVI, Williams Coll., Williamstown, Mass., 1972), Ameraical
mathematical society, Providence, R.I.,pp. (1973), pp. 167-192.




\bibitem{H} Hashizume M., {\it Whittaker models for real reductive
groups}, Japan J. Math., Vol 5, No. 2, (1979),pp (349-401).

\bibitem{I} Ireland K., Rosen M., {\it A classical introduction to the modern number
theory}, Springer, (1998).

\bibitem{J} Jacquet H.,{\it Fonctions de Whittaker associ\'ees aux groupes de
Chevalley.}(French) Bull. Soc. Math. France 95, {1967}, pages
243--309.

\bibitem{J09} Jacquet H., {\it Archimedean Ranking-Selberg
integrals}, in {\it Automorphic forms and L-function II. Local
aspeects. A workshop in honor of Steve Gelbart on the occasion of
his sixtieth birthday. May 15-19, 2006 Rehovot and Tel Aviv},
Amer. Math. Soc. and Bar-Ilan University, (2009), pp. 57-172.

\bibitem{JPS} Jacquet H., Piatetskii-Shapiro, I. I., Shalika, J. A,
\it{Rankin Selberg convolutions.} Amer. J. Math. 105, no. 2,
(1983). pp. 367--464.

\bibitem {JS} Jiang D., Soudry D., {\it The local converse theorem for
${\rm SO}(2n+1)$ and applications.} Ann. of Math. (2) 157 (2003),
no. 3, pp. 743--806.

\bibitem{JS 07} Jiang D., Soudry D., {\it On the genericity of cuspidal
 automorphic forms of $SO(2n+1)$ II}, Compos. Math. 143 (2007), no. 3, pp 721-748.

\bibitem{KP} Kazhdan, D. A., Patterson, S. J,{\it Metaplectic forms},
 Inst. Hautes etudes Sci. Publ. Math. No. 59,(1984), pp. 35--142.

\bibitem{Kim} Kim H.H, {\it Automorphic L-functions} in
{\it Lectures on automorphic $L$-functions.} Fields Institute
Monographs, 20. American Mathematical Society, Providence, RI,
(2004), pp. 97-201.


\bibitem{K} Knapp A.W, Stein E.M, {\it Singular integrals and the
principal series}, IV Proc. Nat. Acad. U.S.A, 72, (1975), pp.
2459-2461.
\bibitem{Kub} Kubota T., {\it Automorphic functions and the
reciprocity law in a number field.} Kinokuniya book store (1969).

\bibitem{Kud} Kudla S.S, {\it Splitting metaplectic covers of
dual reductive pairs.} Israel Journal of Mathematics, Vol { 87},
(1994), pp. 361-401.
\bibitem{L71}  Langlands, R.P., {\it  Euler products.
A James K. Whittemore Lecture in Mathematics given at Yale
University, 1967.} Yale Mathematical Monographs, 1. Yale
University Press, New Haven, Conn.-London, 1971.


\bibitem{L}  Langlands, R.P., {\it On the functional equations
satisfied by Eisenstein series.} Lecture Notes in Mathematics,
Vol. 544. Springer-Verlag, Berlin-New York, 1976.

\bibitem{M} Matsumoto, H., {\it Sur les sous-groupes arithm$\acute{e}$tiques
des groupes semi-simples d$\acute{e}$ploy$\acute{e}$s.}  Ann. Sci.
$\acute{E}$cole Norm. Sup. (4), 2, 1969,  pp. 1-62.

\bibitem{MW} Moeglin C., Waldspurger J.L, {\it Spectral decomposition and Eisenstein
series.} Cambridge Tracts in Mathematics, 113. Cambridge
University Press, Cambridge, 1995.

\bibitem{MVW} Moeglin C., Vign$\acute{e}$ras M.-F., Waldspurger J.L, {\it
Correspondances de Howe sur un corps p-adique.} Springer-Verlag
(1987).

\bibitem{Mo88} Moen, C., {\it Irreducibility of unitary principal
series for covering groups of $SL(2,k)$.}  Pacific Journal of
Mathematics, Vol. 135, no. 1, (1988), pp. 89-110.

\bibitem{P} Perrin P., {\it Reprsentations de Schrodinger, indice
de Maslov et groupe metaplectique.}  insie {\it Noncommutative
harmonic analysis and Lie groups }(Marseille, 1980), pp. 370-407,
Lecture Notes in Math., 880, Springer, Berlin-New York, (1981).

\bibitem{RV} Ramakrishnan D., Valenza R.J, {\it
Fourier analysis on number fields.} Graduate Texts in Mathematics,
186. Springer-Verlag, New York, (1999).

\bibitem{R}  Rao R. R., {\it On some explicit formulas in the theory of Weil
representation.} Pacific Journal of Mathematics, Vol 157, no 2
(1993).

\bibitem{Rod} Rodier F., {\it Whittaker models for admissible representationes of
reductive p-adic split groups}, Proc. of symposia in pure math,
Vol. XXVI, Amer. Math. Soc. (1973), pp 425-430.


\bibitem{T} Tate J.T, {\it Fourier analysis in number fields, and Hecke's zeta-functions} in
. 1967 {\it Algebraic Number Theory} (Proc. Instructional Conf.,
Brighton, 1965), 1967, pp. 305-347.

\bibitem{Ser} Serre, J.P., {\it A course in arithmetic.} Translated from the French.
Graduate Texts in Mathematics, No. 7. Springer-Verlag, New
York-Heidelberg, (1973).

\bibitem{Sha78}  Shahidi F., {\it Functional equation satisfied by certain
$L$-functions.} Compositio Math. 37, no. 2,(1978), pp. 171-207.


\bibitem {Sha80} Shahidi F., { \it Whittaker models for real groups.} Duke Math. J. 47, no. 1, (1980),pp. 99-125.

\bibitem{Sha 1} Shahidi F., {\it On certain $L$-functions.}
Amer. J. Math. 103 (1981), no. 2, pp. 297-355. MR1168488

\bibitem{Sha83} Shahidi F., {\it Local coefficients and
 normalization of intertwining operators for ${\rm GL}(n)$.}
  Compositio Math. 48, no. 3, (1983), pp. 271-295.

\bibitem{Sha84} Shahidi F., {\it Fourier transforms of
intertwining operators and Plancherel measures for $GL(n)$.} Amer.
J. Math. 106, no. 1, (1984). pp. 67--111



\bibitem{Sha85} Shahidi F., {\it Local coefficients as Artin factors for real groups.},
Duke Math. J. Vol 52, no. 4, (1985), pp. 973-1007.


\bibitem{Sha 2} Shahidi F., {\it A proof of Langlands' conjecture on
Plancherel measures; complementary series for $p$-adic groups},
Ann. of Math. (2) 132 (1990), no. 2, pp. 273-330.

\bibitem{Sha 3} Shahidi, F., {\it Langlands' conjecture on
Plancherel measures for $p$-adic groups} in {\it Harmonic analysis
on reductive groups}, Birkhauser Boston, Boston, MA, 1991, pp
277-295.

\bibitem{ShaOxf} Shahidi, F., {\it $L$-functions and representation
 theory of $p$-adic groups.}, in {\it $p$-adic methods and their applications}
 ,Oxford Sci. Publ., Oxford Univ. Press, New York, 1992, pp.
 91-112.


\bibitem{Sha 4} Shahidi, F., {\it Twisted endoscopy and reducibility of
induced representations for $p$-adic groups}, Duke Math. J. 66
(1992), no. 1, pp. 1-41.

\bibitem{Sha 96} Shahidi, F., {\it Intertwining operators, L-functions and representation
theory.} Lecture notes of the eleventh KAIST Mathematics Workshop,
Lecture notes of the KAIST Mathematics Workshop, (1996).
Avaliable at
http://www.math.rutgers.edu/~sdmiller/l-functions/shahidi-korea.pdf

\bibitem{Sh} Shalika J.A., {\it The multiplicity 1 Theorem for
$GL_n$.} Ann. of Math. 100(1) (1974), pp 172-193.


\bibitem{Sil} Silberger A.J., {\it the Knapp-Stein dimension theorem
for p-adic groups}, Proceeding of the Ameraican mathematical
society, Vol. 68, No. 2, (1978), pp. 243-246.

\bibitem{Sil79} Silberger A.J., {\it Correction: ``The Knapp-Stein dimension theorem for
 $p$-adic groups''}, Proceeding of the Ameraican mathematical
society, Vol. 76, No. 1, (1979), pp. 169-170.

\bibitem{Sil book} Silberger A.J., {\it Introduction to harmonic analysis
on reductive $p$-adic groups.}  Mathematical Notes, 23. Princeton
University Press, Princeton, N.J.; University of Tokyo Press,
Tokyo, 1979.


\bibitem{So1} Soudry D., {\it A uniqueness theorem for representatons
of GSO(6) and the strong multiplicty one Theorem for generic
representations of GSp(4).} Israel Journal of Mathematics, Vol 58,
no. 3, (1987), pp. 257-287.

\bibitem{So93}  Soudry D., {\it  Rankin-Selberg convolutions for
${SO}_{2l+1}\times{GL}_n$: local theory.} Mem. Amer. Math. Soc.
105 , no. 500,(1993).

\bibitem {SB} Sun B.,{\it Contragredients of irreducible representations in theta
 correspondences.} Preprint. arXiv:0903.1418v1, (2009).

\bibitem{Sz} Szpruch D., {\it Uniqueness of Whittaker model for the
metaplectic group.} Pacific Journal of Mathematics, Vol. 232, no.
2, (2007),  pp. 453-469.

\bibitem{Sz09} Szpruch D., {\it Computation of the local coefficients for principal
series representations of the metaplectic double cover of
$SL_2(\F)$.} Journal of Number Theory, Vol. 129, (2009), pp.
2180-2213.

\bibitem{Wa} Waldspurger, J.L. {\it Correspondance de Shimura.}
J. Math. Pures Appl. J. Math. Pures Appl. (9) 59 (1980), no. 1,
pages 1-132.

\bibitem{Wa03} Waldspurger, J.L., {\it La formule de Plancherel pour les groupes
$p$-adiques (d'apr\`{e}s Harish-Chandra)} , J. Inst. Math. Jussieu
2 (2003), no. 2, 2003 , pp. 235--333.

\bibitem{Wang} Wang. C.J, {\it On the existence of cuspidal distinguished
representations of metaplectic groups}, Dissertation, Ohio State
University, (2003).


\bibitem{Weil}Weil A., {\it Sur certains groupes d'operateurs unitaires}, Acta Math Vol { 111}, 1964, pp.
143--211.

\bibitem{Whi}  Whittaker E. T.,  Watson, G. N, {\it
A course of modern analysis.}  Fourth edition, Reprinted Cambridge
University Press, New York (1962).

\bibitem{W}  Weiss E., {\it Algebraic number theory.}  Dover
publications INC, (1998).

\bibitem{Zo09} Zorn C., {\it Reducibility of the principal series
for Mp(2,F) over a p-adic field.} Canadian Journal of Mathematics,
to appear. Avaliable at http://www.math.ohio-state.edu/~czorn/



\end{thebibliography}
\end{document}